\documentclass[a4paper, 10pt]{article}

\pdfoutput=1
\usepackage{amsmath}
\usepackage{graphicx}								
\usepackage{subfig}
\usepackage{amsfonts} 							
\usepackage{dsfont}
\usepackage{algorithm}              
\usepackage{algorithmic}
\usepackage{booktabs} 							
\usepackage{multirow}								
\usepackage{setspace}
\usepackage{mathrsfs}
\usepackage{amssymb}
\usepackage{tikz}
\usepackage{color}
\usepackage{xcolor}
\usepackage{mathtools}
\usepackage{xfrac}
\usepackage{upgreek}

\usepackage{framed}
\usepackage{afterpage}
\usepackage{capt-of}
\usepackage{stackengine}

\usepackage{colortbl}
\definecolor{gainsboro}{rgb}{0.86, 0.86, 0.86}
\definecolor{gray(x11gray)}{rgb}{0.75, 0.75, 0.75}
\definecolor{lightgray}{rgb}{0.83, 0.83, 0.83}
\definecolor{silver}{rgb}{0.75, 0.75, 0.75}

\usepackage{geometry}
\usepackage{color}
\definecolor{deepblue}{rgb}{0,0,0.5}
\definecolor{deepred}{rgb}{0.6,0,0}
\definecolor{deepgreen}{rgb}{0,0.5,0}
\definecolor{dollarbill}{rgb}{0.52, 0.73, 0.4}
\definecolor{grannysmithapple}{rgb}{0.66, 0.89, 0.63}
\definecolor{inchworm}{rgb}{0.7, 0.93, 0.36}
\definecolor{olivine}{rgb}{0.6, 0.73, 0.45}
\definecolor{pistachio}{rgb}{0.58, 0.77, 0.45}
\definecolor{teagreen}{rgb}{0.82, 0.94, 0.75}
\definecolor{yellow-green}{rgb}{0.6, 0.8, 0.2}
\definecolor{junebud}{rgb}{0.74, 0.85, 0.34}
\definecolor{babypink}{rgb}{0.96, 0.76, 0.76}
\definecolor{beaublue}{rgb}{0.74, 0.83, 0.9}
\definecolor{corn}{rgb}{0.98, 0.93, 0.36}
\definecolor{DeepSkyBlue}{rgb}{0,0.75,1}
\definecolor{DeepSkyBlue1}{rgb}{0,0.75,1}
\definecolor{DeepSkyBlue2}{rgb}{0,0.7,0.93}
\definecolor{DeepSkyBlue3}{rgb}{0,0.6,0.8}
\definecolor{DeepSkyBlue4}{rgb}{0,0.41,0.54}
\definecolor{DarkSeaGeen4}{rgb}{0.,0.45,0.}
\definecolor{darkcandyapplered}{rgb}{0.64, 0.0, 0.0}
\definecolor{byzantium}{rgb}{0.58, 0.0, 0.83}
\definecolor{lightblue}{rgb}{0.,0.,0.9}
\definecolor{darkblue}{rgb}{0.,0.,0.45}
\definecolor{verydarkblue}{rgb}{25,25,112}
\definecolor{darkspringgreen}{rgb}{0.09, 0.45, 0.27}
\definecolor{darkcoral}{rgb}{0.8, 0.36, 0.27}
\definecolor{crimsonglory}{rgb}{0.75, 0.0, 0.2}
\definecolor{darkblue}{rgb}{0.0, 0.0, 0.55}
\definecolor{darkelectricblue}{rgb}{0.33, 0.41, 0.47}
\definecolor{applegreen}{rgb}{0.55, 0.71, 0.0}

\usepackage{tikz}	
\usetikzlibrary{arrows,chains,matrix,positioning,scopes}

\makeatletter
\tikzset{join/.code=\tikzset{after node path={%
\ifx\tikzchainprevious\pgfutil@empty\else(\tikzchainprevious)%
edge[every join]#1(\tikzchaincurrent)\fi}}}

\makeatother
\tikzset{>=stealth',every on chain/.append style={join}, every join/.style={->}}
\tikzstyle{labeled}=[execute at begin node=$\scriptstyle,   execute at end node=$]

\tikzset{
>=stealth',
help lines/.style={dashed, thick},
axis/.style={<->},
important line/.style={thick},
connection/.style={thick, dotted},
}

\newcommand {\Sum}   {\sum\limits}

\newcommand {\R}		 {\mathbf{r}}

\newcommand {\Ieff}  {I_{\rm eff}}

\newcommand {\Rd}    {{\mathds{R}}^d}

\def \EnergyNorm#1  {{\mid\!\mid\!\mid #1 \mid\!\mid\!\mid}^2 }   

\def \dvrg       {\mathrm{div}}	

\def \traspose#1 {{#1}^{rm T}}
\def \laplace    {\Delta}

\newcommand {\vectorv}  {\boldsymbol {v}}

\newcommand {\vectorp}  {\boldsymbol {p}}
\newcommand {\vectoru}  {\boldsymbol {u}}
\newcommand {\vectorw}  {\boldsymbol {w}}

\newcommand {\vectortau}  {\boldsymbol {\tau}}

\newcommand {\flux}     {\boldsymbol {y}}

\newcommand {\hhat} {\hat{h}} 
\newcommand {\Khat} {\widehat{K}} 
\newcommand {\Bhat} {\widehat{B}} 
\newcommand {\Rhat} {\widehat{R}} 
\newcommand {\Omegahat} {\widehat{\Omega}} 

\def \dx       {\mathrm{\:d}x}

\def\L#1{L^{#1}}

\def\HD#1#2{H^{#1}_{#2}}

\def\mdI{\overline{\mathrm m}^2_{\mathrm{d}}}

\def\mfI{\overline{\mathrm m}^2_{\mathrm{eq}}}

\def\maj#1{\overline{\mathrm M}_{#1}}

\def\mij#1{\underline{\mathrm M}_{#1}}
\def\mdIK{\overline{\mathrm m}^2_{\mathrm{d}, K}}
\def\mdInosq{\overline{\mathrm m}_{\mathrm{d}}}
\def\mfInosq{\overline{\mathrm m}_{\mathrm{f}}}

\newcommand {\CFriedrichs} {C_{{\rm F}\Omega}}

\newcommand*\rfrac[2]{{}^{#1}\!/_{#2}}
\newcommand{\minus}{\scalebox{0.5}[1.0]{\( - \)}}

\newtheorem{theorem}{Theorem}{\bf}{\it}
{\bf}{\it}
{\bf}{\it}
\newtheorem{example}{Example}{\bf}{}
\newtheorem{remark}{Remark}{\bf}{\it}

\def\ProofBegin{\noindent{\bf Proof:} \:}
\def\ProofEnd{{\hfill $\square$}}

\definecolor{formalshade}{rgb}{0.95,0.95,1}

\usepackage{listings} 
\newcommand\pythonstyle{
\lstset{
		language=Python,
		basicstyle=\footnotesize\ttfamily,
		otherkeywords={self},             
		keywordstyle=\bfseries\color{deepblue}\bfseries,
		commentstyle=\itshape\color{purple!40!black},
		stringstyle=\color{deepgreen},
		frame=single,                         
		showstringspaces=false,            %
		belowcaptionskip=.75\baselineskip
}}
\lstnewenvironment{python}[1][]
{
\pythonstyle
\lstset{#1}
}
{}

\newcommand\pythoninline[1]{{\pythonstyle\lstinline!#1!}}

\usepackage{colortbl}

\usetikzlibrary{%
  arrows,%
  shapes.misc,
  shapes.arrows,%
  chains,%
  matrix,%
  positioning,
  scopes,%
  decorations.pathmorphing,
  shadows,%
  backgrounds,
  fit,positioning,shapes.symbols,chains
}
\usetikzlibrary{shapes.geometric,calc}

\author{S. Matculevich%
\thanks{RICAM Linz, Johann Radon Institute, AT-4040 Linz, 
\texttt{svetlana.matculevich@ricam.oeaw.ac.at}}}
\begin{document}
\title{Functional approach to the error control in adaptive IgA schemes\\
for elliptic boundary value problems}
\maketitle

\begin{abstract}
This work presents a numerical study of {\em functional type} a posteriori error estimates for {\em IgA 
approximation schemes} in the context of elliptic boundary-value problems. Along with the detailed 
discussion of the most crucial properties of such estimates, we present the algorithm of a reliable solution  
approximation together with the scheme of an efficient a posteriori error bound generation. 
In this approach, 
we take advantage of B-(THB-) spline's high smoothness for the auxiliary vector function reconstruction,
which, at the same time, allows to use much coarser meshes and decrease the number of unknowns 
substantially. The most representative numerical results, obtained during a systematic testing of error 
estimates, are presented in the second part of the paper. The efficiency of the obtained error bounds is 
analysed from both the error estimation (indication) and the computational expenses points of view. Several 
examples illustrate that functional error estimates (alternatively referred to as the {\em majorants} and 
{\em minorants} of deviation from an exact solution) perform a much 
sharper error control than, for instance, residual-based error estimates. Simultaneously, assembling and 
solving routines for an auxiliary variables reconstruction, which generate the majorant (or minorant) 
of an error, can be executed several times faster than the routines for a primal unknown. 
\end{abstract}

\section{Introduction}

The investigation of effective adaptive refinement procedures has recently become an active area of research 
in the context of fast and efficient solvers for isogeometric analysis (IgA) \cite{HughesCottrellBazilevs2005, 
HughesCottrellBazilevsBook2005}. The adaptivity scheme is naturally linked with reliable and 
quantitatively efficient a posteriori error estimation tools. The latter ones are expected to identify the parts 
of a considered computational domain with relatively high {discretisation} errors and provide a fully 
automated refinement strategy in order to reach desired accuracy levels for an approximate solution. 

Due to a tensor-product setting of IgA splines, mesh refinement has global effects, which include a large 
percentage of superfluous control points in data analysis, unwanted ripples on the surface, etc. 
These issues produce certain challenges at the design stage as well as complications in handling big 
amounts of data, and, therefore, naturally trigger the development of local refinement strategies for IgA. 
At the moment, four different IgA approaches for adaptive mesh refinement are known, i.e., 
T-splines, hierarchical splines, PHT-splines, and LR splines.

The localised splines of the first type, {\em T-splines}, were introduced in 
\cite{LMR:Sederbergetall2003, LMR:Sederbergetall2004} and analysed in \cite{LMR:Bazilevsetall2010, 
LMR:Veigaetall2011, LMR:Scottetall2011, LMR:Scottetall2012}. They are based on the  T-junctions that 
allow eliminating redundant control points from NURBS model. The thorough study confirmed that this 
approach generates an efficient local refinement algorithm for {\em analysis-suitable} T-splines 
\cite{LMR:Lietall2012} and avoids the excessive propagation of control points. 
In \cite{VeigaBuffaChoSangalli2012, VeigaBuffaSangalliVazquez2013}, it was proposed to characterise 
such splines as dual-compatible T-splines, and in \cite{MorgensternPeterseim2015} a refinement strategy 
with linear complexity was described for the bivariate case. 

The alternative approach that implies the local control of refinement is based on {\em hierarchical B-splines} 
(HB-splines), such that a selected refinement region basis functions are replaced with the finer ones 
of the same type. The procedure of designing a basis for the hierarchical spline space was suggested in 
\cite{ForseyBartels1988, LMR:Kraft1997, GreinerHormann1997} and extended in 
\cite{LMR:Vuongetall2011, EvansScottLiThomas2015, Schillingeretall2012}. Such construction  
guarantees the linear independence of the basis and provides nested approximation spaces. However, since 
the partition of unity is not preserved for these splines, {\em truncated hierarchical B-splines} 
(THB-splines) have been developed (see \cite{LMR:GiannelliJuttlerSpeleers2012}). In addition to 
good stability and approximation properties inherited from HB-splines \cite{LMR:GiannelliJuttler2013, 
LMR:SpeleersManni2016}, THB-splines form a convex partition of unity. Therefore, they are suitable for 
the application in CAD. Various usage of THB-spline for arbitrary topologies can be found, e.g., in
\cite{LMR:WeiZhangHughesScott2015, LMR:ZoreJuttler2014, LMR:ZoreJuttlerKosinka2016}.

The locally defined splines of the third type, namely, {\em polynomial splines over hierarchical T-meshes}, 
are constructed for the entire space of piecewise polynomials with given smoothness on the subdivision of 
considered domain. The corresponding application can be found in \cite{LMR:NguyenThanhetall2011, 
LMR:Wangetall2011}. However, in this case, one must assume the reduced regularity of basis 
\cite{LMR:Dengetall2008} or fulfil a certain constraint on admissible mesh configuration 
\cite{Wuetall2012}.  {In \cite{NguyenThanhZhou2017}, 
the adaptive procedure in based on the recovery-based error estimator, in which discontinuous enrichment 
functions are added to the IgA approximation. }

Finally, {\em locally refined splines} (LR-splines) rely on the idea of splitting basis functions. 
This technique achieves localisation but creates difficulties with linear independence 
\cite{LMR:DokkenLychePettersen2013}, which has been studied in \cite{LMR:Bressan2013, 
LMR:BressanJuttler2015}. The application of such type of splines has been thoroughly investigated in 
\cite{LMR:DokkenLychePettersen2013}. In \cite{LMR:JohannessenRemonatoKvamsdal2015}, one can find 
the summary of a detailed comparison of (T)HB-splines and LR splines with respect to sparsity and 
condition numbers. The study concludes that even though LR splines have smaller support than 
THB-splines, the numerical experiments did not reveal any significant advantages of the first ones with respect 
to the sparsity patterns or condition numbers of mass and stiffness matrices.

The refinement tools of IgA mentioned above were combined with various a posteriori error estimation 
techniques. For instance, the a posteriori error estimates based on hierarchical splines were investigated 
in \cite{LMR:DorfelJuttlerSimeon2010,LMR:Vuongetall2011}.
In \cite{LMR:Johannessen2009, LMR:Wangetall2011, LMR:BuffaGiannelli2015, 
LMR:KumarKvamsdalJohannessen2015}, 
the authors used the residual-based a posteriori error estimates and their modifications in order to construct 
mesh refinement algorithms. The latter ones, in particular, require the computation of constants related to
the Clement-type interpolation operators, which are mesh-dependent and often difficult to compute for 
general element shapes. Moreover, these constants must be re-evaluated every time a new mesh is 
generated. The goal-oriented error estimators, which are rather naturally adapted to practical applications, 
have lately been introduced for IgA approximations and can be found in
\cite{LMR:ZeeVerhoosel2011,LMR:DedeSantos2012,LMR:Kuru:2013,LMR:Kuruetall2014}.

In the current work, the terms \emph{error estimate} and \emph{error indicator} distinguish from 
each other. The first one is considered as the total upper (or lower) bound of true energy error. These are 
very important characteristics related to the approximate solution since they can be used to judge whether 
obtained data are reliable or not. In order to locate the areas of the discretised domain that have the highest error
in the approximation, a quantitively sharp \emph{error indicator} is required. The methods of a posteriori 
error estimation listed above are rather error indicators in this terminology and indeed were successfully 
used for mimicking the approximation error distribution. However, their use in the error control, i.e., a reliable 
estimation of the accuracy of obtained data, is rather heuristic in nature. 

Below we investigate a different {\em functional} method providing fully guaranteed error 
estimates, the upper (and lower) bounds of the exact error in the various weighted norms equivalent to 
the global energy norm. These estimates include only global constants (independent of the mesh 
characteristic $h$) and are valid for any approximation from admissible functional space. One of 
the most advantageous properties of functional error estimates is their independence 
of the numerical method used for calculating approximate solutions. The strongest assumption about 
approximations is that they are conforming in the sense that they belong to a certain natural Sobolev 
space suited for the problem. It is important to emphasise that this is still a rather weak assumption and 
that no further restrictions, such as Galerkin orthogonality, are needed. 

Functional error estimates were 
initially introduced in \cite{LMR:Repin:1997, LMR:Repin:1999} and later applied to different mathematical 
models summarised in monographs \cite{NeittaanmakiRepin2004, LMR:RepinDeGruyterMonograph:2008, 
Malietall2014}. They provide guaranteed, sharp, and fully computable upper and lower bounds of errors. 
A pioneering study on the combination of functional type error estimates with the IgA approximations 
generated by tensor-product splines is presented in \cite{LMR:KleissTomar2015} for elliptic boundary 
value problems (BVP). The extensive numerical tests presented in this work confirmed that majorant 
produces not only good upper bounds of the error but also a quantitatively sharp error indicator. 
Moreover, the authors suggest the heuristic algorithm that allows using the smoothness of B-splines for 
a rather efficient calculation of true error upper bound. 

The current work further extends the ideas used in \cite{LMR:KleissTomar2015} for B-splines (NURBS)
and combines the functional approach to the error control with THB-splines. Moreover, our focus is 
concentrated not only on the qualitative and quantitative performance of error estimates but also on 
the required computation time for their reconstruction. The systematic analysis of majorant's numerical 
properties is based on a collection of extensive tests performed on the problems of different complexity. 
For the error control implemented with the help of tensor-structured B-splines (NURBS) and 
THB-splines, we manage to obtain an impressive speed-up in majorant reconstruction by exploiting 
high smoothness of B-splines to our advantage. However, for the problems with sharp local changes 
or various singularities in the solution, the THB-splines implementation in G+Smo restricts the performance 
speed-up when it comes to solving the optimal system for the error majorant as well as for its element-wise 
evaluation. We restrict this study only to the domains modelled by a single patch, which provides at least 
$C^1$-continuity of the approximate solutions inside the patch. However, the application of studied 
majorants can be extended to a multi-patch domain, since the error estimates for stationary problems are 
flexible enough to handle fully non-conforming approximations (this issue has been in details addressed 
in \cite{LazarovRepinTomar2009, TomarRepin2009, RepinTomar2011}). 

The error control for the problems defined on domains of complicated shapes induces another issue related 
to the estimation of Friedrichs' constant used by functional error estimates not only as the weight but also 
as the geometrical characteristics of the considered problem. When such domains are concerned, one can perform 
their decomposition into a collection of non-overlapping convex sub-domains, such that the global constant 
can be replaced by constants in local embedding inequalities (Poincar\'e and Poincar\'e-type inequalities 
\cite{Poincare1890, Poincare1894}). The reliable estimates of these local constants can be found in 
\cite{PayneWeinberger1960, Bebendorf2003a, NazarovRepin2014, LMR:MatculevichRepin:2015}. 
The derivation of functional error estimates exploiting these ideas is discussed in 
\cite{LMR:RepinDeGruyterMonograph:2008, Repin2015} for the elliptic BVP and in 
\cite{LMR:MatculevichRepinPoincare:2014, LMR:MatculevichNeitaanmakiRepin:2015, 
MatculevichRepinDiffUr2016} for the parabolic initial boundary value problem (I-BVP). In order to use this 
method, one needs to impose a crucial restriction on the multi-patch configuration, namely, each patch 
must be a convex sub-domain. Since in the IgA framework patches are treated as mappings from the 
reference domain $\widehat{\Omega}  = (0, 1)^d$, the estimation of local constants is reduced to the 
analysis of the IgA mapping and calculating the corresponding constant for $\widehat{\Omega}$.

The paper proceeds with the following structure. Section \ref{sec:func-error-estimates} formulates the 
general statement of the considered problem and recalls the definition of functional error estimates 
and their main properties in the context of reliable energy error estimation and efficient 
error-distribution indication. The next section serves as an overview of IgA techniques used in 
the current work, i.e., B-splines, NURBS, and THB-splines. 
In Section \ref{sec:funct-estimates-and-iga}, we focus on the algorithms and details of the functional 
error estimates integration into the IgA framework. Last but not least, Section 
\ref{sec:numerical-examples} presents the systematic selection of most relevant numerical examples 
and obtained results that illustrate numerical properties of studied error estimates and indicators.  

\section{Functional approach to the error control}
\label{sec:func-error-estimates}
In this section, we present a model problem, recall the well-posedness results for linear parabolic PDEs, 
which have been thoroughly studied in \cite{Ladyzhenskaya1985,Zeidler1990A,Wloka1987}. We also 
introduce functional a posteriori error estimates for the stated model and discuss its crucial properties. 

Let $\Omega \subset \mathds{R}^d$, $d = \{ 2, 3\}$, be a bounded domain with Lipschitz boundary 
$\Gamma := \partial \Omega$. The elliptic BVP is formulated as 
\begin{alignat}{3}
	- \dvrg \,( A \nabla u) & =\, f	& \quad \mbox{in} \quad \Omega, \label{eq:elliptic-equation-1}\\
  u & =\, 0	& \quad \mbox{on} \quad \Gamma, \label{eq:dirichlet-bc-1}
\end{alignat}
where $f$ supposed to be in $\L{2}(\Omega)$. Alternatively, the problem 
\eqref{eq:elliptic-equation-1}--\eqref{eq:dirichlet-bc-1} can be viewed a system with two (primal and dual) 
unknowns
\begin{alignat}{3}
	- \dvrg \, \vectorp & =\, f	      
	& \quad \mbox{in} \quad \Omega,\label{eq:elliptic-equation}\\
  \vectorp & =\, A \nabla u & \quad \mbox{in} \quad \Omega,\label{eq:dual-part}\\
  u & =\, 0	& \quad \mbox{on} \quad \Gamma.\label{eq:dirichlet-bc}
\end{alignat}
%
We assume that the operator $A$ is symmetric and satisfies the condition of uniform ellipticity for almost 
all (a.a.) $x \in \Omega$, which reads
\begin{equation}
\underline{\nu}_A |\xi|^2 \leq A(x) \: \xi \cdot \xi \leq \overline{\nu}_{A} |\xi|^2, 
\quad \mbox{for all} \quad 
\xi \in \Rd,
\label{eq:operator-a}
\end{equation}
with $0 < \underline{\nu}_A \leq \overline{\nu}_A < \infty$.
Throughout the paper, the following notation for the norms is used:
\begin{equation*}
  \| \, \vectortau \, \|^2_{A, \Omega} := (A \vectortau, \vectortau)_{\Omega}, \quad 
  \| \, \vectortau \, \|^2_{A^{-1}, \Omega} := (A^{-1} \vectortau, \vectortau)_{\Omega}, 
  \quad 
  \mbox{for all}
  \quad 
  \vectortau \in [\L{2}(\Omega)]^{d}, 
\end{equation*}
where $(A \vectoru, \vectorv)_{\Omega} := \int_\Omega A \vectoru \cdot \vectorv \dx$ 
stands for a weighted $\L{2}$ scalar-product for all 
$\vectoru, \vectorv \in [\L{2}(\Omega)]^{d}$.
%
After multiplying (\ref{eq:elliptic-equation-1}) by the test function 
$$\eta \in \HD{1}{0}(\Omega) := \big\{ u \in \L{2}(\Omega) \; \mid \;
\nabla u \in \L{2}(\Omega), \; u |_{\Gamma} = 0 \, \big\},$$
we arrive at the standard generalised formulation of (\ref{eq:elliptic-equation-1})--(\ref{eq:dirichlet-bc-1}): 
find $u \in \HD{1}{0}(\Omega)$ satisfying the integral identity
\begin{equation}
	a(u, \eta) := (A \nabla{u}, \nabla{\eta})_\Omega 
	= (f, \eta)_\Omega =: l({\eta}), \quad \forall \eta \in \HD{1}{0}(\Omega).
	\label{eq:generalized-statement}
\end{equation}
According to \cite{Ladyzhenskaya1985}, the generalised problem \eqref{eq:generalized-statement} has 
a unique solution in $\HD{1}{0}(\Omega)$ provided that $f \in \L{2}(\Omega)$ and condition 
\eqref{eq:operator-a} holds. 

We consider the functional error estimates, which provide a guaranteed two-sided bound of the distance 
$e := u - v$ between the generalised solution of (\ref{eq:generalized-statement}) $u$ and any 
function $v \in \HD{1}{0}(\Omega)$. It is important to emphasise that the suggested functional approach 
to error estimates' derivation is universal for any numerical method used to discretise bilinear form 
\eqref{eq:generalized-statement}. This fact makes it rather unique in comparison with alternative 
approaches, which are always tailored to the discretised version of the identity $a(u, \eta) = l(\eta)$. Later on, 
the considered $v$ is generated numerically, and the distance to $u$ is evaluated in terms of the total 
energy norm
\begin{equation}
	|\!|\!| e |\!|\!|^2_\Omega := \| \, \nabla e \, \|^2_{A, \Omega}
\label{eq:enery-norm-error}
\end{equation}
as well as its element-wise contributions $\| \, \nabla e \, \|^2_{A, K}$ such that
$$|\!|\!| e |\!|\!|^2_\Omega := \sum_{K \in \mathcal{K}_h} \| \, \nabla e \, \|^2_{A, K}.$$
Here, $K$ represents the elements of the mesh $\mathcal{K}_h$ introduced on 
$\Omega$. Hence, besides providing the guaranteed upper bound of total error \eqref{eq:enery-norm-error}, 
the majorant yields a quantitatively sharp indicator of the local error distribution. 

{
\begin{remark}
It has been proved that the integrand of the majorant 
tends to the distribution of the true error in the sense of the measure, if the majorant 
globally converges (from above) to the true value of the error. In other words, 
the measure of the set, where the local difference between majorant and true error 
is bigger than a given epsilon, tends to zero. For the purposes of the error indication, 
it is enough (\cite[Section 3]{LMR:RepinDeGruyterMonograph:2008}). 
\end{remark}
}
To derive the upper bound, we need to transform \eqref{eq:generalized-statement} by subtracting 
$a(v, \eta)$ from left- (LHS) and right-hand side (RHS) and by setting $\eta = e$. Thus, one obtains the 
{error identity}
\begin{equation}
	|\!|\!| e |\!|\!|^2_\Omega \,
	= \big(f, e \big)_\Omega - \big(A \nabla{v}, \nabla e \big)_\Omega.
	\label{eq:energy-balance-equation}
\end{equation}
%
The main idea of functional approach is the introduction of an auxiliary vector-valued variable 
\begin{alignat*}{2}
	\flux \in H(\Omega, {\rm div}) := 
	\Big\{  \flux \in [\L{2}(\Omega)\big]^{d}\;\big|\; 
	& \dvrg \flux \in \L{2} (\Omega) \Big \}
\end{alignat*} 	
satisfying 
\begin{equation}
(\dvrg \flux, v)_\Omega + (\flux, \nabla v)_\Omega = 0.
\label{eq:flux-identity}
\end{equation}
In further calculations, the above-introduced variable allows an additional degree of freedom for the 
majorant (additional optimisation step), whereas, for instance, the residual error estimates do not have 
this flexibility in improving its values.
Next, we add the identity \eqref{eq:flux-identity} to the RHS of \eqref{eq:energy-balance-equation}, 
which yields
\begin{equation}
	|\!|\!| e |\!|\!|^2_\Omega
	= \big(f + \dvrg \flux, e \big)_\Omega 
	+ \big(\flux - A \nabla{v}, \nabla e \big)_\Omega.
	\label{eq:energy-balance-equation-with-flux}
\end{equation}
%
The {\em equilibrated} and {\em dual} {\em residual-functionals} obtained in the RHS of 
\eqref{eq:energy-balance-equation-with-flux} mimic equations \eqref{eq:elliptic-equation} and 
\eqref{eq:dual-part}, respectively, and are denoted by 
\begin{alignat}{2}
	\R_{\rm eq}  (v, \flux) & 
	:= f + \dvrg \flux \quad \mbox{and} \quad
	\R_{\rm d}  (v, \flux) 
	:= \flux - A \nabla{v}.
	\label{eq:residuals}
\end{alignat}
%
\begin{theorem}
\label{th:theorem-minimum-of-majorant-I}
(a) For any functions $v \in \HD{1}{0}(\Omega)$ and $\flux \in H(\Omega, {\rm div})$, we have the 
estimate
\begin{equation}
	|\!|\!| e |\!|\!|^2_\Omega 
	\leq \overline{\rm M}^2 (v, \flux; \beta)
	:= (1 + \beta)\, \|\,\R_{\rm d} \,\|^2_{A^{- 1}, \Omega} 
	 + (1 + \tfrac{1}{\beta}) \,\tfrac{\CFriedrichs^2}{\,\underline{\nu}_A}
	   \big\| \, \R_{\rm eq} \, \big\|^2_{\Omega},
	\label{eq:majorant-1}
\end{equation}
%
where the residuals $\R_{\rm d}$ and $\R_{\rm eq}$ are defined in \eqref{eq:residuals}, $\beta$ is a 
positive parameter, and $\CFriedrichs$ is the constant in the Friedrichs inequality \cite{Friedrichs1937} 
\begin{equation*}
	\| v \|_{\Omega} \leq \CFriedrichs \| \nabla v \|_{\Omega},
	\quad \forall v \in \HD{1}{0}(\Omega).
	\label{eq:friedrichs-inequality}
\end{equation*}

\noindent
(b) For $\beta >0$, the variational problem
\begin{equation*}
\inf\limits_{
\begin{array}{c}
v\in \HD{1}{0}(\Omega)\\[1pt]
\flux \in H(\Omega, {\rm div})
\end{array}
}
\maj{} (v, \flux; \beta) 
\end{equation*}
has a solution (with the corresponding zero-value for the functional), and its minimum is attained if and 
only if $v = u$ and $\flux = A \nabla u$.
\end{theorem}
\ProofBegin For the detailed proof of this theorem we refer the reader to 
\cite[Section 3.2]{RepinDeGruyter2008}.
\ProofEnd

For the real-life problems, the exact solution is usually not known, therefore evaluation of 
$|\!|\!| e |\!|\!|^2_\Omega$ becomes impossible. For the efficiency verification of 
$\overline{\rm M}^2 (v, \flux; \beta)$, the lower bound of the error (further also referred as minorant) 
has been derived using variational arguments (see Theorem \ref{th:theorem-maximum-of-minorant}).

\begin{theorem}
\label{th:theorem-maximum-of-minorant}
For any functions $v, w \in \HD{1}{0}(\Omega)$, we have the estimate
\begin{equation}
	|\!|\!| e |\!|\!|^2_\Omega 
	\geq \sup_{w \in \HD{1}{0}(\Omega)} \underline{\rm M}^2 (v, w)
	:= 2 \,(J(v) - J(w)), \quad 
	J(v) := (f, v) - \tfrac{1}{2} \, \| \nabla v \|^2_\Omega,
	\label{eq:minorant-1}
\end{equation}
%
where $J(v)$ is the variational functional of the problem \eqref{eq:elliptic-equation}--\eqref{eq:dirichlet-bc}.
\noindent
\end{theorem}
\ProofBegin For the detailed proof of this theorem, we refer the reader to 
\cite[Section 4.1]{RepinDeGruyter2008}.
\ProofEnd

\vskip 10pt
\noindent
Remarks below summarise several essential properties of the error estimate derived in
Theorem \ref{th:theorem-minimum-of-majorant-I}.

\begin{remark}
\rm 
Each term on the RHS of \eqref{eq:majorant-1} serves as the penalty of the error that might occur in  
equations \eqref{eq:elliptic-equation} and \eqref{eq:dual-part}. The positive weight $\beta$ can 
be selected optimally in order to get the best value of the majorant. The constant $\CFriedrichs$ acts 
as a geometric characteristic for the considered domain $\Omega$ (unlike, for instance, in the least-square 
methods, where the weights are selected after the terms $\|\,\R_{\rm d} \,\|^2_{A^{- 1}, \Omega}$ and
$\| \, \R_{\rm eq} \,\|^2_{\Omega}$ were calculated). From the author's point of view, this 
constant is essential and cannot be excluded since it scales proportionally to the diameter of the considered 
$\Omega$. Moreover, 
in order to guarantee the reliability of $\overline{\rm M} (v, \flux; \beta)$, the constant $\CFriedrichs$ must 
be estimated from above in a reliable way. Since in practice the term $\| \, \R_{\rm eq} \, \|^2_{\Omega}$ 
is rather small compared to the dominating term $\| \, \R_{\rm d} \, \|^2_{\Omega}$, the Friedrichs 
constant can be replaced by some penalty constant $C \geq \CFriedrichs$ 
(even though it might affect the ratio of the majorant to the error). 
In what follows, to characterise the efficiency of \eqref{eq:majorant-1}, we use the quantity 
$\Ieff(\overline{\rm M}) := {\overline{\rm M}} / {|\!|\!| e |\!|\!|_\Omega}$ that measures the
gap between $\overline{\rm M} (v, \flux; \beta)$ and $|\!|\!| e |\!|\!|_\Omega$.
\end{remark}

\begin{remark}
\rm 
The functional $\overline{\rm M}(v, \flux; \beta)$ generates the upper bound of the error for any auxiliary 
$\flux \in H(\Omega, {\rm div})$ and $\beta >0$, therefore the choice of  $\flux$ might vary. The first 
and most straightforward way to select this variable is to set $\flux = \mathcal{G}(A \, \nabla v)$, where 
$\mathcal{G}: [L^{2}(\Omega)]^d \rightarrow H(\Omega, {\rm div})$ is a certain gradient-averaging  
operator. For this case, using the IgA framework becomes quite advantageous since for splines of the 
degree $p \geq 2$ an obtained $v$ is a $C^1$-continuous function and $\nabla v$ is already in 
$H(\Omega, {\rm div})$. Therefore, no additional post-processing is needed. On the other hand, due 
to the quadratic structure of the majorant, it is rather obvious that the optimal error estimate value is 
achieved at $y = A \nabla u$, i.e., 
\begin{equation}
\| \nabla e\|^2_{A, \Omega} \leq  \overline{\mathrm M}(v, \nabla u) 
= (1 + \beta) \, \| \nabla e\|^2_{A, \Omega}
+ {\CFriedrichs^2 \, (1 + \tfrac{1}{\beta})\, \| f + { {\rm div} ( \nabla u)} \|^2_\Omega}
= (1 + \beta) \, \| \nabla e\|^2_{A, \Omega}.
\label{eq:optimal-majorant}
\end{equation}
From \eqref{eq:optimal-majorant}, it is easy to see that if the auxiliary $\flux$ is chosen optimally and 
$\beta$ is set to zero (in the RHS of \eqref{eq:optimal-majorant}), there is no gap between 
$\overline{\rm M}$ and $\| \nabla e\|^2_{A, \Omega}$. 

One of the numerical methods providing an efficient reconstruction of both dual and primal variables is a 
mixed method. It generates an efficient approximation of the pair 
$(v, \flux) \in W := \HD{1}{0}(\Omega) \times H(\Omega, {\rm div})$ that can be straightforwardly 
substituted into the majorant $\overline{\rm M}(v, \flux)$. Moreover, if the error is measured in terms 
of the combined norm, i.e., including the norm of the error in primal and in dual variables
$$
|\!|\!| (u, \boldsymbol{p}) - (v, \flux) |\!|\!|_{W} := (\| \nabla (u - v) \|^2_\Omega 
+ \| \boldsymbol{p} - \flux \|^2_\Omega 
+ \| \dvrg (\boldsymbol{p} - \flux) \|^2_\Omega)^{\rfrac{1}{2}},$$
it is controlled by the residuals of the majorant as follows:
$$
\tfrac{1}{\sqrt{3}} \, (\|\,\R_{\rm d} \,\|_{A^{- 1}, \Omega} 
+ \,\tfrac{1}{\,\sqrt{\underline{\nu}_A}} \, \| \, \R_{\rm eq} \, \|_{\Omega}) 
\leq |\!|\!| (u, \boldsymbol{p}) - (v, \flux) |\!|\!|_{W} 
\leq \|\,\R_{\rm d} \,\|_{A^{- 1}, \Omega} 
+ \,(1 + 2 \, \tfrac{\CFriedrichs^2}{\,\underline{\nu}_A})^{\rfrac{1}{2}}
\| \, \R_{\rm eq} \, \|_{\Omega}.$$
We note that the ratio between the majorant
$\|\,\R_{\rm d} \,\|_{A^{- 1}, \Omega} 
+ \,\tfrac{1}{\,\sqrt{\underline{\nu}_A}} \| \, \R_{\rm eq} \, \|_{\Omega}$ (that does not include any 
constants) and the error $|\!|\!| (u, \boldsymbol{p}) - (v, \flux) |\!|\!|_{W}$ is controlled by $\sqrt{3}$, 
which proves the robustness of such error estimate. The series of works (see, e.g., 
\cite{RepinSauterSmolianski2004,RepinSmolianski2005,Repinetall2007}) has confirmed the efficiency 
of combination of mixed methods and functional error estimates.
%

The alternative approach providing an accurate $\flux$-reconstruction follows from the minimisation 
problem
\begin{equation}
\big \{ {\boldsymbol{y}}_{\rm min}, \beta_{\rm min} \big \}
:= {\rm arg} \inf \limits_{\beta > 0} \, \inf \limits_{{\boldsymbol y} \in H(\Omega, {\rm div})} 
\overline{\mathrm M}(v, {\boldsymbol y}; \beta).
\label{eq:minimation-problem-for-maj}
\end{equation}
%
The latter one is equivalent to the variational formulation for the optimal ${\boldsymbol y}_{\rm min}$, i.e., 
\begin{equation*}
	\tfrac{\CFriedrichs^2}{\beta_{\rm min} \,} \,
	(\dvrg \boldsymbol{{y}}_{\rm min}, \dvrg \vectorw)_{\Omega} 
         +  (A^{-1}\,\boldsymbol{{y}}_{\rm min}, \vectorw)_{\Omega}
	= - \tfrac{\CFriedrichs^2}{\beta_{\rm min} \,} 
	   \big(f, \dvrg \vectorw)_{\Omega}
	   + (A \nabla v, \vectorw)_{\Omega}, 
	   \quad \forall \vectorw \in H(\Omega, {\rm div}),
\label{eq:majorant-space-time-variational-formulation}
\end{equation*}
where the optimal $\beta$ is given by $\beta_{\rm min} := \tfrac{\CFriedrichs \, \mfInosq}{ \mdInosq}$ 
with
\begin{equation}
\mfInosq := \big\| \, \R_{\rm eq} \, \big\|_{\Omega} \quad \mbox{and} \quad 
\mdInosq := \|\,\R_{\rm d} \,\|_{A^{- 1}, \Omega}.
\label{eq:md-meq-def}
\end{equation}
\end{remark}
\vskip 10pt

In this work, using the IgA approximation schemes' setting, we apply the second method of the efficient 
$\flux$-reconstruction described in detail in Section \ref{sec:funct-estimates-and-iga}. To compare 
the performance of $\overline{\rm M}$ with alternative error estimates we use the standard residual 
error estimator (applied, e.g., in \cite{Giannellietall2016})
\begin{equation}
{\overline{\eta}}^2 = \sum\limits_{K \in \mathcal{K}_h} \,{\overline{\eta}}^2_K, \quad 
{\overline{\eta}}^2_K := h^2_K \, \| f + {\rm div} (A \nabla u_h) \|^2_{\L{2}(K)},
\label{eq:error-indicator}
\end{equation}
where $h_K$ denotes the diameter of cell $K$, and $u_h$ stands for approximation reconstructed 
by the IgA scheme. The term measuring the jumps across the element edges, which is usually included into 
residual error estimates, vanishes in \eqref{eq:error-indicator}
due to the properties of $u_h$ produced by the IgA schemes. It is provided in the G+Smo package and can 
be accessed by using the class available from the G+Smo library \cite[
{\tt\!stable/src/gsErrEstPoissonResidual.h}]{gismoweb}. 

\section{IgA overview: B-splines, NURBS, and THB-splines}

For the consistency of exposition, we first give an overview of the general IgA framework, the definitions of 
B-splines, NURBS, and THB-splines, their use in the geometrical representation of the computational domain 
$\Omega$ and in the construction of IgA discretisation spaces. 

Let $p \geq 2$ denote the degree of polynomials used for the IgA approximations, 
and let $n$ be the number of basis functions used to construct 
a $B$-spline curve. The {\em knot-vector} in $\mathds{R}$ is a non-decreasing set of coordinates in the 
parameter domain, written as $\Xi = \{ \xi_1, ..., \xi_{n+p+1}\}$, $\xi_i \in \mathds{R}$, where 
$\xi_1 = 0$ and $\xi_{n+p+1} = 1$. The knots can be repeated, and the multiplicity of the $i$-th knot is 
indicated by $m_i$. Throughout the paper, we consider only open knot vectors, i.e., 
$m_1 = m_{n+p+1} = p+1$. 
For the one-dimensional parametric domain $\widehat{\Omega} := (0, 1)$,
${\mathcal{\Khat}}_h := \{ \Khat \}$ denotes a locally quasi-uniform mesh, 
where each element $\Khat \in \mathcal{\Khat}_h$ is constructed by distinct 
neighbouring knots. The global size of $ \mathcal{\Khat}_h$ is denoted by 
$$\hhat := \max_{ \Khat \in \mathcal{\Khat}_h} \{ \hhat_{\Khat}\}, 
\quad \mbox{where} \quad
\hhat_{\Khat} := {\rm diam} (\Khat).$$
Henceforth, we assume locally quasi-uniform meshes, i.e., the ratio of two neighbouring elements $\Khat_i$ 
and $\Khat_j$ satisfies the inequality 
$$c_1 \leq \tfrac{\hhat_{\Khat_i}}{\hhat_{\Khat_j}} \leq c_2,\quad \mbox{where} \quad c_1, c_2 > 0.$$

The \emph{univariate B-spline basis functions $\Bhat_{i, p}: \widehat{\Omega} \rightarrow \mathds{R}$} 
are defined by means of the Cox-de Boor recursion formula
\begin{alignat}{2}
\Bhat_{i, p} (\xi) := \tfrac{\xi - \xi_i}{\xi_{i+p} - \xi_i} \, \Bhat_{i, p-1} (\xi)
                         + \tfrac{\xi_{i+p+1} - \xi}{\xi_{i+p+1} - \xi_{i+1}} \Bhat_{i+1, p-1} (\xi), 
                         \quad 
\Bhat_{i, 0} (\xi) \, := 
\begin{cases} 
1 & \mbox{if} \quad \xi_i \leq \xi < \xi_{i+1}  \\
0 & \mbox{otherwise}
\end{cases}
,
%
\end{alignat}
where a division by zero is defined to be zero. The B-splines are $(p-m_i)$-times continuously 
differentiable across the $i$-th knot with multiplicity $m_i$. Hence, if $m_i = 1$ for inner knots, the B-splines 
of the degree e.o.c. are $C^{p-1}$ continuous across them. 

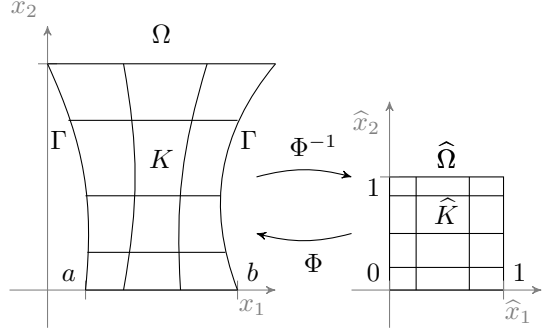
\begin{figure}
\centering
\begin{tikzpicture}[scale=0.5]
\def \ticksize {0.25};
\def \T {6.0};

\node[above left] at (1.0, 0.0) {$a$ };
\node[above right] at (5.0, 0.0) {$b$ };

\coordinate [label={above left:$$}] (a) at (1.0, 0.0);
\coordinate [label={above right:$$}] (b) at (5.0, 0.0);
\coordinate [label={above right:$$}] (aT) at (0.0, \T);
\coordinate [label={above right:$$}] (bT) at (6.0, \T);

\node[left] at (0.0, \T) {$$ };
\coordinate [label={left:$$}] (T) at (0.0, \T);
\node[above left] at (0.0, 0.0) {$$ };
\coordinate [label={above left:$$}] (O) at (0.0, 0.0);
\node[above] at (3.0, \T + 0.3) {$\Omega$ };
\node[above] at (3.0, \T/2) {${K}$ };
\node[above] at (0.3, 3.5) {$\Gamma$ };
\coordinate [label={above:$\Gamma$}] (sigma) at (5.3, 3.5);

\coordinate [label={above left:$$}] (xi10) at (2.0, 0.0);
\coordinate [label={above left:$$}] (xi20) at (3.5, 0.0);

\coordinate [label={above left:$$}] (xi1T) at (2.0, \T);
\coordinate [label={above left:$$}] (xi2T) at (4.2, \T);

\coordinate [label={above left:$$}] (yi10) at (1.05, 1.0);
\coordinate [label={above left:$$}] (yi20) at (4.68, 1.0);

\coordinate [label={above left:$$}] (yi12) at (1.0, 2.5);
\coordinate [label={above left:$$}] (yi22) at (4.55, 2.5);

\coordinate [label={above left:$$}] (yi13) at (0.55, 4.5);
\coordinate [label={above left:$$}] (yi23) at (5.0, 4.5);

\draw[->,thin,gray] (-1,0) --++(7,0)node[below left]{$x_1$};
\draw[->,thin,gray] (0,-1) --++(0,8)node[above left]{$x_2$};

\draw[thin, gray] (1, \ticksize) -- (1, -\ticksize);
\draw[thin, gray] (5.0, \ticksize) -- (5.0, -\ticksize);
\draw[thin, gray] (0.0-\ticksize, 6.0) -- (0.0 +\ticksize, 6.0);

\draw[black] (a) -- (b);
\draw[black] (aT) -- (bT);
\draw[black] (a) to [bend left=-15] (aT);
\draw[black] (b) to [bend left=30] (bT);

\draw[black] (xi10) to [bend left=-10] (xi1T);
\draw[black] (xi20) to [bend left=10] (xi2T);

\draw[black] (yi10) -- (yi20);
\draw[black] (yi12) -- (yi22);
\draw[black] (yi13) -- (yi23);

\draw[->,thin,gray] (8,0) --(13,0)node[below left]{$\widehat{x}_1$};
\draw[->,thin,gray] (9,-1) --(9,5)node[below left]{$\widehat{x}_2$};

\draw[black] (9, 0) -- (12, 0) -- (12, 3) -- (9, 3) -- (9, 0);
\node[above] at (10.5, \T/2) {$\widehat{\Omega}$ };
\node[above] at (10.5, \T/4) {$\widehat{{K}}$ };

\coordinate [label={above:$\widehat{\Omega}$}] (Q) at (10.5, \T/2);

\coordinate [label={above left:$$}] (xi10hat) at (9.7, 0.0);
\coordinate [label={above left:$$}] (xi20hat) at (9.7, 3.0);

\coordinate [label={above left:$$}] (xi1That) at (11.1, 0.0);
\coordinate [label={above left:$$}] (xi2That) at (11.1, 3.0);

\coordinate [label={above left:$$}] (yi10hat) at (9.00, 0.6);
\coordinate [label={above left:$$}] (yi20hat) at (12.00, 0.6);

\coordinate [label={above left:$$}] (yi12hat) at (9.00, 1.5);
\coordinate [label={above left:$$}] (yi22hat) at (12.00, 1.5);

\coordinate [label={above left:$$}] (yi13hat) at (9.00, 2.5);
\coordinate [label={above left:$$}] (yi23hat) at (12.00, 2.5);

\draw[black] (xi10hat) -- (xi20hat);
\draw[black] (xi1That) -- (xi2That);

\draw[black] (yi10hat) -- (yi20hat);
\draw[black] (yi12hat) -- (yi22hat);
\draw[black] (yi13hat) -- (yi23hat);

\draw[black, ->] (5.5, \T/2) to [bend left=15] (8.0, \T/2);
\draw[black, <-] (5.5, \T/4) to [bend left=-15] (8.0, \T/4);

\node[above] at (7, \T/2 + 0.3) {$\Phi^{-1}$ };
\node[below] at (7, \T/4 - 0.3) {$\Phi$ };


\node[above left] at (9.0, 0.0) {$0$ };
\node[above right] at (12.0, 0.0) {$1$ };
\node[below left] at (9.0, 3.2) {$1$ };

\draw[thin, gray] (9, \ticksize) -- (9, -\ticksize);
\draw[thin, gray] (12, \ticksize) -- (12, -\ticksize);
\draw[thin, gray] (-\ticksize + 9.0, 3.0) -- (9.0+\ticksize, 3.0);
\end{tikzpicture}
\caption{Mapping of $\widehat{\Omega}$ to $\Omega$. }
\end{figure}
The {\emph{multivariate B-splines}} on the parameter domain $\Omegahat := (0, 1)^{d}$, $d = \{1, 2, 3\}$, 
are defined as tensor products of the corresponding {univariate} ones. 
In the multidimensional case, we define a knot-vector dependent on the coordinate direction 
$\Xi^\alpha = \{ \xi^\alpha_1, ..., \xi^\alpha_{n^\alpha+p^\alpha+1}\}$, 
$\xi^\alpha_i \in \mathds{R}$, where $\alpha = 1, ..., d$ indicates the direction (in 
space or time). 
Furthermore, we introduce a set of multi-indices
$${\mathcal{I}} = \big\{\, \boldsymbol{i} = (i_1,  ..., i_{d}): i_\alpha = 1, ..., n_\alpha, \quad  
\,\alpha = 1, ..., d \big\}$$ 
and a multi-index $\boldsymbol{p} := (p_1, ..., p_{d})$ indicating the order of polynomials. 
The tensor-product of univariate B-spline basis functions generates multivariate B-spline basis 
functions 
\begin{equation}
\Bhat_{\boldsymbol{i}, \boldsymbol{p}} ({\boldsymbol \xi}) 
:= \prod_{\alpha = 1}^{d} \Bhat_{i_\alpha, p_\alpha} (\xi^\alpha), \quad  
\mbox{where} \quad {\boldsymbol \xi} = (\xi^1, ..., \xi^{d}) \in \Omegahat.
\label{eq:mutlivariable}
\end{equation}
%
The \emph{univariate and multivariate NURBS basis} functions are defined in a parametric domain by 
means of B-spline basis functions, 
i.e., for a given $\boldsymbol{p}$ and any $\boldsymbol{i} \in {\mathcal{I}}$, NURBS basis functions are 
defined as $\Rhat_{\boldsymbol{i}, \boldsymbol{p}}: \Omegahat \rightarrow \mathds{R}$ 
{
\begin{equation}
\Rhat_{\boldsymbol{i}, \boldsymbol{p}} ({\boldsymbol \xi}) 
:= \tfrac{w_{\boldsymbol{i}} \, \Bhat_{\boldsymbol{i}, \boldsymbol{p}} ({\boldsymbol \xi})}{\sum_{\boldsymbol{i} \in {\mathcal{I}}} w_{\boldsymbol{i}} \, \Bhat_{\boldsymbol{i}, \boldsymbol{p}} ({\boldsymbol \xi})},
\end{equation}
}
where $w_i \in \mathds{R}^+$.
%
To recall basic definitions related to THB-splines, we follow the structure outlined in 
\cite{Giannellietall2016} and consider a finite sequence of nested $d$-variate tensor-product spline spaces 
$\widehat{V}^0 \subset \widehat{V}^1 \subset ... \subset \widehat{V}^N$ defined on the axis aligned box-domain 
$\widehat{\Omega}^0 \subset \mathds{R}^d$. To each space $V^\ell$ we assign a  
tensor-product B-spline basis of degree ${\boldsymbol p}$

$$\big\{ \Bhat^\ell_{{\boldsymbol i}, {\boldsymbol p}} \}_{\boldsymbol{i} \in \mathcal{I}^\ell}, \quad 
\mathcal{I}^\ell := \{ {\boldsymbol i} = (i_1, ..., i_d), 
\quad  i_k = 1, ..., n_k^\ell \quad \mbox{for} \quad k = 1, ..., d \, \big\},$$ 
%
where $\mathcal{I}^\ell$ is a set of multi-indices for each level, and $n_k^\ell$ denotes the number of 
univariate B-spline basis functions in the $k$-th coordinate direction. After assuming  
that $\mathcal{I}^\ell$ has a fixed ordering and rewriting the basis as \linebreak
${\boldsymbol{\rm \Bhat}}^\ell({\boldsymbol \xi}) 
= (\Bhat^\ell_{{\boldsymbol i}, {\boldsymbol p}}({\boldsymbol \xi}))_{\boldsymbol{i} \in \mathcal{I}^\ell},$
it can be considered as a column-vector of basis functions. Then, a spline function 
$s: \widehat{\Omega}^0 \rightarrow \mathds{R}^m$ is defined by $\boldsymbol{\rm \Bhat}^\ell({\boldsymbol \xi})$
and a coefficient matrix $C^\ell$, i.e.,  
$$s({\boldsymbol \xi}) 
= \sum_{\boldsymbol{i} \in \mathcal{I}^\ell } \Bhat^\ell_{\boldsymbol{i}, \boldsymbol{p}} ({\boldsymbol \xi}) c^\ell_i 
= \boldsymbol{\rm \Bhat}^\ell({\boldsymbol \xi})^{\rm T}\,  C^\ell ,$$
where $c^\ell_i \in \mathds{R}^m$ are row-coefficients of $C^\ell$.

Since $\widehat{V}^\ell \subset \widehat{V}^{\ell+1}$, the basis $\boldsymbol{\rm \Bhat}^\ell$ can be 
represented by the linear combination of $\boldsymbol{\rm \Bhat}^{\ell+1}$, namely, 
$$s({\boldsymbol \xi}) 
= \boldsymbol{ \rm \Bhat}^\ell({\boldsymbol \xi})^{\rm T}\, C^\ell 
= \boldsymbol{ \rm \Bhat}^{\ell+1}({\boldsymbol \xi})^{\rm T}\, R^{\ell+1} \, C^\ell,$$ 
where $R^{\ell+1}$ is a refinement matrix. Its entries can be obtained from B-splines refinement rules 
(see \cite{PieglTiller1997}). Along with nested space, a corresponding sequence of nested domains is 
considered 
\begin{equation}
\widehat{\Omega}^0 \supseteq \widehat{\Omega}^1 \supseteq ... \supseteq \widehat{\Omega}^N, 
\label{eq:omega-hierarchy}
\end{equation}
where each $\widehat{\Omega}^\ell \in \mathds{R}^d$ is covered by a collection of cells with respect to 
the tensor-product grid of level $l$. In this work,  we focus on dyadic cell refinement for the bi- and 
trivariate cases with uniform degrees $p_{\alpha} = p$ for all levels and coordinate directions, therefore,
$\boldsymbol{p} = p$ in further exposition. 

%

Let the {\em characteristic matrix} $X^\ell$ of $\boldsymbol{\Bhat}^\ell(\boldsymbol{\xi})$ w.r.t. domains
$\Omega^\ell$ and $\Omega^{\ell+1}$ is defined as
$$X^\ell := {\rm diag} (x^\ell_{\boldsymbol{i}})_{\boldsymbol{i} \in \mathcal{I}^\ell}, \quad 
x^\ell_{\boldsymbol{i}} := \begin{cases}
1, \quad \mbox{if} \quad 
{\rm supp} \Bhat^\ell_{\boldsymbol{i}, \boldsymbol{p}} \subseteq \Omega^\ell \wedge {\rm supp} \Bhat^\ell_{\boldsymbol{i}, \boldsymbol{p}} \nsubseteq \Omega^{\ell+1} \\
0, \quad \mbox{otherwise}.
\end{cases}$$
%
Next, for each level $\ell$, the set of the indices of {\em active functions} can be defined with 
$\mathcal{I}^\ell_* := \{\mathcal{I}^\ell : x^\ell_{\boldsymbol{i}} = 1\}.$ To store the indices of all 
active functions at all hierarchical levels, we define an {\em index set}
$$\boldsymbol{\mathcal{I}} := \{ (\ell, \boldsymbol{i}) : \ell \in \{0, ..., N\}, \boldsymbol{i} \in \mathcal{I}^\ell_*\}.$$
{
The initial hierarchical data structure defined by the tensor-product mesh $\widehat{\Omega}^0$ (see 
Figure \ref{fig:thb-domains}). In particular, we illustrate the knot lines of the spaces 
$V^0 \subset V^1 \subset V^2$ (where levels increase form the left to the right). 
By means of the insertion operation, new subdomains can be added to obtain new representations 
$\widehat{\Omega}^0 \supseteq \widehat{\Omega}^1 \supseteq \widehat{\Omega}^2$ .
The sets of active basis functions $\mathcal{I}^\ell_*$ as well as the characteristic matrices $X^\ell$ 
for all levels $\ell$ are extracted simultaneously with the new box insertion and initialisation of the basis.
Figure \label{fig:thb-domains-and-levels} illustrates meshes at the refinement levels $\ell = 0, 1, 2$. 

\begin{figure}[!h]
	\centering
	\subfloat[Knot lines and inserted subdomain (shaded)]{
	\includegraphics[scale=0.3]{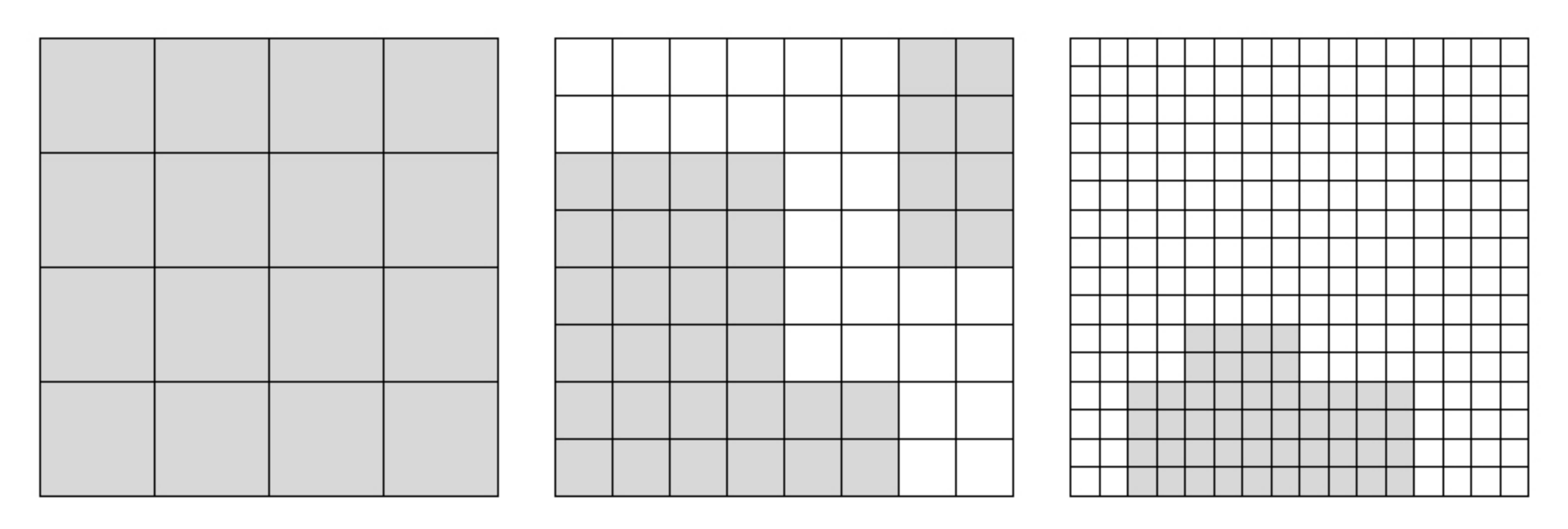}
	\label{fig:thb-domains}
	} \\
	\subfloat[{Hierarchical meshes for refinement levels}]{
	\includegraphics[scale=0.3]{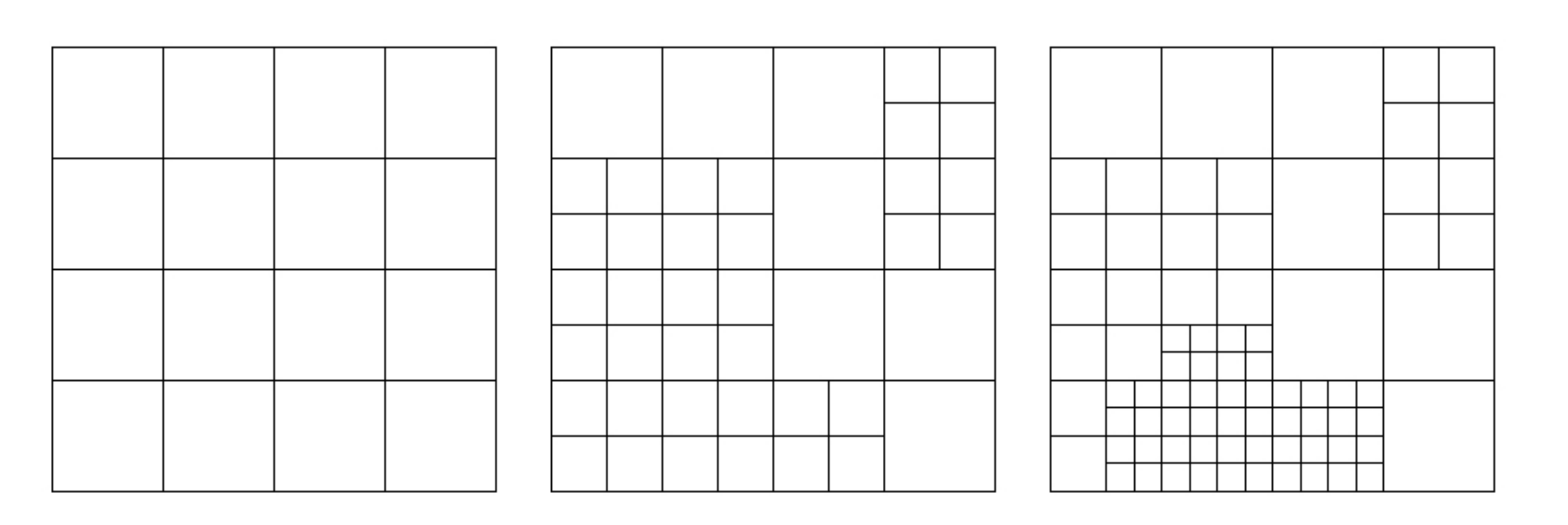}}
	\label{fig:thb-domains-and-levels}
	\caption{Bivariate configuration for hierarchical levels $\ell = 0, 1, 2$ 
	(the picture is taken from \cite{Giannellietall2016}).}
\end{figure}
}

The THB-spline basis related to the hierarchical domains is defined as
$$\widehat{\boldsymbol{{\rm T}}} (\boldsymbol{\xi}) 
= (\mathcal{K}^\ell_{\boldsymbol{i}}(\boldsymbol{\xi}))_{(l, \boldsymbol{i}) \in \boldsymbol{\mathcal{I}}}, \quad 
\mathcal{K}^\ell_{\boldsymbol{i}}(\boldsymbol{\xi}) 
= {\rm trunc}^N({\rm trunc}^{N-1}(... {\rm trunc}^{\ell+1}(\Bhat^\ell_{\boldsymbol{i}, \boldsymbol{p}}(\boldsymbol{\xi})))),$$
where the {\em truncation} of any function $s({\boldsymbol \xi}) \in \widehat{V}^\ell$ w.r.t. level $\ell+1$ 
is defined by $${\rm trunc}^{\ell+1} (s({\boldsymbol \xi})) 
= \boldsymbol{ \rm \Bhat}^{\ell+1}({\boldsymbol \xi})^{\rm T}\, (I^{\ell+1} - X^{\ell+1}) \, R^{\ell+1}\, C^\ell.$$
Here, $I^{\ell+1}$ denotes an identity matrix $I^{\ell+1}$ of size $|I^{\ell+1}| \times |I^{\ell+1}|,$
the multiplication of $R^{\ell+1}$ by $C^\ell$ represents $s({\boldsymbol \xi})$ w.r.t. to the level $\ell+1$, 
and additional multiplication by $(I^{\ell+1} - X^{\ell+1})$ performs the truncation operation. For the 
detailed discussion of truncation operation, we refer the reader to \cite{LMR:GiannelliJuttlerSpeleers2012, 
LMR:GiannelliJuttlerSpeleers2014, Giannellietall2016}.
{
We illustrate an effect of the truncation for the case of univariate quadratic spline basis functions 
(see Figure \ref{fig:hb-thb-splines}). Figure \ref{fig:1} presents the HB-splines for the hierarchical levels
$\ell = 2, 3, 4$, whereas Figure \ref{fig:2} shows the same levels for THB-splines. On the last raw, 
basis functions influenced by the truncation on coarser levels are exposed, i.e., 
THB-splines before and after being truncated are denoted by the grey and black marker, respectively.

\begin{figure}[!h]
	\centering
	\subfloat[HB-splines on the levels 2, 3, 4]{
	\includegraphics[scale=0.045]{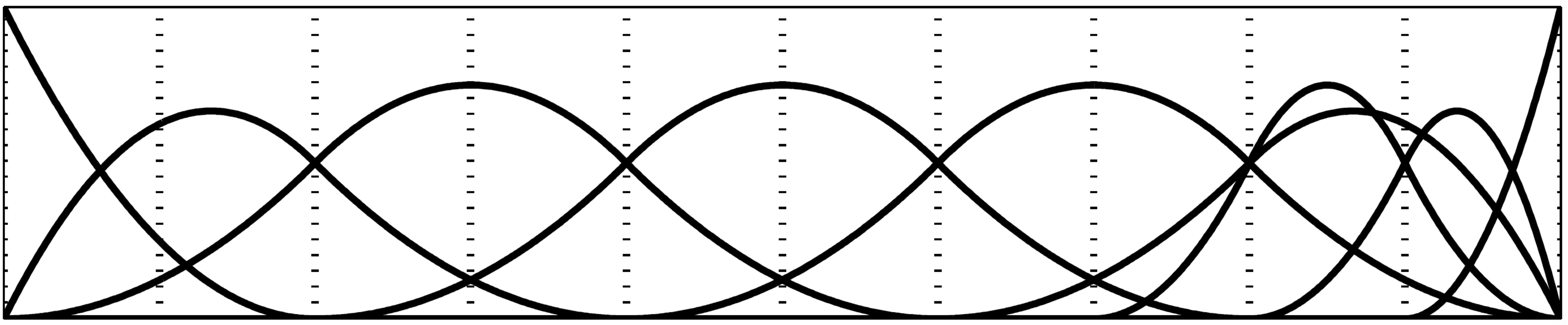}
	\quad
	\includegraphics[scale=0.045]{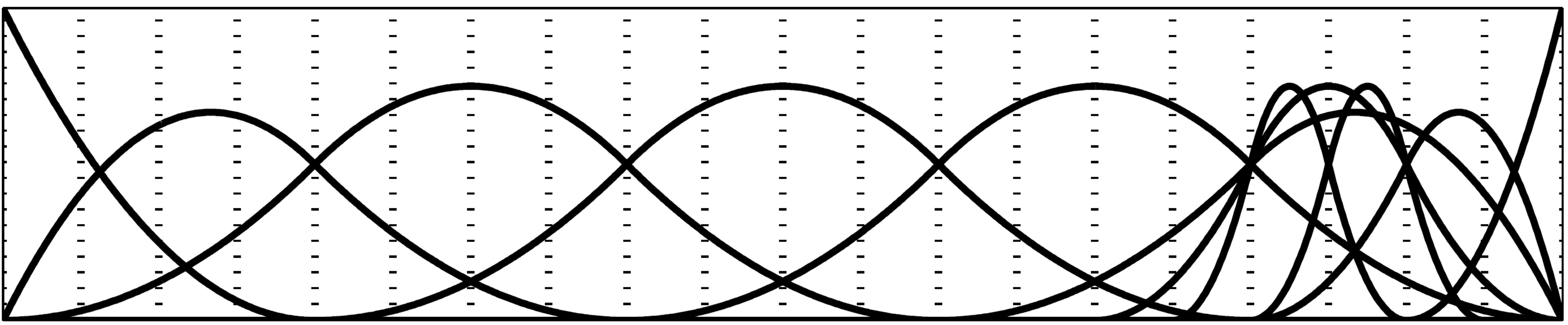}
	\quad 
	\includegraphics[scale=0.045]{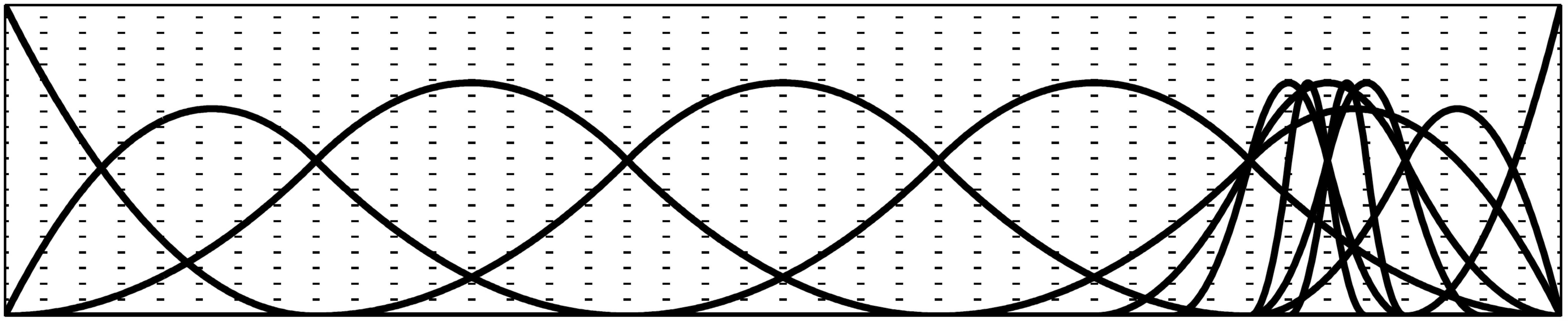}
	\label{fig:1}
	} \quad
	\subfloat[THB-splines on the levels 2, 3, 4]{
	\includegraphics[scale=0.045]{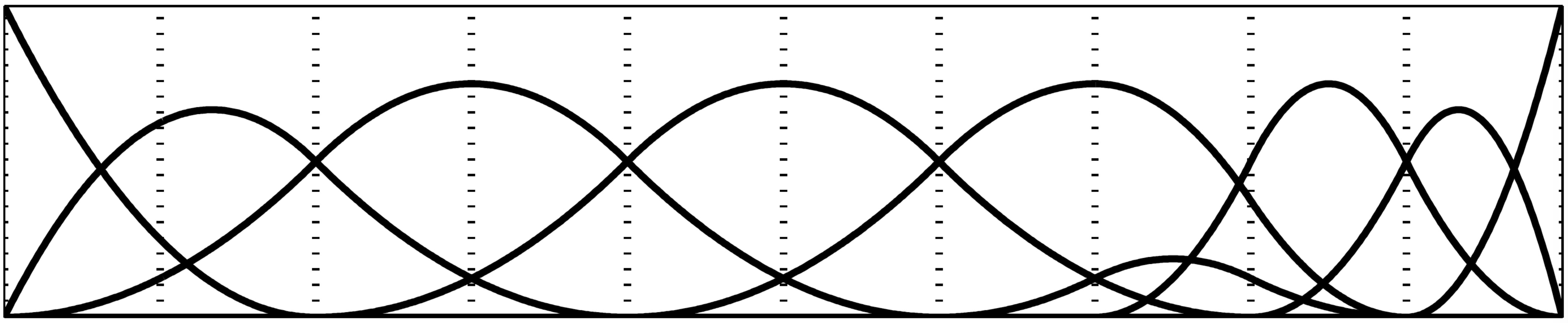}
	\quad
	\includegraphics[scale=0.045]{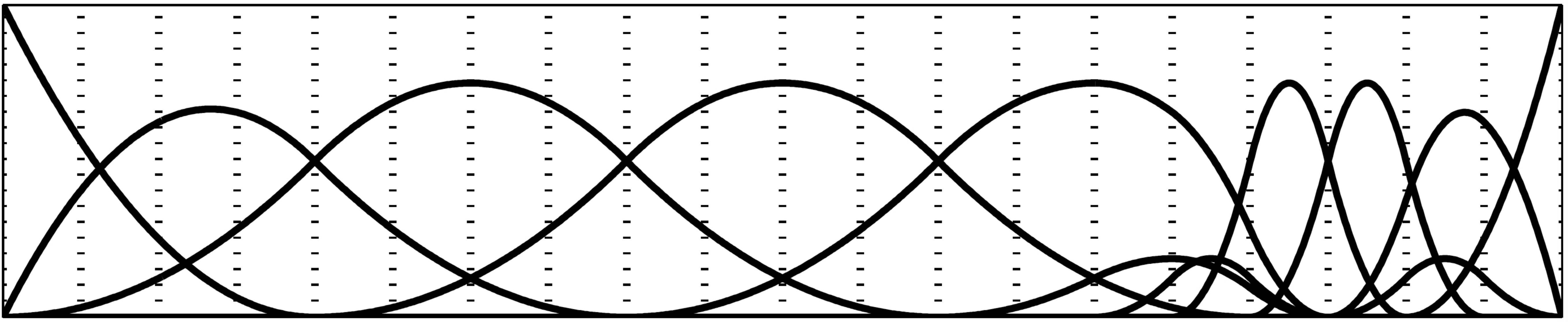}
	\quad
	\includegraphics[scale=0.045]{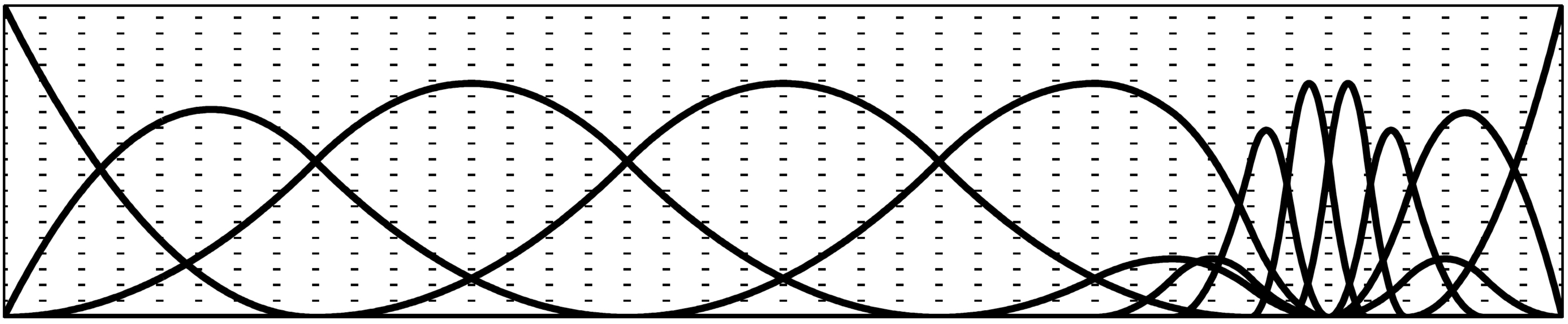}
	\label{fig:2}
	} \quad
	\subfloat[Influenced by truncation THB splines on the levels 1, 2, 3]{
	\includegraphics[scale=0.045]{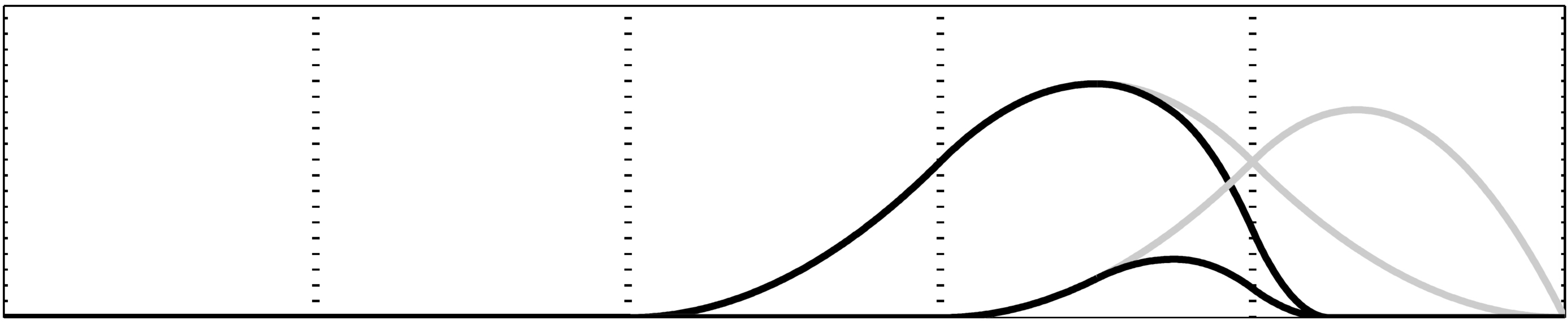}
	\quad
	\includegraphics[scale=0.045]{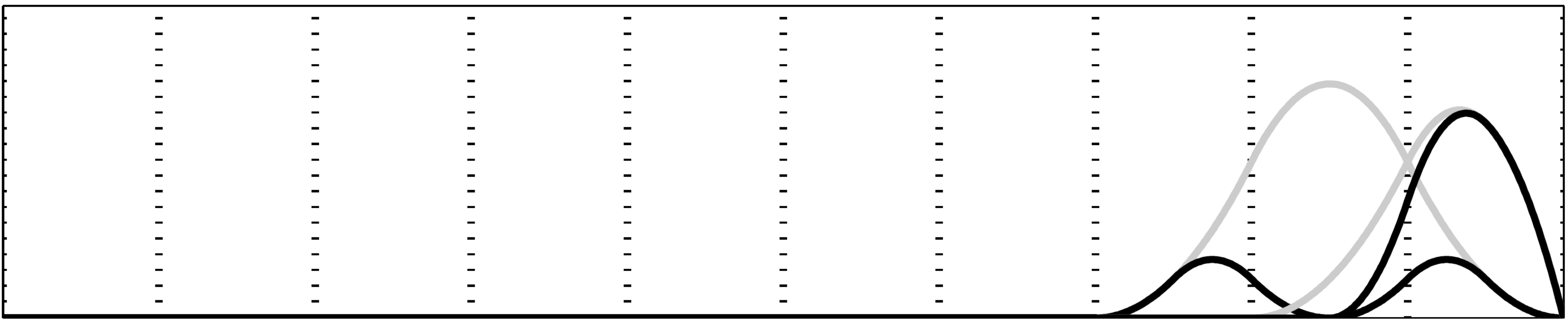}
	\quad
	\includegraphics[scale=0.045]{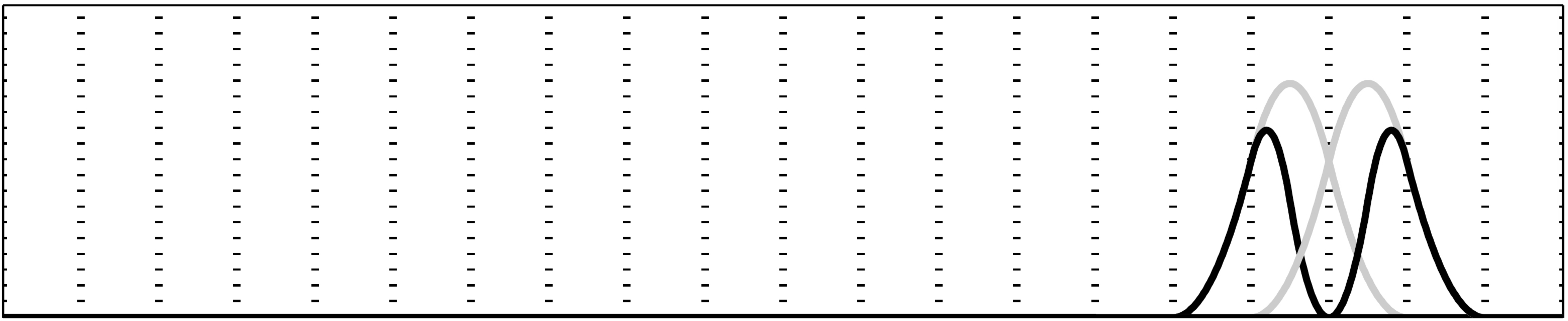}
	} 
	\caption{Comparison of HB- and THB-splines on the different refinement levels 
	(the picture is taken from \cite{Giannellietall2016}).}
	\label{fig:hb-thb-splines}
\end{figure}
}

The physical domain $\Omega \subset \mathds{R}^{d}$ is defined by the geometrical mapping  
of the parametric domain $\widehat{\Omega} := (0, 1)^{d}$:
\begin{equation}
\Phi: \widehat{\Omega} \rightarrow \Omega := \Phi(\widehat{\Omega}) \subset \mathds{R}^{d}, 
\quad 
\Phi({\boldsymbol \xi}) := \sum_{\boldsymbol{i} \in {\mathcal{I}}} {\widehat{\mathcal{B}}}_{\boldsymbol{i}, \boldsymbol{p}}({\boldsymbol \xi}) \, {c}_i,
\label{eq:geom-mapping}
\end{equation}
where ${c}_i \in \mathds{R}^{d}$ are control points, and 
${\widehat{\mathcal{B}}}_{\boldsymbol{i}, \boldsymbol{p}}$ stands for
either B-splines, NURBS, or THB-basis functions . 
The mesh $\mathcal{K}_h$ discretising $\Omega$ consists of elements $K \in \mathcal{K}_h$ that are 
the images of $\Khat \in \mathcal{\Khat}_h$, 
i.e., 
$$\mathcal{K}_h := \big\{K = \Phi(\Khat) : \Khat \in \mathcal{\Khat}_h \big\}.$$
The global mesh-size is denoted by 
\begin{equation}
h := \max\limits_{K \in \mathcal{K}_h} \{ \, h_{K}\,\}, \quad  
h_{K} := \| \nabla \Phi \|_{\L{\infty} (K)} \, \hhat_{\Khat}.
\label{eq:global-mesh-size}
\end{equation}
Moreover, we assume that $\mathcal{K}_h$ is a quasi-uniform mesh, i.e., there exists a positive constant 
$C_u$ independent of $h$, such that $h_{K} \leq h \leq C_u \, h_{K}.$

\section{Functional error estimates within the IgA framework}
\label{sec:funct-estimates-and-iga}

In this section, we present the algorithms used for general reliable computations and functional-type error 
estimates reconstruction. Then we proceed with commenting on the implementation of these error 
estimates in G+Smo and their integration into the library's structure. Finally, we present a series of 
examples demonstrating numerical properties of derived error majorants. 

\subsection{Reliable reconstruction of IgA approximations. Algorithms}

In order to keep the presentation concise, we restrict (\ref{eq:elliptic-equation})--(\ref{eq:dirichlet-bc}) 
to the Dirichlet--Poisson problem
\begin{equation}
- \laplace u = f \;\;\; {\rm in} \;\;\; \Omega := (0, 1)^d \in \mathds{R}^d, \; d = \{2, 3\}, \
\qquad u = 0 \;\;\; {\rm on} \;\;\; {\Gamma} = \partial \Omega. 
\label{eq:poisson}
\end{equation}
Let the approximation 
$$u_h \in V_{0h} := V_h \cap H^1_{0}(\Omega),  \quad \mbox{where} \quad
V_{h} \equiv {\mathcal{S}}^{p, p}_h :=\big\{ \phi^{(p)}_{h,i} := \widehat{V}_h \circ \Phi^{-1} \big\}.$$ 
Here, 
$\widehat{V}_h \equiv \widehat{\mathcal{S}}^{p, p}_h$ is generated with NURBS of degree $p$, i.e., 
$\widehat{V}_h:= {\rm span} \,\big\{ {\widehat{\mathcal{B}}}_{\boldsymbol{i}, \boldsymbol{p}} \big\}_{\boldsymbol{i} \in \mathcal{I}}$. 
Due to the one-patch setting and restriction on the knots' multiplicity of $\widehat{\mathcal{S}}^{p, p}_h$, the 
smoothness $u_h \in C^{p-1}$ is automatically provided. Since no numerical algorithms specific to the 
hierarchical levels of the localised splines will be discussed below, we use the same notation for spaces 
generated by THB-splines.
%
Therefore, the constructed approximation can be written as
$$u_h(x) = u_h(x_1, . . . , x_d) := \sum_{i \in \mathcal{I}} \underline{ \rm u}_{h, i} \, \phi^{(p)}_{h, i}(x),$$ 
where 
$\underline{\rm u}_h 
:= \big[ \underline{ \rm u}_{h,i}\big]_{i \in \mathcal{I}} \in {\mathds{R}}^{|\mathcal{I}|}$ 
is a vector of degrees of freedom (d.o.f.) defined by the system 

\begin{equation}
{\rm K}^{(p)}_h \, \underline{ \rm  u}_h = {\rm f}^{(p)}_h, \quad 
{\rm K}^{(p)}_h := \big[(\nabla \phi^{(p)}_{h, i}, \nabla \phi^{(p)}_{h, j})_\Omega \big]_{i, j \in \mathcal{I}}, \quad
{\rm f}^{(p)}_h :=\big[(f, \phi^{(p)}_{h, i})_\Omega \big]_{i \in \mathcal{I}}. 
\label{eq:system-uh}
\end{equation}
%
The majorant corresponding to the problem \eqref{eq:poisson} reads as
\begin{equation}
\overline{\rm M}^2(u_h, \flux_h) 
:= (1 + \beta) \, \mdI
    + (1 + \tfrac{1}{\beta}) \, \CFriedrichs^2 \, \mfI 
 = (1 + \beta) \, \| \flux_h - \nabla u_h \|^2_\Omega
    + (1 + \tfrac{1}{\beta}) \, \CFriedrichs^2 \, \| \dvrg \flux_h + f \|^2_\Omega, 
\label{eq:majorant-elliptic} 
\end{equation}
where $\mfInosq$ and $\mdInosq$ are defined in \eqref{eq:md-meq-def}, 
$\beta > 0$ and ${\flux_h} \in Y_h \subset H(\Omega, {\rm div})$. Here, the approximation space for 
%
$$\flux_h \in Y_h 
\equiv \oplus^{d} \mathcal{S}^{q, q}_h 
\equiv {\mathcal{S}}^{q, q}_h \oplus ... \oplus {\mathcal{S}}^{q, q}_h
:= \big\{ \widehat{Y}_h \circ \Phi^{-1} \big\}$$ 
is generated by the {push-forward} of a corresponding space in the parametric domain
$$\widehat{Y}_h := \oplus^{d} \widehat{\mathcal{S}}^{q, q}_h 
\equiv \widehat{\mathcal{S}}^{q, q}_h \oplus ... \oplus \widehat{\mathcal{S}}^{q, q}_h.$$ 
%
Here, $\widehat{\mathcal{S}}^{q, q}_h$ is a space of NURBS with the degree $q$ for each of $d$ 
components of $\flux_h = (y_h^{(1)}, ..., y_h^{(d)})^{\rm T}.$ The details of the numerical 
reconstruction of \eqref{eq:majorant-elliptic} were thoroughly studied in \cite{LMR:KleissTomar2015}. 
The best estimate follows from the optimisation of $\overline{\rm M}(u_h, \flux_h)$ w.r.t. function 
$$\flux_h(x) := \sum_{i \in \mathcal{I}} \underline{ \rm \bf y}_{h,i} \, {{\boldsymbol \psi}}_{h,i}(x).$$ 
The basis functions ${{\boldsymbol \psi}}_{h,i}$ generate the space $Y_h$, whereas 
$\underline{ \rm \bf y}_{h} 
:= \big[ \underline{ \rm \bf y}_{h,i}\big]_{i \in \mathcal{I} \times d} \in {\mathds{R}}^{d|\mathcal{I}|}$ 
is a vector of d.o.f of $\flux_h$ defined by a system 
\begin{equation}
\left( {\CFriedrichs^2} \, {\rm Div}_h + {\beta} \, {\rm M}_h \right)\, 
\underline{ \rm \bf y}_{h} = 
- {\CFriedrichs^2} \, {\rm z}_h + {\beta} \, {\rm g}_h, 
\label{eq:system-fluxh}
\end{equation}
where 
\begin{alignat*}{2} 
	{{\rm {Div}}_h} & 
	:= 
	\big[(\dvrg {\boldsymbol \psi}_i, \dvrg {\boldsymbol \psi}_j)_\Omega
	\big]_{i, j=1}^{d|\mathcal{I}|} 
, \qquad 
	{\rm z}_{h} := 
	\big[
	\big(f,  \dvrg {\boldsymbol \psi}_j \big)_\Omega  
	\big]_{j=1}^{d|\mathcal{I}|}  
	,  
	\\
	{{\rm {M}}_h} & := 
	\big[ 
	({\boldsymbol \psi}_i, {\boldsymbol \psi}_j)_\Omega
	\big]_{i, j=1}^{d|\mathcal{I}|}
	, 
	\qquad \qquad \qquad \quad
	{\rm g}_h := 
	\big[ \big(\nabla v, {\boldsymbol \psi}_j \big)_\Omega 
	\big]_{j=1}^{d|\mathcal{I}|} 
	.
\end{alignat*}
%
According to the numerical results obtained in \cite{LMR:KleissTomar2015}, the most efficient
majorant reconstruction (with uniform refinement) is obtained when $q$ is set substantially higher 
than $p$. Let us assume that $q = p + m$, $m \in \mathds{N}^+$. At the same time, when $u_h$ is 
reconstructed on the mesh $\mathcal{K}_h$, we use a coarser one $\mathcal{K}_{Mh}$,  
$M \in \mathds{N}^+$ 
in order to recover $\flux_h$. For the reader's convenience, all used notation is summarised 
in Table \ref{tab:table-of-notation}. The initial mesh $\mathcal{K}^0_{h}$ and the basis functions defined on 
it are assumed to be given via the geometry representation of the computational domain. The exact 
representation of geometry on the initial (the coarsest) level is preserved in the process of mesh refinement.

For the reconstruction of $\underline{\rm M}(v, w)$, let the approximation 
$$w_h \in {W}_{0h} := W_h \cap H^1_{0}(\Omega),  \quad \mbox{where} \quad
W_{h} \equiv {\mathcal{S}}^{r, r}_h :=\big\{ \varphi^{r}_{h,i} := \widehat{W}_h \circ \Phi^{-1} \big\}.$$ 
Here, 
$\widehat{W}_h \equiv \widehat{\mathcal{S}}^{r, r}_h$ is approximation space generated with NURBS 
of degree $r$ on the parameter domain, i.e., $\widehat{W}_h 
:= {\rm span} \,\big\{ {\widehat{\mathcal{B}}}_{\boldsymbol{i}, \boldsymbol{r}} \big\}_{\boldsymbol{i} \in \mathcal{I}}$. 
Then, the auxiliary approximation can be written as
$$w_h(x) = w_h(x_1, . . . , x_d) := \sum_{i \in \mathcal{I}} \underline{ \rm w}_{h, i} \, \phi^{(r)}_{h, i},$$ 
where 
$\underline{\rm w}_h 
:= \big[ \underline{ \rm w}_{h,i}\big]_{i \in \mathcal{I}} \in {\mathds{R}}^{|\mathcal{I}|}$ 
is a vector of degrees of freedom (d.o.f.) defined by the system 

\begin{equation}
{\rm K}^{(r)}_h \, \underline{ \rm  u}_h = {\rm f}^{(r)}_h, \quad 
{\rm K}^{(r)}_h := \big[(\nabla \varphi^{(r)}_{h, i}, \nabla \varphi^{(r)}_{h, j})_\Omega \big]_{i, j \in \mathcal{I}},
\quad
{\rm f}^{(r)}_h :=\big[(f, \varphi^{(r)}_{h, i})_\Omega \big]_{i \in \mathcal{I}}. 
\label{eq:system-wh}
\end{equation}
%
Analogously to the selection of the $q$ for the space $Y_h$, we let $r = p + l$, $l \in \mathds{N}^+$. 
At the same time, we use a coarser mesh $\mathcal{K}_{Lh}$,  $L \in \mathds{N}^+$ 
for the $w_h$ approximation. 

\begin{table}[htbp]
\footnotesize
\begin{center}
\begin{tabular}{r c p{12cm} }
\toprule
$p$ & $$ & degree of the splines used for $u_h$ approximation \\[2pt]
$q$ & $$ & degree of the splines used for $\flux_h$ approximation \\[2pt]
$r$ & $$ & degree of the splines used for $w_h$ approximation \\[2pt]
$m$ & $ $ & $q - p$ \\[2pt]
$l$ & $ $ & $r - p$ \\[2pt]
$S_h^{p, p}$ ($S_h^{r, r}$) & $$ & approximation space for the scalar-functions generated by splines \\[2pt]
$\oplus^d S_h^{q, q}$ & $$ & approximation space for the $d$-dimensional vector-functions generated by splines \\[2pt]
$S_h^{q, q} \oplus S_h^{q, q}$ & $$ & approximation space for the two-dimensional vector-functions generated by splines \\[2pt]
$M$  & $ $ & coarsening ratio of the global size of the mesh for $u_h$ 
		 approximation to the global size of the mesh for $\flux_h$ reconstruction\\[2pt]
$L$  & $ $ & coarsening ratio of the global size of the mesh for $u_h$ 
		approximation to the global size of the mesh for $w_h$ reconstruction\\[2pt]
$\mathcal{K}_{h}$ ($\mathcal{K}^{u_h}_h$) & $ $ & mesh used for $u_h$ approximation \\[2pt]
$\mathcal{K}_{Mh}$  ($\mathcal{K}^{\flux_h}_h$, $M = 1$)& $ $ & mesh used for $\flux_h$ reconstruction\\[2pt] 
$\mathcal{K}_{Lh}$  ($\mathcal{K}^{w_h}_h$, $L = 1$)& $ $ & mesh used for $w_h$ reconstruction\\[2pt] 
$N_{\rm ref}$ & $$ & number of uniform or adaptive refinement steps \\[2pt]
$N_{\rm ref, 0}$ & $$ & number of initial refinement steps performed before testing \\[2pt]
${\mathds{M}}_{*}({\theta})$ & $$ & marking criterion $*$ with the parameter $\theta$\\[2pt]
\bottomrule
\end{tabular}
\end{center}
\vskip -10pt
\caption{\small  Table of notations.}
\label{tab:table-of-notation}
\end{table}

The classical strategy of the reliable $u_h$-approximation is summarised in Algorithm 
\ref{alg:reliable-uh-reconstruction}. Let us assume that the problem data such as $f$, $u_0$, and 
$\Omega$ of (\ref{eq:elliptic-equation})--(\ref{eq:dirichlet-bc}) are provided.
The Input of Algorithm \ref{alg:reliable-uh-reconstruction} is the initial mesh $\mathcal{K}_h$ 
(or the one obtained on the previous refinement step). It provides the refined version of $\mathcal{K}_h$ 
denoted by $\mathcal{K}_{h_{\rm ref}}$ as an output. The process of new mesh generation can be 
divided into classical block-chain, i.e.,
$${\rm APPROXIMATE} 
\rightarrow {\rm ESTIMATE} 
\rightarrow {\rm MARK} 
\rightarrow {\rm REFINE}.$$ 

On the {\rm APPROXIMATE} step, we construct the system that provides the d.o.f. of $u_h$, i.e., 
we assemble the matrix ${{\rm {K}}^{(p)}_h}$ and RHS ${{\rm {f}}^{(p)}_h}$ defined in \eqref{eq:system-uh}, and 
solve it with a direct sparse ${\rm LDL^{\rm T}}$ Cholesky factorisations for $d = 2$ and conjugate 
gradient (CG) method for $d = 3$. In the follow-up report, we will investigate how the selection of the  
initial guess enhances the performance of the iterative solver. In particular,  we use the work 
\cite{Deuflhard1994, BornemannDeuflhard1996} that studies the 
so-called cascadic preconditioned conjugate gradient (CPCG) method. The latter one has an improved 
speed of convergence due to initial guess chosen as an interpolation of the approximation obtained 
on the previous refinement (hierarchical) level. It appears that such a cascadic structure of the 
meshes by itself realises some kind of preconditioning. The time spent on assembling 
and solving sub-procedures for $u_h$ is tracked and saved in 
vectors ${t_{\rm as}(u_h)}$ and ${t_{\rm sol}(u_h)}$, respectively. This notation is used in the 
upcoming examples to analyse the efficiency of Algorithm \ref{alg:reliable-uh-reconstruction} and 
compare the computational costs for its blocks.

The next {\rm ESTIMATE} step is first and foremost responsible for the reconstruction of global estimate 
$\maj{} (u_h, \flux_h)$ as well as the element-wise error indicator distribution $\mdI{} (u_h, \flux_h)$ 
(see \eqref{eq:md-meq-def}) that follows. 
The time spent for this is measured by ${t_{\rm as}(\flux_h)} + {t_{\rm sol}(\flux_h)}$. 
Simultaneously with the upper bound, we reconstruct minorant $\mij{} (u_h, w_h)$, whereas the time 
spent for its reconstruction is tracked by ${t_{\rm as}(w_h)} + {t_{\rm sol}(w_h)}$. Their detailed 
description of latter estimates generation is presented in Algorithms \ref{alg:estimate-step-maj} and 
\ref{alg:estimate-step-min}.

In the chain-block {\rm MARK}, we apply a marking criterion denoted by ${\mathds{M}}_{*}({\theta})$. 
It provides an algorithm for defining the threshold  $\Theta_* $ for selecting those $K \in \mathcal{K}_h$ 
for further refinement that satisfies the criterion
$${\mdIK \geq \Theta_* ({\mathds{M}}_{*}({\theta})), \quad K \in \mathcal{K}_h}.$$
In the G+smo library \cite{gismoweb}, several marking strategies are considered. The first criterion 
defines an `absolute threshold', and it is denoted as {\rm \bf GARU} (an abbreviation for `greatest 
appearing residual utilisation'). The corresponding threshold reads as
$$\Theta_{\rm GARU} 
:= \theta \, \max\limits_{K \in \mathcal{K}_h} \{ \mdIK \}, \quad \theta \in (0, 1).$$
The percentage of marked elements (dictated by this criterion) varies at each refinement step since 
$\Theta_{\rm GARU}$ considers only the absolute value of the largest local error, without 
taking into account the element-wise distribution of the error.

The second marking criterion defining the {`relative threshold'} is denoted as 
${\mathds{M}}_{\rm \bf PUCA}$, where {\rm \bf PUCA} stands for `percent-utilising cutoff 
ascertainment'. The corresponding amount of elements selected for the refinement can be 
approximated as follows:
$$|\{K: \mdIK > \Theta_{\rm PUCA} \}_{K \in \mathcal{K}_h}| 
\approx (1 - \theta) \cdot | \{K\}_{K \in \mathcal{K}_h}|, 
\quad \theta \in (0, 1).$$
For instance, if we let $\theta = 0.7$,  $\Theta_{\rm PUCA}$ is chosen such that 
$\mdIK \geq \Theta_{\rm PUCA}$ holds for $30\%$ of elements. 

Last and most widely used criterion is called bulk marking (also known as the D\"orfler marking 
\cite{Dorfler1996})
and is denoted as ${\mathds{M}}_{{\rm \bf BULK}}({\theta})$. According to this marking strategy, 
we select the subset of elements from the collection ${\mathcal{K}}_h$ that has been sorted w.r.t. 
element-wise contributions $\mdIK$, i.e., 
$\mathcal{K}^\prime_h \xleftarrow[\text{K}]{} \,
\mathcal{K}^{\rm sort}_h : = {\rm sort}_{\mdIK} \{ {\mathcal{K}}_h \},$
until we satisfy 
$$\Sum_{K \in \mathcal{K}^\prime_h } \mdIK \geq \Theta_{\rm BULK} 
:= (1 - \theta) \, \Sum_{K \in \mathcal{K}_h } \mdIK, \quad  \theta \in (0, 1). $$
This way, we form a subset of elements which contains the highest indicated errors. The selection 
process stops when the error accumulated on previous steps exceeds the `bulk' level (threshold) defined 
by $\theta$.  In the case of uniform refinement, all elements of $\mathcal{K}_h$ are marked for 
refinement (i.e, $\theta = 0$). If the numerical IgA scheme is implemented correctly, the error is supposed 
to decrease at least as $O(h^p)$ (which is verified throughout the numerical tests in Section 
\ref{sec:numerical-examples}).

Finally, on the last {\rm REFINE} step, we apply the refinement algorithm $\mathcal{R}$ to those 
elements that have been selected on the {\rm MARK} level. Since the THB-splines are based on the 
subdomains of different hierarchical levels, the procedure $\mathcal{R}$ increases  the level of 
subdomains that have been selected by ${\mathds{M}}_{*}({\theta})$. For $\mathcal{R}$, 
a dyadic cell refinement is applied. 
{To prevent the cases of refinement, when the inserted box is not 
aligned with the current hierarchical mesh (occurrence of the L-shaped cells), `affected' cells of lower 
levels are locally subdivided to adapt to the inserted box. For that, in further examples, we specify the 
extension of the refined box by one cell (see, e.g., Figure \ref{fig:thb-extension}).
Here, Figure \ref{fig:thb-extension-a} illustrates the box insertion (yellow area) in the second hierarchical 
level of THB-spline. In Figure \ref{fig:thb-extension-b}, blue cells around the inserted box are 
the `one-cell' extension of the yellow area. Green cells of the first level are the so-called 
`affected' cells of zero level that have been locally subdivided to adapt to the inserted box.
}

\begin{algorithm}[!t]
\begin{algorithmic} 
\small
\STATE {{\bf Input:}}
$\mathcal{K}_h$\, 
\COMMENT{discretisation of $\Omega$} \\
\qquad \quad \, ${\rm span} \,\big\{ \phi^{(p)}_{h, i} \big\}$, $i = 1, ..., |\mathcal{I}|$ \,
\COMMENT{$V_{h}$-basis} \\[8pt]
\STATE {{\bf APPROXIMATE}}: \\[4pt]
\begin{itemize}
\item ASSEMBLE the matrix ${{\rm {K}}^{(p)}_h}$ and RHS ${{\rm {f}}^{(p)}_h}$
\hfill :${\bf t_{\rm as}(u_h)}$
\\[4pt]
\item SOLVE ${\rm K}^{(p)}_h \, \underline{ \rm  u}_h = {\rm f}^{(p)}_h$
\hfill :${\bf t_{\rm sol}(u_h)}$
\\[4pt]
\item Reconstruct $u_h = \Sum_{i \in \mathcal{I}} \underline{\rm u}_i \, \phi^{(p)}_{h,i}(x)$ \\[8pt]
\end{itemize}
\STATE {{\bf {ESTIMATE}}}: Reconstruct $\maj{} (u_h, \flux_h)$ and $\mdI{} (u_h, \flux_h)$ 
\hfill :${\bf t_{\rm as}(\boldsymbol{{y}}_h) + t_{\rm sol}(\boldsymbol{{y}}_h)}$
\\[8pt]
\STATE \qquad \qquad \qquad \; Reconstruct $\mij{} (u_h, w_h)$
\hfill :${\bf t_{\rm as}(w_h) + t_{\rm sol}(w_h)}$\\[8pt]
\STATE{{\bf MARK}}: Using the marking criteria $\mathds{M}_*(\theta)$, select the elements $K$ of mesh 
$\mathcal{K}_h$ that must be refined \\[8pt]
%
\STATE {{\bf REFINE}}: Execute the refinement strategy:  
$\mathcal{K}_{h_{\rm ref}} = \mathcal{R}(\mathcal{K}_h)$ \\[8pt]
\STATE {{\bf Output:}} $\mathcal{K}_{h_{\rm ref}}$ \COMMENT{refined discretisation of $\Omega$}
\end{algorithmic}
\caption{\small Reliable reconstruction of $u_h$ (a single refinement step)}
\label{alg:reliable-uh-reconstruction}
\end{algorithm}

{\begin{figure}[!t]
	\centering
	{
	\subfloat[]{
	\includegraphics[width=4cm, clip]{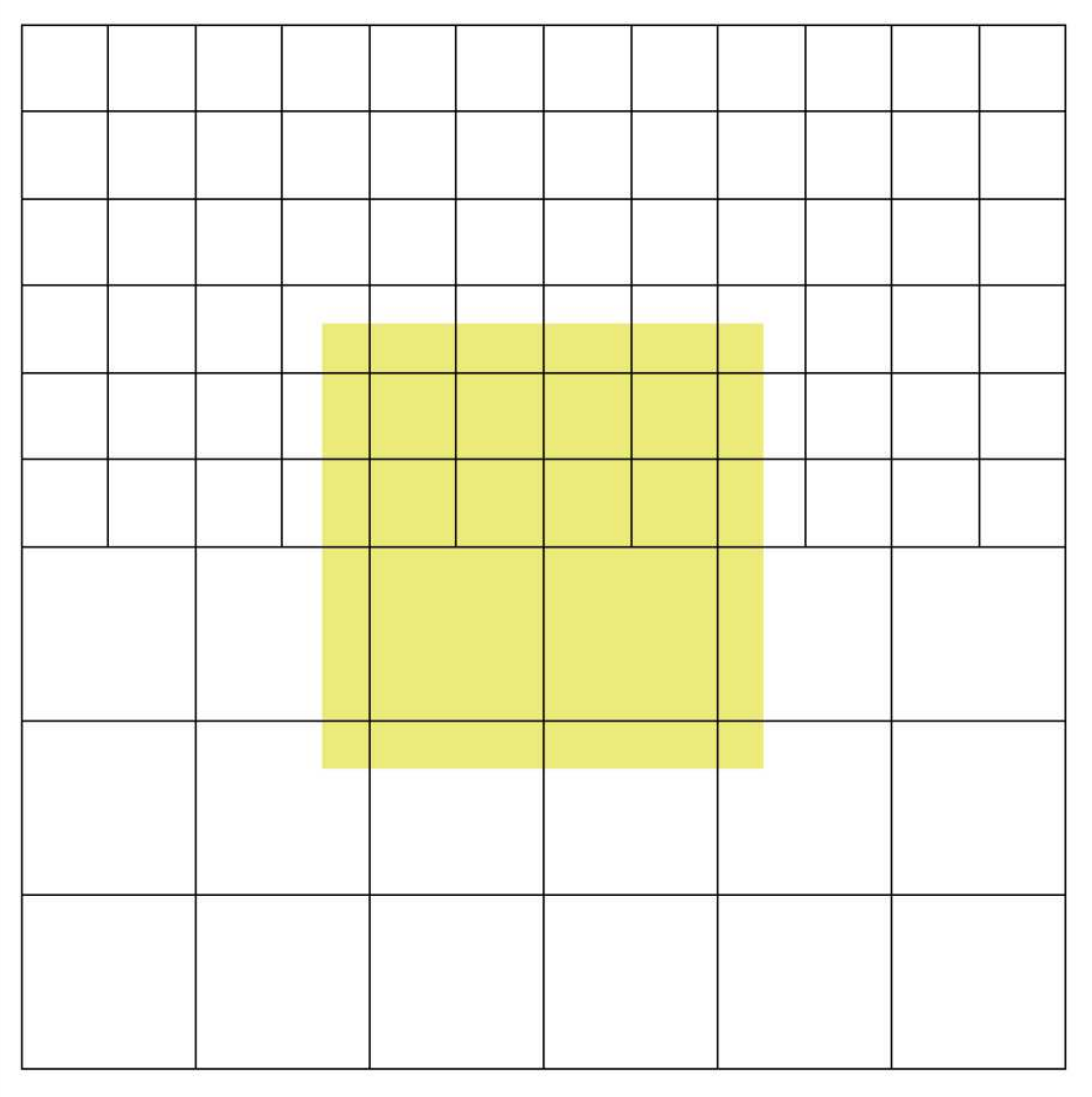}
	\label{fig:thb-extension-a}} \qquad
	\subfloat[]{
	\includegraphics[width=4cm, clip]{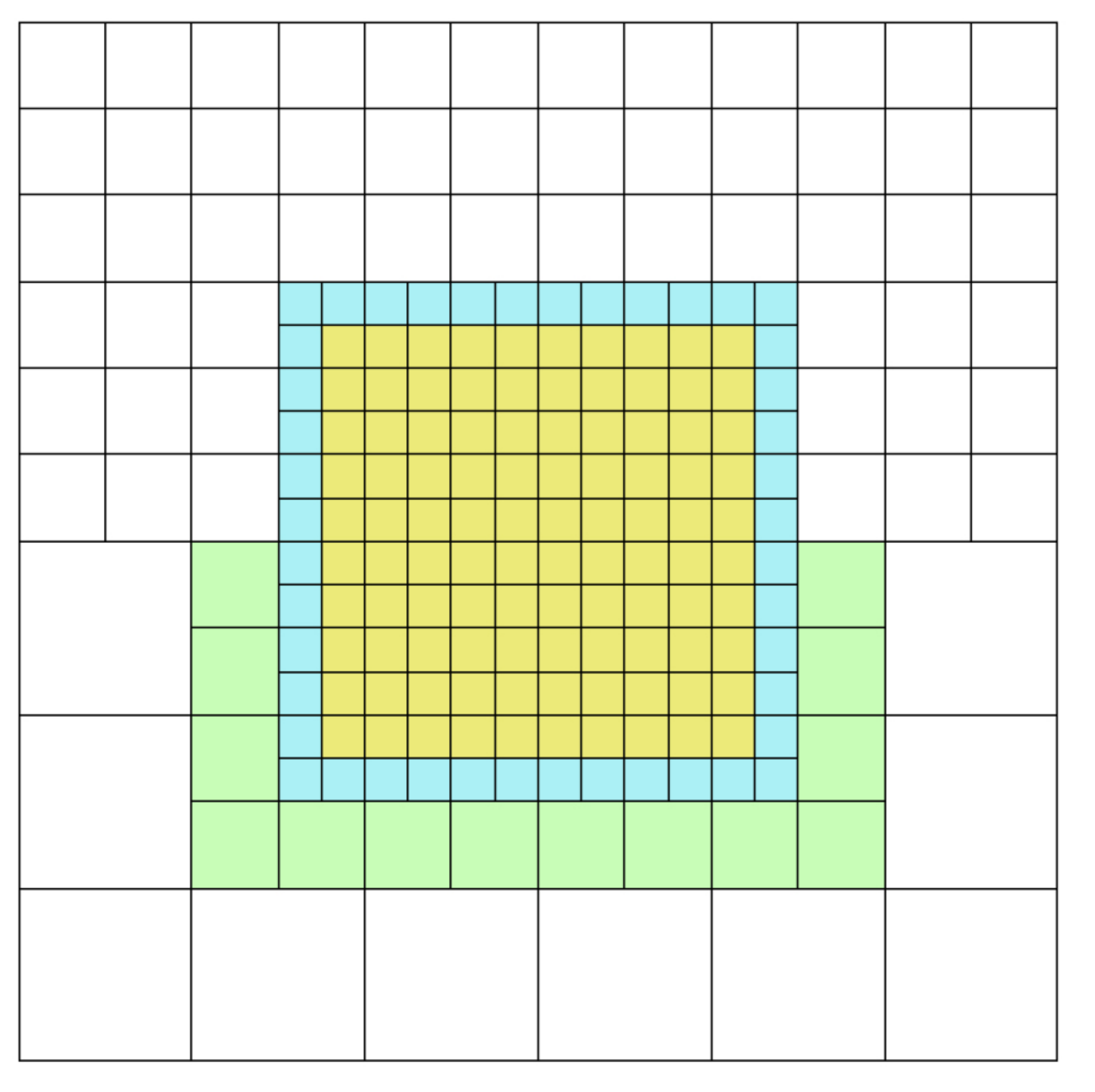}
	\label{fig:thb-extension-b}}
	}
	\caption{Example of the box insertion in the second hierarchical level of THB-spline (the picture is taken from \cite{Giannellietall2016}).}
	\label{fig:thb-extension}	
\end{figure}
}

\begin{algorithm}[!t]
\small
\caption{\small \quad \small {\bf ESTIMATE} step (majorant reconstruction)}
\begin{algorithmic} 
\STATE  {{\bf Input:}} 
 $u_h$ \,\COMMENT{approximation}\\[2pt] 
\qquad \quad \, $\mathcal{K}_{Mh}$ \, \COMMENT{disctretisation of $\Omega$}, \\[2pt]
\qquad \quad \, ${\rm span} \,\big\{ {\boldsymbol \psi}_{h, i} \big\}$, $i = 1, ..., d |\mathcal{I}|$ \, \COMMENT{$Y_{h}$-basis}, \\[2pt]
\qquad \quad \, $N^{\rm it}_{\rm maj}$ \, \COMMENT{number of optimisation iterations} \\[4pt]
\STATE  {\bf ASSEMBLE}  
${{\rm {Div}}_h}, {{\rm {M}}_h} \in \mathds{R}^{d |\mathcal{I}| \times d |\mathcal{I}|}$ 
and ${\rm z}_h$, ${\rm g}_h \in {\mathds R}^{d |\mathcal{I}|}$
\hfill :${\bf t_{\rm as}(\boldsymbol{{y}}_h)}$
\\[4pt]
\STATE  {Set} $\beta^{(0)} = 1$ \\[4pt]
\FOR {$n = 1$ {\bf to} $N^{\rm it}_{\rm maj}$} 
\STATE {\bf SOLVE} \quad 
$\big( \rfrac{\CFriedrichs^2}{\beta^{(n-1)} \, } {{\rm {Div}}_h} + {{\rm {M}}_h} \big)\, \underline{ \rm \bf y}^{(n)}_{h} = 
	-\rfrac{\CFriedrichs^2}{\beta^{(n-1)} \,} {\rm z}_h + {\rm g}_h$ 
\hfill :${t_{\rm sol}(\boldsymbol{{y}}_h)}$ \\[6pt]
\STATE  {Reconstruct}
$\flux^{(n)}_h 
:= \sum_{i \in \mathcal{I} \times d} \underline{ \rm \bf y}^{(n)}_{h,i} \, {{\boldsymbol \psi}}_{h,i}$ \\[6pt]
\STATE  {Compute } 
$\overline{\mathrm m}^{2, (n)}_{\mathrm{d}} := \| \, f + \dvrg \flux^{(n)}_h \, \|^2_{\Omega}$ and 
$\overline{\mathrm m}^{2, (n)}_{\mathrm{eq}} := \| \flux^{(n)}_h - \nabla u_h \,\|^2_{\Omega}$ \\[6pt]
\STATE  {Compute }
$
\beta^{(n)} = 
\tfrac{\CFriedrichs \, \mfInosq^{(n)}}{ \mdInosq^{(n)}}$ \\[6pt]
\ENDFOR 
\vskip 8pt
\STATE {Compute}
%
$\overline{\rm M}^2 (u_h, \flux^{(n)}_h; \beta^{(n)}) 
	:= (1 + \beta^{(n)}) \, \overline{\mathrm m}^{2, (n)}_{\mathrm{eq}} + 
   	    (1 + \tfrac{1}{\beta^{(n)}}) \, {\CFriedrichs^2}
		\overline{\mathrm m}^{2, (n)}_{\mathrm{d}}$ \\[2pt]
\STATE {{\bf Output:}} $\maj{}$ \qquad \quad \; \COMMENT{total error majorant on $\Omega$}, \\
\qquad \quad \;\;\;\, $\mdInosq = \mdInosq^{(n)}$ \COMMENT{indicator of error distribution over $\mathcal{K}_h$}
\end{algorithmic}
\label{alg:estimate-step-maj}
\end{algorithm}

Let us now consider the structure of Algorithm \ref{alg:estimate-step-maj}, which clarifies the 
{\rm ESTIMATE} step of Algorithm \ref{alg:reliable-uh-reconstruction} in the context of functional type 
error estimates. As the first Input argument, the algorithm receives the approximate solution $u_h$ 
reconstructed with the IgA scheme. Since the majorant is minimised with respect to a vector-valued 
variable $\flux_h \in Y_{h}$, the algorithm is also provided with the collection of basis functions generating the 
space $Y_{h} := {\rm span} \,\big\{ {\boldsymbol \psi}_{h, i} \big\}$, $i = 1, ..., d |\mathcal{I}|$. The 
last input parameter $N^{\rm it}_{\rm maj}$ defines the number of the optimisation loops executed 
to obtain a good enough minimiser of $\maj{}$. According to the tests performed in 
\cite{RepinSauterSmolianski2003} as well as experience of the author, 
one or two iterations are usually rather sufficient in order to achieve reasonable accuracy of error majorant. 
Technically, if the ratio between $\mfI$ and $\mdI$ is 
small enough, the loop can be exited even if $n < N^{\rm it}_{\rm maj}$. This condition might cut the 
computational costs for the error control. However, for the consistency of exposition this is not incorporated 
into Algorithm \ref{alg:estimate-step-maj} but only noted here as a remark.

It is crucial to emphasise that both matrices ${{\rm {Div}}_h}, {{\rm {M}}_h}$ and vectors ${\rm z}_h$, 
${\rm g}_h$ are assembled only once and remain unchanged in the minimisation procedure. The loop is 
iterated $N^{\rm it}_{\rm maj}$ times, where on each step the optimal $\underline{ \rm \bf y}^{(n)}_{h}$ 
and $\beta^{(n)}$ are reconstructed. In our implementation, 
the optimality system for the flux (cf. \eqref{eq:system-fluxh}) is 
solved by direct sparse ${\rm LDL^{\rm T}}$ Cholesky factorisations for $d=2$ and by a conjugate 
gradient method for $d=3$ (again, the initial guess is reconstructed from the approximation obtained in 
the earlier refinement). The time spent on {\rm ASSEMBLE} and {\rm SOLVE} steps with regard to the 
system \eqref{eq:system-fluxh} is measured by ${\bf t_{\rm as}(\flux_h)}$ and ${\bf t_{\rm sol}(\flux_h)}$, 
respectively, and compared to values ${\bf t_{\rm as}}(u_h)$ and ${\bf t_{\rm sol}}(u_h)$ in forthcoming 
numerical examples. 

Algorithm \ref{alg:estimate-step-min} illustrates the sequence of steps for lower error bound 
reconstruction. Both assembling and solving are analogous to Algorithm \ref{alg:reliable-uh-reconstruction} 
for the primal approximation, but use the basis $\phi^r_{h, i}$ of higher regularity. At the same 
time, we use the mesh $\mathcal{K}_{L h}$ that is up to $L$ times coarser than the mesh used for $u_h$. 
The time-efficiency of the minorant reconstruction is tracked by ${\bf t_{\rm as}(w_h)}$ and 
${\bf t_{\rm sol}(w_h)}$ measurements and later compared to those related to $u_h$ and $\flux_h$ 
approximation.

Besides the computational costs related to the assembling and solving  
of \eqref{eq:system-uh} and \eqref{eq:system-fluxh}, we measure the time spent on the element-wise 
(e/w) evaluation of error, majorant, minorant, and the residual error estimator. They are denoted by 
$t_{\rm e/w}(\| \nabla e \|)$, $t_{\rm e/w}(\overline{\rm M})$, $t_{\rm e/w}(\underline{\rm M})$, and 
$t_{\rm e/w}(\overline{\eta})$, respectively. 
\begin{algorithm}[!t]
\small
\caption{\small \quad \small {\bf ESTIMATE} step (minorant reconstruction)}
\begin{algorithmic} 
\STATE  {{\bf Input:}} 
 $u_h$ \COMMENT{approximation}\\[2pt] 
\qquad \quad \, $\mathcal{K}_{L h}$ \COMMENT{disctretisation of $\Omega$}, \\[2pt]
\qquad \quad \, ${\rm span} \,\big\{  \phi^{(r)}_{h, i} \big\}$, $i = 1, ..., |\mathcal{I}|$ \COMMENT{$W_{h}$-basis} \\[2pt]
\STATE  {\bf ASSEMBLE}  
${\rm K}^{(r)}_h \in \mathds{R}^{|\mathcal{I}| \times |\mathcal{I}|}$ 
and ${\rm f}^{(r)}_h \in {\mathds R}^{|\mathcal{I}|}$
\hfill :${\bf t_{\rm as}(w_h)}$
\\[4pt]
\STATE {\bf SOLVE} \quad 
${\rm K}^{(r)}_h \, \underline{ \rm  w}_h = {\rm f}^{(r)}_h$ 
\hfill :${\bf {t_{\rm sol}(w_h)}}$ \\[6pt]
\STATE  {Reconstruct}
$w_h(x) = w_h(x_1, . . . , x_d) := \sum_{i \in \mathcal{I}} \underline{ \rm w}_{h, i} \, \phi^{(r)}_{h, i}$ \\[6pt]
\STATE {Compute}
%
$\mij{}^2 (u_h, w_h) 
:= 2 \, (f, u_h - w_h) - (\| \nabla u_h \|^2_{\Omega} + \, \| \nabla w_h \|^2_{\Omega})$ \\[2pt]
\STATE {{\bf Output:}} $\mij{}$ \COMMENT{total error minorant on $\Omega$} \\
\end{algorithmic}
\label{alg:estimate-step-min}
\end{algorithm}

\section{Numerical examples}
\label{sec:numerical-examples}

In the current section, we present a series of examples demonstrating the numerical properties of 
the error majorants discussed above. We start with relatively simple examples, in which we aim to introduce 
the main properties of majorant and, at the same time, familiarise the reader with the structure of 
performed numerical tests. This approach is intended to bring the focus to analysis in more complicated 
examples discussed further.


\begin{example}
\label{ex:unit-domain-example-2}
\rm
First, we consider a basic example with 
%
$$u = (1 - x_1) \, x_1^2 \,  (1 - x_2) \, x_2, \quad 
f = - \big(2\, (1 - 3\, x_1)\,(1 - x_2)\,x_2 - 2\, (1 - x_1)\,x_1^2\big) \quad \mbox{in} \quad \Omega,$$ 
%
and homogeneous Dirichlet boundary condition (BC). 
\begin{table}[!t]
\footnotesize
\centering
\newcolumntype{g}{>{\columncolor{gainsboro}}c} 	
\begin{tabular}{c|c|ccc|gc|c}
\quad \# ref. \qquad & 
\quad  $\| \nabla e \|_\Omega$ \qquad & 	  
\quad $\overline{\rm M}$ \qquad &    
\quad $\mdInosq$ \qquad & 	       
\quad  $\mfInosq$ \qquad &  
\quad $\Ieff (\overline{\rm M})$ \qquad & 
\quad $\Ieff (\overline{\rm \eta})$ \qquad & 
\quad e.o.c \qquad \\
\midrule
\multicolumn{8}{c}{ (a) $\flux_h \in S_{3h}^{5, 5} \oplus S_{3h}^{5, 5}$ ($m = 3$, $M = 3$)} \\
\midrule
   3 &   2.5648e-03 &   3.2806e-03 &   3.1546e-03 &   5.5974e-04 &       1.2791 &      11.0113 &   3.4565 \\
   5 &   1.5952e-04 &   1.9770e-04 &   1.9084e-04 &   3.0441e-05 &       1.2393 &      10.9580 &   2.3602 \\
   7 &   9.9673e-06 &   1.1974e-05 &   1.1921e-05 &   2.3799e-07 &       1.2013 &      10.9546 &   2.0901 \\
   9 &   6.2294e-07 &   7.4549e-07 &   7.4502e-07 &   2.0851e-09 &       1.1967 &      10.9545 &   2.0225 \\
  11 &   3.8934e-08 &   4.6571e-08 &   4.6564e-08 &   3.2185e-11 &       1.1962 &      10.9545 &   2.0056 \\
\midrule
\multicolumn{8}{c}{ (b) $\flux_h \in S_{7h}^{9, 9} \oplus S_{7h}^{9, 9}$ ($m = 7$, $M = 7$)} \\
\midrule
   3 &   2.5648e-03 &   2.6756e-03 &   2.5800e-03 &   4.2495e-04 &       1.0432 &      11.0113 &   3.4565 \\
   5 &   1.5952e-04 &   1.7737e-04 &   1.6869e-04 &   3.8537e-05 &       1.1118 &      10.9580 &   2.3602 \\
   7 &   9.9673e-06 &   1.0215e-05 &   1.0035e-05 &   7.9975e-07 &       1.0248 &      10.9546 &   2.0901 \\
   9 &   6.2294e-07 &   6.9080e-07 &   6.3274e-07 &   2.5797e-07 &       1.1089 &      10.9545 &   2.0225 \\
  11 &   3.8934e-08 &   4.0932e-08 &   3.9140e-08 &   7.9608e-09 &       1.0513 &      10.9545 &   2.0056 \\

\end{tabular}
\caption{\small {\em Ex. \ref{ex:unit-domain-example-2}}. Error, majorant (with dual and equilibrated terms), 
efficiency indices, 
and e.o.c. w.r.t. unif. ref. steps.}
\label{tab:unit-domain-example-2-error-majorant-v-2-y-5-y-9-uniform-ref}
\end{table}
%
\begin{table}[!t]
\footnotesize
\centering
\newcolumntype{g}{>{\columncolor{gainsboro}}c} 	
\begin{tabular}{c|cc|cg|cg|cgc}
\# ref & 
\# d.o.f.($u_h$) &  \# d.o.f.($\flux_h$) &  
\; $t_{\rm as}(u_h)$ \; & 
\; $t_{\rm as}(\flux_h)$ \; & 
\; $t_{\rm sol}(u_h)$ \; & 
\; $t_{\rm sol}(\flux_h)$ \; &
$t_{\rm e/w}(\| \nabla e \|)$ & 
$t_{\rm e/w}(\overline{\rm M})$ & 
$t_{\rm e/w}(\overline{\eta})$ \\
\midrule
\multicolumn{10}{c}{ (a) $\flux_h \in S_{3h}^{5, 5} \oplus S_{3h}^{5, 5}$ ($q = 5$, $m = 3$, $M = 3$)} \\
\midrule
   1 &          9 &         36 &     0.0013 &     0.0023 &           0.0001 &           0.0017 &       0.0000 &       0.0008 &       0.0006 \\
   3 &         36 &         36 &     0.0010 &     0.0025 &           0.0001 &           0.0018 &       0.0005 &       0.0025 &       0.0022 \\
   5 &        324 &         81 &     0.0094 &     0.0216 &           0.0008 &           0.0094 &       0.0087 &       0.0159 &       0.0276 \\
   7 &       4356 &        441 &     0.0729 &     0.2506 &           0.0439 &           0.0730 &       0.1830 &       0.1571 &       0.2858 \\
   9 &      66564 &       4761 &     1.2661 &     4.0725 &           3.6962 &           8.3926 &       2.5329 &       2.8220 &       4.1740 \\
  11 &    1052676 &      68121 &    22.4621 &    68.2723 &         211.1700 &         570.4293 &      37.8862 &      37.9696 &      65.6940 \\
\midrule
\multicolumn{10}{c}{ (b) $\flux_h \in S_{7h}^{9, 9} \oplus S_{7h}^{9, 9}$ ($q = 9$, $m = 7$, $M = 7$)} \\
\midrule
   1 &          9 &        100 &     0.0008 &     0.0234 &           0.0001 &           0.0122 &       0.0002 &       0.0025 &       0.0004 \\
   3 &         36 &        100 &     0.0006 &     0.0167 &           0.0001 &           0.0132 &       0.0003 &       0.0038 &       0.0012 \\
   5 &        324 &        100 &     0.0089 &     0.0258 &           0.0010 &           0.0057 &       0.0048 &       0.0269 &       0.0167 \\
   7 &       4356 &        100 &     0.0750 &     0.0140 &           0.0401 &           0.0093 &       0.1564 &       0.5749 &       0.3073 \\
   9 &      66564 &        169 &     1.1129 &     0.1967 &           3.2580 &           0.0763 &       2.5923 &       6.2473 &       4.2985 \\
  11 &    1052676 &        625 &    17.6219 &     3.9372 &         196.0170 &           1.2941 &      35.1466 &      99.9845 &      61.1072 \\
\end{tabular}
\caption{\small {\em Ex. \ref{ex:unit-domain-example-2}}. 
Time for assembling and solving the systems that generate $u_h$ and $\flux_h$, 
time for e/w evaluation of error, majorant, and residual error estimator w.r.t. unif. ref. steps. 
}
\label{tab:unit-domain-example-2-times-v-2-y-5-y-9-uniform-ref}
\end{table}

Let the primal variable be approximated by the splines of degree $p = 2$, i.e., the discretisation space 
$S^{p, p}_h$. For the uniform refinement (unif. ref.), we first test the idea introduced in \cite{LMR:KleissTomar2015} 
and compare two different settings for spaces approximating auxiliary dual variable $\flux_h \in S^{q, q}_{Mh}$:
\begin{alignat}{2}
(a) \; & q = 5, \; m = 3, \; M = 3, \quad  \mbox{and} \quad (b) \; q = 9, \; m = 7, \; M = 7.
\end{alignat}
After performing $N_{\rm ref} = 11$ unif. ref. steps, we present the 
obtained numerical results {in} Tables 
\ref{tab:unit-domain-example-2-error-majorant-v-2-y-5-y-9-uniform-ref}--\ref{tab:unit-domain-example-2-times-v-2-y-5-y-9-uniform-ref}
(where the upper and the lower parts correspond to the cases (a) and (b), respectively). 
The efficiency of functional error majorant is confirmed by corresponding indices, i.e., 
{${ \Ieff (\overline{{\rm M}}) ={1.1961}}$ } for the case (a) and 
{${ \Ieff (\overline{{\rm M}}) ={1.0024}}$ } for the case (b) (see the shaded column of Table 
\ref{tab:unit-domain-example-2-error-majorant-v-2-y-5-y-9-uniform-ref}).
The expected error order of convergence (e.o.c.) $p = 2$ is confirmed by the last column of Table 
\ref{tab:unit-domain-example-2-error-majorant-v-2-y-5-y-9-uniform-ref}. 
{The time needed for the residual-based estimate evaluation is dependent only 
on the approximation $u_h$. Since \# d.o.f.($u_h$) stays the same, 
$t_{\rm e/w}(\overline{\eta})$ stays approximately the same. 
The time required for the evaluation of the majorant, in turn, is more complex, since it is dependent on 
both $u_h$ and $\flux_h$. In the case (b), it is larger than in the case (a), because majorant is evaluated 
on the mesh $\mathcal{K}_h$ (the same one that is used for $u_h$ approximation). In the same time, 
the degree of splines used to approximate $y_h$ is higher ($\flux_h \in S_{7h}^{9, 9} \oplus S_{7h}^{9, 9}$). 
That is the reason we obtain an overhead in the evaluation time for the majorant. 
Therefore, 
 the fact the time for evaluation of the majorant  is larger than that of the residual error estimator (when 
 higher order approximation is used) can be explained exactly by the increase of the degree of B-splines.
}

\newpage
When the computational costs are considered, the time spent on the reconstruction of $\flux_h$ \linebreak
(i.e., $t_{\rm as}(\flux_h)$ + $t_{\rm sol}(\flux_h)$) is about {${2-3}$ times} higher than the time 
$t_{\rm as}(u_h)$ + $t_{\rm sol}(u_h)$ in the setting (a). However, for the case (b), the assembling 
time of ${{\rm {Div}}_h}$ and ${{\rm {M}}_h}$ denoted by $t_{\rm as}(\flux_h)$ takes 
approximately {${{1}/{4}}$}-th of the assembling time for ${\rm {K}_h}$ denoted by $t_{\rm as}(u_h)$. 
Similarly, solving the system \eqref{eq:system-fluxh} denoted by $t_{\rm sol}(\flux_h)$ requires only 
{${{1}/{150}}$-th} of time spent on solving \eqref{eq:system-uh}, i.e., $t_{\rm sol}(u_h)$. 
%

Due to the smoothness of the exact solution in this example, we can even use splines of lower 
degree for the flux approximation, e.g., $q = 3$, but at the same time reconstruct it on a much coarser 
mesh than for $u_h$, e.g., mesh $M = 8$ times coarser. The resulting efficiency indices are illustrated in Table 
\ref{tab:unit-domain-example-2-times-v-2-y-3-y-8-uniform-ref}, and corresponding times spent on 
the reconstruction of $u_h$ and $\flux_h$ (i.e., $\overline{\rm M} (u_h, \flux_h)$ and 
$\mdInosq (u_h, \flux_h)$) are presented in 
Table \ref{tab:unit-domain-example-2-error-majorant-v-2-y-3-y-8-uniform-ref}. By looking at 
Table \ref{tab:unit-domain-example-2-times-v-2-y-3-y-8-uniform-ref}, one can see the considerable 
speed-up in the time required for reconstruction of $\flux_h$ in comparison to $u_h$: 
$$
\tfrac{t_{\rm as}(u_h)}{t_{\rm as}(\flux_h)} \approx \tfrac{19.9547}{0.0114} \approx 1750 \quad  
\mbox{and} \quad
\tfrac{t_{\rm sol}(u_h)}{t_{\rm sol}(\flux_h)} \approx \tfrac{107.0756}{0.0020} \approx 53538.$$

\begin{table}[!t]
\footnotesize
\centering
\newcolumntype{g}{>{\columncolor{gainsboro}}c} 	
\begin{tabular}{c|c|ccc|gc|c}
\quad \# ref. \qquad & 
\quad  $\| \nabla e \|_\Omega$ \qquad & 	  
\quad $\overline{\rm M}$ \qquad &    
\quad $\mdInosq$ \qquad & 	       
\quad  $\mfInosq$ \qquad &  
\quad $\Ieff (\overline{\rm M})$ \qquad & 
\quad $\Ieff (\overline{\rm \eta})$ \qquad & 
\quad e.o.c. \qquad \\
\midrule
   3 &   2.5648e-03 &   3.0154e-03 &   3.0000e-03 &   6.8225e-05 &       1.1757 &      11.0115 &   3.4566 \\
   5 &   1.5952e-04 &   1.7571e-04 &   1.6338e-04 &   5.4779e-05 &       1.1015 &      10.9580 &   2.3602 \\
   7 &   9.9672e-06 &   1.1959e-05 &   1.0675e-05 &   5.7051e-06 &       1.1998 &      10.9547 &   2.0901 \\
   9 &   6.2294e-07 &   6.4308e-07 &   6.3365e-07 &   4.1905e-08 &       1.0323 &      10.9545 &   2.0225 \\
  11 &   3.8934e-08 &   4.0029e-08 &   3.9480e-08 &   2.4369e-09 &       1.0281 &      10.9545 &   2.0056 \\
\end{tabular}
\caption{\small {\em Ex. \ref{ex:unit-domain-example-2}}. Error, majorant (with dual and reliability terms), 
efficiency indices, and e.o.c.
for $\flux_h \in S_{8h}^{3, 3} \oplus S_{8h}^{3, 3}$ w.r.t. uniform ref. steps.}
\label{tab:unit-domain-example-2-error-majorant-v-2-y-3-y-8-uniform-ref}
\end{table}
\begin{table}[!t]
 \footnotesize
\centering
\newcolumntype{g}{>{\columncolor{gainsboro}}c} 	
\begin{tabular}{c|cc|cg|cg|cgc}
\# ref. & 
\# d.o.f.($u_h$) &  \# d.o.f.($\flux_h$) &  
\; $t_{\rm as}(u_h)$ \; & 
\; $t_{\rm as}(\flux_h)$ \; & 
\; $t_{\rm sol}(u_h)$ \; & 
\; $t_{\rm sol}(\flux_h)$ \; &
$t_{\rm e/w}(\| \nabla e \|)$ & 
$t_{\rm e/w}(\overline{\rm M})$ & 
$t_{\rm e/w}(\overline{\eta})$ \\
\midrule
  1 &          9 &         16 &     0.0009 &     0.0015 &     0.0000 &     0.0001 &       0.0001 &       0.0003 &       0.0003 \\
   3 &         36 &         16 &     0.0008 &     0.0005 &     0.0000 &     0.0001 &       0.0006 &       0.0007 &       0.0018 \\
   5 &        324 &         16 &     0.0081 &     0.0006 &     0.0005 &     0.0001 &       0.0184 &       0.0112 &       0.0285 \\
   7 &       4356 &         16 &     0.0753 &     0.0004 &     0.0173 &     0.0001 &       0.1391 &       0.1071 &       0.2534 \\
   9 &      66564 &         25 &     1.1899 &     0.0009 &     1.3832 &     0.0001 &       2.2776 &       1.6354 &       4.0632 \\
  11 &    1052676 &        121 &    19.9547 &     0.0114 &   107.0756 &     0.0020 &      36.0268 &      26.0721 &      63.6307 \\
\end{tabular}
\caption{\small {\em Ex. \ref{ex:unit-domain-example-2}}. 
Time for assembling and solving the systems that generate $u_h$ and $\flux_h$, 
time for e/w evaluation of error, majorant, and residual error estimator 
for $\flux_h \in S_{8h}^{3, 3} \oplus S_{8h}^{3, 3}$ w.r.t. unif. ref. steps.}
\label{tab:unit-domain-example-2-times-v-2-y-3-y-8-uniform-ref}
\end{table}

For the cases when the exact solution is not provided, the quality of the majorant can be 
verified by comparison of its values to the lower bound of the error. We assume that for the space 
approximating $w_h \in S^{r, r}_h$, we choose $r = 3$, and for the mesh ${{K}_{Lh}}$, the 
coarsening parameter is taken $L = 8$. To study the efficiency of the estimates, values of 
the majorant and minorant are compared to the error as well as to each other in Table
\ref{tab:unit-domain-example-2-error-majorant-minorant-eta-v-2-y-3-w-3-uniform-ref}. According to 
efficiency index $\Ieff (\underline{\rm M})$, minorant remains sharp w.r.t the increasing number of refinement steps 
(see column seven of Table \ref{tab:unit-domain-example-2-error-majorant-minorant-eta-v-2-y-3-w-3-uniform-ref}). 
As a result, the ratio of the upper and lower bounds ${\overline{\rm M}} / {\underline{\rm M}}$ 
is very close to one. This fact confirms that 
provided by $\overline{\rm M}$ and $\mij{}$ two-sided bound is guaranteed and the error of reconstructed 
approximation is contained inside of the interval $[\underline{\rm M}, \overline{\rm M}]$. 
In addition to the efficiency of the error estimates, 
we compare the time spent for assembling and solving systems 
\eqref{eq:system-uh},
\eqref{eq:system-fluxh}, and
\eqref{eq:system-wh} in Table 
\ref{tab:unit-domain-example-2-error-majorant-minorant-eta-v-2-y-3-w-3-uniform-ref}.
The last rows with ratios between time spend on the $u_h$-approximation w.r.t. $\flux_h$ and $w_h$
show that `expenses' related to the approximation of latter two variables are thousand times cheaper 
than the time dedicated to the primal variable. 

\begin{table}[!t]
\footnotesize
\centering
\newcolumntype{g}{>{\columncolor{gainsboro}}c} 	
\begin{tabular}{c|c|cg|cg|g|c}
\quad \# ref. \qquad & 
\quad  $\| \nabla e \|_\Omega$ \qquad & 	  
\quad $\overline{\rm M}$ \qquad &    
\quad $\Ieff (\overline{\rm M})$ \qquad & 
\quad $\underline{\rm M}$ \qquad &    
\quad $\Ieff (\underline{\rm M})$ \qquad & 
\quad ${\overline{\rm M}}/{\underline{\rm M}}$ \qquad & 
\quad e.o.c. \qquad \\
\midrule
   3 &     2.5648e-03 &   3.0154e-03 &       1.1757 &   2.5648e-03 &       1.0000 &       1.1757 &   3.4566 \\
   5 &     1.5952e-04 &   1.7571e-04 &       1.1015 &   1.5952e-04 &       1.0000 &       1.1015 &   2.3602 \\
   7 &     9.9672e-06 &   1.1959e-05 &       1.1998 &   9.9672e-06 &       1.0000 &       1.1998 &   2.0901 \\
   9 &     6.2294e-07 &   6.4308e-07 &       1.0323 &   6.2294e-07 &       1.0000 &       1.0323 &   2.0225 \\
\end{tabular}
\caption{{\em Ex. \ref{ex:unit-domain-example-2}}. 
Error, majorant, minorant, residual based error indicator 
with corresponding efficiency indices, and e.o.c. for 
$\flux_h \in S_{8h}^{3, 3} \oplus S_{8h}^{3, 3}$ and $w_h \in S_{8h}^{3, 3}$ 
w.r.t. unif. ref. steps.}
\label{tab:unit-domain-example-2-error-majorant-minorant-eta-v-2-y-3-w-3-uniform-ref}
\end{table}

\begin{table}[!t]
\footnotesize
\centering
\newcolumntype{g}{>{\columncolor{gainsboro}}c} 	
\begin{tabular}{c|ccc|ccc|ccc}
\parbox[c]{0.7cm}{\#ref.} & 
\parbox[c]{1.4cm}{\#d.o.f.($u_h$)} &  
\parbox[c]{1.4cm}{\#d.o.f.($\flux_h$)} &   
\parbox[c]{1.4cm}{\#d.o.f.($w_h$)} &
\parbox[c]{1.3cm}{$t_{\rm as}(u_h)$} & 
\parbox[c]{1.3cm}{$t_{\rm as}(\flux_h)$} & 
\parbox[c]{1.3cm}{$t_{\rm as}(w_h)$} & 
\parbox[c]{1.3cm}{$t_{\rm sol}(u_h)$} & 
\parbox[c]{1.3cm}{$t_{\rm sol}(\flux_h)$} & 
\parbox[c]{1.3cm}{$t_{\rm sol}(w_h)$} \\ 
\midrule
%
   3 &         36 &         16 &         16 &   8.01e-04 &   4.70e-03 &   2.39e-03 &         6.20e-05 &         1.72e-04 &         3.80e-05 \\ 
   5 &        324 &         16 &         16 &   8.25e-03 &   4.95e-03 &   8.54e-04 &         1.98e-04 &         1.32e-04 &         7.00e-06 \\ 
   7 &       4356 &         16 &         16 &   8.73e-02 &   1.84e-03 &   2.58e-03 &         1.46e-02 &         2.08e-04 &         2.10e-05 \\ 
   9 &      66564 &         25 &         25 &   1.26e+00 &   1.08e-02 &   2.30e-03 &         1.38e+00 &         5.32e-04 &         1.00e-05 \\ 
   \midrule
    &
    \multicolumn{3}{c|}{ }  &    
    \multicolumn{3}{c|}{547.82 \quad : \quad 4.69 \quad : \qquad 1 \quad } &
    \multicolumn{3}{c}{138000 \quad : \quad 53.2 \quad : \qquad 1 \quad } 
   \end{tabular}
\caption{{\em Ex. \ref{ex:unit-domain-example-2}}. 
Time for solving the systems that generate $u_h$, $\flux_h$, and $w_h$,
with direct and iterative methods for $\flux_h \in S_{8h}^{3, 3} \oplus S_{8h}^{3, 3}$ and 
$w_h \in S_{8h}^{3, 3}$ w.r.t. unif. ref. steps.}
\label{tab:unit-domain-example-2-times-v-2-y-3-w-3-uniform-ref}
\end{table}

\begin{figure}[!t]
	\centering
	{
	\subfloat[ref. \# 5: \quad ${\mathcal{K}}_h({\mathds{M}}_{{\rm \bf BULK}}(0.4))$]{
	\includegraphics[width=5.7cm, trim={2cm 2cm 2cm 2cm}, clip]{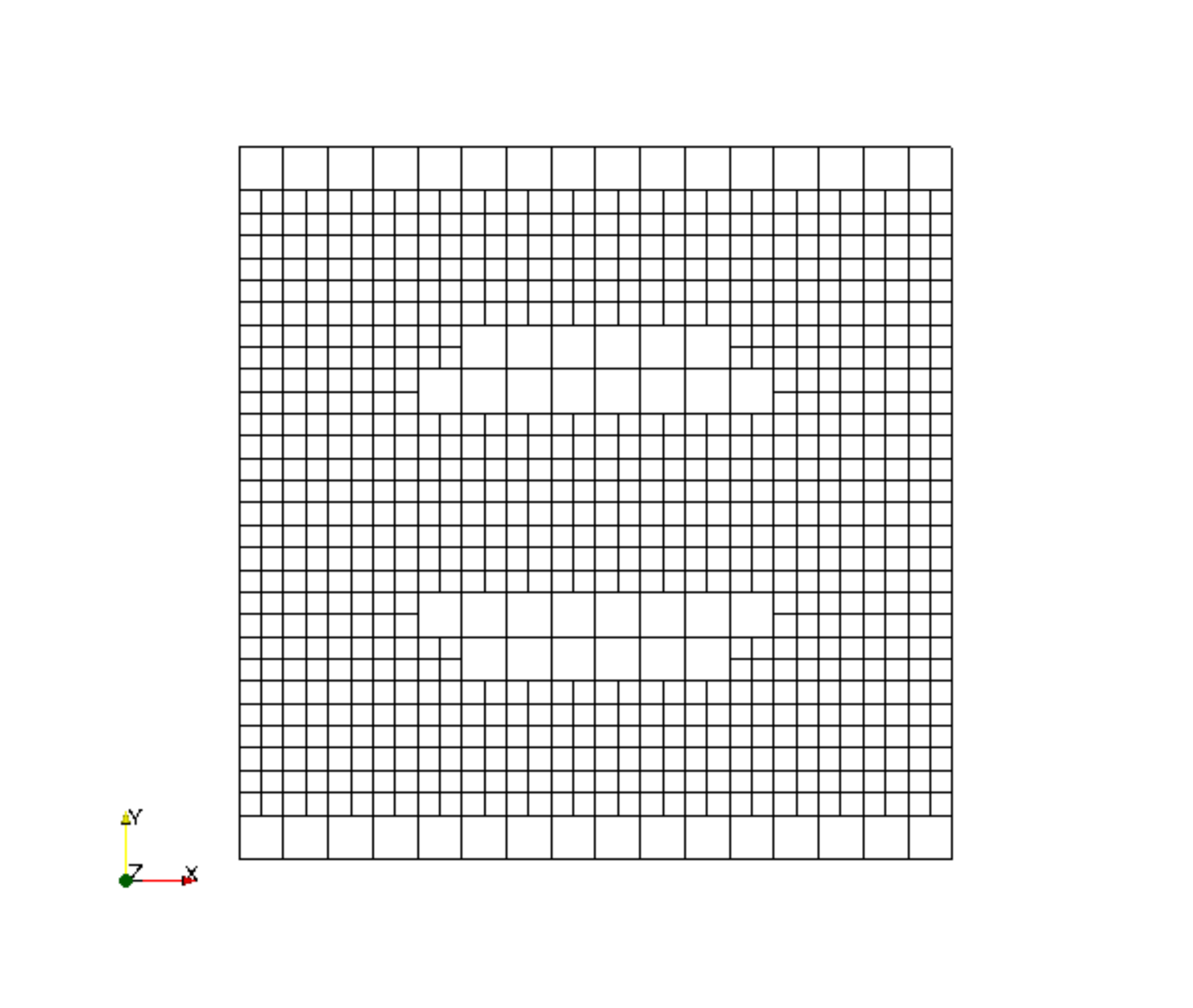}}
	\subfloat[ref. \# 5: \quad ${\mathcal{K}}_h({\mathds{M}}_{{\rm \bf BULK}}(0.2))$]{
	\includegraphics[width=5.7cm,  trim={2cm 2cm 2cm 2cm}, clip]{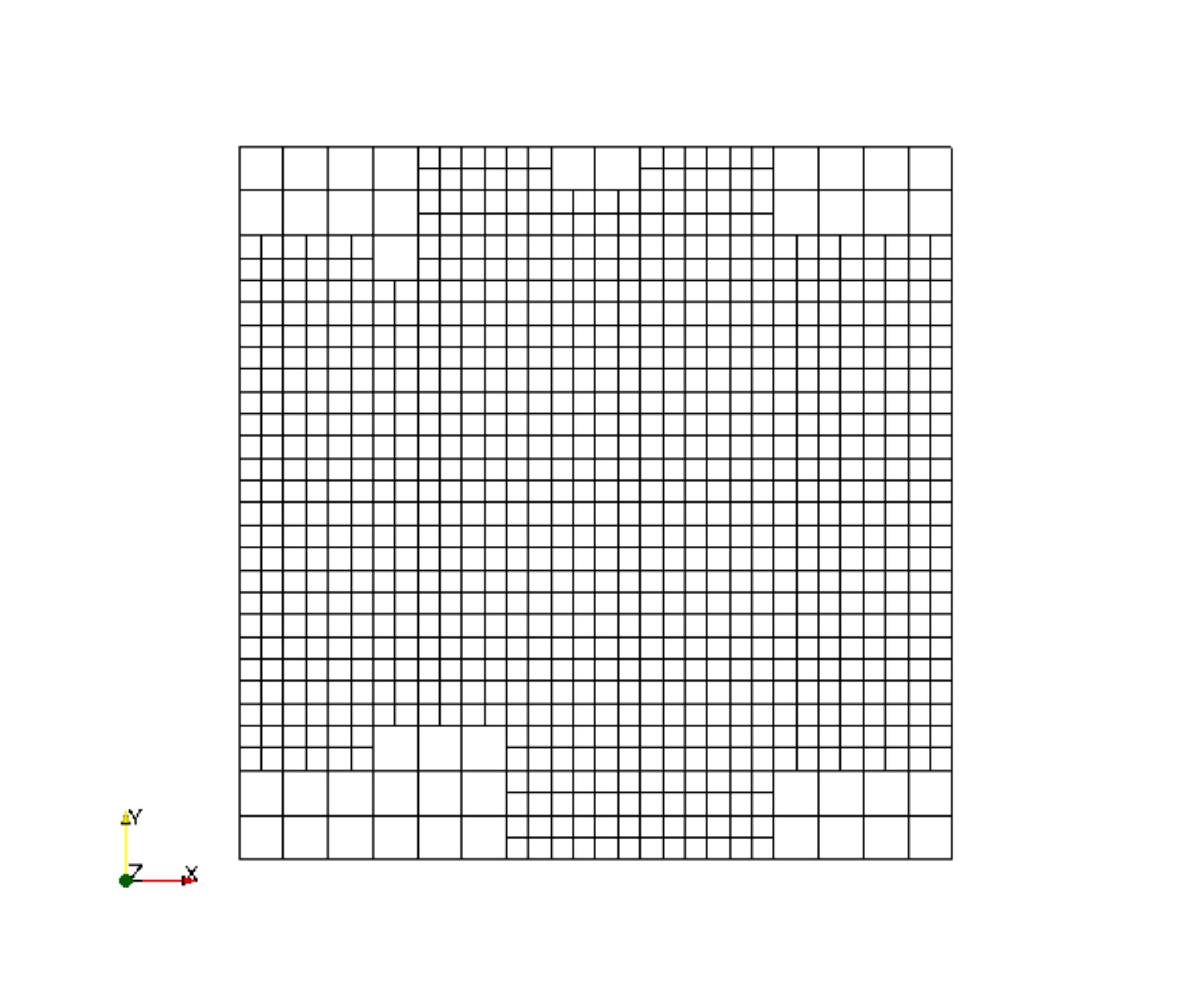}}
	}\\
	{
	\subfloat[ref. \# 6: \quad  ${\mathcal{K}}_h({\mathds{M}}_{{\rm \bf BULK}}(0.4))$]{
	\includegraphics[width=5.7cm,  trim={2cm 2cm 2cm 2cm}, clip]{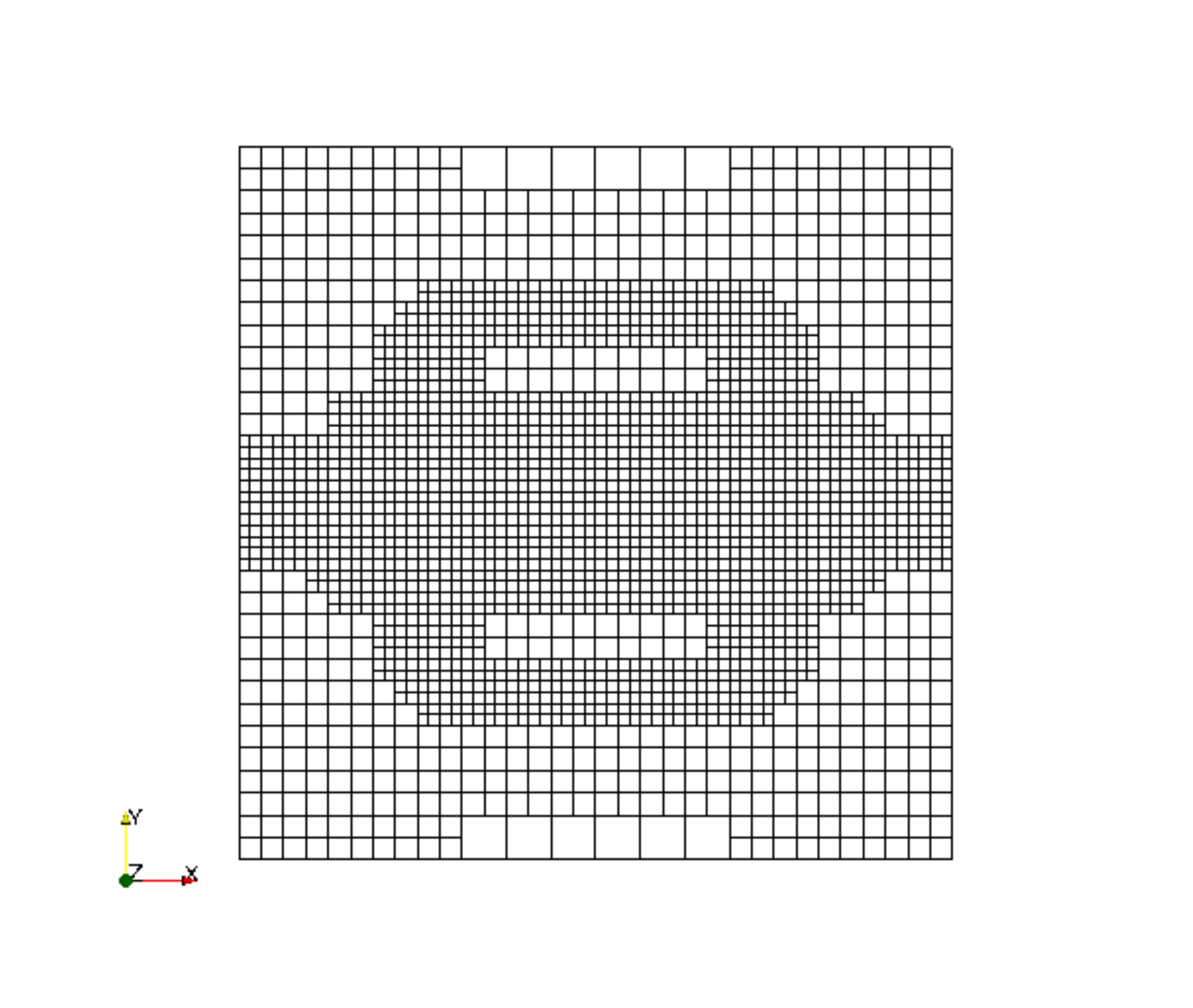}}
	\subfloat[ref. \# 6: \quad ${\mathcal{K}}_h({\mathds{M}}_{{\rm \bf BULK}}(0.2))$]{
	\includegraphics[width=5.7cm,  trim={2cm 2cm 2cm 2cm}, clip]{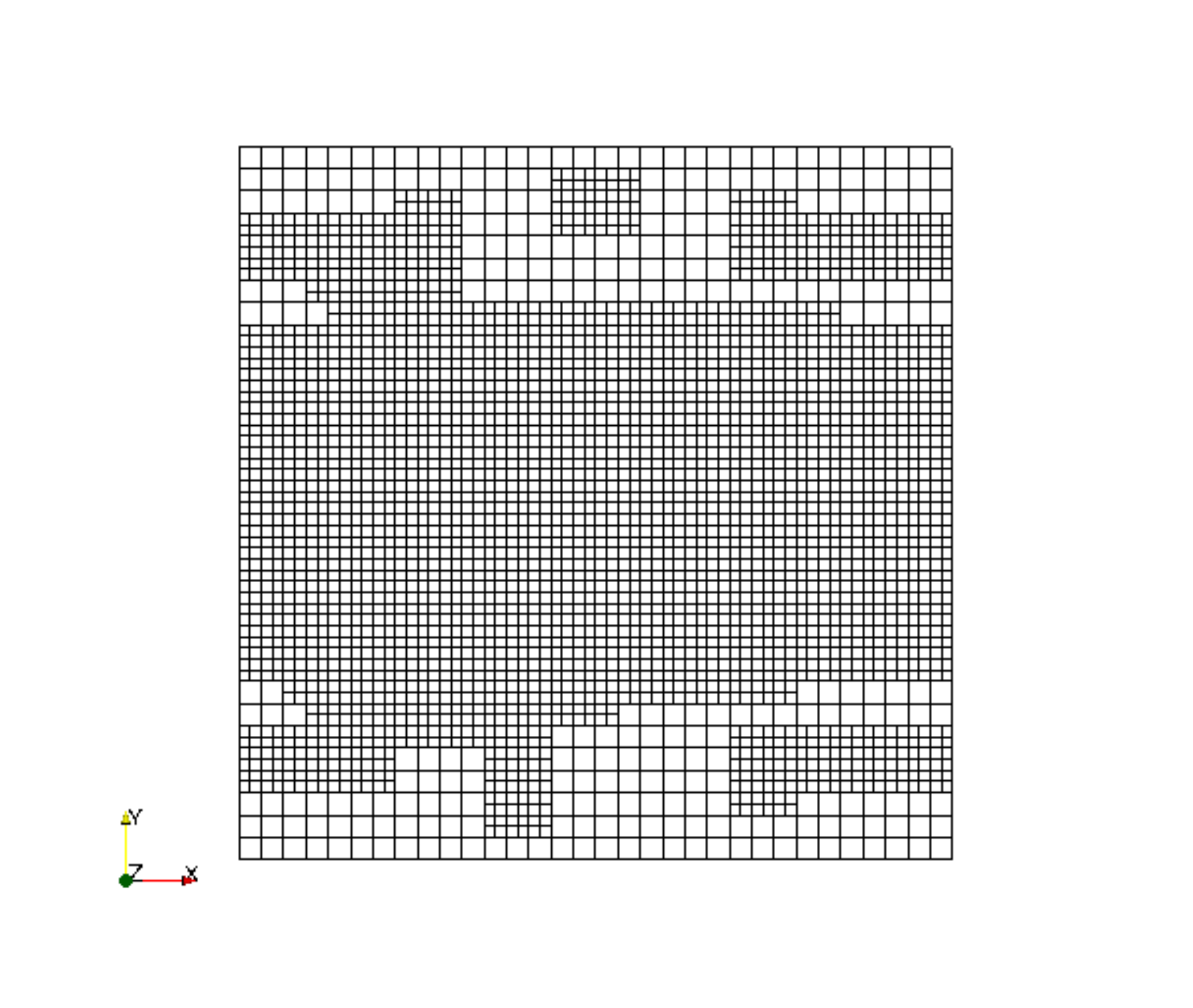}}
	}\\
	{
	\subfloat[ref. \# 7: \quad ${\mathcal{K}}_h({\mathds{M}}_{{\rm \bf BULK}}(0.4))$]{
	\includegraphics[width=5.7cm,  trim={2cm 2cm 2cm 2cm}, clip]{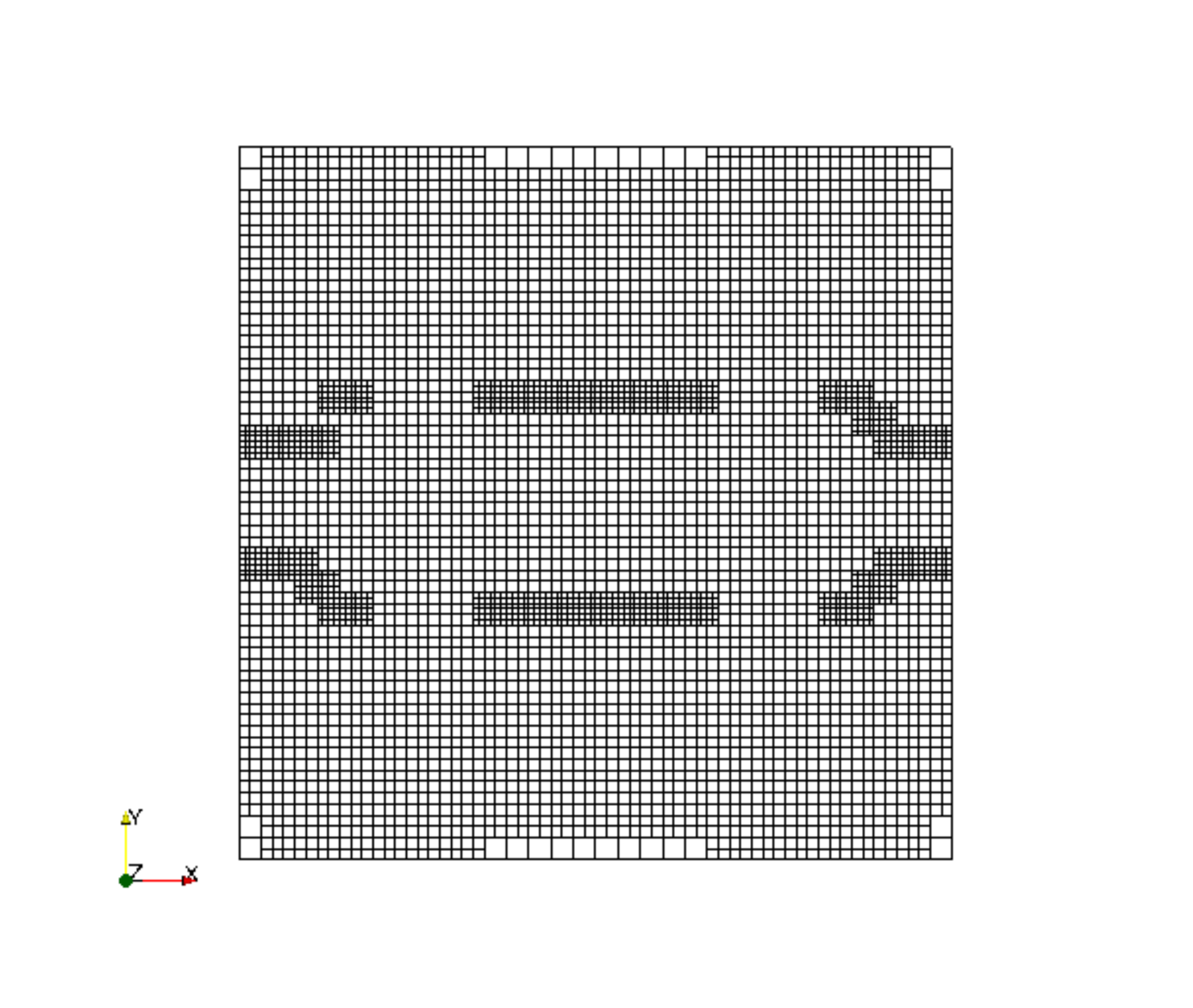}}
	\subfloat[ref. \# 7: \quad ${\mathcal{K}}_h({\mathds{M}}_{{\rm \bf BULK}}(0.2))$]{
	\includegraphics[width=5.7cm,  trim={2cm 2cm 2cm 2cm}, clip]{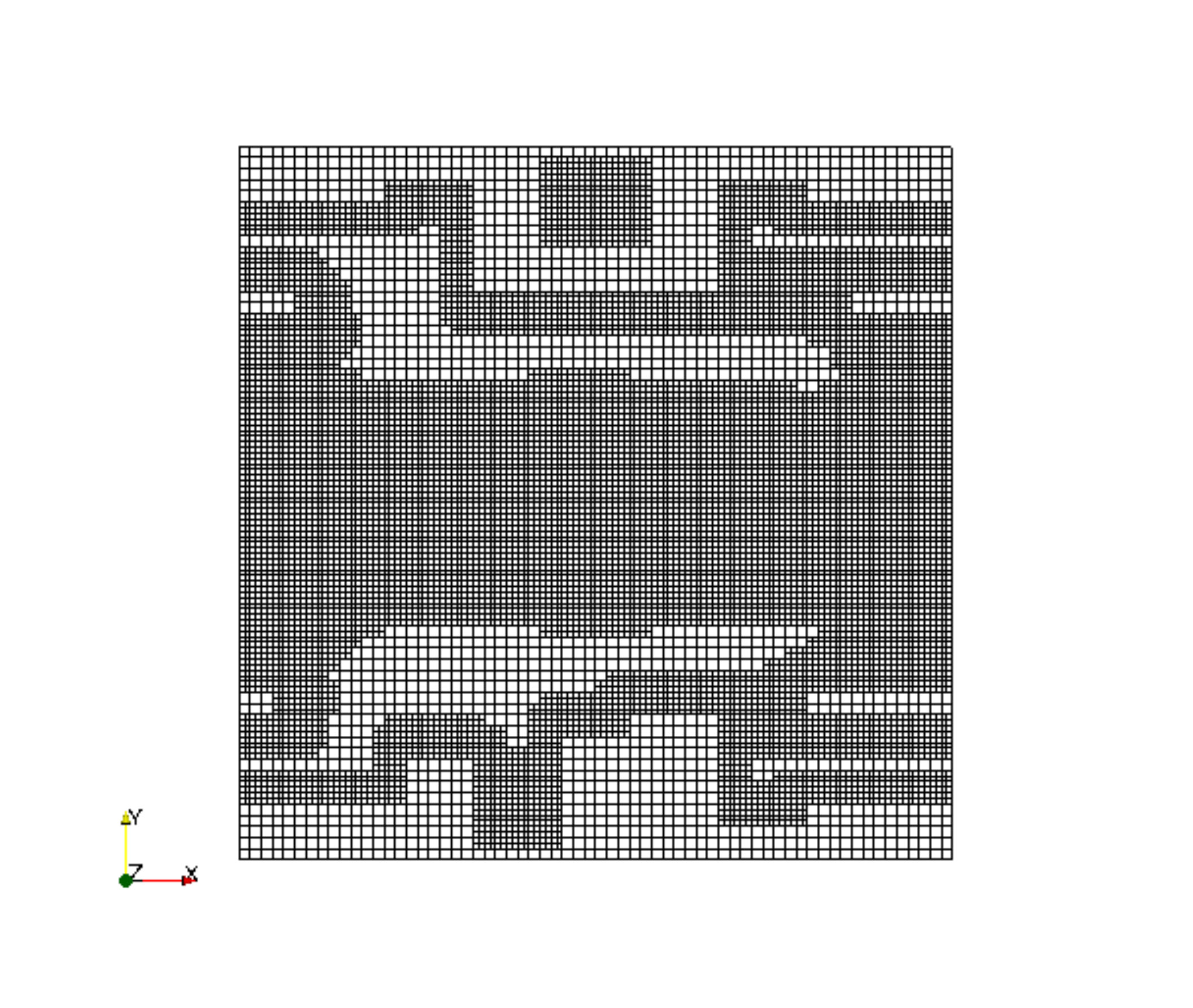}}
	}\\
	{
	\subfloat[ref. \# 8: \quad ${\mathcal{K}}_h({\mathds{M}}_{{\rm \bf BULK}}(0.4))$]{
	\includegraphics[width=5.7cm, trim={2cm 2cm 2cm 2cm}, clip]{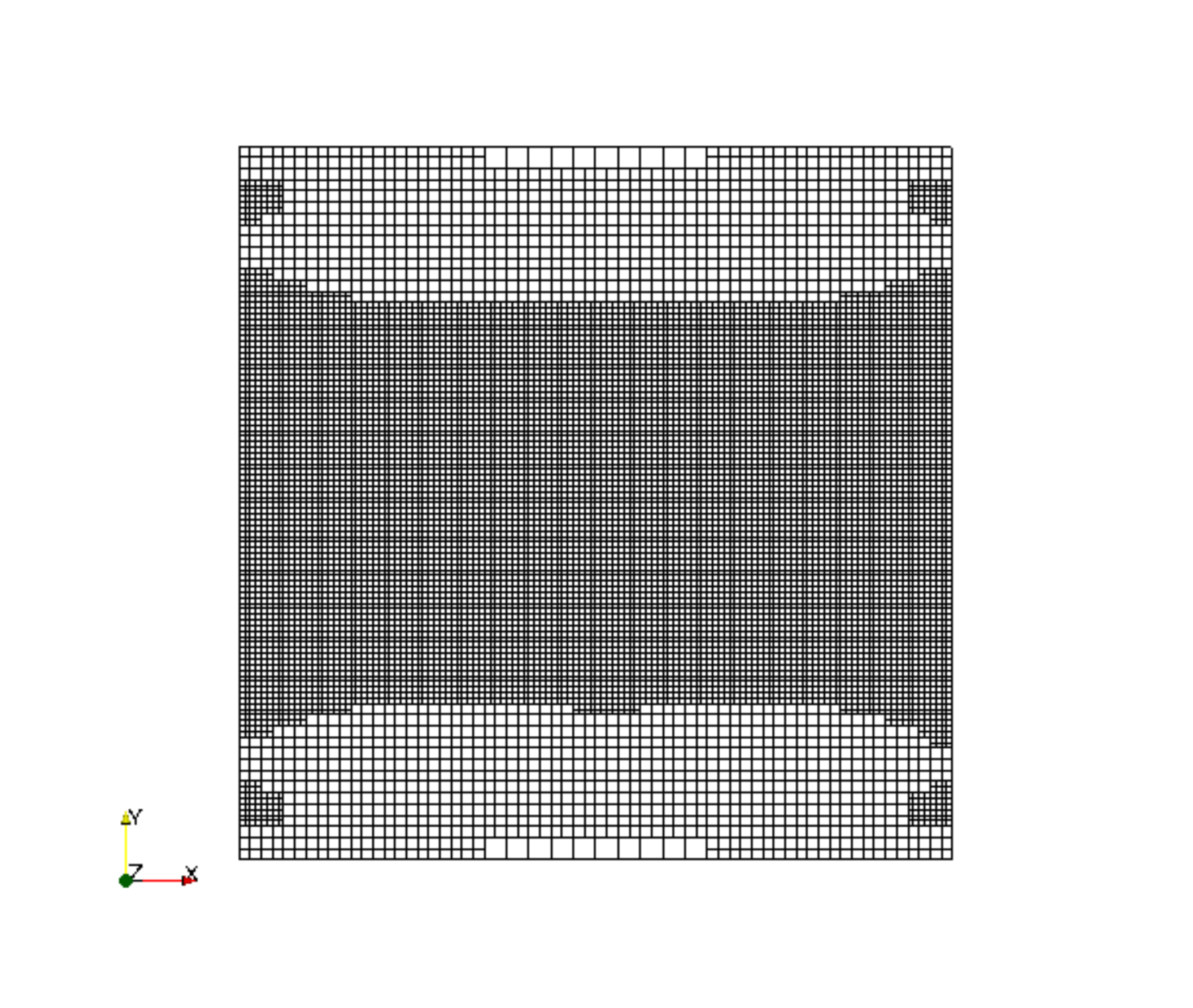}}
	\subfloat[ref. \# 8: \quad ${\mathcal{K}}_h({\mathds{M}}_{{\rm \bf BULK}}(0.2))$]{
	\includegraphics[width=5.7cm,  trim={2cm 2cm 2cm 2cm}, clip]{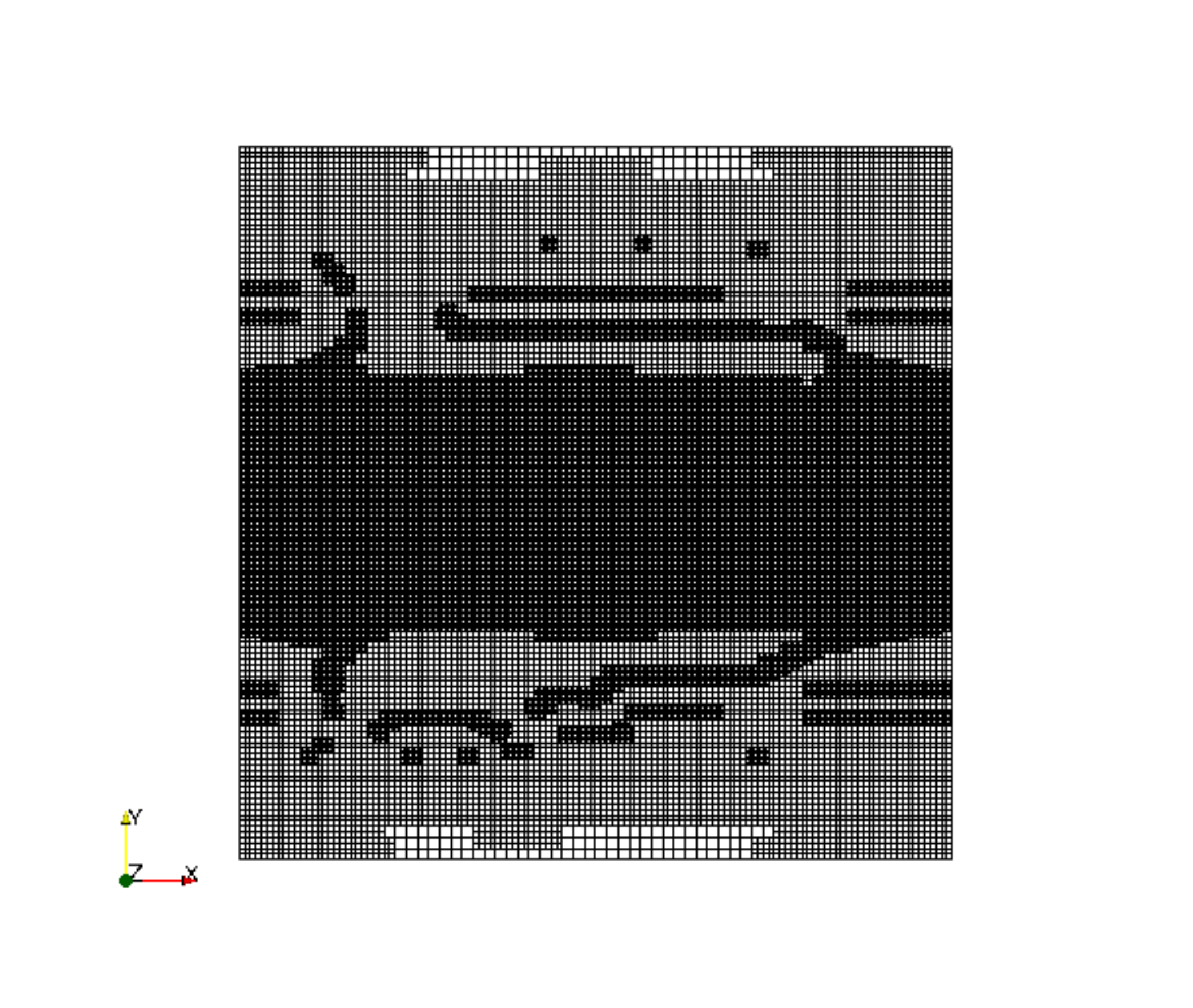}}
	}
	\caption{\small {\em Ex. \ref{ex:unit-domain-example-2}}. 
	Evolution of adaptive meshes obtained with the marking criteria
	${\mathds{M}}_{{\rm \bf BULK}}(0.4)$ and ${\mathds{M}}_{{\rm \bf BULK}}(0.2)$ w.r.t. adaptive ref. steps.}
	\label{fig:unit-domain-example-2-times-v-2-y-3-adaptive-ref}	
\end{figure}

\begin{table}[!t]
 \footnotesize
\centering
\newcolumntype{g}{>{\columncolor{gainsboro}}c} 	
\begin{tabular}{c|c|ccc|gc|c}
\quad \# ref.  \quad & 
\quad  $\| \nabla e \|_\Omega$ \quad   & 	  
\quad $\overline{\rm M}$ \quad &    
\quad $\mdI$ \quad & 	       
\quad $\mfI$ \quad  &  
\quad $\Ieff (\overline{\rm M})$ \quad & 
\quad $\Ieff (\overline{\rm \eta})$ \quad & 
e.o.c. \\
\midrule
   3 &   7.9113e-03 &   8.9838e-03 &   8.1878e-03 &   3.5362e-03 &       1.1356 &       8.7164 &   0.9252  \\
   5 &   5.5188e-04 &   7.5068e-04 &   7.3580e-04 &   6.6131e-05 &       1.3602 &       9.9802 &   2.9781  \\
   7 &   4.8373e-05 &   6.1113e-05 &   5.9155e-05 &   8.7003e-06 &       1.2634 &      10.0270 &   2.4156  \\
   9 &   6.1176e-06 &   8.6725e-06 &   8.4757e-06 &   8.7451e-07 &       1.4176 &      10.3944 &   1.6446   \\
  11 &   6.1657e-07 &   6.3268e-07 &   6.2654e-07 &   2.7268e-08 &       1.0261 &      10.7801 &   2.5543  \\
\end{tabular}
\caption{\small {\em Ex. \ref{ex:unit-domain-example-2}}. Error, majorant (with dual and reliability terms), 
efficiency indices, and e.o.c. w.r.t. adaptive ref. steps with the marking ${\mathds{M}}_{\rm BULK}({0.2})$.}
\label{tab:unit-domain-example-2-error-majorant-v-2-y-3-adaptive-ref}
\end{table}
\begin{table}[!t]
 \footnotesize
\centering
\newcolumntype{g}{>{\columncolor{gainsboro}}c} 	
\begin{tabular}{c|cc|cg|cg|cgc}
\parbox[c]{0.7cm}{\#ref.}  & 
\parbox[c]{1.2cm}{\#d.o.f.($u_h$)} &  
\parbox[c]{1.3cm}{\#d.o.f.($\flux_h$)} &  
\parbox[c]{1.3cm}{$t_{\rm as}(u_h)$} & 
\parbox[c]{1.3cm}{$t_{\rm as}(\flux_h)$} & 
\parbox[c]{1.3cm}{$t_{\rm sol}(u_h)$} & 
\parbox[c]{1.3cm}{$t_{\rm sol}(\flux_h)$} &
\parbox[c]{1.3cm}{$t_{\rm e/w}(\| \nabla e \|)$} & 
\parbox[c]{1.3cm}{$t_{\rm e/w}(\overline{\rm M})$} & 
\parbox[c]{1.3cm}{$t_{\rm e/w}(\overline{\eta})$} \\
\midrule
   3 &         36 &         16 &     0.0055 &     0.0021 &           0.0001 &           0.0002 &       0.0032 &       0.0146 &       0.0170 \\
   5 &        305 &         16 &     0.0839 &     0.0023 &           0.0008 &           0.0002 &       0.1223 &       0.1920 &       0.2382 \\
   7 &       3224 &         16 &     1.4683 &     0.0035 &           0.0581 &           0.0002 &       1.4490 &       2.8517 &       2.8234 \\
   9 &      38276 &         16 &    27.1005 &     0.0021 &           1.9923 &           0.0002 &      22.0243 &      30.7559 &      38.1258 \\
  11 &     396360 &         49 &  3153.3647 &     0.0495 &          73.2963 &           0.0017 &     218.8799 &     328.0449 &     410.5585 \\
\end{tabular}
\caption{\small {\em Ex. \ref{ex:unit-domain-example-2}}. 
Time for assembling and solving the systems that generate $u_h$ and $\flux_h$, 
time for e/w evaluation of error, majorant, and residual error estimator 
w.r.t. adaptive ref. steps with the marking ${\mathds{M}}_{\rm BULK}({0.2})$.}
\label{tab:unit-domain-example-2-times-v-2-y-3-adaptive-ref}
\end{table}
%
%
For an adaptive refinement strategy, we combine the THB-splines 
\cite{LMR:Kraft1997, LMR:Vuongetall2011,LMR:GiannelliJuttlerSpeleers2012}, which support local 
refinement, and functional error estimates. We use bulk marking criterion with 
parameter $\theta = 0.4$. Let us start with the following setting: $u_h \in S^{2, 2}_h$, where 
$S^{2, 2}_h$ is generated by THB-splines,
and $\flux_h \in S^{3, 3}_{8h} \oplus S^{3, 3}_{8h}$, where 
$S^{3, 3}_{8h} \oplus S^{3, 3}_{8h}$ is generated by the basis of THB-splines as well.
Overall, $N_{\rm ref} = 11$ refinements  are executed to obtain the error 
illustrated in Table \ref{tab:unit-domain-example-2-error-majorant-v-2-y-3-adaptive-ref}. The time spent to 
generate $u_h$ and $\flux_h$ is illustrated in Table 
\ref{tab:unit-domain-example-2-times-v-2-y-3-adaptive-ref}. By using a mesh that is up to $8$ times coarser 
than the one for $u_h$, we manage to spare computational time for reconstructing 
the optimal $\flux_h$ and speed up the overall reconstruction of majorant. In the current configuration, 
we obtain the following ratios: 
$$
\tfrac{t_{\rm as}(u_h)}{t_{\rm as}(\flux_h)} \approx \tfrac{3153.3647}{0.0495} \approx 63704 \quad 
\mbox{and} \quad 
\tfrac{t_{\rm sol}(u_h)}{t_{\rm sol}(\flux_h)} \approx \tfrac{73.2963}{0.0017} \approx 43115. 
$$
The comparison of meshes obtained while refining with different parameters can be found on 
Figure \ref{fig:unit-domain-example-2-times-v-2-y-3-adaptive-ref}, i.e., $\theta = 0.4$ (left column) and 
$\theta = 0.2$ (right column).
It is obvious from the plots that the smaller bulk parameter $\theta$ is, the higher the percentage 
of refined elements in the mesh is. 

\end{example}


\begin{example}
\label{ex:unit-domain-kleis-tomar-paper}
\rm 

Next, we consider an example with the exact solution, such that its gradient growth is controlled by the 
parameters. This way, we can study properties of the majorant on the subdomains of $\Omega$, where 
$u_h$ has fast-growing gradients. Namely, we let $\Omega$ be a unit square, and let the exact solution 
and RHS be chosen as follows:
\begin{alignat*}{3}
u 	& = \sin (k_1\,\pi\,x_1)\, \sin(k_2\,\pi\,x_2)  			& \quad \mbox{in} \quad \Omega, \\
f 	& = (k_1^2 + k_2^2) \, \pi^2\, \sin(k_1\,\pi\,x_1)\, \sin(k_2 \, \pi \, x_2) 	& \quad \mbox{in} \quad \Omega, \\
u_D 	& = 0								         & \quad \mbox{on} \quad \Gamma.
\end{alignat*}

\begin{figure}[!t]
	\centering
	\subfloat[$k_1 = k_2 = 1$]{
	\includegraphics[scale=0.7]{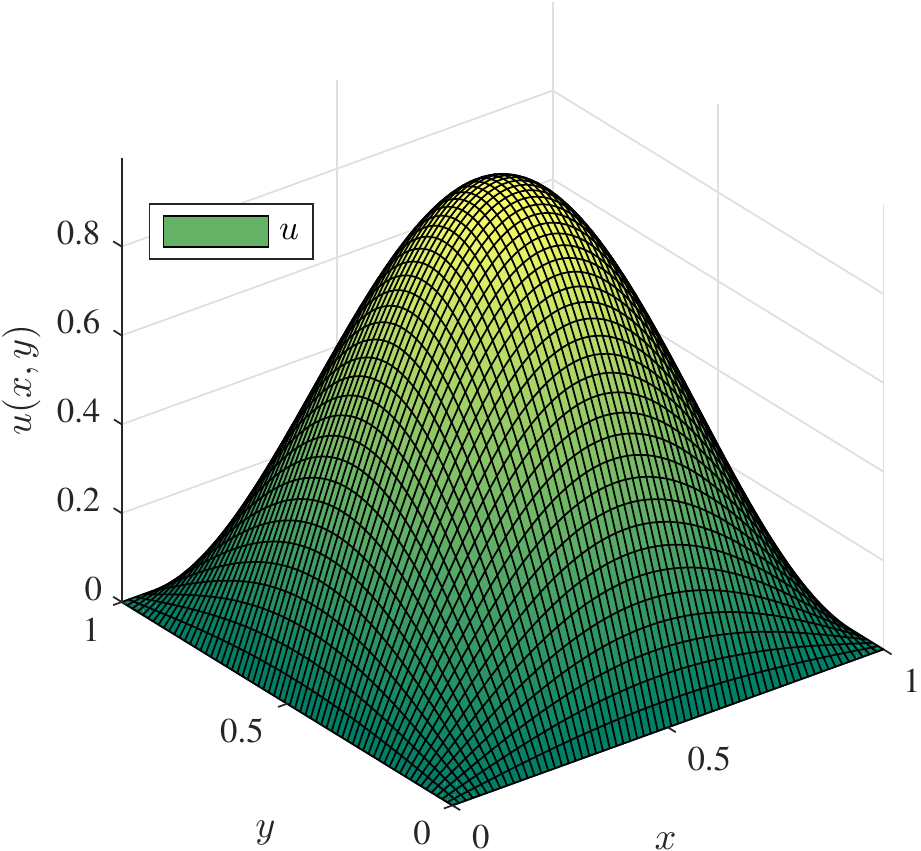}
	\label{fig:example-4-exact-solution-a}}
	\quad
	\subfloat[$k_1 = 6, k_2 = 3$]{
	\includegraphics[scale=0.7]{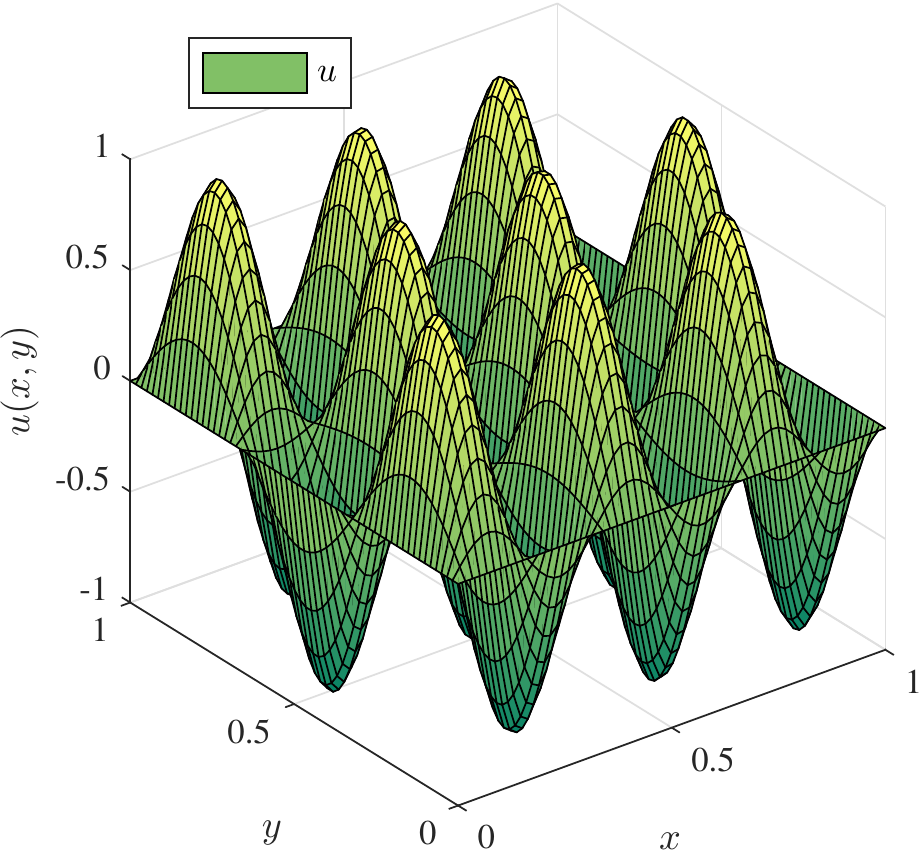}
	\label{fig:example-4-exact-solution-b}}
	\caption{\small {\em Ex. \ref{ex:unit-domain-kleis-tomar-paper}}. 
	Exact solution $u = \sin (k_1\,\pi\,x_1)\, \sin(k_2\,\pi\,x_2)$.}
	\label{fig:example-4-exact-solution}
\end{figure}

\begin{table}[!t]
 \footnotesize
\centering
\newcolumntype{g}{>{\columncolor{gainsboro}}c} 	
\begin{tabular}{c|c|ccc|gc|c}
\quad \# ref. \quad & 
\quad  $\| \nabla e \|_\Omega$ \qquad   & 	  
\quad \quad \quad \quad $\overline{\rm M}$ \quad \quad \quad \quad &    
\quad \quad   $\mdI$ \qquad \quad & 	       
\quad \qquad $\mfI$ \qquad \quad  &  
\qquad $\Ieff (\overline{\rm M})$ \qquad & 
\qquad $\Ieff (\overline{\rm \eta})$ \qquad & 
\qquad e.o.c. \qquad \\
\midrule
   2 &   5.5286e-02 &   6.3291e-02 &   5.7322e-02 &   2.6518e-02 &       1.1448 &      10.3894 &   3.9940  \\
   4 &   3.2077e-03 &   4.0140e-03 &   3.4919e-03 &   2.3195e-03 &       1.2514 &      10.9176 &   2.3839 \\
  5 &   7.9894e-04 &   1.4534e-03 &   1.4273e-03 &   1.1597e-04 &       1.8191 &      10.9451 &   2.1856  \\
   6 &   1.9955e-04 &   1.2390e-03 &   1.1931e-03 &   2.0405e-04 &       6.2091 &      10.9521 &   2.0914 \\
   8 &   1.2468e-05 &   9.8611e-05 &   3.7673e-05 &   2.7074e-04 &       7.9091 &      10.9543 &   2.0226 \\
  10 &   7.7924e-07 &   8.4668e-07 &   7.7970e-07 &   2.9758e-07 &       1.0865 &      10.9544 &   2.0056 \\
\end{tabular}
\caption{\small {\em Ex. \ref{ex:unit-domain-kleis-tomar-paper}}, $k_1 = k_2 = 1$.
Error, majorant (with dual and reliability terms), 
efficiency indices, and e.o.c. w.r.t. unif. ref. steps.}
\label{tab:unit-domain-example-2-error-majorant-v-2-y-5-uniform-ref}
\end{table}
\begin{table}[!t]
\footnotesize
\centering
\newcolumntype{g}{>{\columncolor{gainsboro}}c} 	
\begin{tabular}{c|cc|cg|cg|cg}
\# ref.   & 
\# d.o.f.($u_h$) &  \# d.o.f.($\flux_h$) &  
\; $t_{\rm as}(u_h)$ \; & 
\; $t_{\rm as}(\flux_h)$ \; & 
\; $t_{\rm sol, dir}(u_h)$ \; & 
\; $t_{\rm sol, dir}(\flux_h)$ \; &
\; $t_{\rm sol, iter}(u_h)$ \; & 
\; $t_{\rm sol, iter}(\flux_h)$ \; \\
%
\midrule
   2 &         36 &         36 &     0.0007 &     0.0013 &           0.0001 &           0.0014 &    0.0000 &           0.0010   \\
   4 &        324 &         36 &     0.0091 &     0.0016 &           0.0013 &           0.0010 &       0.0001 &           0.0007  \\
   5 &       1156 &         36 &     0.0289 &     0.0015 &           0.0057 &           0.0008 &       0.0598 &       0.1015 \\
   6 &       4356 &         36 &     0.0723 &     0.0017 &           0.0342 &           0.0007 &        0.0067 &           0.0003 \\
   8 &      66564 &         81 &     1.5561 &     0.0141 &           3.1404 &           0.0036 &       0.6299 &           0.0045 \\
  10 &    1052676 &        441 &    20.2540 &     0.2224 &         166.5165 &           0.1367 &       33.6298 &           0.1121 \\
\end{tabular}
\caption{\small {\em Ex. \ref{ex:unit-domain-kleis-tomar-paper}}, $k_1 = k_2 = 1$.
Time for assembling and solving the systems that generate $u_h$ and $\flux_h$ (with direct and iterative solvers).}
\label{tab:unit-domain-example-2-times-v-2-y-5-uniform-ref}
\end{table}
%
First, let $k_1 = k_2 = 1$. For such parameters, the exact solution is illustrated in Figure 
\ref{fig:example-4-exact-solution-a}. The function $u_h$ is approximated by $S_h^{2, 2}$, whereas 
$\flux_h \in S_{6h}^{5, 5} \oplus S_{6h}^{5, 5}$, and $N_{\rm ref} = 11$ unif. ref. steps are considered. The 
resulting performance of majorant is presented in Table 
\ref{tab:unit-domain-example-2-error-majorant-v-2-y-5-uniform-ref}. At the same time, the computational effort spent 
on $\flux_h$-reconstruction is several times lower than for $u_h$, i.e., 
$$
\tfrac{t_{\rm as}(u_h)}{t_{\rm as}(\flux_h)} \approx \tfrac{20.2540}{0.2224} \approx 91
\quad \mbox{and} \quad 
\tfrac{t_{\rm sol}(u_h)}{t_{\rm sol}(\flux_h)} \approx \tfrac{166.5165}{0.1367} \approx 1218,$$  
which can be observed from Table \ref{tab:unit-domain-example-2-times-v-2-y-5-uniform-ref}.

Let us assume again that the exact solution is not given a priori. Then, reconstruction of the minorant can help to evaluate 
the reliability of the error majorant. The performance of the latter one is illustrated in 
Table \ref{tab:unit-domain-example-3-error-majorant-minorant-eta-v-2-y-6-w-7-adaptive-ref}, whereas 
Table \ref{tab:unit-domain-example-2-times-v-2-y-6-w-7-adaptive-ref} provides the comparison of the time one 
spends on the reconstruction of $u_h$ as well as $\flux_h$ and $w_h$. It is easy to observe from the column providing 
the ratios ${\overline{\rm M}} / {\underline{\rm M}}$ that values of $\overline{\rm M}$ are quite reliable and 
efficient. Moreover, Table \ref{tab:unit-domain-example-2-times-v-2-y-6-w-7-adaptive-ref} confirms that the expenses 
on the reconstruction of both $\overline{\rm M}$ and $\underline{\rm M}$ are rather modest in comparison to the 
expenses on the primal variable. 
 
\begin{table}[!t]
\footnotesize
\centering
\newcolumntype{g}{>{\columncolor{gainsboro}}c} 	
\begin{tabular}{c|c|cgl|cg|g|c}
\quad \# ref. \qquad & 
\quad  $\| \nabla e \|_\Omega$ \qquad & 	  
\quad $\overline{\rm M}$ \qquad &    
\quad $\Ieff (\overline{\rm M})$ \qquad & 
\quad $\underline{\rm M}$ \qquad &    
\quad $\Ieff (\underline{\rm M})$ \qquad & 
\quad ${\overline{\rm M}} / {\underline{\rm M}}$ \qquad & 
\quad e.o.c. \qquad \\
\midrule
   3 &     5.5285e-02 &   5.7554e-02 &       1.0411 &   5.5115e-02 &       0.9969 &       1.0443 &   3.9940 \\
   4 &     1.3024e-02 &   1.6545e-02 &       1.2704 &   1.3012e-02 &       0.9991 &       1.2715 &   2.8302 \\
   6 &     1.6936e-03 &   2.2063e-03 &       1.3027 &   1.6885e-03 &       0.9970 &       1.3067 &   1.6060 \\
   7 &     4.5250e-04 &   1.2058e-03 &       2.6648 &   4.3221e-04 &       0.9552 &       2.7898 &   2.1530 \\
   8 &     1.6513e-04 &   2.2145e-04 &       1.3411 &   1.6487e-04 &       0.9984 &       1.3432 &   1.9801 \\
   9 &     5.4927e-05 &   1.2373e-04 &       2.2527 &   5.4114e-05 &       0.9852 &       2.2865 &   3.8077 \\
\end{tabular}
\caption{{\em Ex. \ref{ex:unit-domain-kleis-tomar-paper}}. 
Error, majorant, minorant, residual based error indicator 
with corresponding efficiency indices, and e.o.c. for 
$\flux_h \in S_{6h}^{6, 6} \oplus S_{6h}^{6, 6}$ and $w_h \in S_{6h}^{7, 7}$ 
w.r.t. adaptive ref. steps.}
\label{tab:unit-domain-example-3-error-majorant-minorant-eta-v-2-y-6-w-7-adaptive-ref}
\end{table}

\begin{table}[!t]
\footnotesize
\centering
\newcolumntype{g}{>{\columncolor{gainsboro}}c} 	
\begin{tabular}{c|ccc|ccc|ccc}
\parbox[c]{0.7cm}{\#ref.} & 
\parbox[c]{1.4cm}{\#d.o.f.($u_h$)} &  
\parbox[c]{1.4cm}{\#d.o.f.($\flux_h$)} &   
\parbox[c]{1.4cm}{\#d.o.f.($w_h$)} &
\parbox[c]{1.3cm}{$t_{\rm as}(u_h)$} & 
\parbox[c]{1.3cm}{$t_{\rm as}(\flux_h)$} & 
\parbox[c]{1.3cm}{$t_{\rm as}(w_h)$} & 
\parbox[c]{1.3cm}{$t_{\rm sol}(u_h)$} & 
\parbox[c]{1.3cm}{$t_{\rm sol}(\flux_h)$} & 
\parbox[c]{1.3cm}{$t_{\rm sol}(w_h)$} \\ 
\midrule
   1 &          9 &         49 &         64 &   2.14e-03 &   3.80e-02 &   2.88e-02 &         3.60e-05 &         1.11e-03 &         1.70e-05 \\
   3 &         36 &         49 &         64 &   8.72e-03 &   3.96e-02 &   2.61e-02 &         5.50e-05 &         2.23e-03 &         1.29e-04 \\
   4 &        100 &         49 &         64 &   2.89e-02 &   3.10e-02 &   2.53e-02 &         1.40e-05 &         2.50e-03 &         3.60e-05 \\
   6 &        952 &         49 &         64 &   5.65e-01 &   3.56e-02 &   2.25e-02 &         2.74e-03 &         2.03e-03 &         5.80e-05 \\
   7 &       3244 &         49 &         64 &   1.34e+00 &   1.93e-02 &   2.02e-02 &         8.11e-03 &         1.21e-03 &         1.02e-04 \\
   8 &       8980 &         64 &         81 &   5.68e+00 &   1.19e-01 &   6.92e-02 &         1.01e-01 &         3.46e-03 &         7.00e-05 \\
   9 &      16009 &         64 &         81 &   5.88e+00 &   8.24e-02 &   6.07e-02 &         2.34e-01 &         2.78e-03 &         9.70e-05 \\
   \midrule
    &
    \multicolumn{3}{c|}{ }  &    
    \multicolumn{3}{c|}{96.86 \quad : \quad 1.35 \qquad : \qquad 1 \quad } &
    \multicolumn{3}{c}{2412.37 \quad : \quad 28.66 \qquad : \qquad 1 \qquad } 
   \end{tabular}
\caption{{\em Ex. \ref{ex:unit-domain-kleis-tomar-paper}}. 
Time for solving the systems that generate $u_h$, $\flux_h$, and $w_h$, with direct and iterative methods for 
$\flux_h \in S_{6h}^{6, 6} \oplus S_{6h}^{6, 6}$ and $w_h \in S_{6h}^{7, 7}$ w.r.t. adaptive ref. steps.}
\label{tab:unit-domain-example-2-times-v-2-y-6-w-7-adaptive-ref}
\end{table}
%

Let us consider now a more complicated case with $k_1 = 6$ and $k_2 = 3$ (see Figure 
\ref{fig:example-4-exact-solution-b}).
For an efficient  flux reconstruction, we apply the same strategy as discussed in {\em Ex. \ref{ex:unit-domain-example-2}}, 
i.e., we increase the degree of B-splines used for the space approximating $\flux_h$, but at the same time, 
we use $M = 8$ times coarser mesh, i.e., $S_{8h}^{9, 9} \oplus S_{8h}^{9, 9}$.
First, we analyse the results obtained by global refinement; they are presented in Tables 
\ref{tab:unit-domain-example-3-error-majorant-v-2-y-9-uniform-ref} and 
\ref{tab:unit-domain-example-3-times-v-2-y-9-uniform-ref}. We consider $N_{\rm ref} = 8$ unif. ref. steps (starting 
from a rather fine initial mesh generated by $N_{\rm ref, 0} = 4$ initial ref. steps of original geometry and the 
basis assigned for it). In column six of Table \ref{tab:unit-domain-example-3-error-majorant-v-2-y-9-uniform-ref}, 
one can see that $\Ieff$ takes values up to 6.2091 but decreases back to 1.0865 once we start refining the basis for the 
variable $\flux_h$ as well. In particular, at the refinements steps 5 and 6, the initial mesh of 36 d.o.f. or $\flux_h$ 
becomes relatively coarse in comparison to the basis for $u_h$ and must be refined in order to obtain efficient 
values of $\overline{\rm M}$. Concerning the time spent on assembling and solving the systems in 
\eqref{eq:system-uh} and \eqref{eq:system-fluxh}, we obtain the following ratios taken from 
Table \ref{tab:unit-domain-example-3-times-v-2-y-9-uniform-ref}, namely,
$$\tfrac{t_{\rm as}(u_h)}{t_{\rm as}(\flux_h)} \approx \tfrac{17.3623}{3.2302} \approx  2 
\quad \mbox{and} \quad
\tfrac{t_{\rm sol}(u_h)}{t_{\rm sol}(\flux_h)} \approx \tfrac{144.7056}{1.3482} \approx 107.$$ 
%
%
In the case of adaptive refinement, we also use the space $S_{7h}^{9, 9} \oplus S_{7h}^{9, 9}$ generated by 
THB-splines. Let the bulk threshold be defined by parameter $\theta = 0.4$, which causes the refinement of  
approximately $60\%$ of all elements for the primal variable $u_h$. The obtained numerical results are presented in 
Tables \ref{tab:unit-domain-example-3-error-majorant-v-2-y-9-adaptive-ref}--\ref{tab:unit-domain-example-3-times-v-2-y-9-adaptive-ref}. 
\begin{table}[!t]
 \footnotesize
\centering
\newcolumntype{g}{>{\columncolor{gainsboro}}c} 	
\begin{tabular}{c|c|ccc|gc|c}
\quad \# ref.  \quad & 
\quad  $\| \nabla e \|_\Omega$ \qquad   & 	  
\quad \quad \quad \quad $\overline{\rm M}$ \quad \quad \quad \quad &    
\quad \quad   $\mdI$ \qquad \quad & 	       
\quad \qquad $\mfI$ \qquad \quad  &  
\qquad $\Ieff (\overline{\rm M})$ \qquad & 
\qquad $\Ieff (\overline{\rm \eta})$ \qquad & 
\qquad \quad e.o.c. \qquad \quad \\
\midrule
   3 &   3.1030e-02 &   3.1818e-02 &   3.1057e-02 &   3.3805e-03 &       1.0254 &      10.8671 &   2.1371 \\
   5 &   1.9203e-03 &   1.9909e-03 &   1.9303e-03 &   2.6939e-04 &       1.0367 &      10.9490 &   2.0254 \\
   7 &   1.1995e-04 &   1.8194e-04 &   1.2028e-04 &   2.7397e-04 &       1.5168 &      10.9541 &   2.0058 \\
\end{tabular}
\caption{\small {\em Ex. \ref{ex:unit-domain-kleis-tomar-paper}}, $k_1 = 6, k_2 = 3$. 
Error, majorant (with dual and reliability terms), 
efficiency indices, and e.o.c. w.r.t. 
unif. ref. steps.}
\label{tab:unit-domain-example-3-error-majorant-v-2-y-9-uniform-ref}
\end{table}
\begin{table}[!th]
 \footnotesize
\centering
\newcolumntype{g}{>{\columncolor{gainsboro}}c} 	
\begin{tabular}{c|cc|cg|cg|cgc}
\# ref.   & 
\# d.o.f.($u_h$) &  \# d.o.f.($\flux_h$) &  
\; $t_{\rm as}(u_h)$ \; & 
\; $t_{\rm as}(\flux_h)$ \; & 
\; $t_{\rm sol}(u_h)$ \; & 
\; $t_{\rm sol}(\flux_h)$ \; &
$t_{\rm e/w}(\| \nabla e \|)$ & 
$t_{\rm e/w}(\overline{\rm M})$ & 
$t_{\rm e/w}(\overline{\eta})$ \\
\midrule
   1 &        324 &        625 &     0.0053 &     2.9646 &           0.0007 &           0.2622 &       0.0047 &       0.1177 &       0.0158 \\
\rowcolor{gainsboro}   
   3 &       4356 &        625 &     0.0831 &     3.4472 &           0.0396 &           0.6094 &       0.2142 &       0.4602 &       0.3388 \\
   5 &      66564 &        625 &     1.1135 &     2.9721 &           2.3809 &           1.5243 &       2.4025 &       6.4769 &       3.9421 \\
   7 &    1052676 &        625 &    17.3623 &     3.2302 &         144.7056 &           1.3482 &      45.4342 &     102.9160 &      71.1602 \\
\end{tabular}
\caption{\small {\em Ex. \ref{ex:unit-domain-kleis-tomar-paper}}, $k_1 = 6, k_2 = 3$. 
Time for assembling and solving the systems that generate $u_h$ and $\flux_h$ 
as well as the time spent on e/w evaluation of error, majorant, and 
residual error estimator w.r.t. unif. ref. steps.}
\label{tab:unit-domain-example-3-times-v-2-y-9-uniform-ref}
\end{table}
%

\begin{table}[!t]
\footnotesize
\centering
\newcolumntype{g}{>{\columncolor{gainsboro}}c} 	
\begin{tabular}{c|c|ccc|gc|c}
\quad $N_{\rm ref}$ \quad & 
\quad  $\| \nabla e \|_\Omega$ \qquad   & 	  
\quad \quad \quad \quad $\overline{\rm M}$ \quad \quad \quad \quad &    
\quad \quad   $\mdI$ \qquad \quad & 	       
\quad \qquad $\mfI$ \qquad \quad  &  
\qquad $\Ieff (\overline{\rm M})$ \qquad & 
\qquad $\Ieff (\overline{\rm \eta})$  \qquad & 
e.o.c.
\\
\midrule
   3 &   4.2892e-02 &   4.3740e-02 &   4.2910e-02 &   3.6902e-03 &       1.0198 &       9.6688 &   2.1266 \\
   5 &   5.3723e-03 &   5.5714e-03 &   5.4777e-03 &   4.1637e-04 &       1.0371 &       9.9389 &   1.8082 \\
   7 &   6.4564e-04 &   7.2116e-04 &   6.5719e-04 &   2.8420e-04 &       1.1170 &      10.5034 &   2.3521 \\
\end{tabular}
\caption{\small {\em Ex. \ref{ex:unit-domain-kleis-tomar-paper}}, $k_1 = 6, k_2 = 3$. 
Error, majorant (with dual and reliability terms), efficiency indices w.r.t. adaptive ref. steps.}
\label{tab:unit-domain-example-3-error-majorant-v-2-y-9-adaptive-ref}
\end{table}
\begin{table}[!t]
 \footnotesize
\centering
\newcolumntype{g}{>{\columncolor{gainsboro}}c} 	
\begin{tabular}{c|cc|cg|cg|cgc}
\# ref.  & 
\# d.o.f.($u_h$) &  \# d.o.f.($\flux_h$) &  
\; $t_{\rm as}(u_h)$ \; & 
\; $t_{\rm as}(\flux_h)$ \; & 
\; $t_{\rm sol}(u_h)$ \; & 
\; $t_{\rm sol}(\flux_h)$ \; &
$t_{\rm e/w}(\| \nabla e \|)$ & 
$t_{\rm e/w}(\overline{\rm M})$ & 
$t_{\rm e/w}(\overline{\eta})$ \\
\midrule
   1 &        324 &        625 &     0.1125 &    23.3092 &           0.0010 &           0.2418 &       0.0918 &       7.0919 &       0.1938 \\
\rowcolor{gainsboro}
   3 &       3468 &        625 &     1.2291 &    22.0679 &           0.0313 &           0.6313 &       1.5198 &      12.0444 &       2.9476 \\
   5 &      31640 &        625 &    44.5650 &    22.1864 &           0.9953 &           1.2794 &      18.0067 &     107.7677 &      32.6181 \\
   7 &     205060 &        625 &  1135.8096 &    21.3158 &          17.4630 &           1.3050 &     105.0698 &     583.0759 &     190.6939 \\
%
\end{tabular}
\caption{\small {\em Ex. \ref{ex:unit-domain-kleis-tomar-paper}}, $k_1 = 6, k_2 = 3$. Time for assembling and 
solving the systems generating d.o.f. of $u_h$ and $\flux_h$ as well as the time spent on e/w evaluation of error, 
majorant, and residual error estimator w.r.t. adaptive ref. steps.}
\label{tab:unit-domain-example-3-times-v-2-y-9-adaptive-ref}
\end{table}

Let us compare the performance of majorant in the uniform and adaptive refinement strategies. Due 
to the implementation of THB-splines evaluation on G+Smo \cite{gismoweb}, the assembling of matrices both for 
$\flux_h$ and $u_h$ is slower w.r.t. B-splines (compare the third and fourth columns of Table 
\ref{tab:unit-domain-example-3-times-v-2-y-9-adaptive-ref} to the third and fourth columns of Table 
\ref{tab:unit-domain-example-3-times-v-2-y-9-uniform-ref}). For ${\rm d.o.f.}(u_h) \approx 4000$, in the first case 
we spend $t_{\rm as}(u_h) = 0.0813$ secs (second row is highlighted with grey background in Table 
\ref{tab:unit-domain-example-3-times-v-2-y-9-uniform-ref}) in comparison to $t_{\rm as}(u_h) = 3.1285$ secs for 
the THB-splines (row with grey background  in Table 
\ref{tab:unit-domain-example-3-times-v-2-y-9-adaptive-ref}), 
which is about 45 times slower. Moreover, these ratios 
grow as ${\rm d.o.f.}(u_h)$ increases. For the auxiliary variable $\flux_h$, the assembling time for THB-splines is 4--5 
times slower than when using B-splines. A similar increase in time can be observed for the element-wise evaluation of 
the error, majorant, and residual error estimator illustrated in the last three columns of Table
\ref{tab:unit-domain-example-3-times-v-2-y-9-adaptive-ref} (in comparison to Table 
\ref{tab:unit-domain-example-3-times-v-2-y-9-uniform-ref}). This slowdown can be explained by a  
bottleneck of G+Smo library when the evaluation of THB-splines is concerned.

Analogously to the previous example, we demonstrate the evaluation of adaptive meshes for different marking criteria, 
i.e., marking ${\mathds{M}}_{{\rm \bf BULK}}(0.4)$ (left column of Figure 
\ref{fig:unit-domain-example-3-times-v-2-y-3-adaptive-ref}) 
and ${\mathds{M}}_{{\rm \bf BULK}}(0.6)$	(right column of Figure \ref{fig:unit-domain-example-3-times-v-2-y-3-adaptive-ref}). 
It resembles the patterns obtained in \cite[Example 1]{LMR:KleissTomar2015}, however, in the current case, due to the local structure 
of THB-splines, many superfluous d.o.f. are eliminated.
\begin{figure}[!t]
	\centering
	{
	\subfloat[ref. \# 2:  \quad ${\mathcal{K}}_h({\mathds{M}}_{{\rm \bf BULK}}(0.4))$]{
	\includegraphics[width=5.9cm, trim={2cm 2cm 2cm 2cm}, clip]{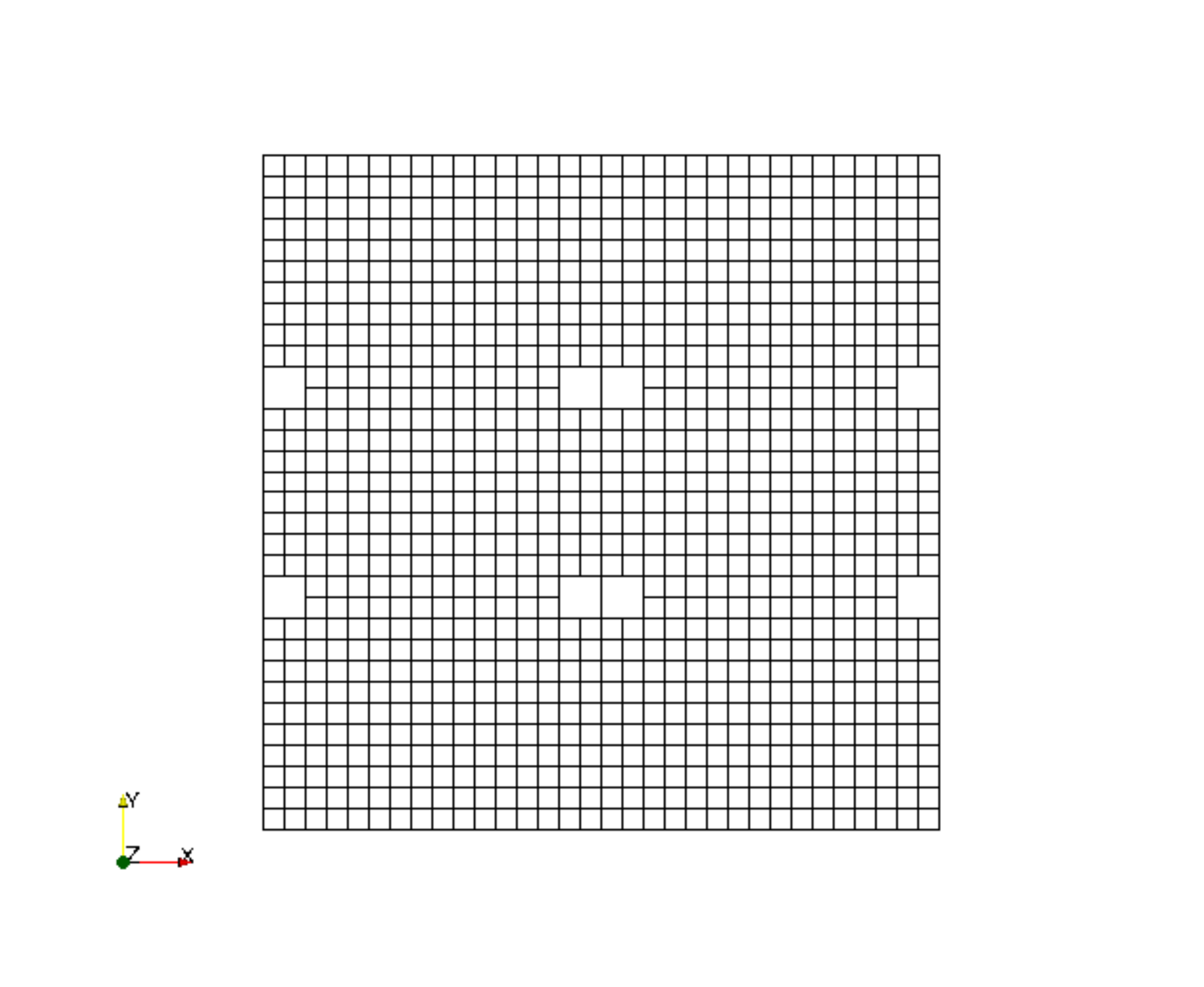}} \qquad 
	\subfloat[ref. \# 2:  \quad ${\mathcal{K}}_h({\mathds{M}}_{{\rm \bf BULK}}(0.6))$]{
	\includegraphics[width=5.9cm, trim={2cm 2cm 2cm 2cm}, clip]{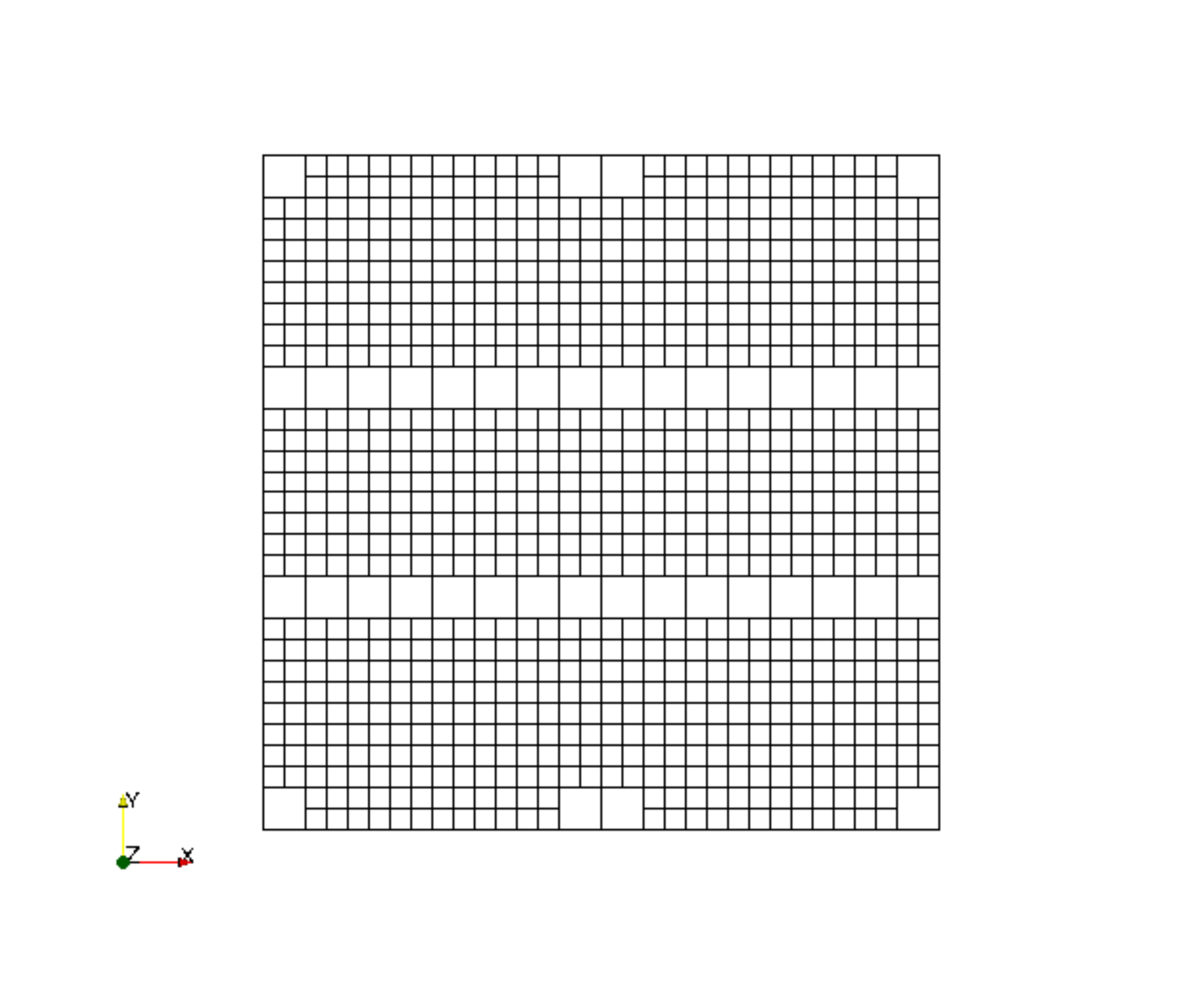}}
	}\\
	{
	\subfloat[ref. \# 3:   \quad ${\mathcal{K}}_h({\mathds{M}}_{{\rm \bf BULK}}(0.4))$]{
	\includegraphics[width=5.9cm, trim={2cm 2cm 2cm 2cm}, clip]{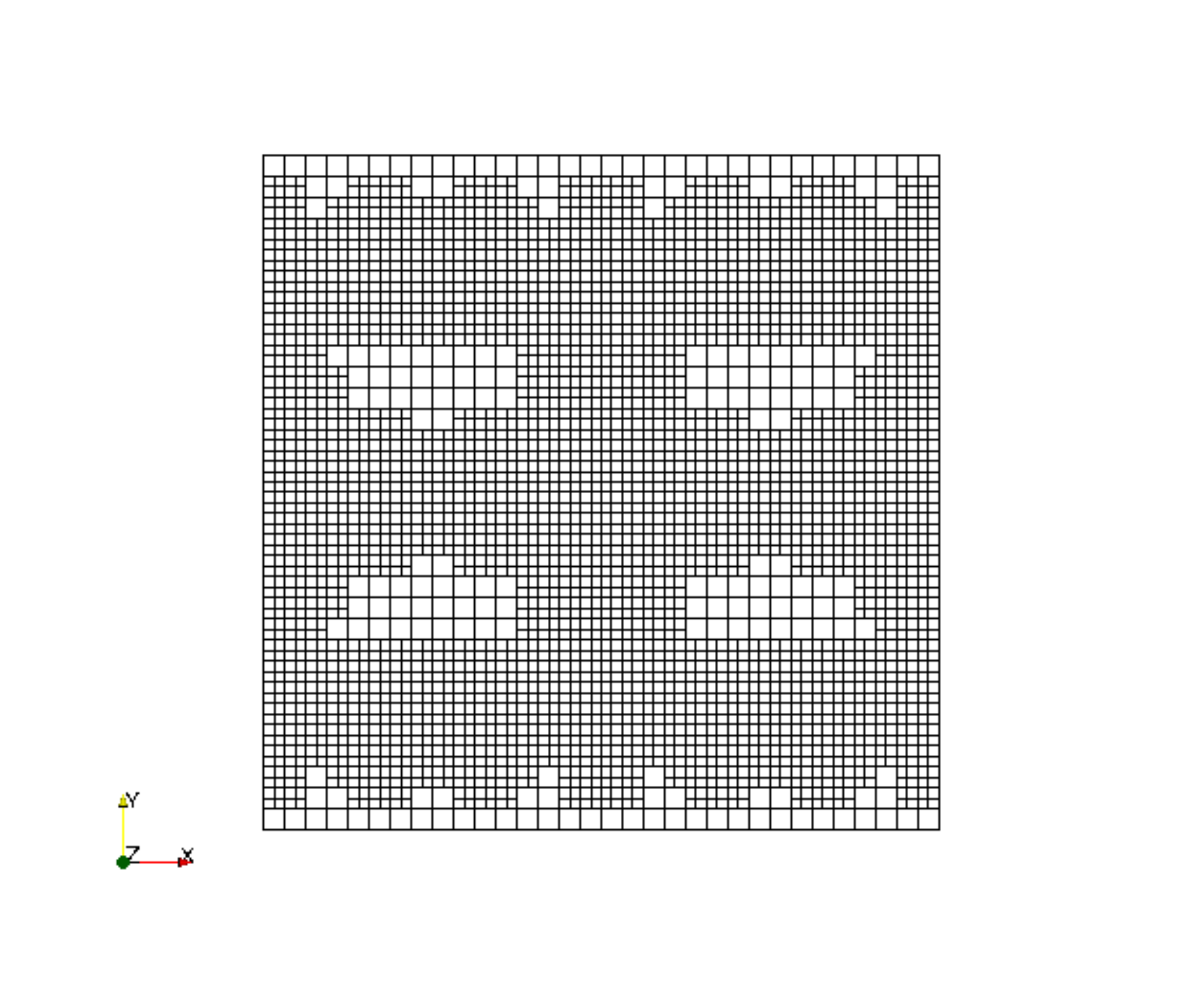}} 
	\quad 
	\subfloat[ref. \# 3:   \quad ${\mathcal{K}}_h({\mathds{M}}_{{\rm \bf BULK}}(0.6))$]{
	\includegraphics[width=5.9cm, trim={2cm 2cm 2cm 2cm}, clip]{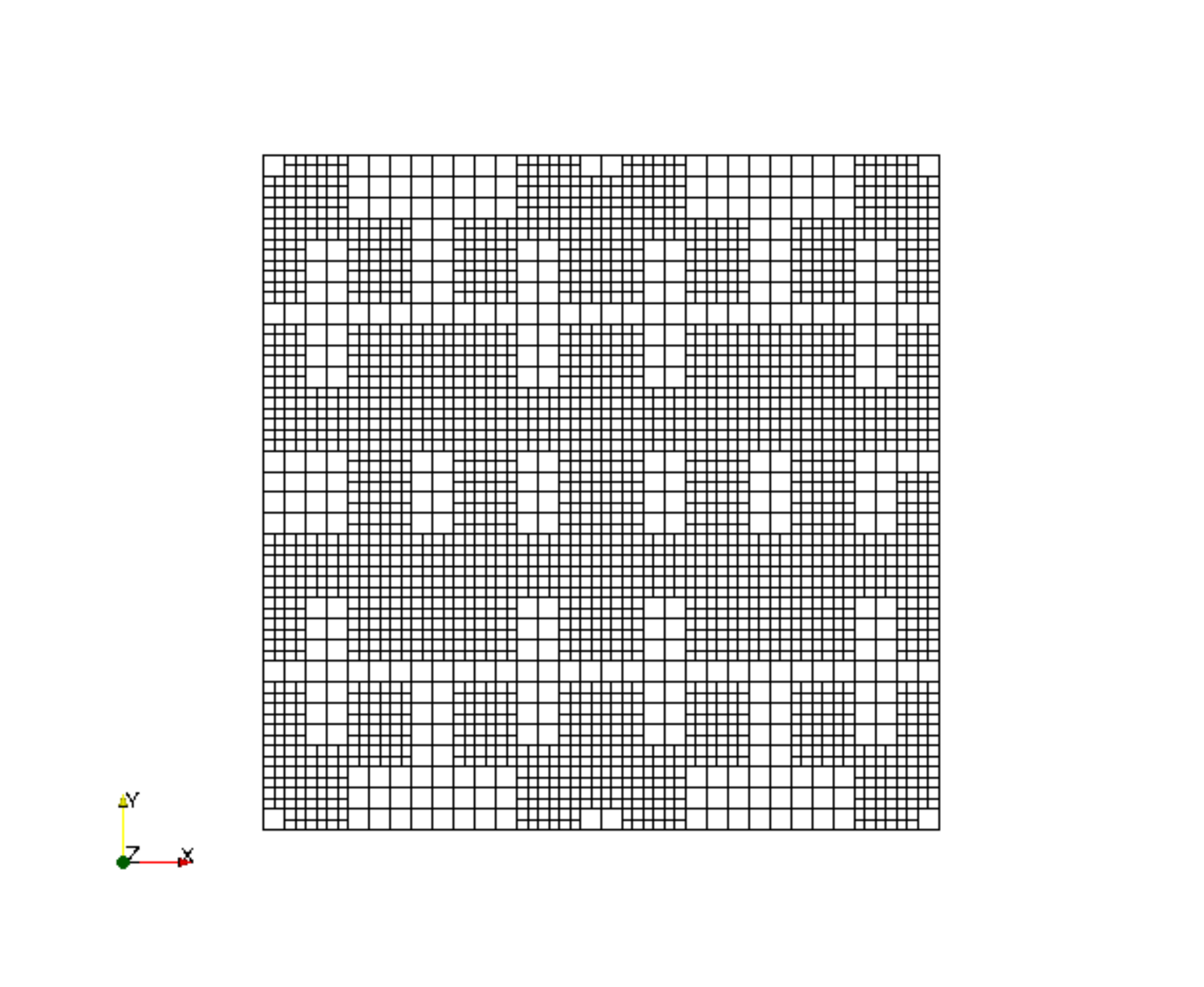}}
	}
	{
	\subfloat[ref. \# 4:  \quad ${\mathcal{K}}_h({\mathds{M}}_{{\rm \bf BULK}}(0.4))$]{
	\includegraphics[width=5.9cm, trim={2cm 2cm 2cm 2cm}, clip]{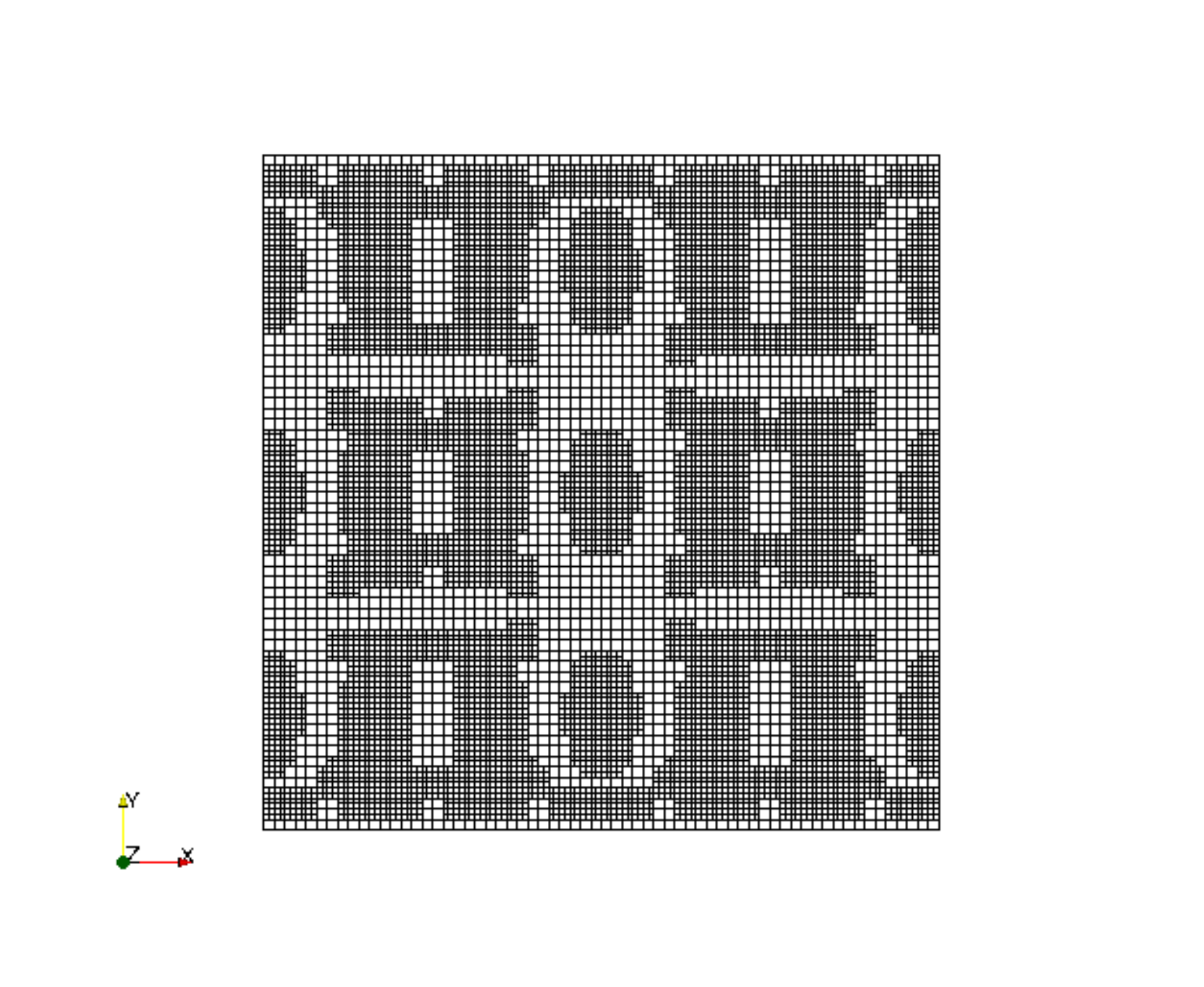}} \qquad 
	\subfloat[ref. \# 4:   \quad ${\mathcal{K}}_h({\mathds{M}}_{{\rm \bf BULK}}(0.6))$]{
	\includegraphics[width=5.9cm, trim={2cm 2cm 2cm 2cm}, clip]{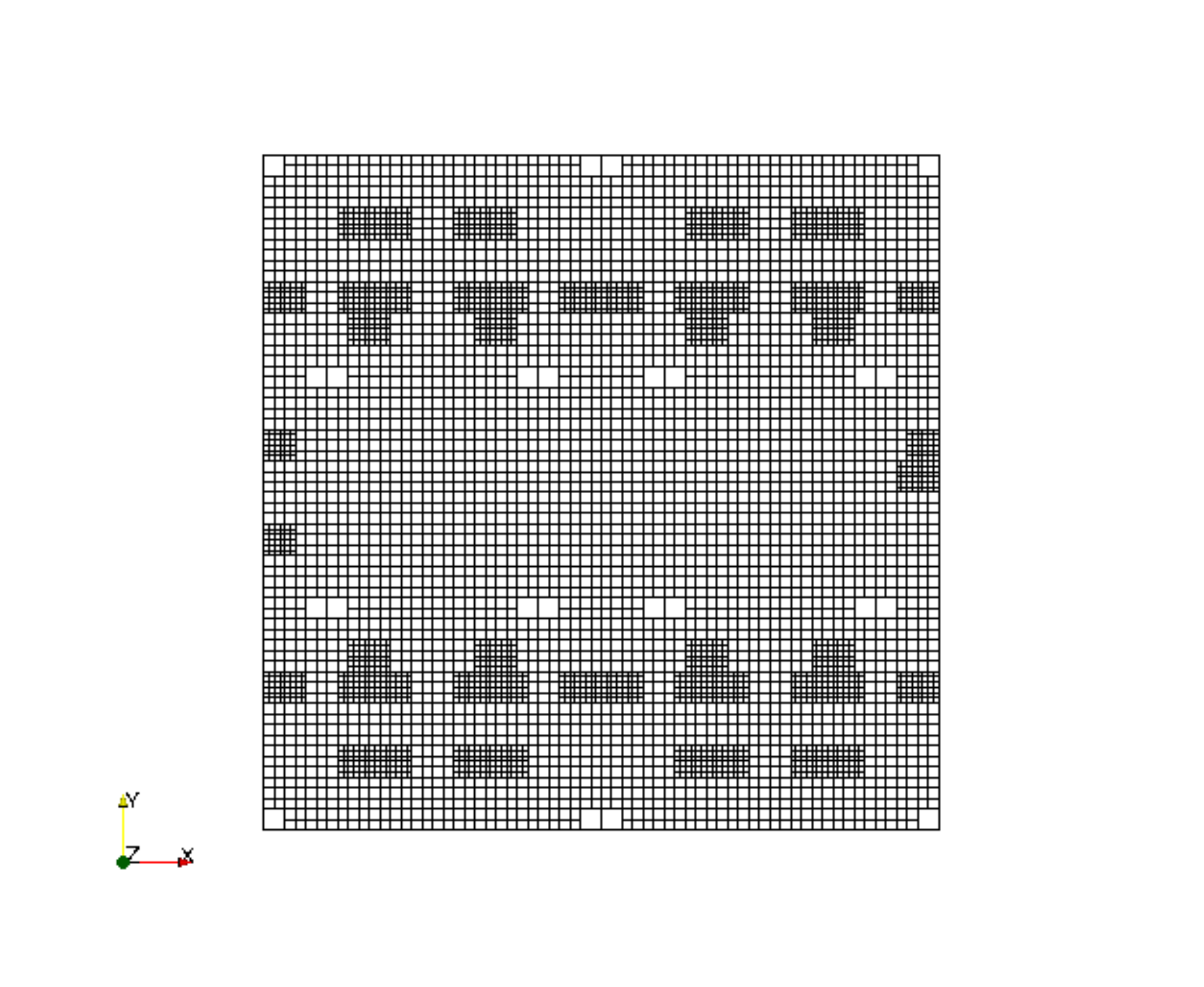}}
	}\\
	{
	\subfloat[ref. \# 5:  \quad ${\mathcal{K}}_h({\mathds{M}}_{{\rm \bf BULK}}(0.4))$]{
	\includegraphics[width=5.9cm, trim={2cm 2cm 2cm 2cm}, clip]{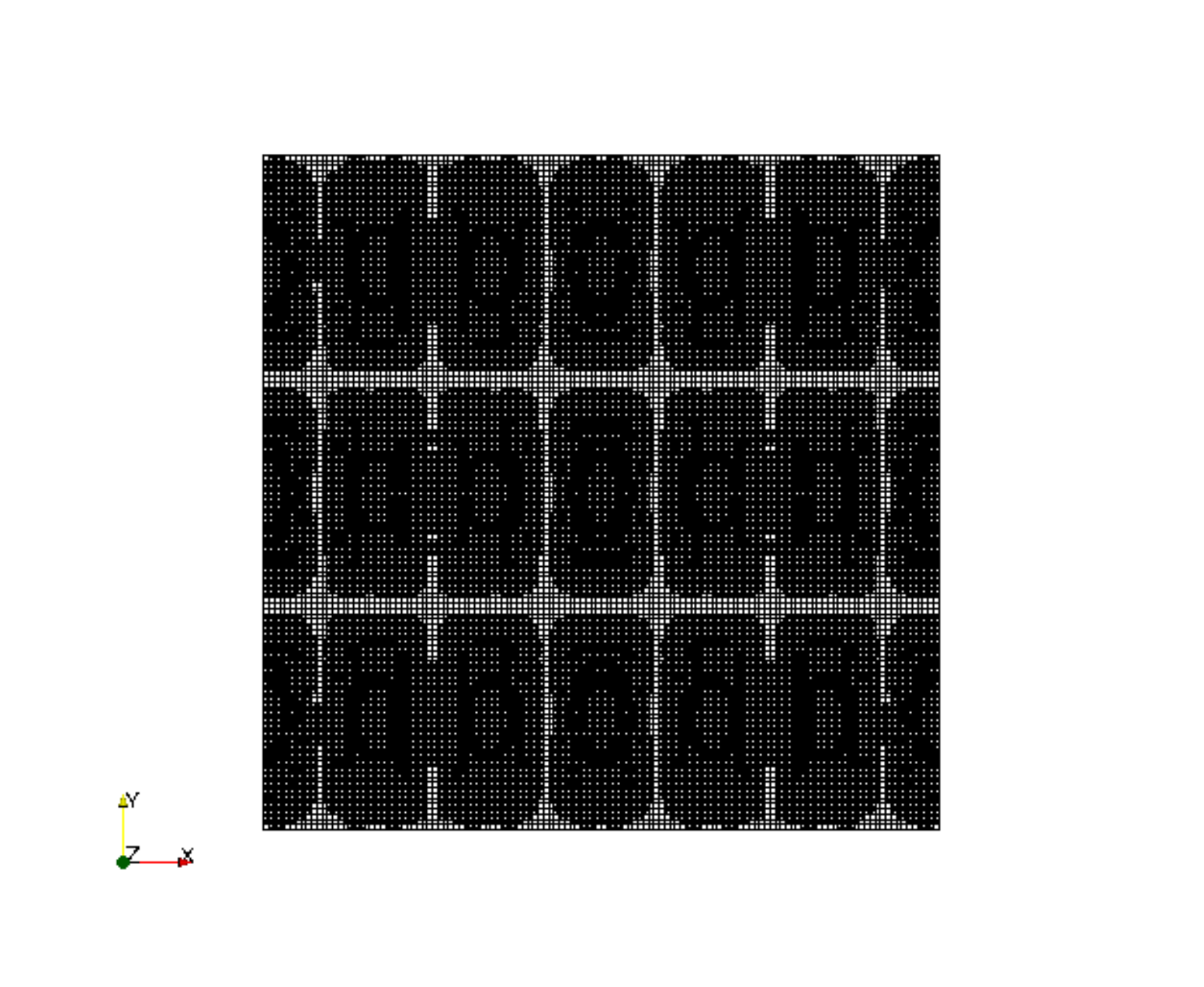}} \qquad 
	\subfloat[ref. \# 5:   \quad ${\mathcal{K}}_h({\mathds{M}}_{{\rm \bf BULK}}(0.6))$]{
	\includegraphics[width=5.9cm, trim={2cm 2cm 2cm 2cm}, clip]{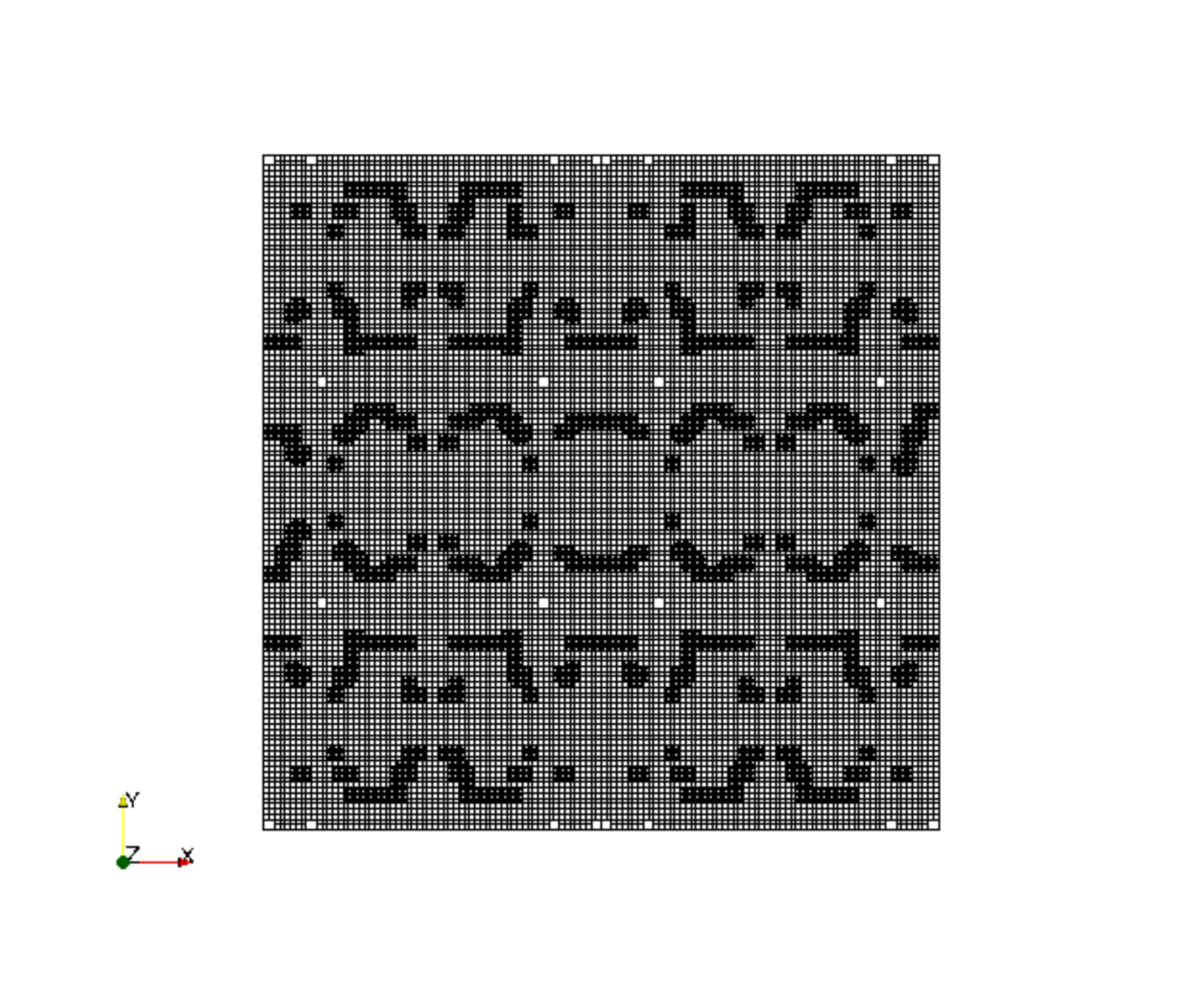}}
	}\\
	\caption{\small {\em Ex. \ref{ex:unit-domain-kleis-tomar-paper}}, 
	$k_1 = 6, k_2 = 3$. 
	Evolution of adaptive meshes obtained with the marking criteria
	${\mathds{M}}_{{\rm \bf BULK}}(0.4)$ (left) and ${\mathds{M}}_{{\rm \bf BULK}}(0.6)$	(right) 
	w.r.t. adaptive ref. steps.}
	\label{fig:unit-domain-example-3-times-v-2-y-3-adaptive-ref}
	\end{figure}
\end{example}

%
\begin{example}
\label{ex:rectangular-domain-kick-jump}
\rm 
Next, we consider an example with a sharp local jump in the exact solution. Let $\Omega := (0, 2) \times (0, 1)$, 
\begin{alignat*}{3}
u 	& = (x_1^2 - 2 \, x_1) \, (x_2^2 - x_2) \, e^{-100 \,|(x_1, x_2) - (1.4, 0.95)|}  & \quad \mbox{in} \quad \Omega,
\label{eq:example-8-exact-solution}
\end{alignat*}
where the jump is located in the point $(x_1, x_2) = (1.4, 0.95)$ (see Figure \ref{fig:example-8-exact-solution}), 
$f$ is calculated by substituting $u$ into \eqref{eq:poisson}, and the Dirichlet BCs are homogenous.
\begin{figure}[!t]
	\centering
	\includegraphics[scale=0.7]{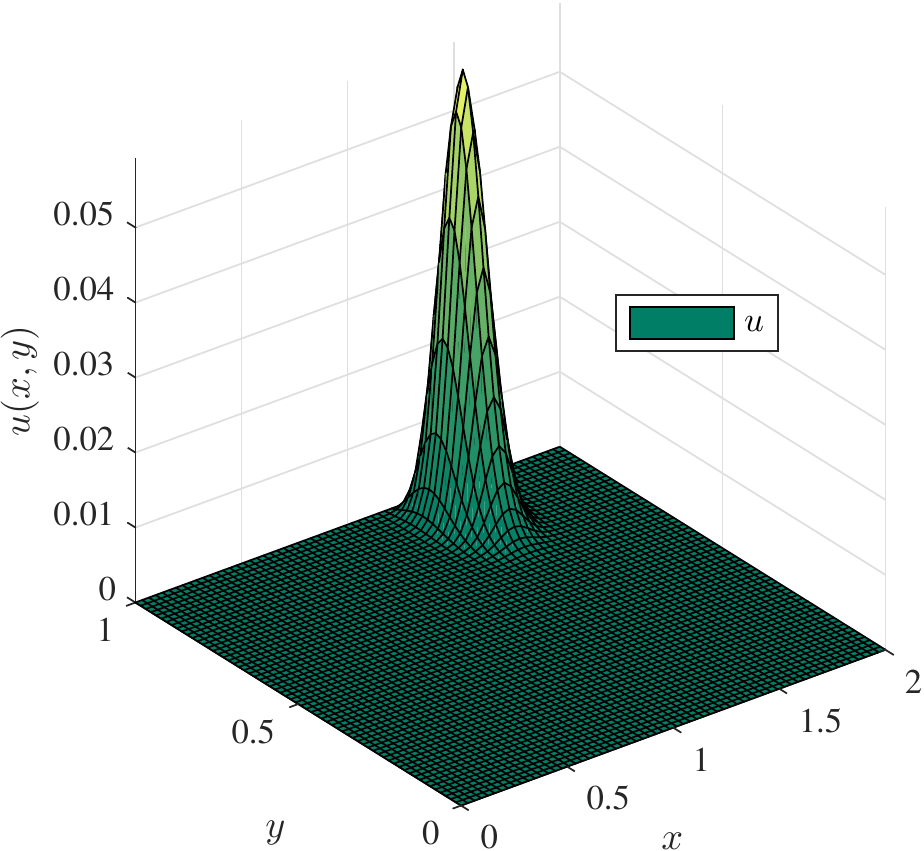}
	\caption{\small {\em Ex. \ref{ex:rectangular-domain-kick-jump}}. 
	Exact solution $u = (x_1^2 - 2 \, x_1) \, (x_2^2 - x_2) \, e^{-100 \,|(x_1, x_2) - (1.4, 0.95)|} $.}
	\label{fig:example-8-exact-solution}
\end{figure}
%
First, we run the test with the global ref. strategy. The obtained results are summarised in Tables 
\ref{tab:rect-domain-example-8-error-majorant-v-2-y-4-uniform-ref}--\ref{tab:rect-domain-example-8-times-v-2-y-4-uniform-ref}. 
Several systematically performed tests demonstrated that in order to perform a reliable estimation of the error in 
$u_h \in S_{h}^{2, 2}$, it is optimal to take $\flux_h \in S_{3h}^{4, 4} \oplus S_{3h}^{4, 4}$, i.e., we obtain 
efficient error bounds with the minimal computation effort spent on assembling and solving \eqref{eq:system-fluxh}.
\begin{table}[!t]
 \footnotesize
\centering
\newcolumntype{g}{>{\columncolor{gainsboro}}c} 	
\begin{tabular}{c|c|ccc|gc|c}
\quad \# ref. \quad & 
\quad  $\| \nabla e \|_\Omega$ \qquad   & 	  
\quad \quad \quad \quad $\overline{\rm M}$ \quad \quad \quad \quad &    
\quad \quad   $\mdI$ \qquad \quad & 	       
\quad \qquad $\mfI$ \qquad \quad  &  
\qquad $\Ieff (\overline{\rm M})$ \qquad & 
\qquad $\Ieff (\overline{\rm \eta})$ \qquad & 
\qquad \quad e.o.c. \qquad \quad \\
\midrule
   2 &   5.5665e-03 &   1.7635e-01 &   5.6126e-02 &   4.2226e-01 &      31.6798 &       8.5283 &   2.9617 \\
   4 &   2.6655e-04 &   1.1913e-02 &   4.0094e-03 &   2.7762e-02 &      44.6942 &      10.6943 &   2.1266 \\
   6 &   1.6374e-05 &   3.0360e-05 &   1.7856e-05 &   4.3919e-05 &       1.8541 &      10.8469 &   2.0163 \\
   8 &   1.0223e-06 &   1.1654e-06 &   1.1146e-06 &   1.7861e-07 &       1.1400 &      10.8565 &   2.0031 \\
\end{tabular}
\caption{\small {\em Ex. \ref{ex:rectangular-domain-kick-jump}}. 
Error, majorant (with dual and reliability terms), 
efficiency indices, and e.o.c. w.r.t. unif. ref. steps.}
\label{tab:rect-domain-example-8-error-majorant-v-2-y-4-uniform-ref}
\end{table}
\begin{table}[!t]
 \footnotesize
\centering
\newcolumntype{g}{>{\columncolor{gainsboro}}c} 	
\begin{tabular}{c|cc|cg|cg|cgc}
\# ref.  & 
\# d.o.f.($u_h$) &  \# d.o.f.($\flux_h$) &  
\; $t_{\rm as}(u_h)$ \; & 
\; $t_{\rm as}(\flux_h)$ \; & 
\; $t_{\rm sol}(u_h)$ \; & 
\; $t_{\rm sol}(\flux_h)$ \; &
$t_{\rm e/w}(\| \nabla e \|)$ & 
$t_{\rm e/w}(\overline{\rm M})$ & 
$t_{\rm e/w}(\overline{\eta})$ \\
\midrule
   2 &       1156 &        400 &     0.0299 &     0.1359 &     0.0058 &     0.0254 &       0.0705 &       0.0665 &       0.1066 \\
   4 &      16900 &       1296 &     0.5347 &     0.5352 &     0.4255 &     0.1728 &       0.9429 &       0.6846 &       1.4733 \\
   6 &     264196 &      17424 &     7.3424 &     7.9156 &    23.5576 &    16.7765 &      14.1433 &      10.1504 &      22.3347 \\
   8 &    4202500 &     266256 &   107.9652 &   121.7985 &  1516.5229 &   970.6061 &     238.5717 &     155.4827 &     370.6186 \\
\end{tabular}
\caption{\small {\em Ex. \ref{ex:rectangular-domain-kick-jump}}. 
Time for assembling and solving the systems that generate $u_h$ and $\flux_h$ 
as well as the time spent on e/w evaluation of error, majorant, and 
residual error estimator w.r.t. unif. ref. steps.}
\label{tab:rect-domain-example-8-times-v-2-y-4-uniform-ref}
\end{table}

Still, the most interesting test-case is the one that checks the performance of the majorant for the adaptive refinement. 
A series of tests showed that the optimal setting (in terms of quality of the error bounds and computational time 
spent on its reconstruction) is the approximation of $\flux_h$ with THB-basis functions of the degree $4$, i.e., 
$\flux_h \in S_{h}^{4, 4} \oplus S_{h}^{4, 4}$. At the same time, we consider the same mesh $\mathcal{K}_h$ 
that is used for approximation of $u_h$.
%
The obtained decrease of error and majorant with the marking criteria ${\mathds{M}}_{{\rm \bf BULK}}(0.4)$ and 
${\mathds{M}}_{{\rm \bf BULK}}(0.6)$ is illustrated in 
Table \ref{tab:rect-domain-example-8-error-majorant-v-2-y-3-adaptive-ref}, the corresponding time expenses are 
summarised in Table \ref{tab:rect-domain-example-8-times-v-2-y-3-adaptive-ref}. The most efficient error decrease
is obtained for the bulk parameter $\theta = 0.4$, which can be detected from Figure 
\ref{fig:convergence-majorant-bulk}. 

Figure \ref{fig:rect-domain-example-8-times-v-2-y-3-adaptive-ref} presents the evolution of physical meshes obtained 
during the refinement steps with different marking criteria ${\mathds{M}}_{{\rm \bf BULK}}(0.4)$ and 
${\mathds{M}}_{{\rm \bf BULK}}(0.2)$. Again, it is easy to observe from the graphics that the percentage of the refined elements on the right is higher than the percentage of such elements on the left. 

\begin{figure}[!t]
	\centering
	\includegraphics[scale=0.7]{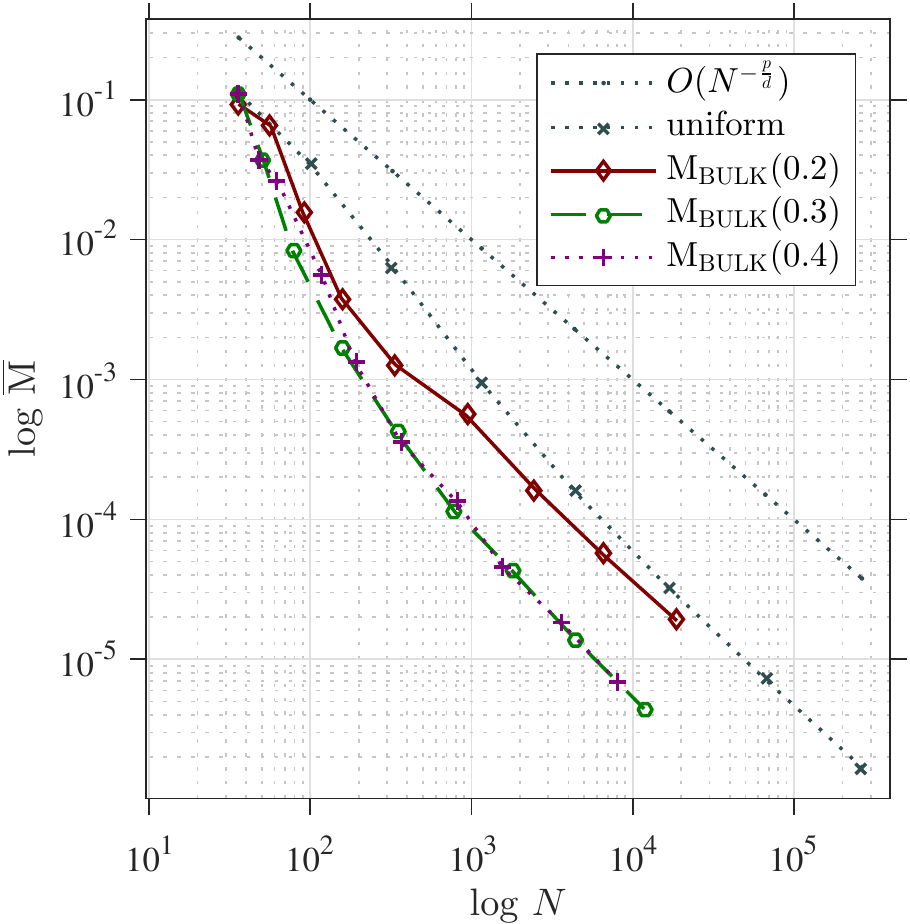}
	\caption{Convergence of the majorant for different marking criteria.}
	\label{fig:convergence-majorant-bulk}
\end{figure}

When the exact solution contains large local changes in the gradient (such as the one in the current example), 
the assembling and solving the system \eqref{eq:system-fluxh} becomes computationally heavier than respective 
procedures for \eqref{eq:system-uh}. This can be 
explained by the size of the generated optimal system \eqref{eq:system-fluxh} providing the reconstruction of 
vector-valued $\flux_h$. This drawback can be possibly eliminated by introducing multi-threading techniques 
(e.g., OpenMP, MPI) into the implementation of THB-splines. However, this matter stays beyond the focus of current 
paper and will be addressed in the upcoming technical report. 

\begin{table}[!t]
 \footnotesize
\centering
\newcolumntype{g}{>{\columncolor{gainsboro}}c} 	
\begin{tabular}{c|c|ccc|gc|c}
\quad \# ref. \quad & 
\quad  $\| \nabla e \|_\Omega$ \qquad   & 	  
\quad \quad \quad \quad $\overline{\rm M}$ \quad \quad \quad \quad &    
\quad \quad   $\mdI$ \qquad \quad & 	       
\quad \qquad $\mfI$ \qquad \quad  &  
\qquad $\Ieff (\overline{\rm M})$ \qquad & 
\qquad $\Ieff (\overline{\rm \eta})$ \qquad & 
\qquad \quad e.o.c. \qquad \quad \\
\midrule
\multicolumn{8}{c}{ (a) $\theta = 0.4$} \\
\midrule
   2 &   4.0502e-02 &   2.8670e-01 &   1.0725e-01 &   6.3030e-01 &       7.0786 &       5.5395 &   3.8159 \\
   4 &   5.2223e-03 &   3.7058e-02 &   1.5518e-02 &   7.5657e-02 &       7.0960 &       8.9835 &   4.7759 \\
   6 &   8.7993e-04 &   2.1154e-03 &   9.8705e-04 &   3.9631e-03 &       2.4040 &       8.6122 &   3.3554 \\
   8 &   1.1156e-04 &   1.9665e-04 &   1.1451e-04 &   2.8852e-04 &       1.7627 &       9.5564 &   2.8809 \\
\midrule
\multicolumn{8}{c}{ (b) $\theta = 0.2$} \\
\midrule
   2 &   3.6612e-02 &   2.0753e-01 &   8.3292e-02 &   4.3637e-01 &       5.6683 &       5.7646 &   3.0397 \\
   4 &   1.3527e-03 &   4.0090e-03 &   1.6115e-03 &   8.4210e-03 &       2.9637 &       9.3503 &   5.0006 \\
   6 &   1.5416e-04 &   3.1033e-04 &   1.6976e-04 &   4.9375e-04 &       2.0130 &      10.0080 &   2.0152 \\
   8 &   1.7351e-05 &   2.2095e-05 &   1.7587e-05 &   1.5832e-05 &       1.2734 &      10.4611 &   2.1194 \\
\end{tabular}
\caption{\small {\em Ex. \ref{ex:rectangular-domain-kick-jump}}. 
Error, majorant (with dual and reliability terms), efficiency indices, error, e.o.c. w.r.t. adaptive ref. steps.}
\label{tab:rect-domain-example-8-error-majorant-v-2-y-3-adaptive-ref}
\end{table}
\begin{table}[!t]
 \footnotesize
\centering
\newcolumntype{g}{>{\columncolor{gainsboro}}c} 	
\begin{tabular}{c|cc|cg|cg|cgc}
\# ref.  & 
\# d.o.f.($u_h$) &  \# d.o.f.($\flux_h$) &  
\; $t_{\rm as}(u_h)$ \; & 
\; $t_{\rm as}(\flux_h)$ \; & 
\; $t_{\rm sol}(u_h)$ \; & 
\; $t_{\rm sol}(\flux_h)$ \; &
$t_{\rm e/w}(\| \nabla e \|)$ & 
$t_{\rm e/w}(\overline{\rm M})$ & 
$t_{\rm e/w}(\overline{\eta})$ \\
\midrule
\multicolumn{10}{c}{ (a) $\theta = 0.4$} \\
\midrule
   2 &        124 &        145 &     0.0547 &     0.4156 &     0.0002 &     0.0025 &       0.0532 &       0.1738 & 0.1418\\ 
   4 &        243 &        245 &     0.2481 &     2.6372 &     0.0008 &     0.0063 &       0.3529 &       0.9892 & 0.9102 \\
   6 &        736 &        633 &     0.7903 &    10.7018 &     0.0052 &     0.0393 &       0.9605 &       3.5833 & 2.2969 \\
   8 &       2460 &       2231 &     2.3106 &    33.6222 &     0.0349 &     0.4035 &       2.9163 &      11.5602 & 4.4633\\
\midrule
\multicolumn{10}{c}{ (b) $\theta = 0.2$} \\
\midrule
   2 &        140 &        160 &     0.0449 &     0.3684 &     0.0002 &     0.0028 &       0.0651 &       0.1682 &       0.1111 \\
   4 &        366 &        366 &     0.2435 &     2.7071 &     0.0015 &     0.0124 &       0.3104 &       1.0138 &       0.5248 \\
   6 &       2043 &       1883 &     1.8328 &    25.3480 &     0.0246 &     0.3000 &       2.0493 &       8.0682 &       3.4391 \\
   8 &      15373 &      13974 &    17.4622 &   264.3344 &     0.4683 &     8.7014 &      15.9445 &      58.1277 &      25.9898 \\
\end{tabular}
\caption{\small {\em Ex. \ref{ex:rectangular-domain-kick-jump}}. 
Time for assembling and solving the systems generating d.o.f. of $u_h$ and $\flux_h$ as well as the time spent 
on e/w evaluation of error, majorant, and residual error estimator w.r.t. adaptive ref. steps.}
\label{tab:rect-domain-example-8-times-v-2-y-3-adaptive-ref}
\end{table}

\begin{figure}[!t]
	\centering
	{
	\subfloat[ref. \# 3: \quad ${\mathcal{K}}_h({\mathds{M}}_{{\rm \bf BULK}}(0.4))$]{
	\includegraphics[width=5.2cm, trim={3cm 4cm 3cm 4cm}, clip]{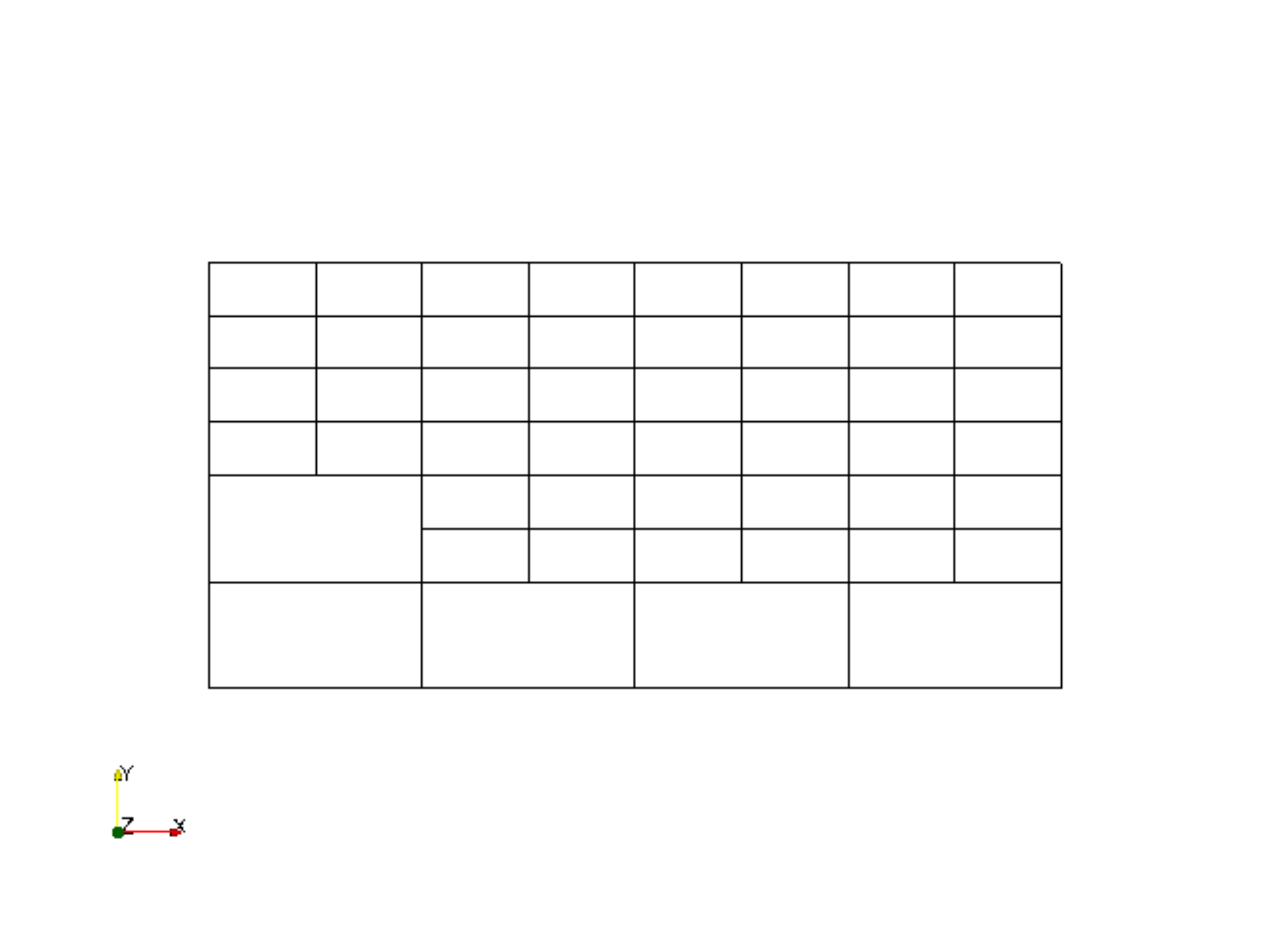}} 
	\quad 
	\subfloat[ref. \# 3:  \quad ${\mathcal{K}}_h({\mathds{M}}_{{\rm \bf BULK}}(0.2))$]{
	\includegraphics[width=5.2cm, trim={3cm 4cm 3cm 4cm}, clip]{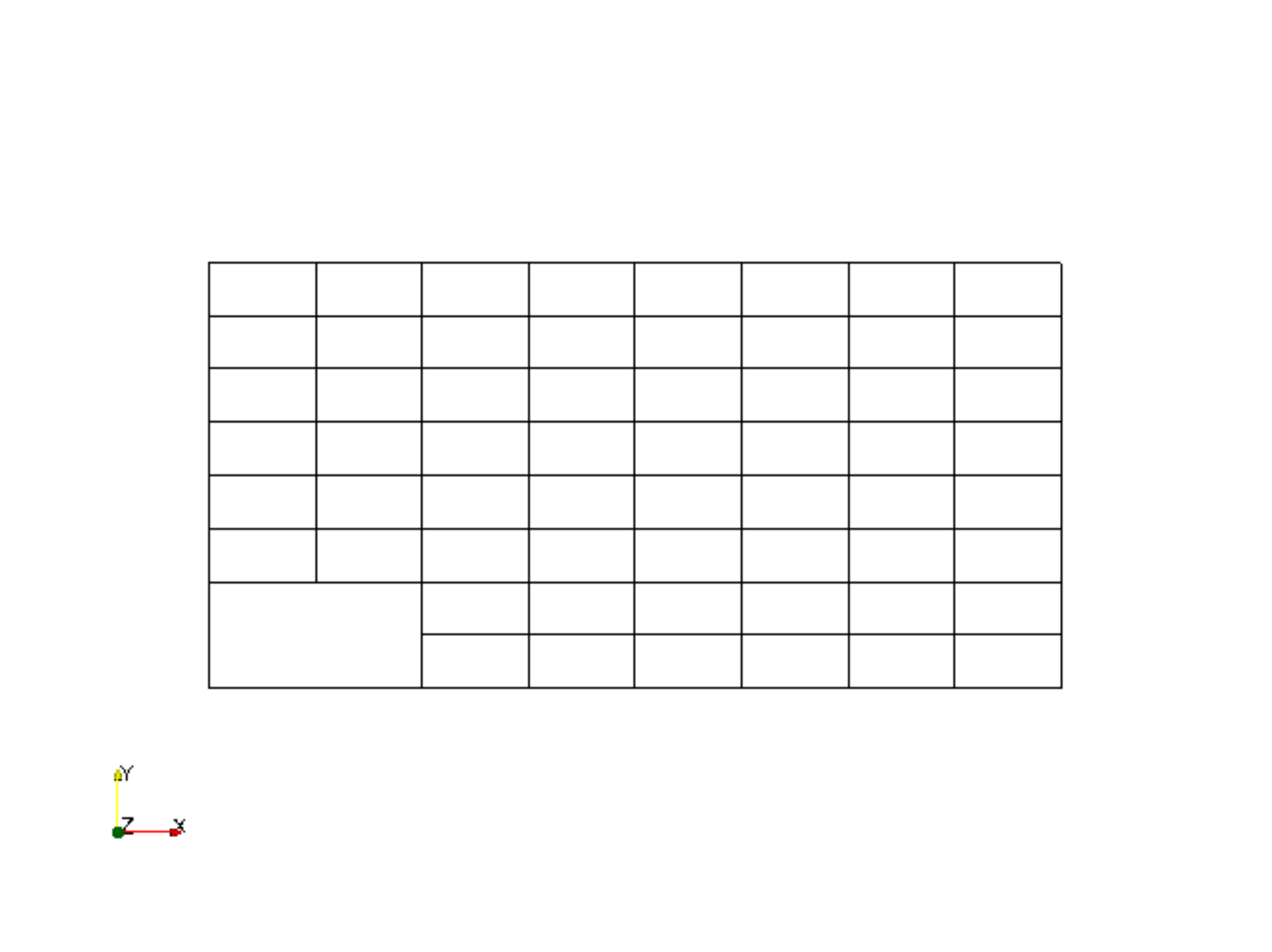}}
	}\\[2pt]
	{
	\subfloat[ref. \# 4: \quad ${\mathcal{K}}_h({\mathds{M}}_{{\rm \bf BULK}}(0.4))$]{
	\includegraphics[width=5.2cm, trim={3cm 4cm 3cm 4cm}, clip]{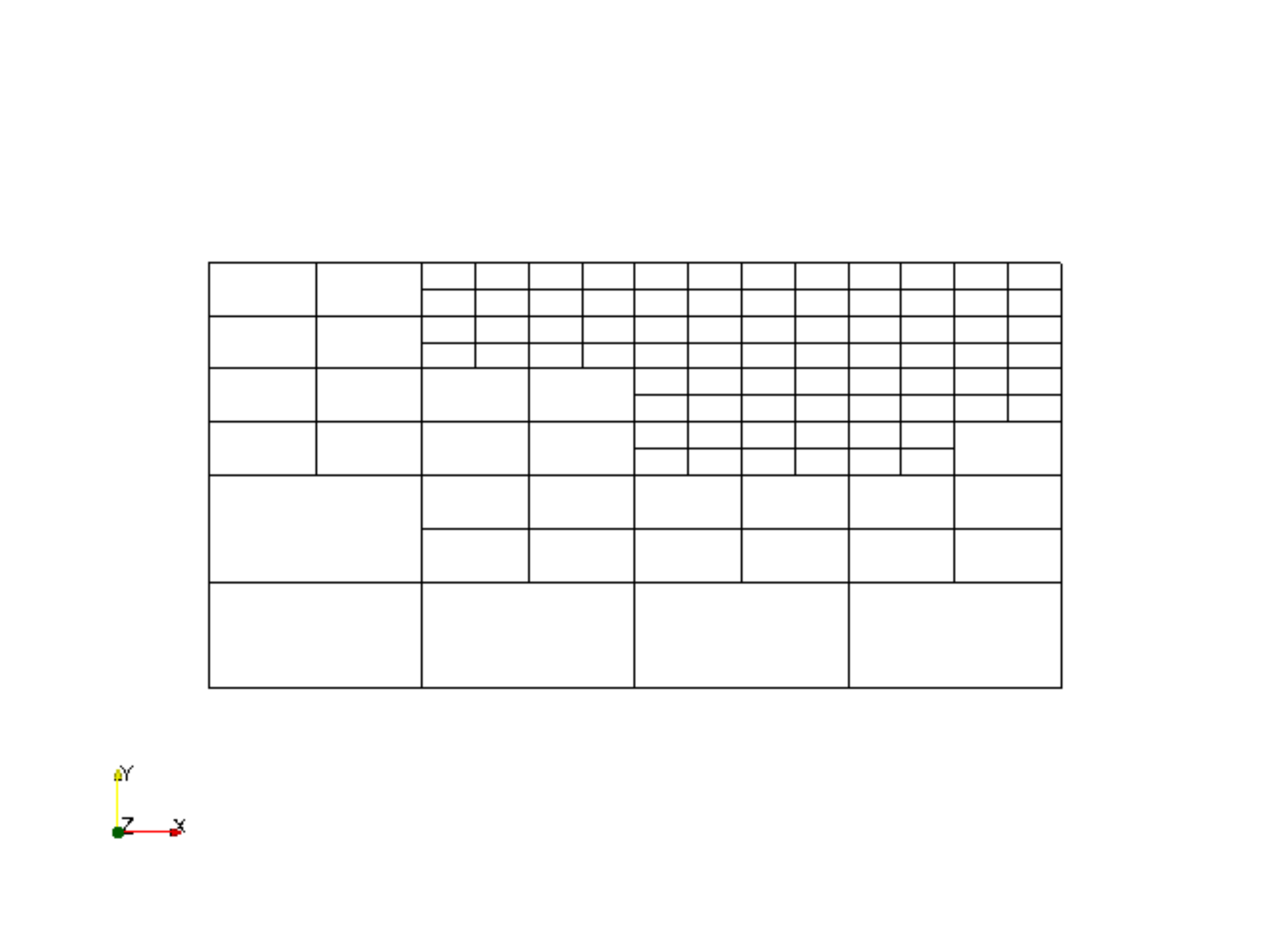}} 
	\quad 
	\subfloat[ref. \# 4:  \quad ${\mathcal{K}}_h({\mathds{M}}_{{\rm \bf BULK}}(0.2))$]{
	\includegraphics[width=5.2cm, trim={3cm 4cm 3cm 4cm}, clip]{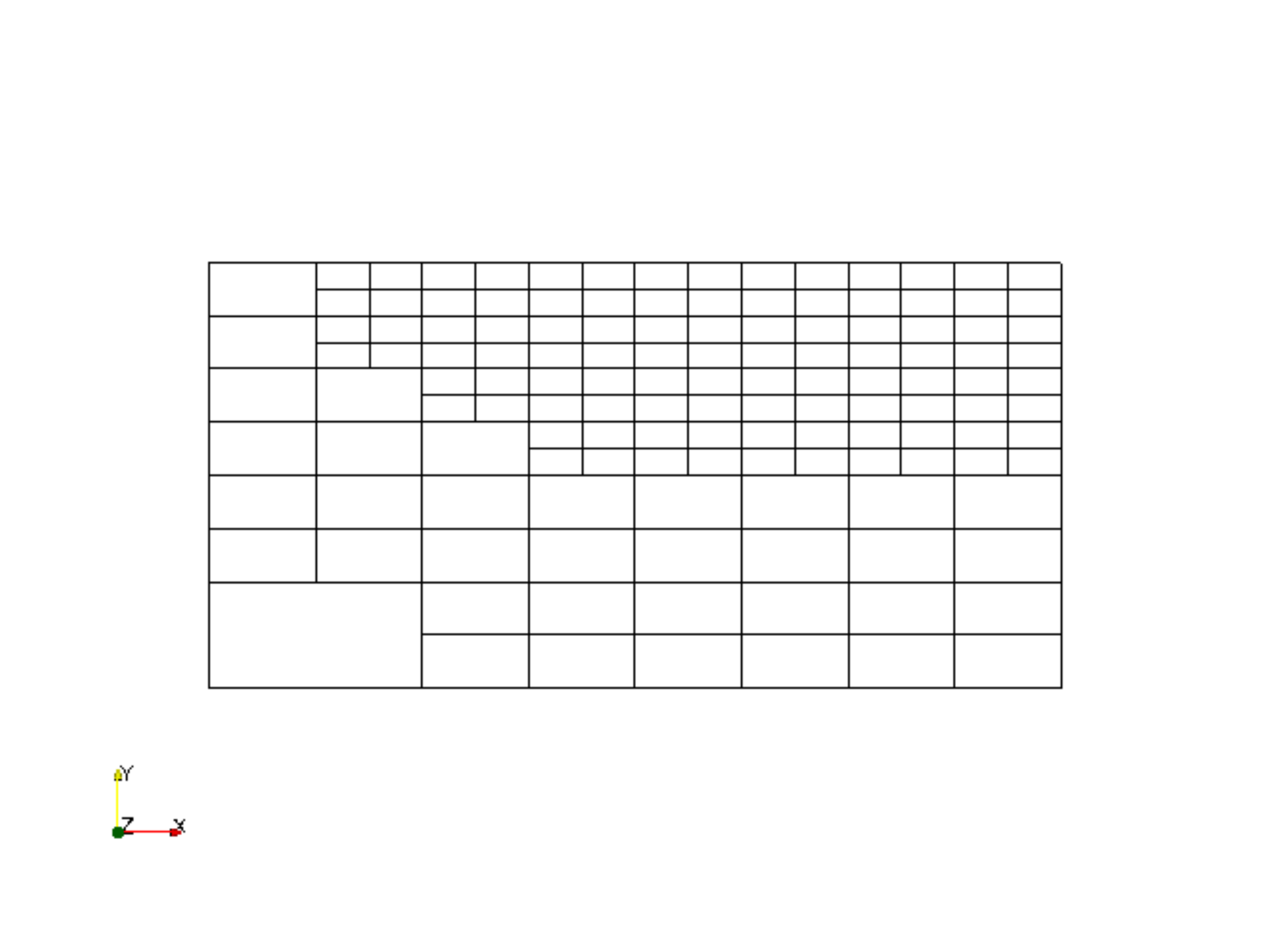}}
	}\\[2pt]
	{
	\subfloat[ref. \# 5:  \quad ${\mathcal{K}}_h({\mathds{M}}_{{\rm \bf BULK}}(0.4))$]{
	\includegraphics[width=5.2cm, trim={3cm 4cm 3cm 4cm}, clip]{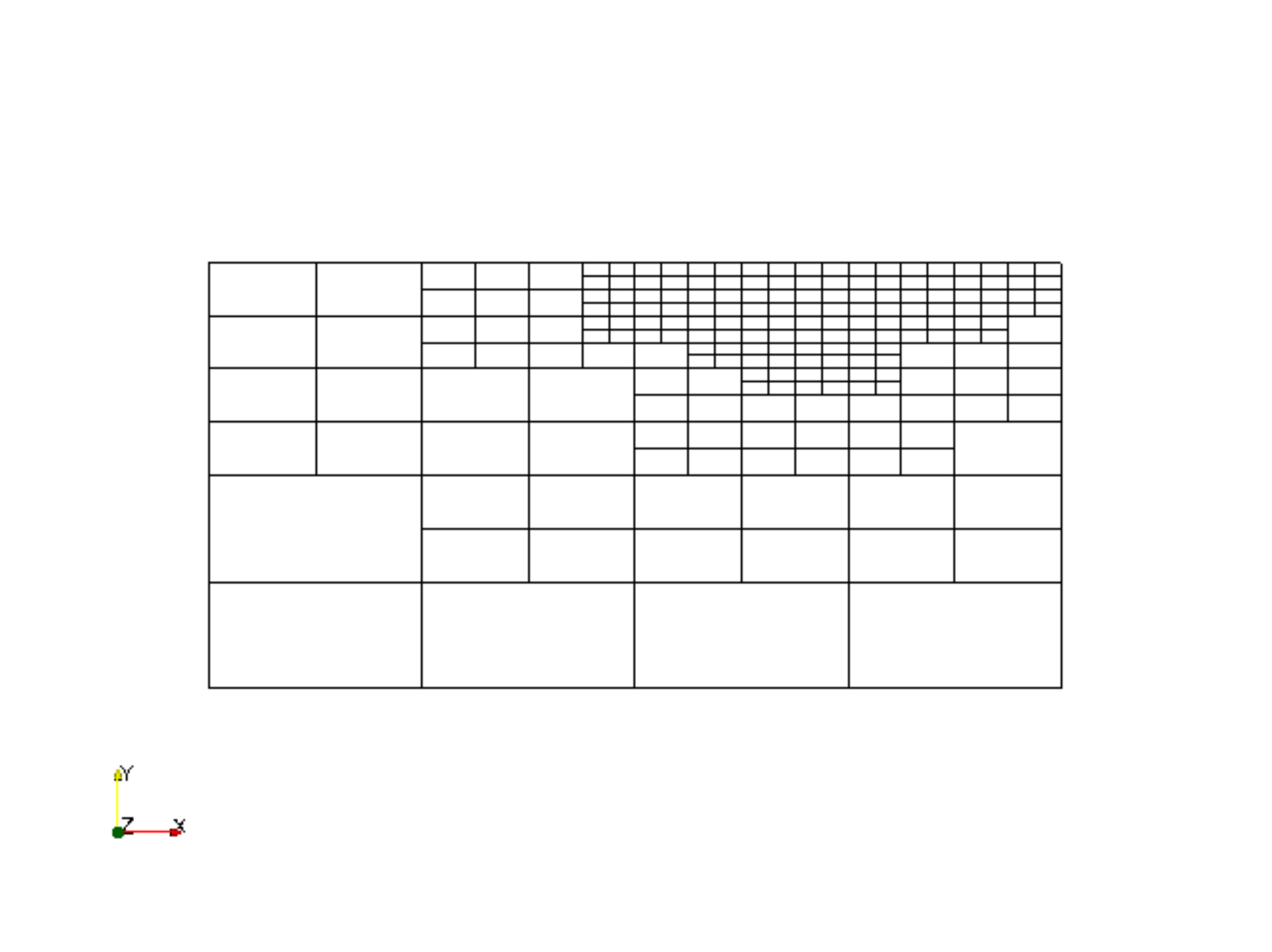}} 
	\quad
	\subfloat[ref. \# 5:  \quad ${\mathcal{K}}_h({\mathds{M}}_{{\rm \bf BULK}}(0.2))$]{
	\includegraphics[width=5.2cm, trim={3cm 4cm 3cm 4cm}, clip]{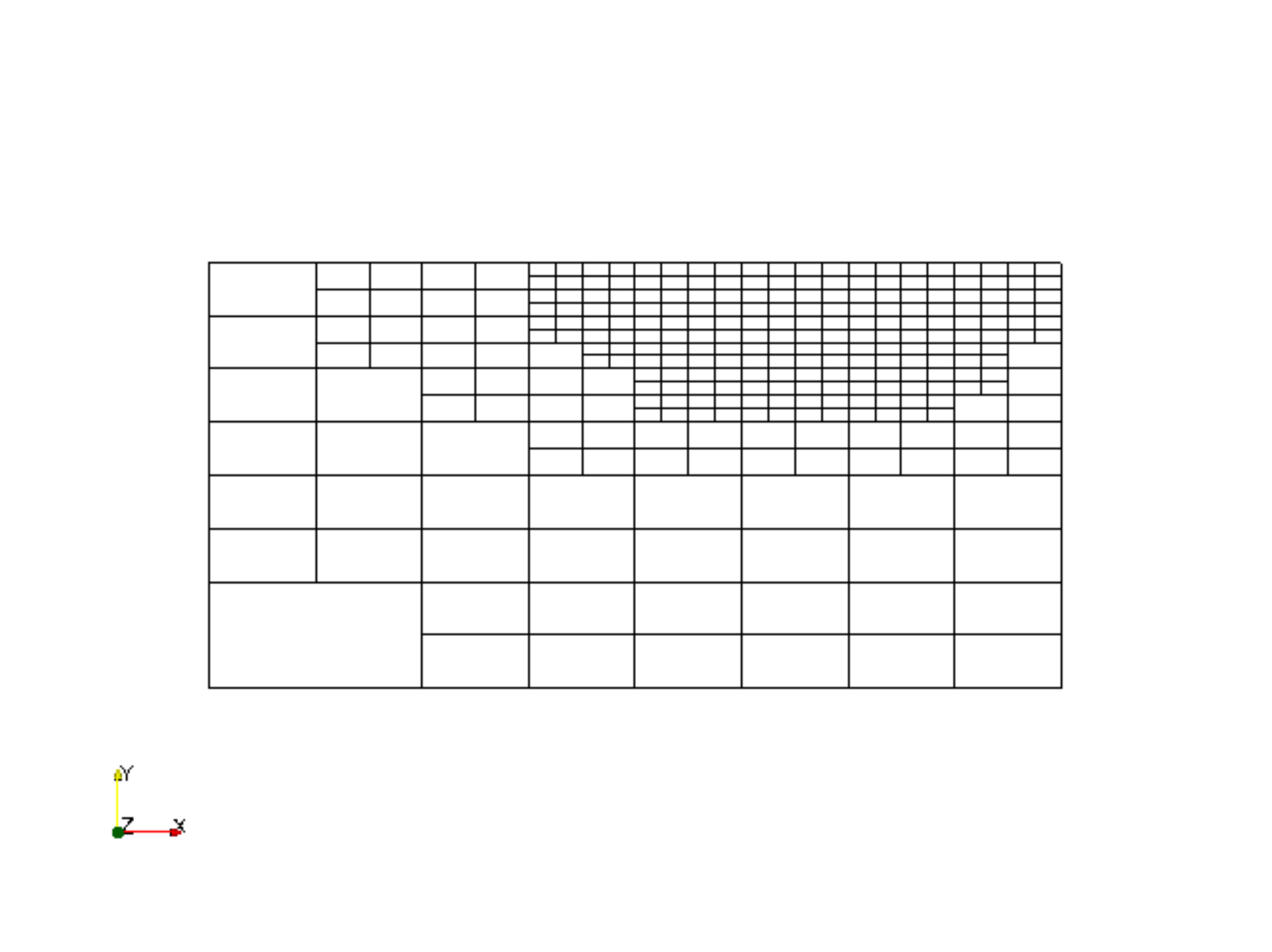}}
	}\\[2pt]
	{
	\subfloat[ref. \#6:  \quad ${\mathcal{K}}_h({\mathds{M}}_{{\rm \bf BULK}}(0.4))$]{
	\includegraphics[width=5.2cm, trim={3cm 4cm 3cm 4cm}, clip]{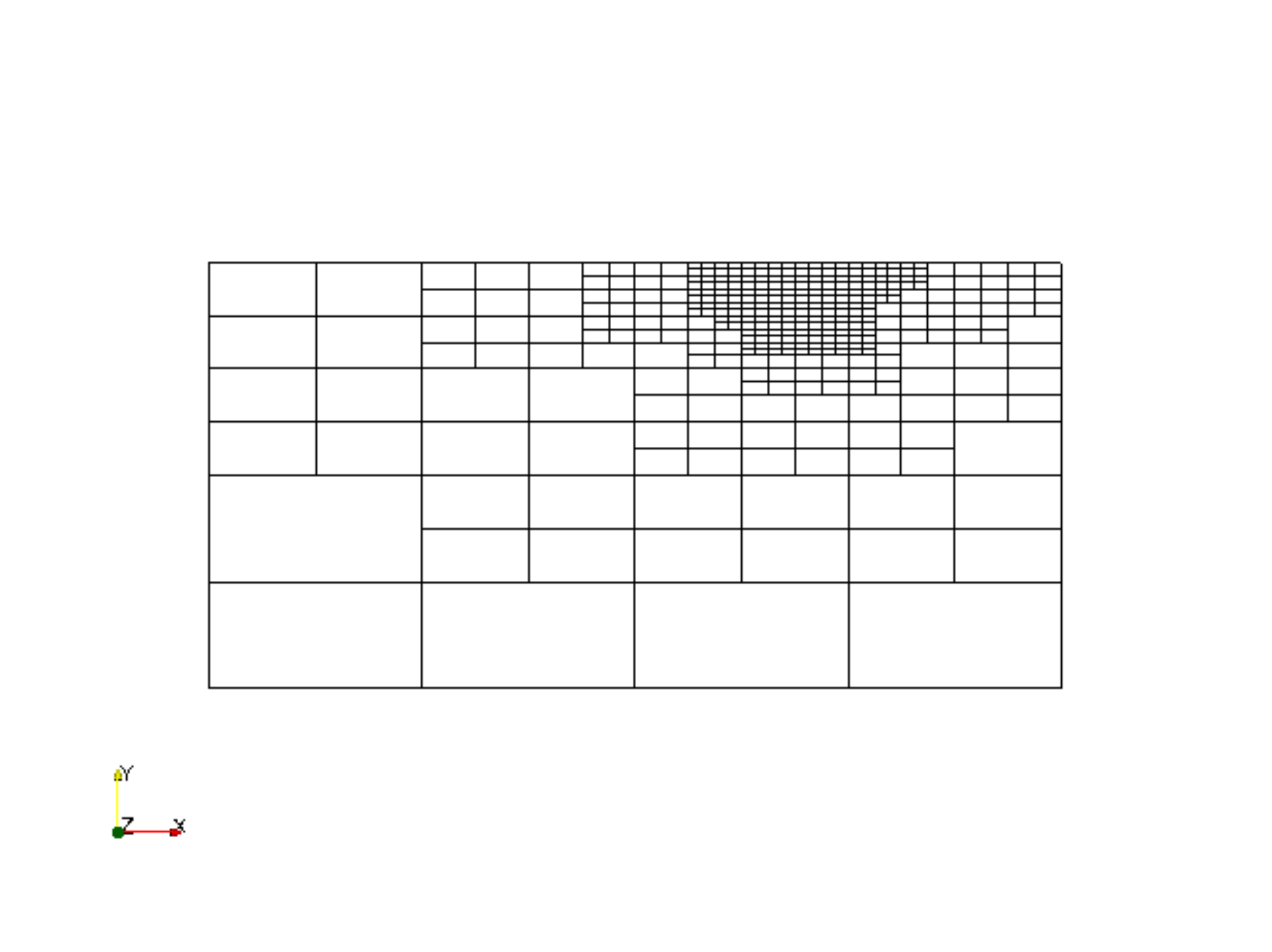}} 
	\quad
	\subfloat[ref. \# 6:  \quad ${\mathcal{K}}_h({\mathds{M}}_{{\rm \bf BULK}}(0.2))$]{
	\includegraphics[width=5.2cm, trim={3cm 4cm 3cm 4cm}, clip]{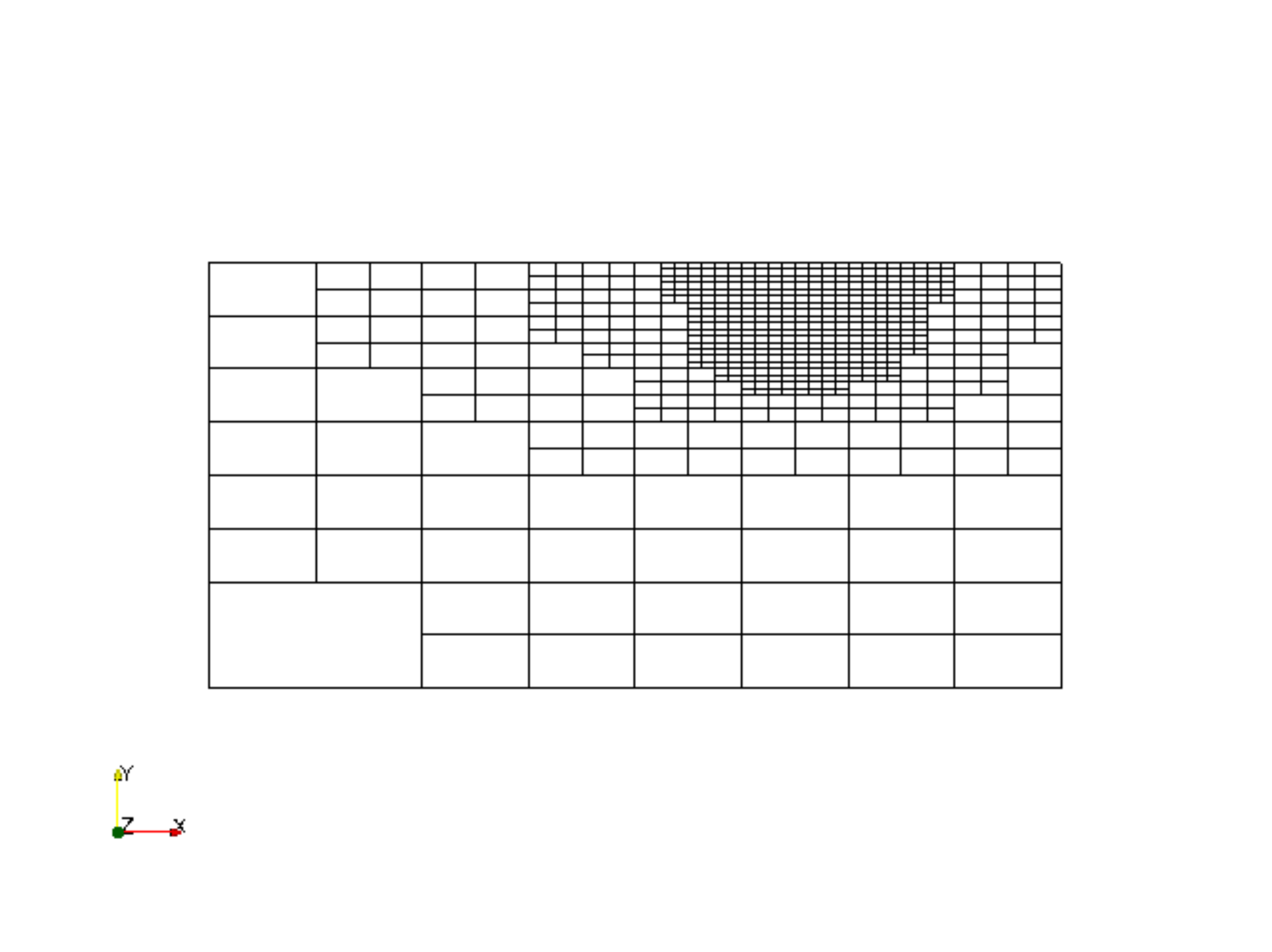}}
	}\\[2pt]
	{
	\subfloat[ref. \# 8:  \quad ${\mathcal{K}}_h({\mathds{M}}_{{\rm \bf BULK}}(0.4))$]{
	\includegraphics[width=5.2cm, trim={3cm 4cm 3cm 4cm}, clip]{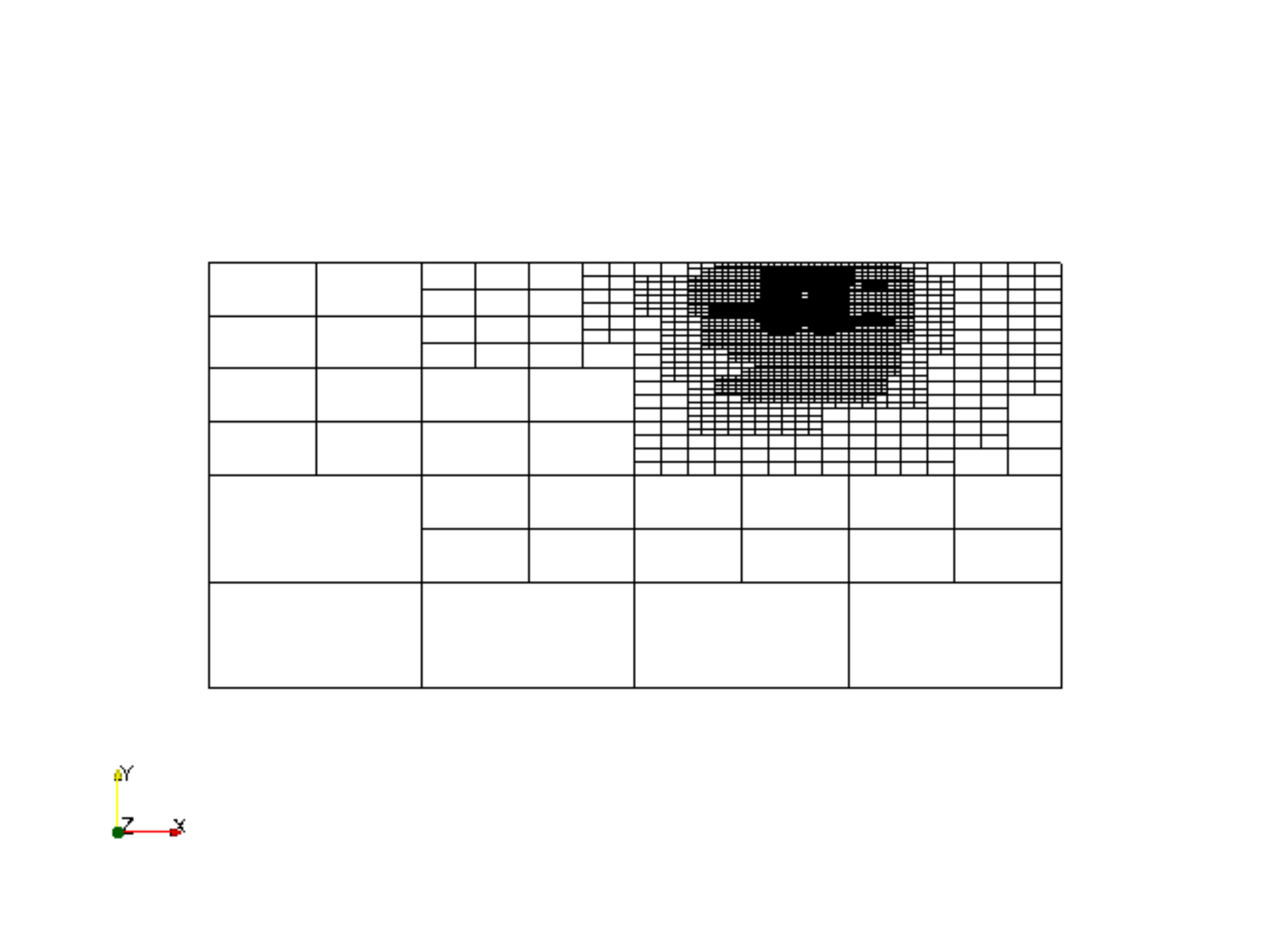}} 
	\quad
	\subfloat[ref. \# 8:  \quad ${\mathcal{K}}_h({\mathds{M}}_{{\rm \bf BULK}}(0.2))$]{
	\includegraphics[width=5.2cm, trim={3cm 4cm 3cm 4cm}, clip]{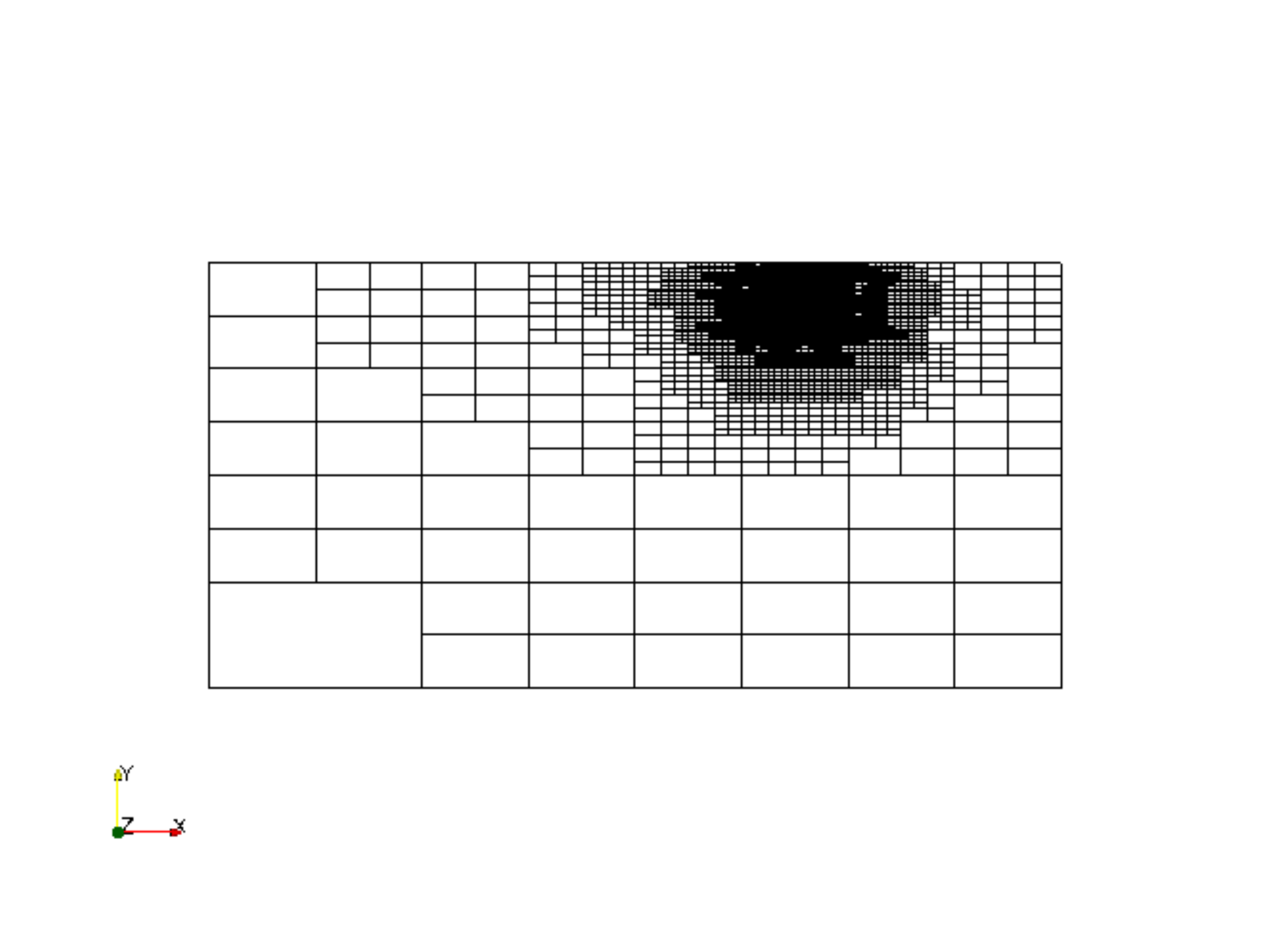}}
	}\\[2pt]
	{
	\subfloat[ref. \# 9:  \quad ${\mathcal{K}}_h({\mathds{M}}_{{\rm \bf BULK}}(0.4))$]{
	\includegraphics[width=5.2cm, trim={3cm 4cm 3cm 4cm}, clip]{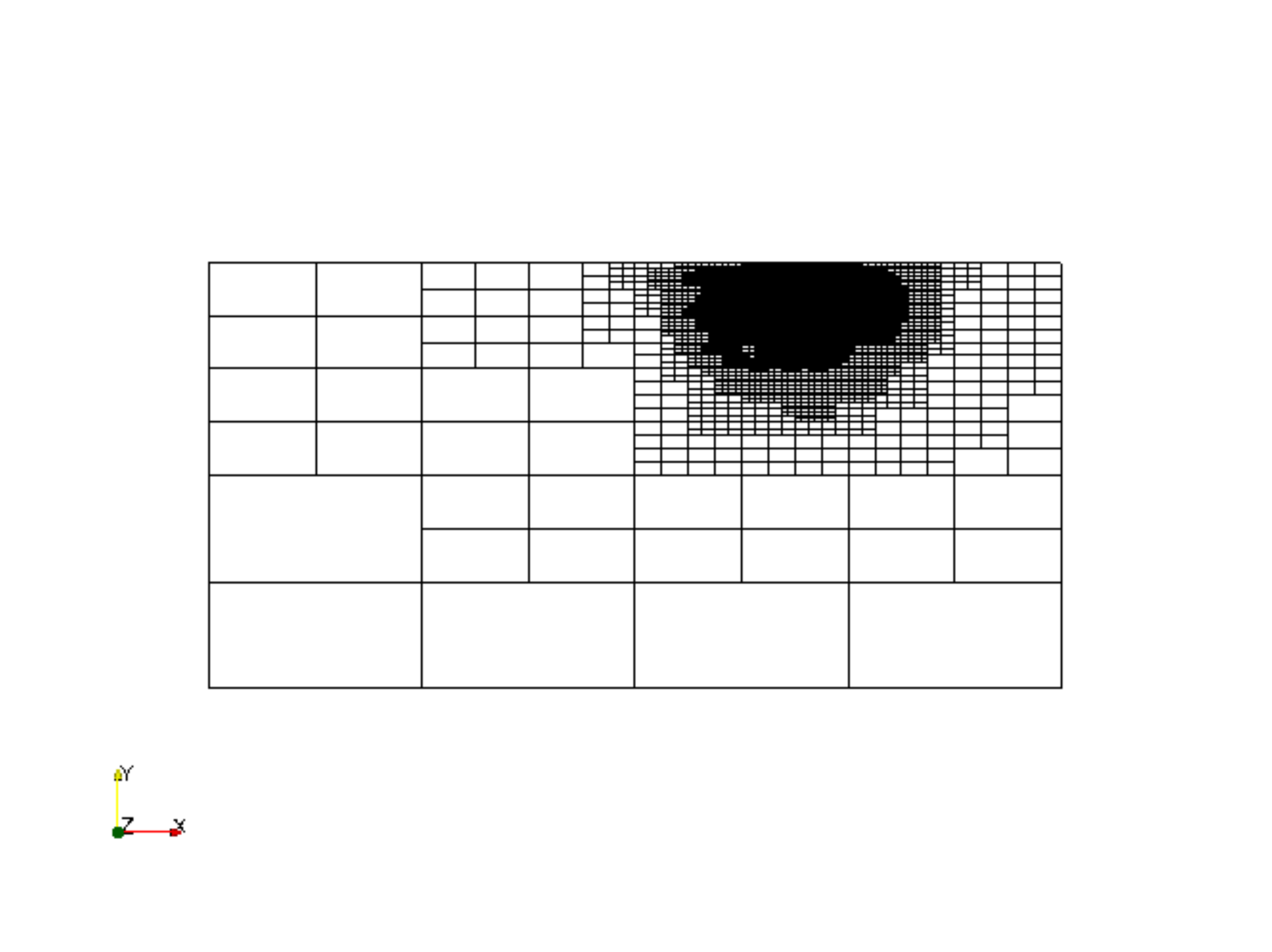}} 
	\quad
	\subfloat[ref. \# 9:  \quad ${\mathcal{K}}_h({\mathds{M}}_{{\rm \bf BULK}}(0.2))$]{
	\includegraphics[width=5.2cm, trim={3cm 4cm 3cm 4cm}, clip]{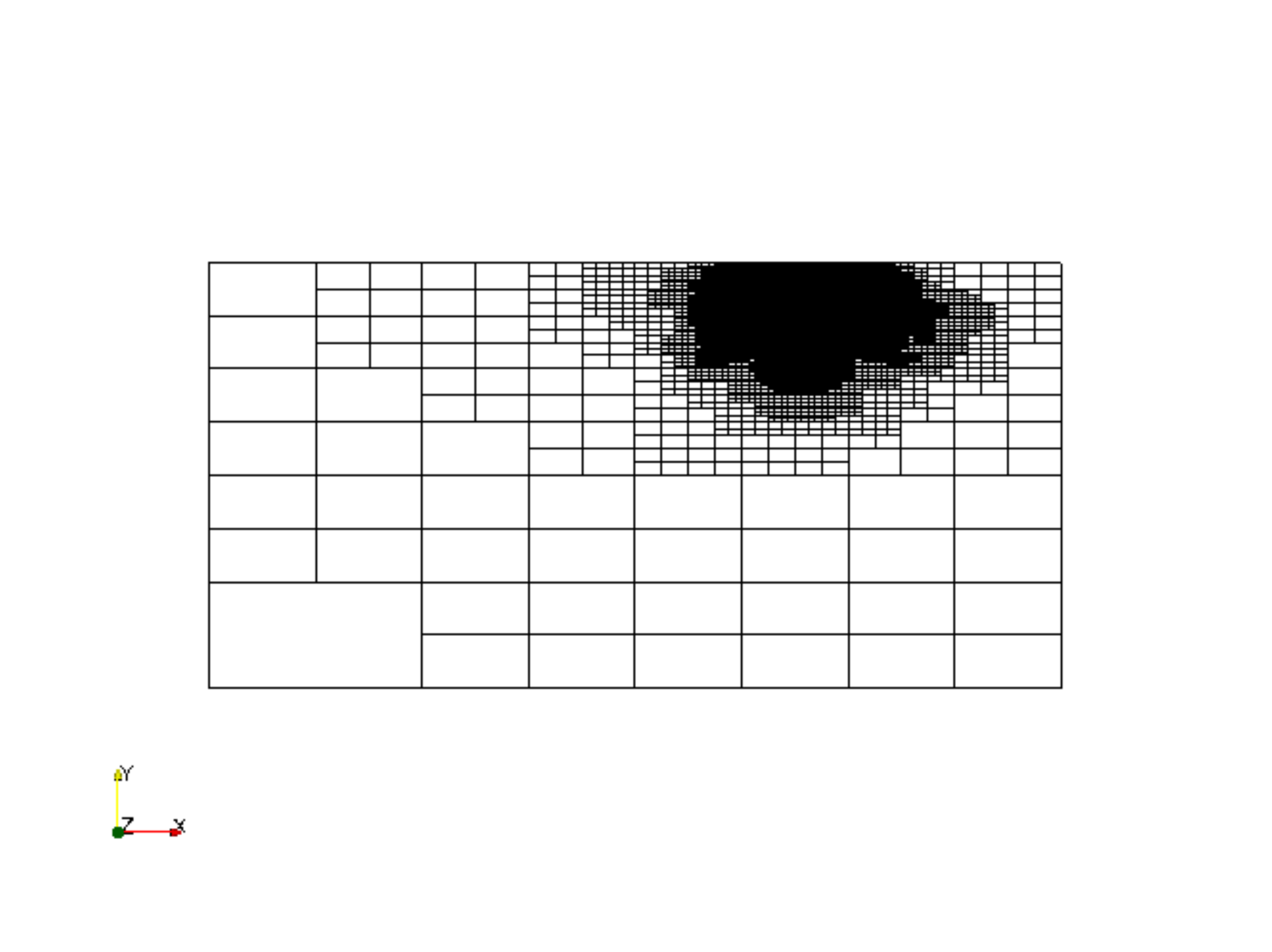}}
	}
	\caption{\small {\em Ex. \ref{ex:rectangular-domain-kick-jump}}. 
	Adaptive meshes obtained for the bulk parameters {$\theta = 0.4$ (left) and $\theta = 0.2$ (right)}.}
	\label{fig:rect-domain-example-8-times-v-2-y-3-adaptive-ref}
	\end{figure}
\end{example}

\newpage
\begin{example}
\label{ex:l-shape-domain-example-12}
\rm 
One of the classical benchmark examples (containing a singularity in the exact solution) is the problem with 
a $L$-shaped domain $\Omega := (\minus 1, 1) \times (\minus 1, 1) \backslash [0, 1) \times [0, 1)$. The Dirichlet BC 
are defined on $\Gamma$ by the load $u_D = r^{1/3} \, \sin (\theta)$, where  
$$r = (x_1^2 + x_2^2) \quad \mbox{and} \quad 
\theta = \begin{cases}
\tfrac13 \, (2\, {\rm atan2} (x_2, x_1) - \pi), \quad \;\;  x_2 > 0,\\ 
\tfrac13 \, (2\,{\rm atan2} (x_2, x_1) + 3\,\pi), \quad x_2 \leq 0.
\end{cases}
$$ 
The corresponding exact solution $u = u_D$ has the singularity in the 
point $(r, \theta) =(0, 0)$ (see also Figure \ref{fig:example-12-exact-solution}).
\begin{figure}[!t]
	\centering
	\subfloat[$u(x_1, x_2)$]{
	\includegraphics[scale=0.6]{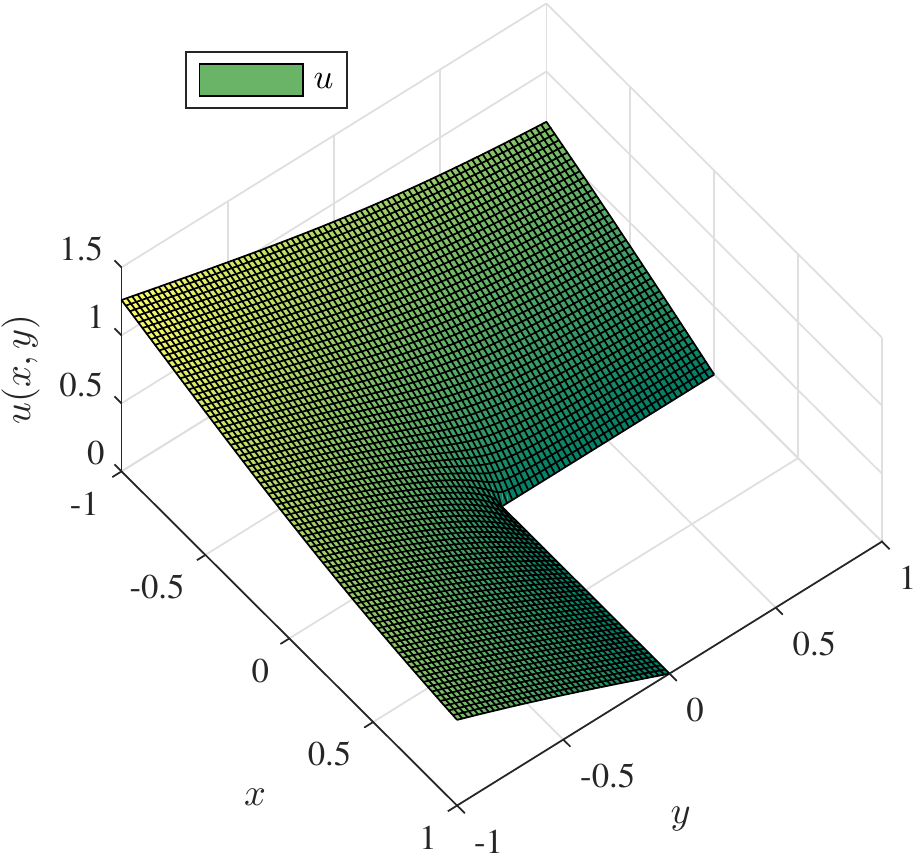}
	\label{fig:example-12-exact-solution}}
	\quad
	\subfloat[$\mathcal{K}_h$ defined on $\Omega$]{
	\includegraphics[width=5.0cm, trim={2cm 2cm 4cm 2cm}, clip]{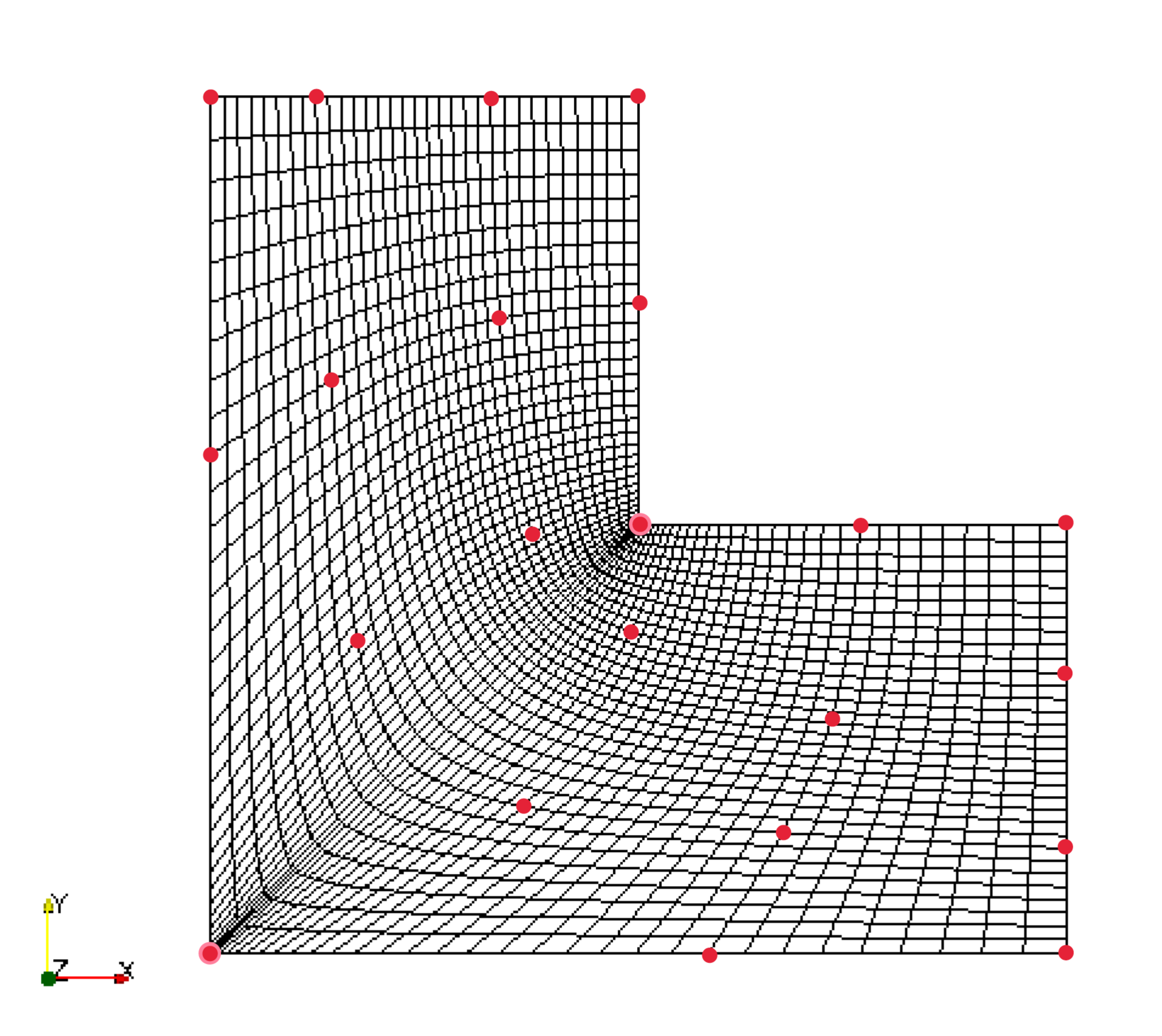}
	\label{fig:initial-mesh-2}}	
	\caption{\small {\em Ex. \ref{ex:l-shape-domain-example-12}}. 
	(a) Exact solution $u = r^{1/3} \, \sin (\theta)$.
	(b) Initial geometry data with Greville's points with double control points at the corners and 
	a corresponding mesh generated with $C^1$-continuous geometrical mapping.}
\end{figure}

The initial geometry data and the mesh are defined by Greville's points with double control points in the corners
marked with red circled markers in Figure \ref{fig:initial-mesh-2}, and considered knots-vectors 
{$\kappa = \{0,0,0, 0.25,0.5,0.75, 1,1,1\}$} and ${s} = \{ 0,0,0, 0.5, 1,1,1 \}$. 
Due to the doubled control points in the corners of the L-shape domain, only the re-entrant corner and 
its counterpart on the other side are singular (instead of the whole diagonal), i.e., the Jacobian of the geometry map in 
these two points is not regular. Since no integration points are placed at these corners, computational evaluation of the 
integrals make sense. The downside of such a setting is that with the increase of refinement steps the cells near these 
corners become rather thin and loose shape-regularity. In addition, since on the functional level the requirements on the 
regularity of $\flux$ is not fulfilled, the global error estimate is not reliable and has rather heuristic character, therefore, 
we only consider its performance from error indication point of view.

For these settings, 
the performance of the error majorant is compared to the performance of the residual error indicator in Table 
\ref{tab:lshape-domain-example-12-error-majorant-v-2-y-3-adaptive-ref}, where the first one is constructed with the 
help of fluxes $\flux_{h} \in S_{h}^{{3, 3}} \oplus S_h^{{3, 3}}$. By increasing the degree of splines that approximate
 $\flux_{h}$, one could reconstruct a sharper indicator from $\overline{\rm M}$. Whereas the residual error 
estimate (dependent only on $u_h$ and local $h_K$) always stays on the same `accuracy level', i.e., 
$\Ieff(\overline{\eta}) \approx 183.3358$. So the advantage of using the error majorant instead of residual-based error 
estimates is rather obvious for such kind of problems.  
However, the time required for reconstruction of $\overline{\rm M}$, increases as well (see Table 
\ref{tab:lshape-domain-example-12-times-v-2-y-3-adaptive-ref}). Hence, the selection of space for the 
dual variable $\flux_h$ is always dependent on the smoothness of the exact solution (or RHS) and on the 
allocated time for the a posteriori error estimates control.

In Figure \ref{fig:l-shape-domain-example-12-theta-10-domains-comparison}, we illustrate the evolution of 
the adaptive meshes discretising parametric domain $\widehat{\Omega}$ (left) and corresponding to them meshes 
discretising physical domain $\Omega$ (right). From presented graphics, one can see that the refinement is localised 
in the area close to 
the singularity point and no superfluous refinement is performed. We also perform the test, where, starting with the 
same initial mesh Figure \ref{fig:initial-mesh-2}, we compare $\mathcal{K}_h$ generated by the refinement based on 
the majorant (error indicator $\mdIK$) and by the refinement based on true error distribution $\| \nabla e\|^2_K$ 
(with bulk parameter for the marking chosen to be $\theta = 0.2$).
These meshes are illustrated in Figure \ref{fig:l-shape-domain-example-12-theta-20-error-majorant-refinement}, and 
it is obvious that the majorant provides an adequate strategy for the adaptive refinement.
The efficiency of the 
studied error bounds is also confirmed by comparing of majorant decrease for different marking criteria with 
respect to the \# d.o.f.($u_h$). In particular, Figure \ref{fig:l-shape-domain-example-9-convergence-majorant} shows 
that using majorant in combination with different markers does not only 
improve e.o.c. (in comparison to the one provided by the uniform refinement) but provides an even better one 
than $p = 2$. 
{In the limit with $h \rightarrow 0$, the e.o.c. for the refinement strategy using any marking criteria is similar. 
However, on the initial refinement steps GARU (greatest appearing residual utilisation) marker locate those elements, 
where the error at least $90\%$ of the highest error (over all the elements discretising the domain). For the considered 
L-shape domain, the maximum error is very high on first refinement steps and is mainly localised in the area close 
to the singular corner. Taking that into account, GARU$90\%$ selects small number of $K \in \mathcal{K}_h$ with very 
high errors. After the refinement, new elements contribute with smaller values into the absolute maximum error and 
eventually e.o.c. becomes similar to the other marking strategies. PUCA$10\%$ marking criterion, in turn, always 
refines $90\%$ of total number of the elements (such that the one with the highest errors come first). Such a strategy 
is not connected to the values of the distributed error. Therefore, the error decay using PUCA$10\%$ is rather 
latent w.r.t to the number of refinements.
}
%
%
%

 %
 \begin{table}[!t]
 \footnotesize
\centering
\newcolumntype{g}{>{\columncolor{gainsboro}}c} 	
\begin{tabular}{c|c|ccc|gc|c}
\quad \# ref. \quad & 
\quad  $\| \nabla e \|_\Omega$ \qquad   & 	  
\quad \quad \quad \quad $\overline{\rm M}$ \quad \quad \quad \quad &    
\quad \quad   $\mdI$ \qquad \quad & 	       
\quad \qquad $\mfI$ \qquad \quad  &  
\qquad $\Ieff (\overline{\rm M})$ \qquad & 
\qquad $\Ieff (\overline{\rm \eta})$ \qquad & 
\qquad \quad e.o.c. \qquad \quad \\
\midrule
   2 &   1.3359e-02 &   3.6140e-02 &   2.7595e-02 &   1.8983e-02 &       2.7053 &      42.6770 &   5.5339 \\
   3 &   8.7388e-03 &   2.3661e-02 &   1.8061e-02 &   1.2442e-02 &       2.7076 &      55.7909 &   6.0339 \\
   4 &   5.7148e-03 &   1.5862e-02 &   1.1998e-02 &   8.5837e-03 &       2.7755 &      74.0819 &   4.5481 \\
   5 &   3.7337e-03 &   1.0751e-02 &   8.0370e-03 &   6.0286e-03 &       2.8794 &      99.4303 &   2.4711 \\
   6 &   2.4386e-03 &   7.3896e-03 &   5.4326e-03 &   4.3474e-03 &       3.0303 &     134.3125 &   1.7454 \\
   7 &   1.5774e-03 &   4.8596e-03 &   3.6979e-03 &   2.5806e-03 &       3.0807 &     183.3358 &   1.3410 \\
\end{tabular}
\caption{\small {\em Ex. \ref{ex:l-shape-domain-example-12}}. 
Error, majorant (with dual and reliability terms), 
efficiency indices, error, e.o.c. 
w.r.t. adaptive ref. steps, $\flux_{h} \in S_{h}^{{3, 3}} \oplus S_h^{{3, 3}}$.}
\label{tab:lshape-domain-example-12-error-majorant-v-2-y-3-adaptive-ref}
\end{table}
\begin{table}[!t]
 \footnotesize
\centering
\newcolumntype{g}{>{\columncolor{gainsboro}}c} 	
\begin{tabular}{c|cc|cg|cg|cgc}
\# ref.  & 
\# d.o.f.($u_h$) &  \# d.o.f.($\flux_h$) &  
\; $t_{\rm as}(u_h)$ \; & 
\; $t_{\rm as}(\flux_h)$ \; & 
\; $t_{\rm sol}(u_h)$ \; & 
\; $t_{\rm sol}(\flux_h)$ \; &
$t_{\rm e/w}(\| \nabla e \|)$ & 
$t_{\rm e/w}(\overline{\rm M})$ & 
$t_{\rm e/w}(\overline{\eta})$ \\
\midrule
   1 &        612 &        760 &     0.2051 &     1.6904 &          0.0001 &          0.0437 &       0.3510 &       1.2626 &       0.6734 \\
   3 &        823 &        971 &     0.2290 &     2.1511 &          0.0013 &          0.0559 &       0.3223 &       1.8207 &       0.6083 \\
   5 &       1400 &       1543 &     0.5463 &     7.1518 &          0.0069 &          0.1361 &       0.7253 &       5.1051 &       1.2739 \\
   7 &       4368 &       4468 &     2.2633 &    43.5814 &          0.1138 &          0.7445 &       2.4629 &      18.6119 &       4.2903 \\
\end{tabular}
\caption{\small {\em Ex. \ref{ex:l-shape-domain-example-12}}. 
Time for assembling and solving of systems generating d.o.f. of $u_h$ and $\flux_h$ 
as well as the time spent on e/w evaluation of error, majorant, and 
residual error estimator w.r.t. adaptive ref. steps, $\flux_{h} \in S_{h}^{{3, 3}} \oplus S_h^{{3, 3}}$.}
\label{tab:lshape-domain-example-12-times-v-2-y-3-adaptive-ref}
\end{table}

%
%
\begin{figure}[!t]
	\centering
	{
	\subfloat[ref. \# 2: \quad $\widehat{\mathcal{K}}_h$ $\leftarrow$ $\| \nabla e\|^2_K$]{
	\includegraphics[width=4.5cm, trim={4cm 2cm 4cm 2cm}, clip]{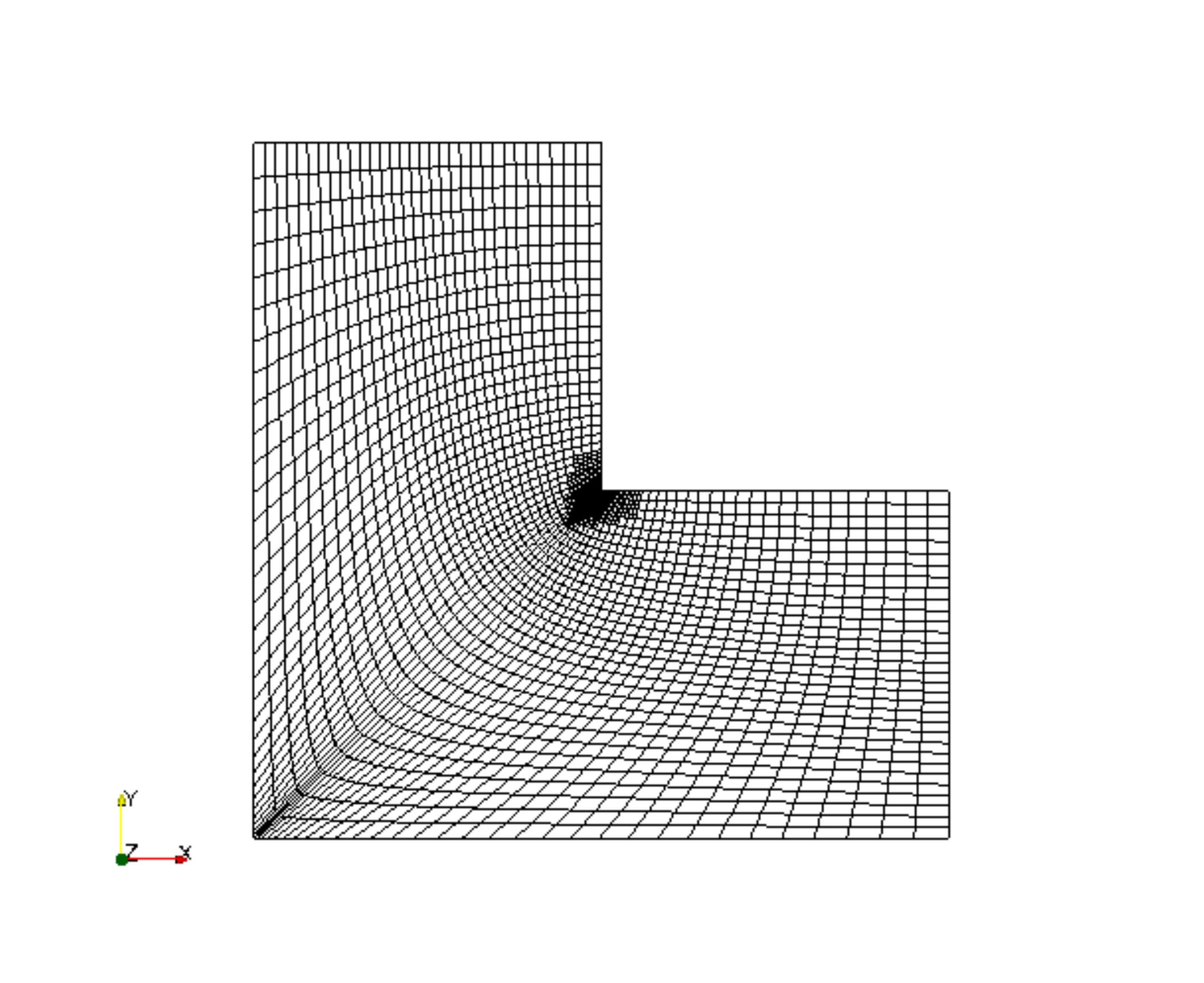}} \qquad
	\subfloat[ref. \# 2: \quad $\widehat{\mathcal{K}}_h$ $\leftarrow$ $\mdIK$]{
	\includegraphics[width=4.5cm, trim={4cm 2cm 4cm 2cm}, clip]{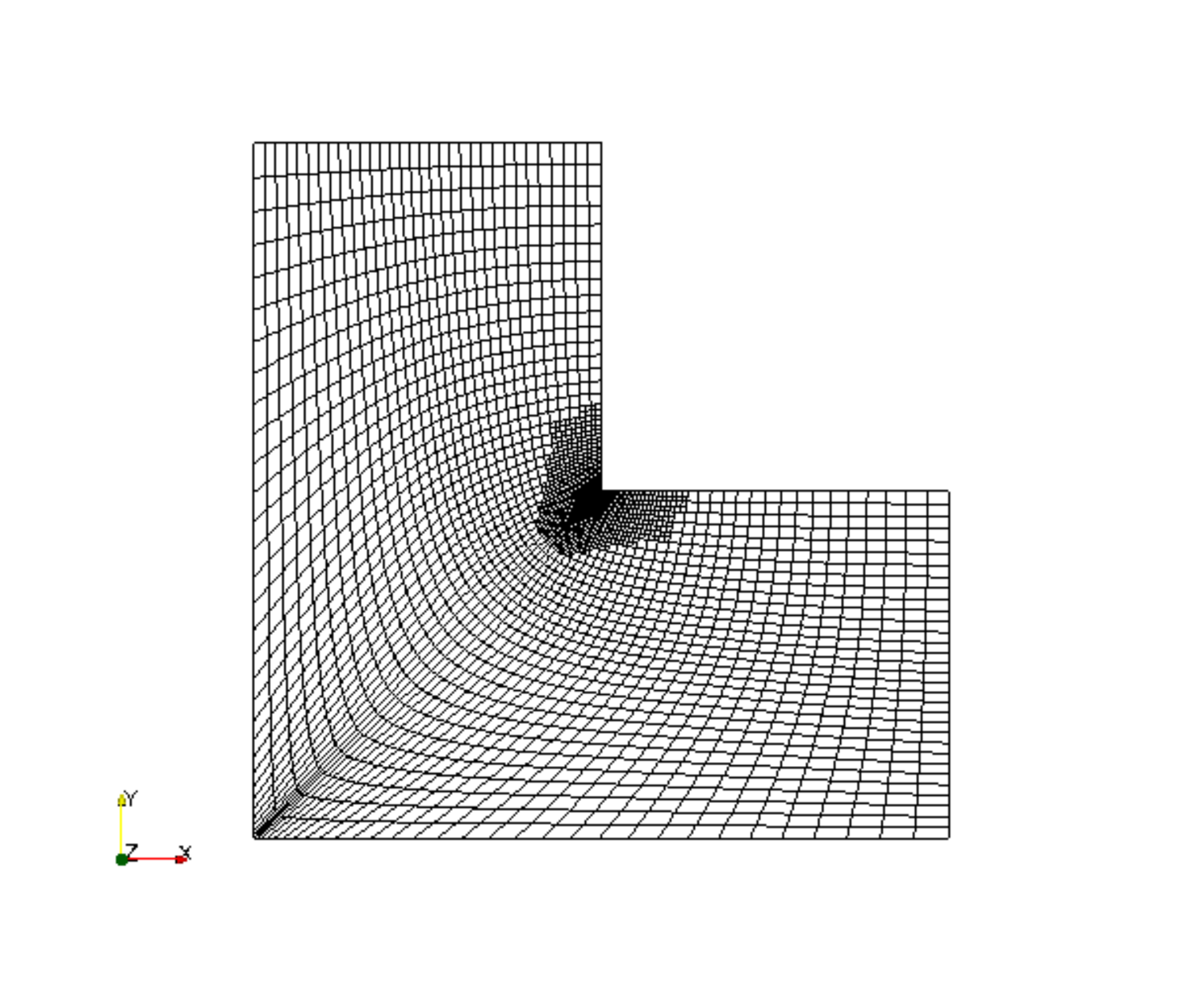}}  
	%
	}\\
	{
	\subfloat[ref. \# 4: \quad $\widehat{\mathcal{K}}_h$ $\leftarrow$  $\| \nabla e\|^2_K$]{
	\includegraphics[width=4.5cm, trim={4cm 2cm 4cm 2cm}, clip]{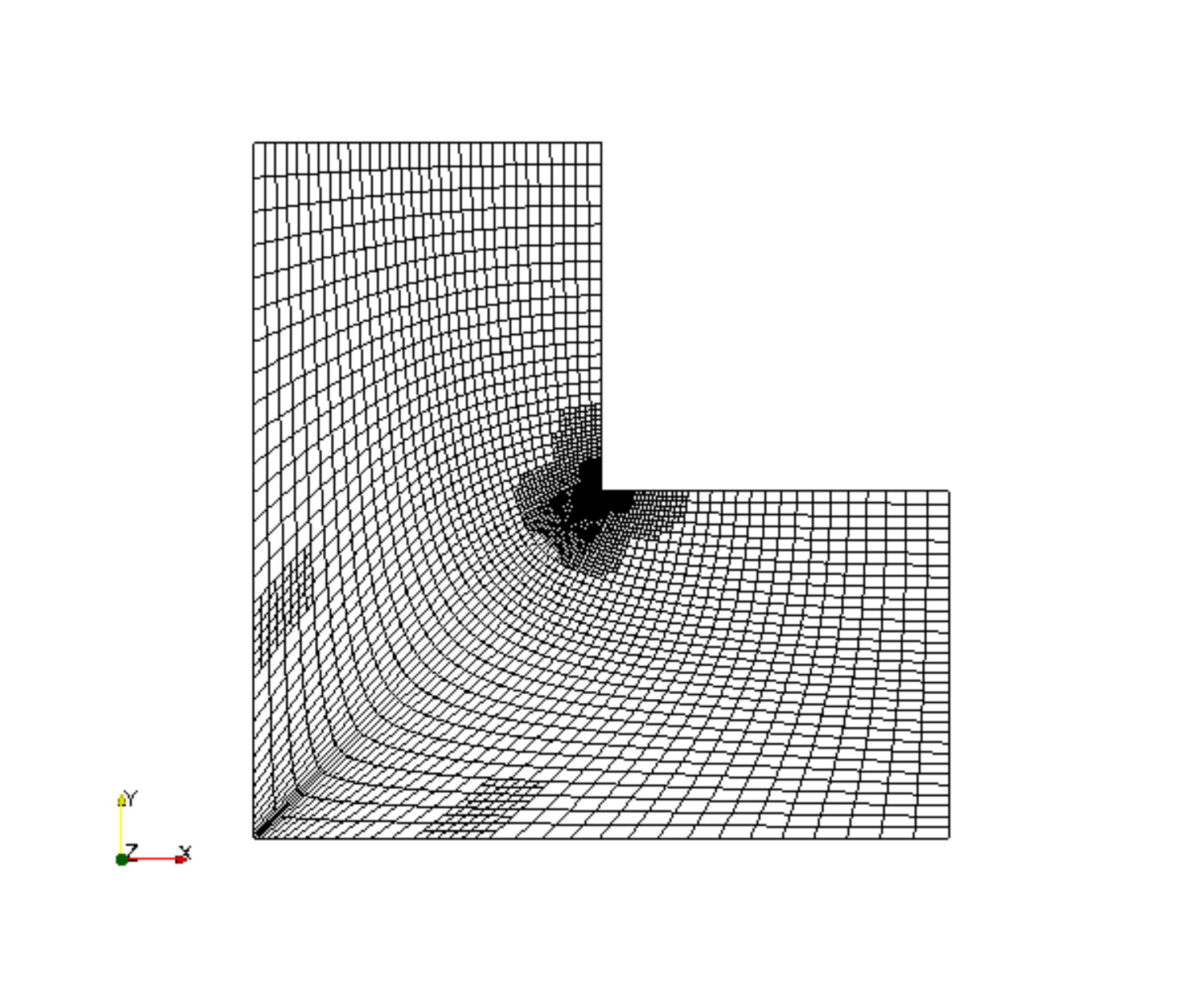}} \qquad
	\subfloat[ref. \# 4: \quad $\widehat{\mathcal{K}}_h$ $\leftarrow$  $\mdIK$]{
	\includegraphics[width=4.5cm, trim={4cm 2cm 4cm 2cm}, clip]{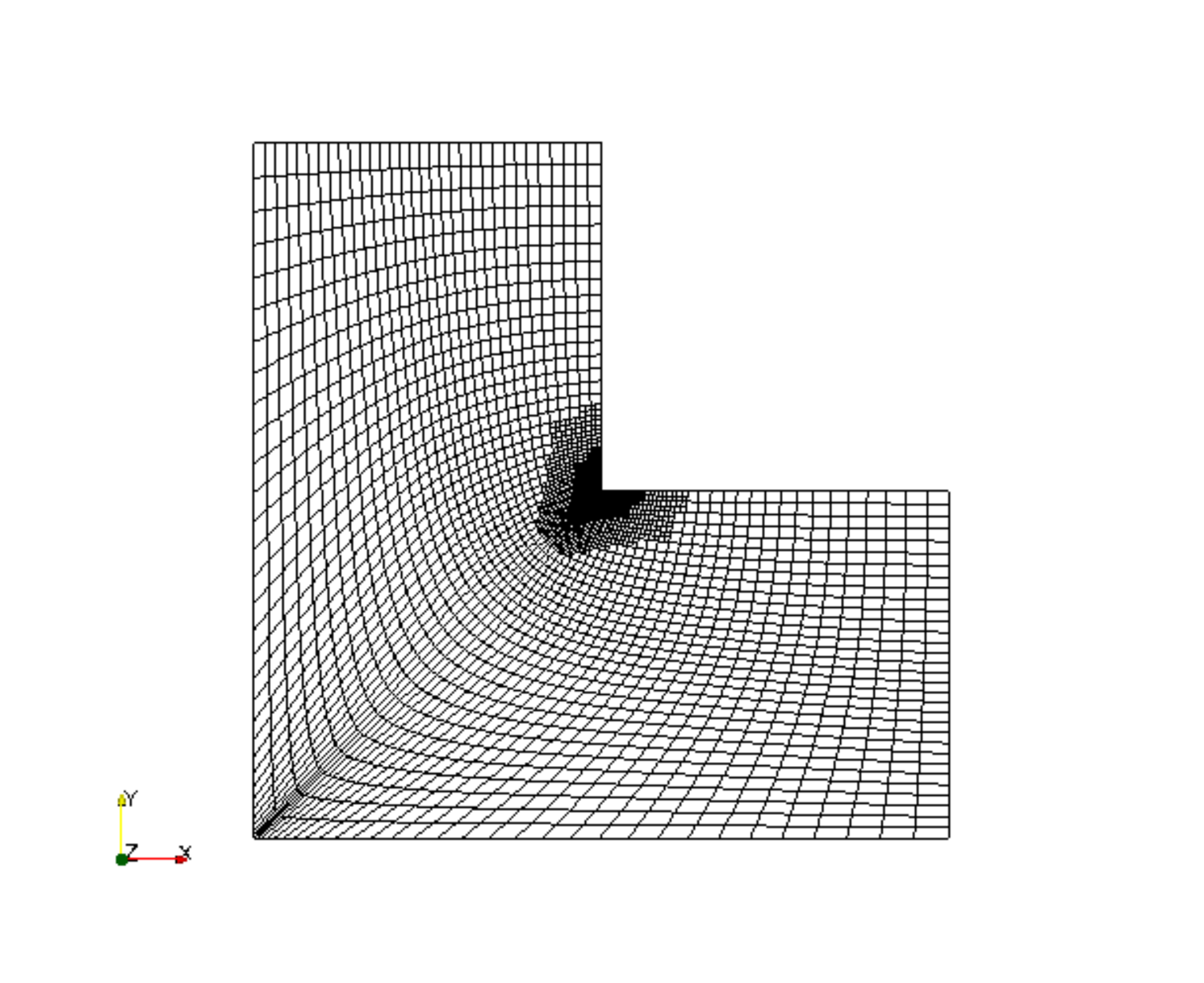}}  
	%
	}\\
	{\subfloat[ref. \# 6: \quad $\widehat{\mathcal{K}}_h$ $\leftarrow$  $\| \nabla e\|^2_K$]{
	\includegraphics[width=4.5cm, trim={4cm 2cm 4cm 2cm}, clip]{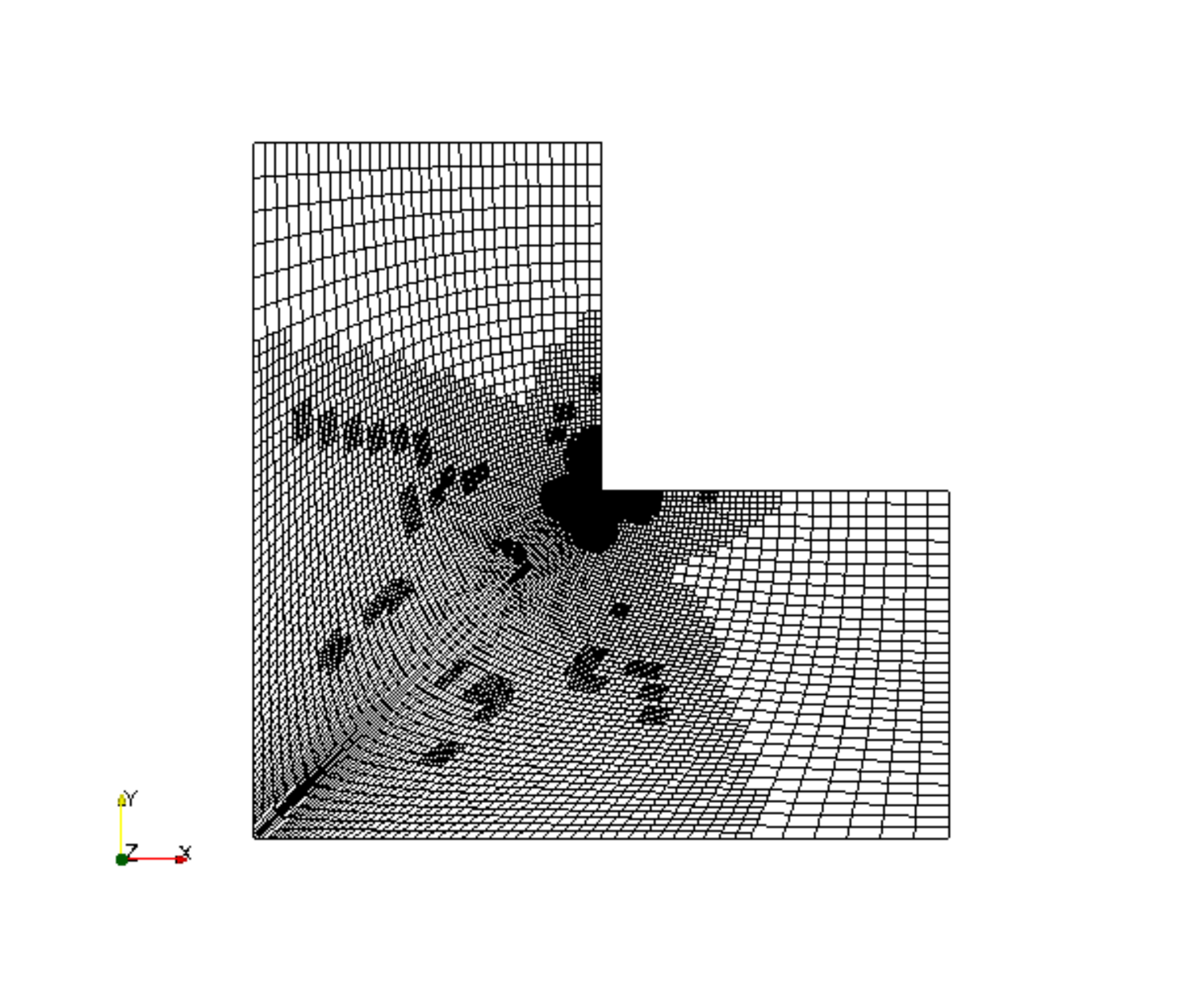}} \qquad
	\subfloat[ref. \# 6: \quad  $\widehat{\mathcal{K}}_h$ $\leftarrow$  $\mdIK$]{
	\includegraphics[width=4.5cm, trim={4cm 2cm 4cm 2cm}, clip]{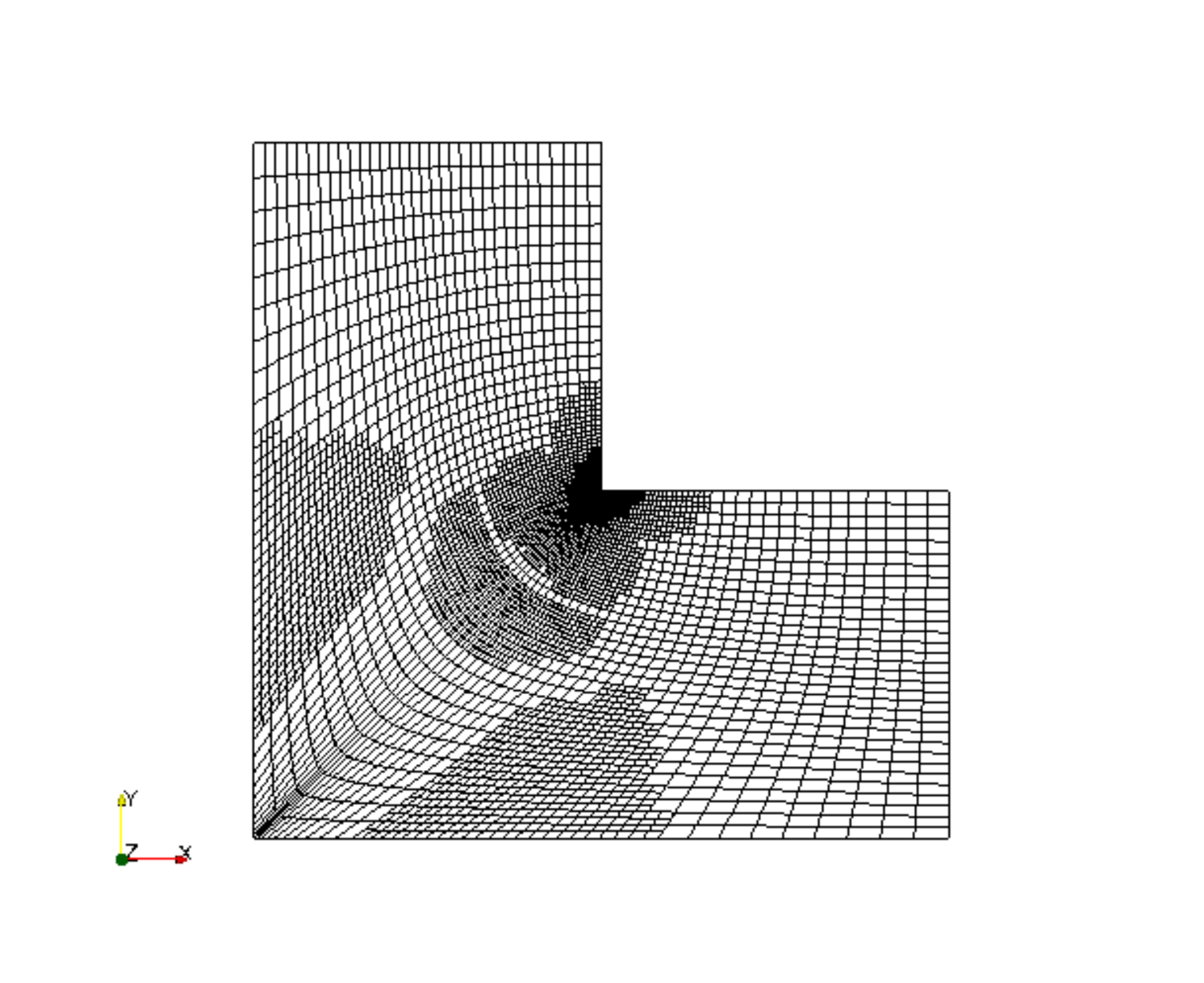}}  
	%
	}\\
	{\subfloat[ref. \# 7: \quad $\widehat{\mathcal{K}}_h$ $\leftarrow$  $\| \nabla e\|^2_K$]{
	\includegraphics[width=4.5cm, trim={4cm 2cm 4cm 2cm}, clip]{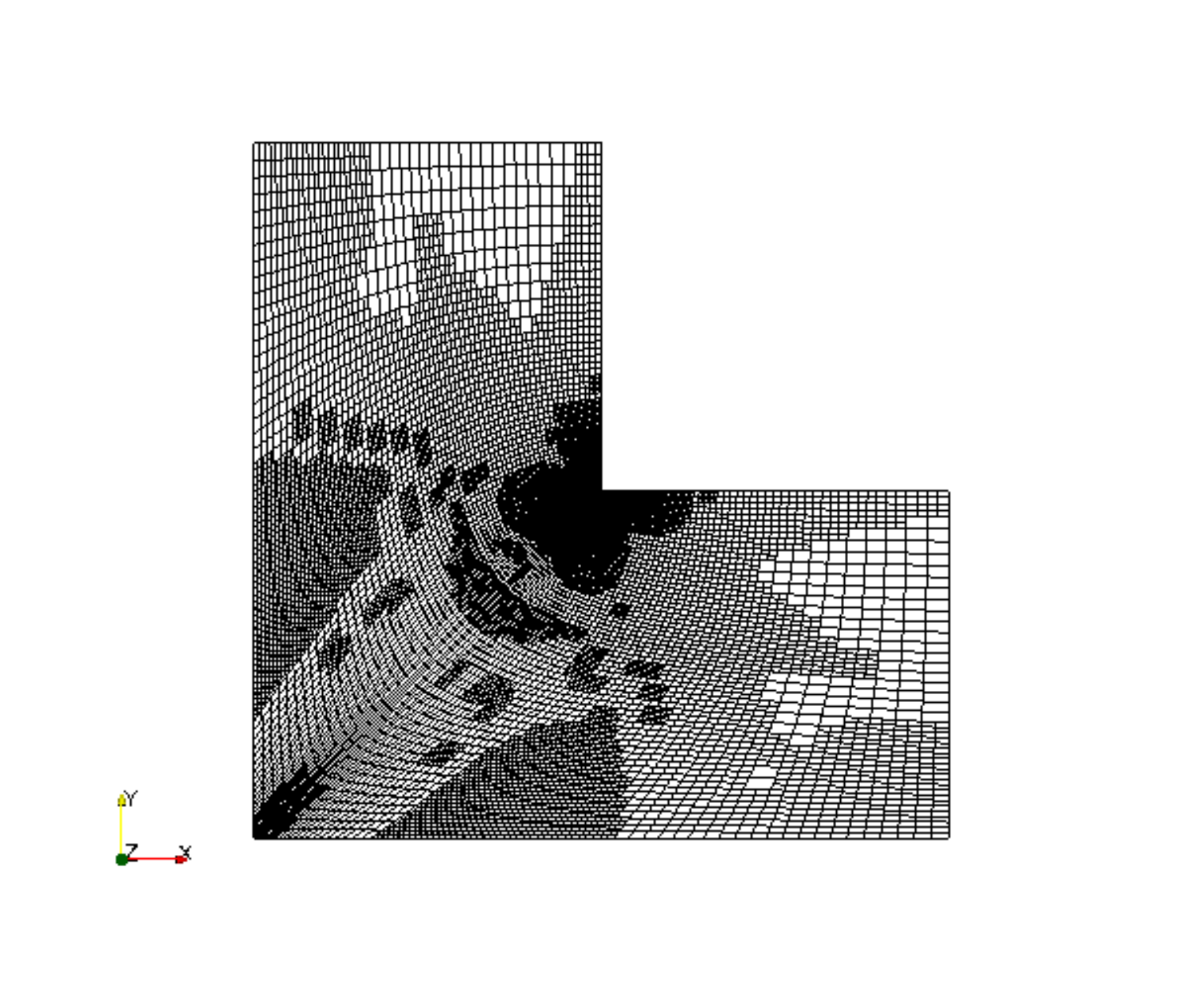}} \qquad
	\subfloat[ref. \# 7: \quad  $\widehat{\mathcal{K}}_h$ $\leftarrow$  $\mdIK$]{
	\includegraphics[width=4.5cm, trim={4cm 2cm 4cm 2cm}, clip]{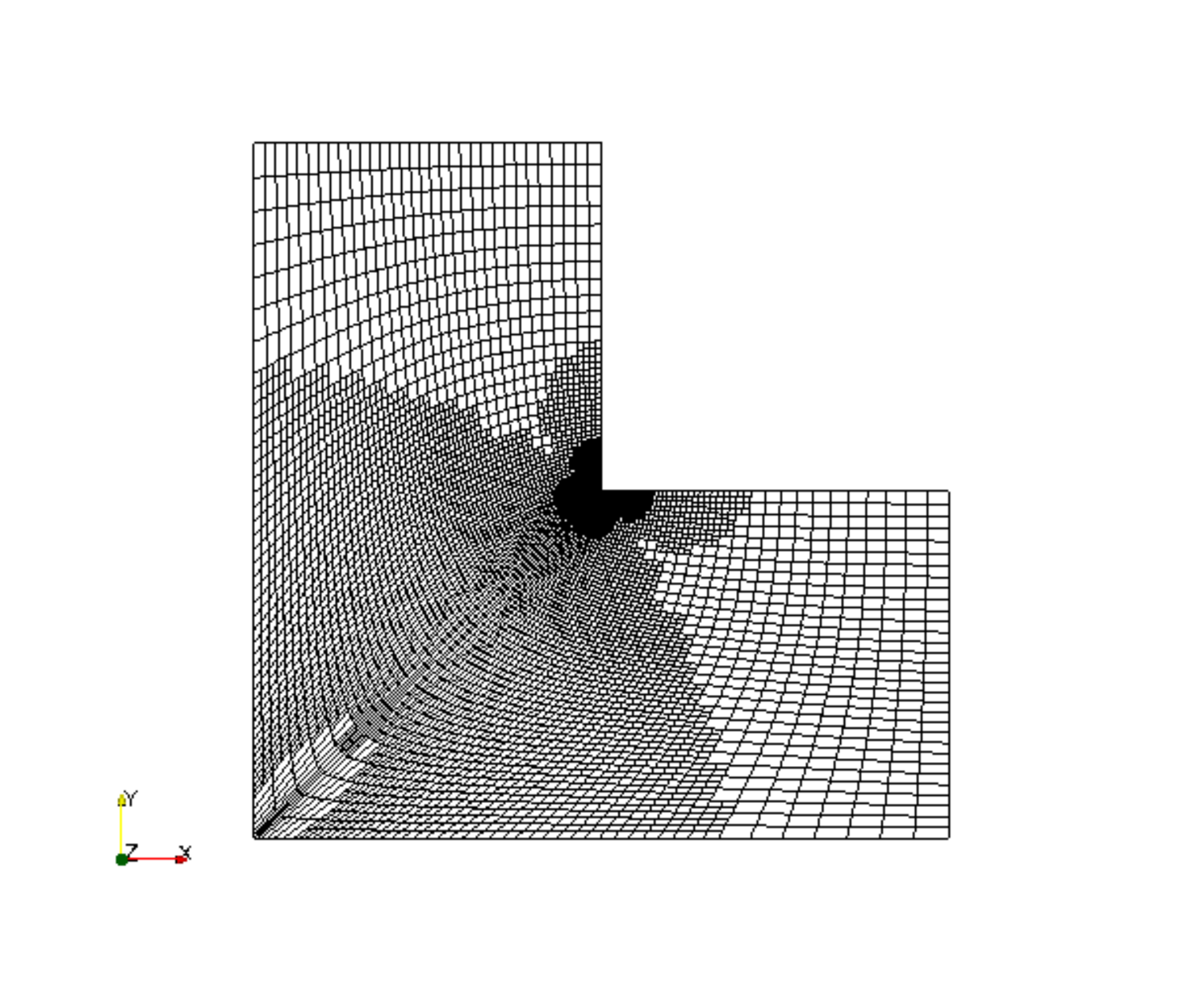}}}
	\caption{Comparison of the meshes generated by two refinement strategies, i.e., the true error (left) and 
	error indicator provided by the majorant (right), with the marking criterion ${\mathds{M}}_{{\rm \bf BULK}}(0.2)$.}
	\label{fig:l-shape-domain-example-12-theta-20-error-majorant-refinement}
\end{figure}

\begin{figure}[!t] 
	\centering
	{
	\subfloat[ref. \# 4: \quad $\widehat{\Omega}$ and $\widehat{\mathcal{K}}_h$]{
	\includegraphics[width=4.5cm, trim={8cm 2cm 8cm 2cm}, clip]{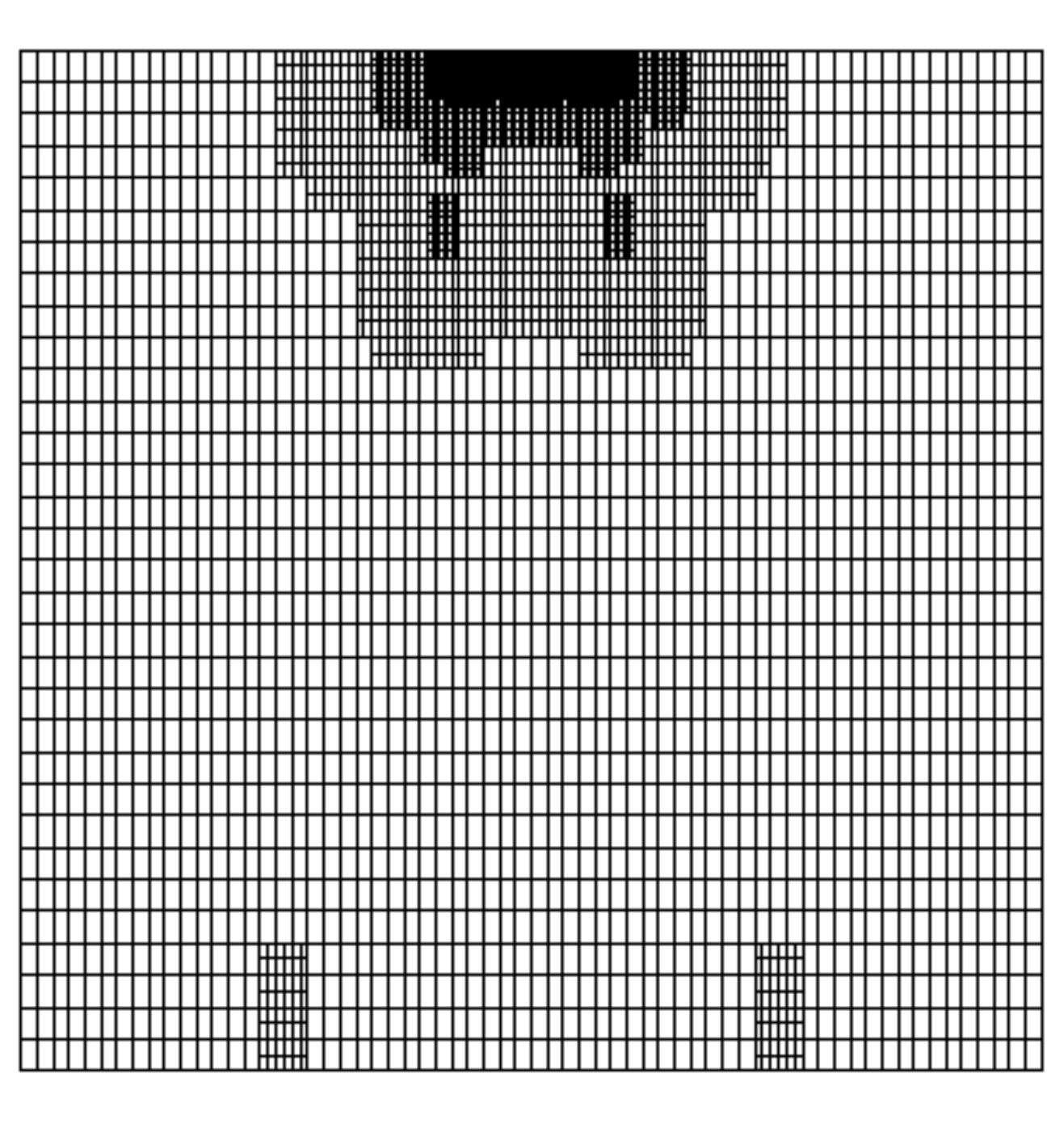}}
	\qquad
	\subfloat[ref. \# 4: \quad  $\Omega$ and $\mathcal{K}_h$]{
	\includegraphics[width=4.5cm, trim={8cm 2cm 8cm 2cm}, clip]{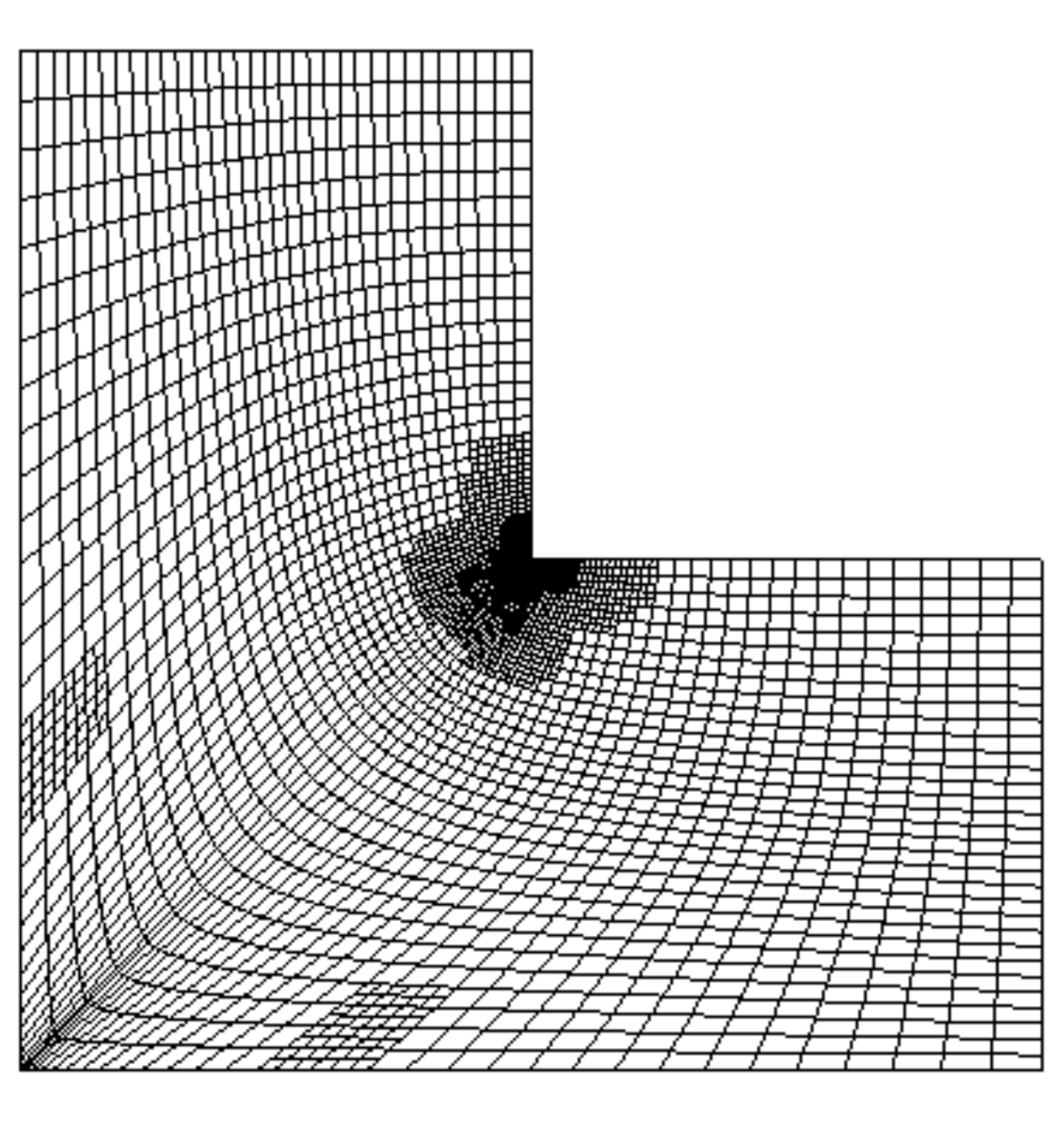}} 
	}\\
	{
	\subfloat[ref. \# 5: \quad $\widehat{\Omega}$ and $\widehat{\mathcal{K}}_h$]{
	\includegraphics[width=4.5cm, trim={8cm 2cm 8cm 2cm}, clip]{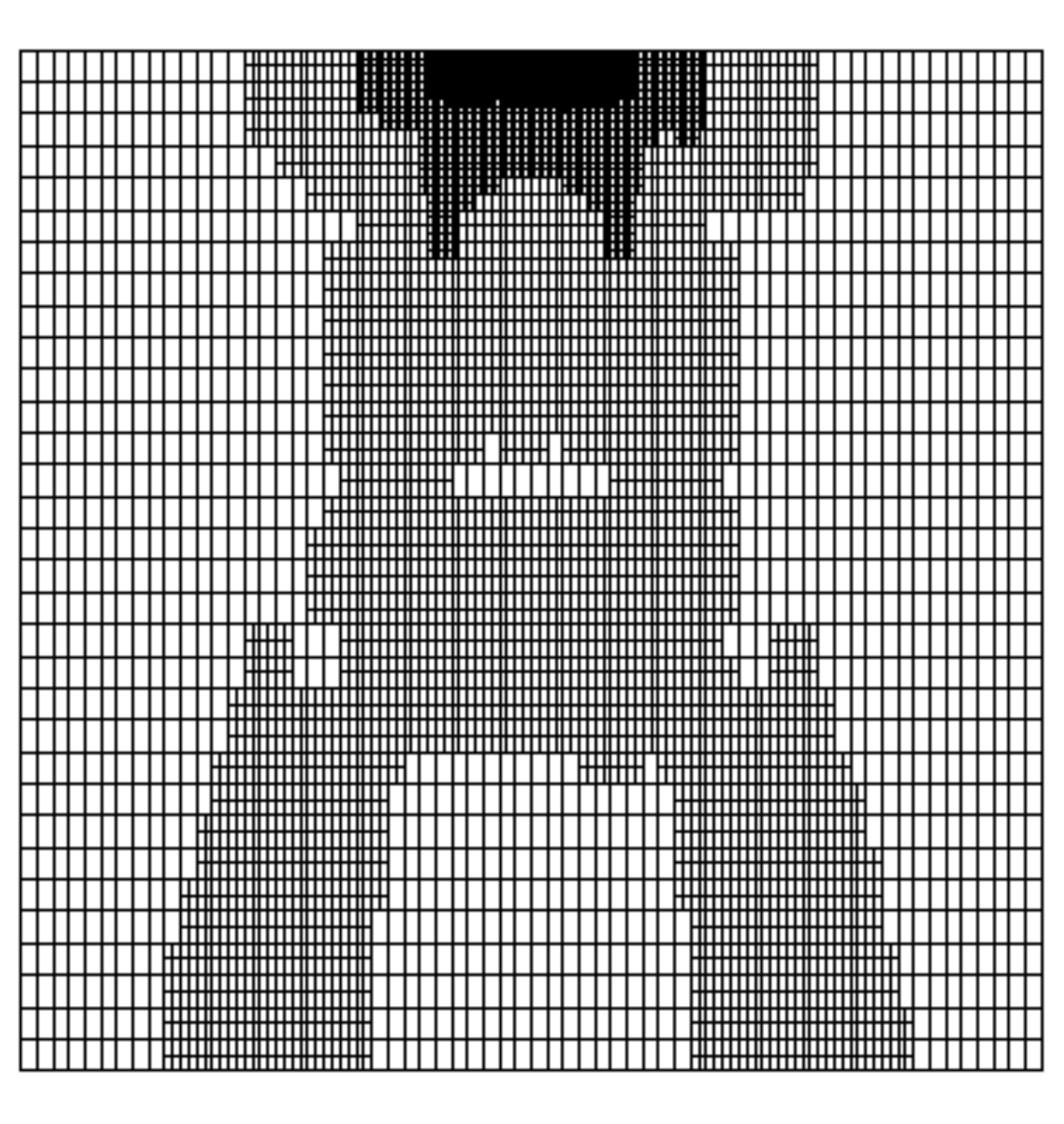}}
	\qquad
	\subfloat[ref. \# 5: \quad  $\Omega$ and $\mathcal{K}_h$]{
	\includegraphics[width=4.5cm, trim={8cm 2cm 8cm 2cm}, clip]{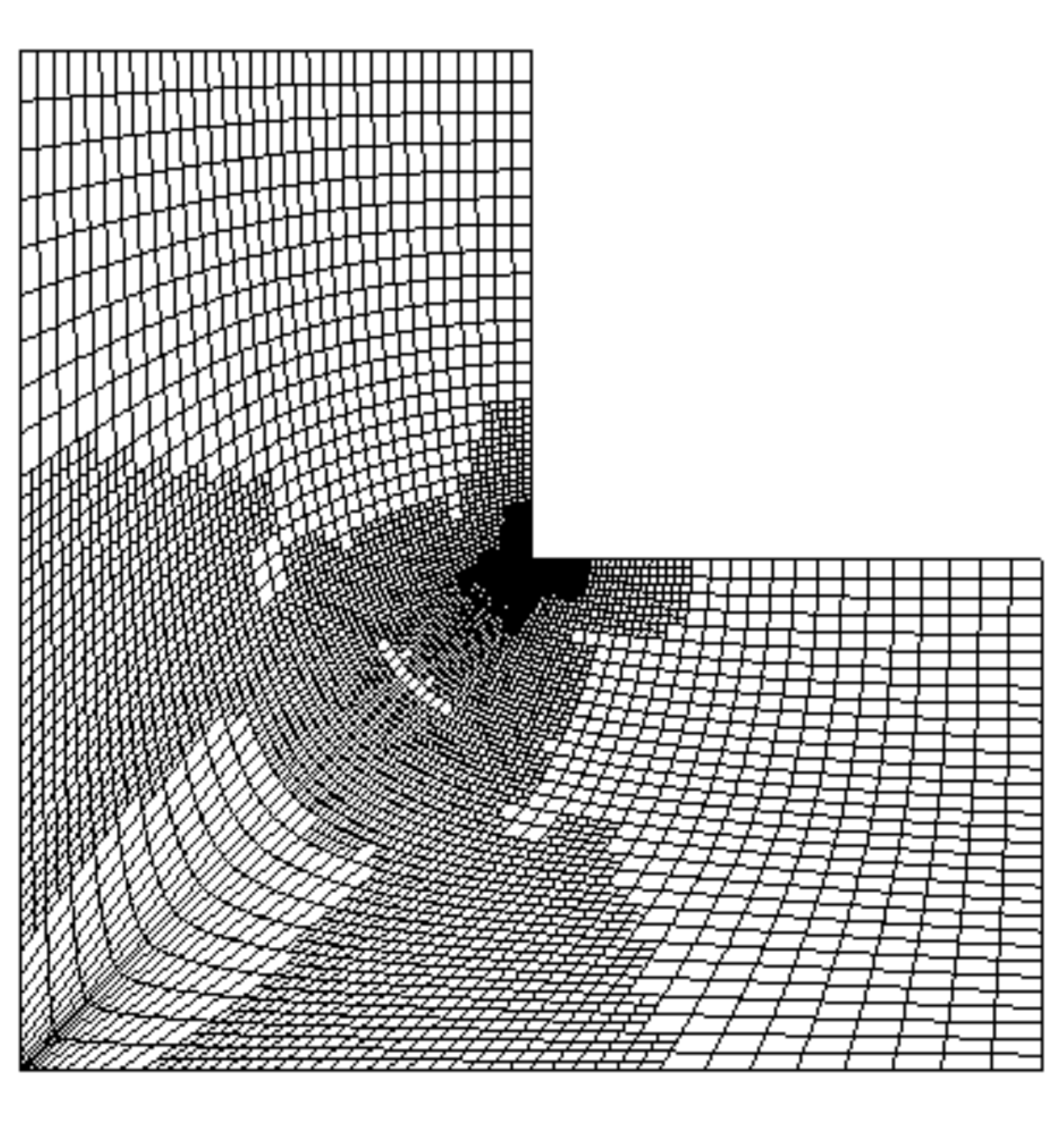}} 
	}\\
	{
	\subfloat[ref. \# 6: \quad $\widehat{\Omega}$ and $\widehat{\mathcal{K}}_h$]{
	\includegraphics[width=4.5cm, trim={8cm 2cm 8cm 2cm}, clip]{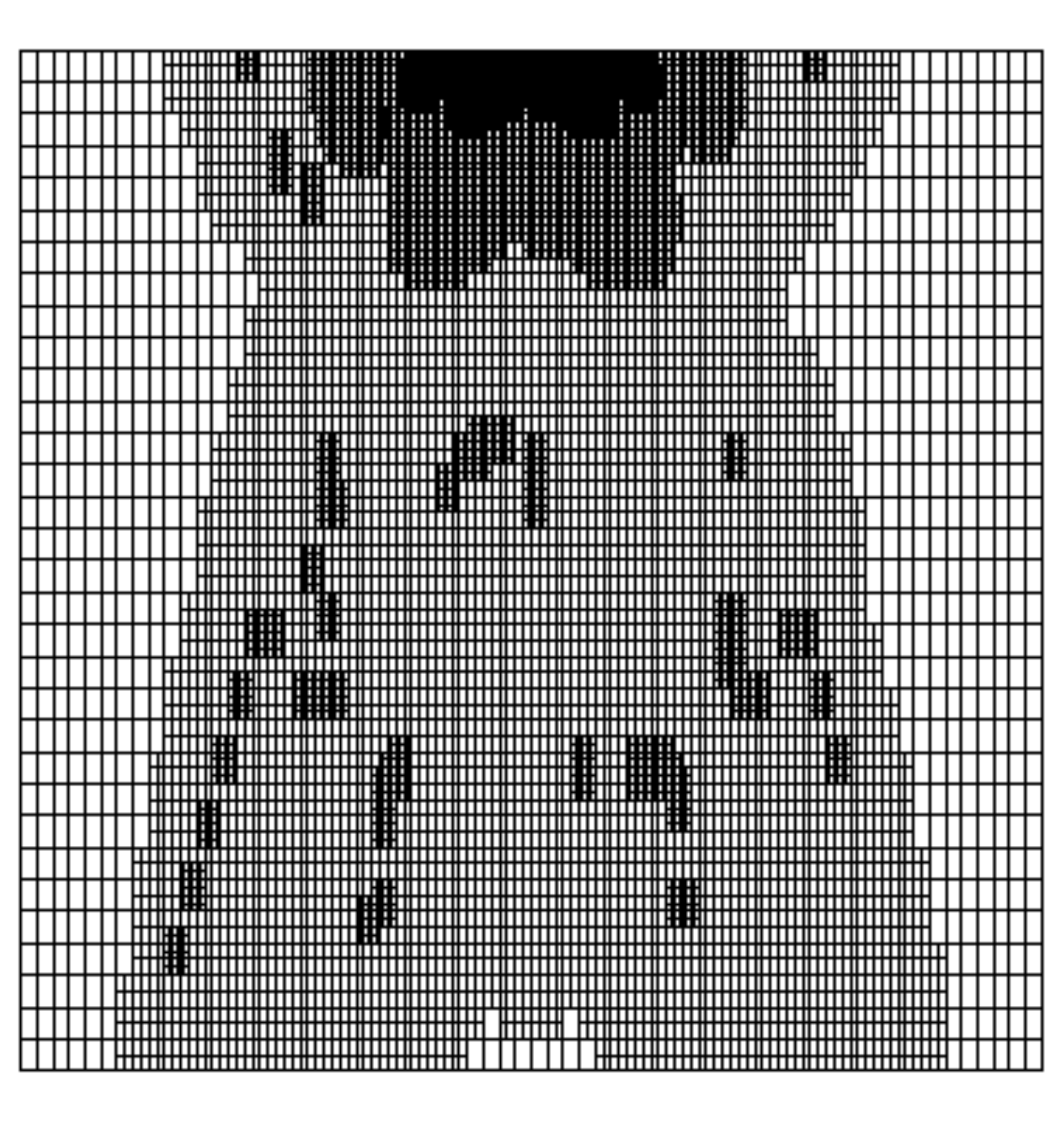}}
	\qquad
	\subfloat[ref. \# 6: \quad $\Omega$ and $\mathcal{K}_h$]{
	\includegraphics[width=4.5cm, trim={8cm 2cm 8cm 2cm}, clip]{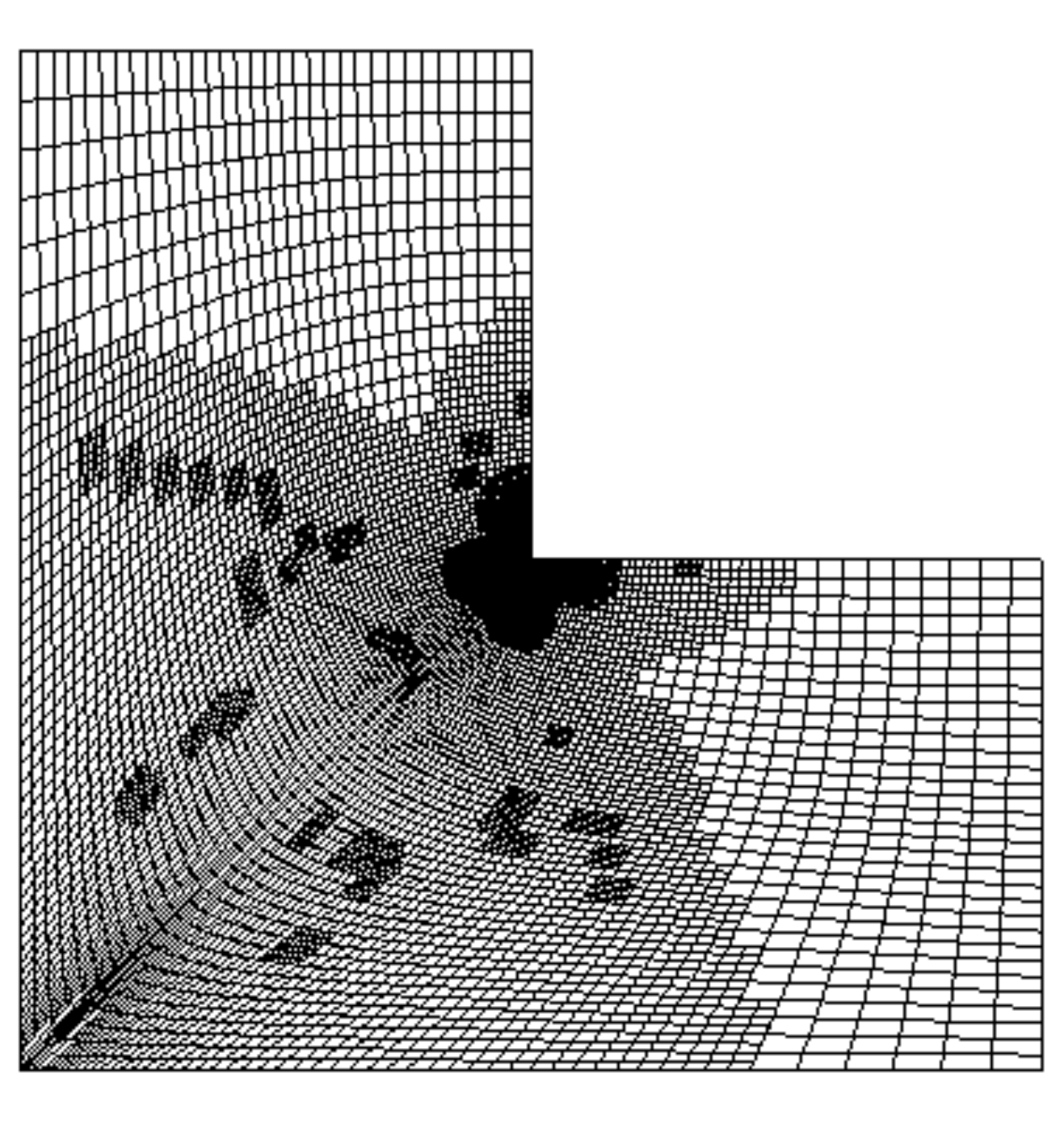}} 
	}\\
	{
	\subfloat[ref. \# 7: \quad $\widehat{\Omega}$ and $\widehat{\mathcal{K}}_h$]{
	\includegraphics[width=4.5cm, trim={8cm 2cm 8cm 2cm}, clip]{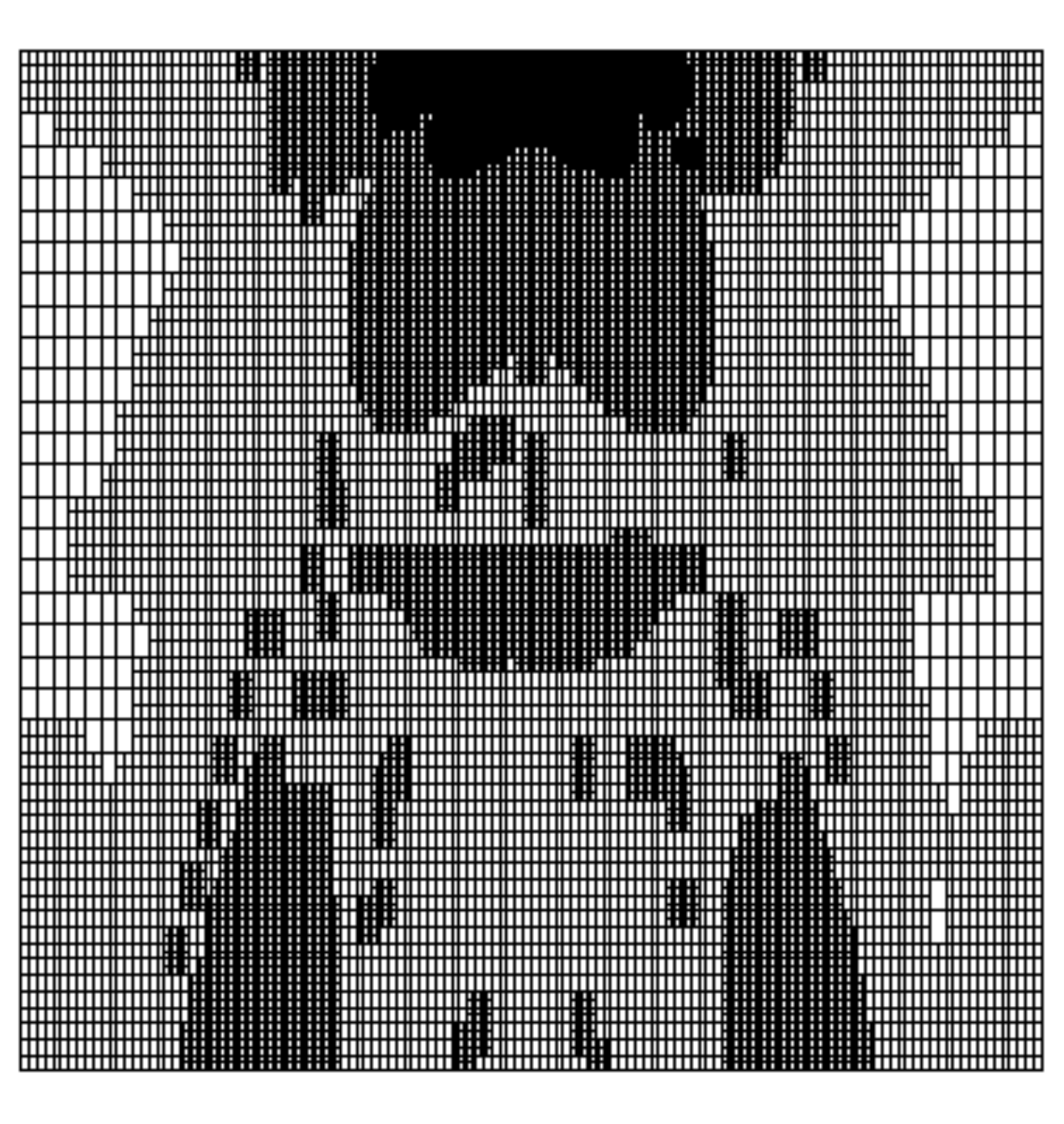}}
	\qquad
	\subfloat[ref. \# 7: \quad $\Omega$ and $\mathcal{K}_h$]{
	\includegraphics[width=4.5cm, trim={8cm 2cm 8cm 2cm}, clip]{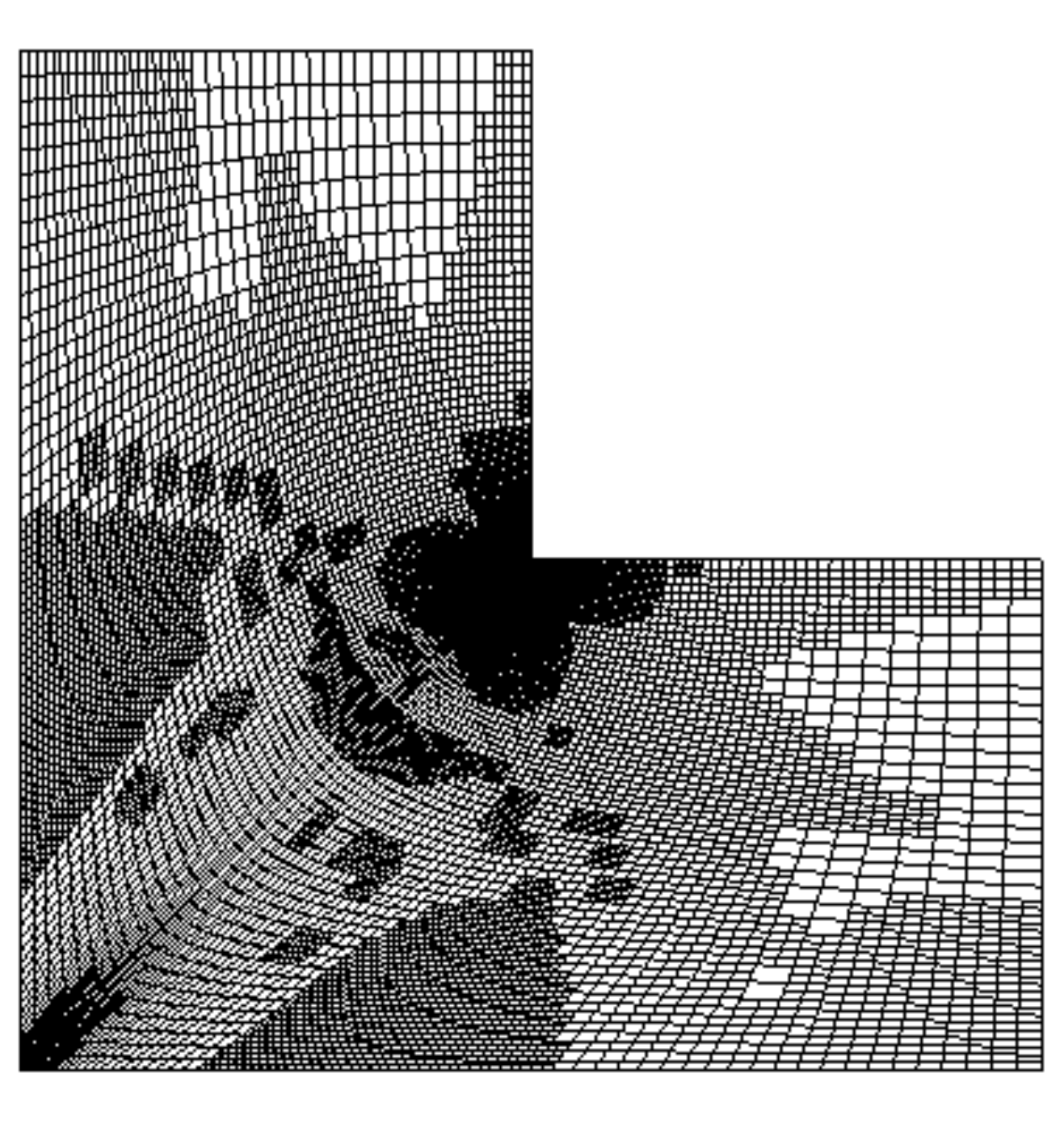}} 
	}
	\caption{\small {\em Ex. \ref{ex:l-shape-domain-example-12}}. 
	Comparison of meshes on the physical and parametrical domains w.r.t. adaptive ref. steps, 
	${\mathds{M}}_{{\rm \bf BULK}}(0.1)$.}
	\label{fig:l-shape-domain-example-12-theta-10-domains-comparison}
\end{figure}
\begin{figure}[!t]
	\centering
	\includegraphics[scale=0.7]{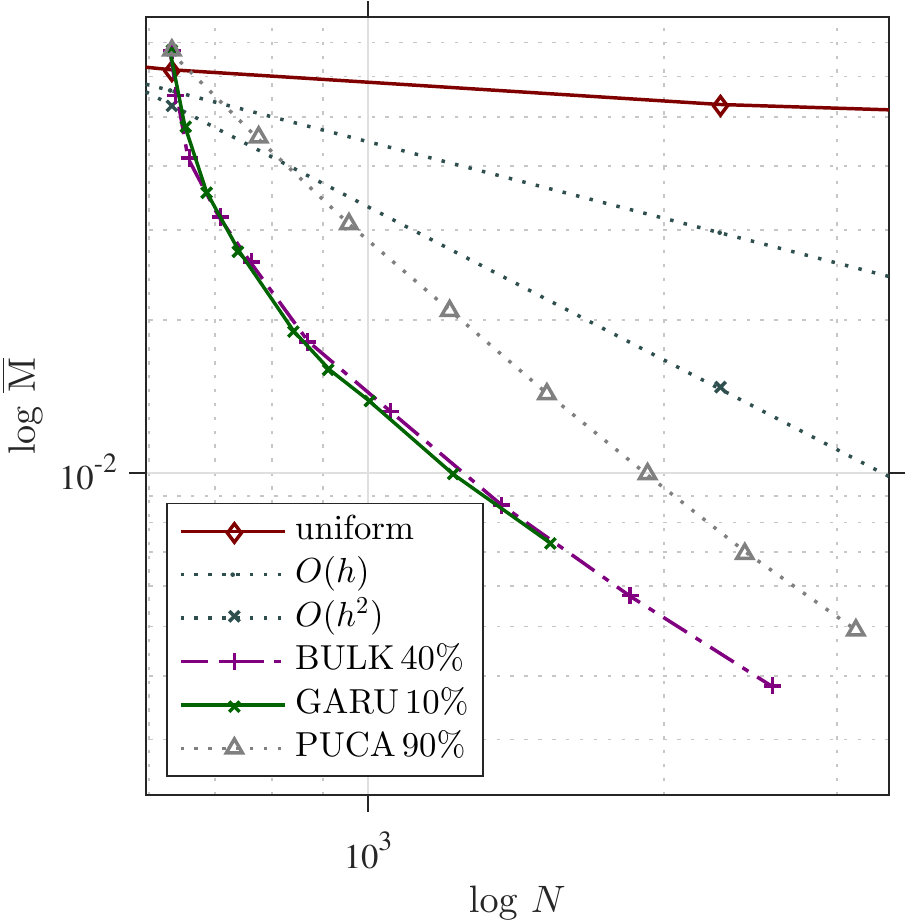}
	\caption{Convergence of majorant.}
	\label{fig:l-shape-domain-example-9-convergence-majorant}
\end{figure}

\end{example}


\newpage
\begin{example}
\label{ex:quater-annulus-example-10}
\rm 
\begin{figure}[!h]
	\centering
	\includegraphics[scale=0.7]{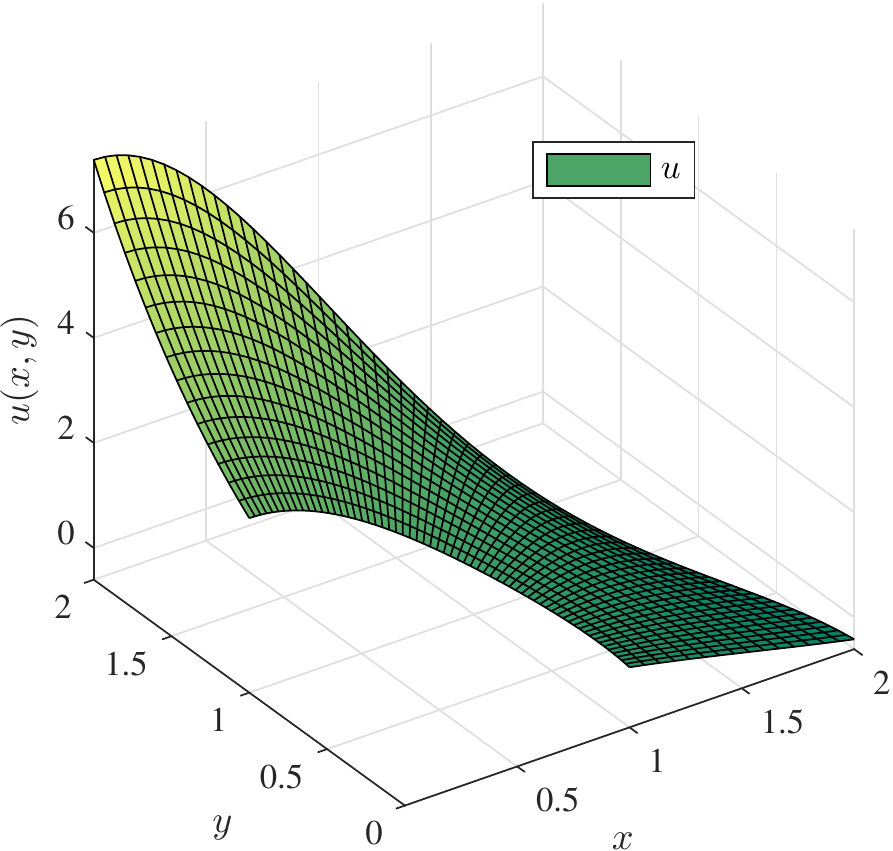}
	\caption{\small {\em Ex. \ref{ex:quater-annulus-example-10}}. 
	Exact solution $cos (x_1) \, e^{x_2}$.}
	\label{fig:example-10-exact-solution}
\end{figure}
                       
The final example with two-dimensional domain is defined on a quoter-annulus with the following exact solution 
and RHS:
\begin{alignat*}{2}
u 	& = \cos x_1 \, e^{x_2}			\quad \;\;\mbox{in} \quad \Omega, \\
f 	& = 0						\qquad\qquad\quad \; \mbox{in} \quad \Omega, \\
u_D 	& = \cos x_1 \, e^{x_2}			\quad \;\;\; \mbox{on} \quad \Gamma.
\end{alignat*}
Due to the IgA framework, the Dirichlet BC for the approximation $u_h$ are fully satisfied, therefore, functional 
error estimates can be applied for the domains with curved boundaries providing a fully reliable (non-heuristic) error 
control. 
The results obtained for the adaptive refinement for $u_h \in S_h^{2, 2}$, 
$\flux \in S_{4h}^{4, 4} \oplus S_{4h}^{4, 4}$, and bulk parameter $\theta = 0.4$, are illustrated in Tables 
\ref{tab:quater-annulus-example-10-error-majorant-v-2-y-4-adaptive-ref}--\ref{tab:quater-annulus-example-10-times-v-2-y-4-adaptive-ref}.
Again, even for such specific geometry, the generation of the optimal $\flux_h$ (assembling and solving the 
\eqref{eq:system-fluxh}) requires several times less computational effort than $u_h$: 
$
\tfrac{t_{\rm as}(u_h)}{t_{\rm as}(\flux_h)} \approx \tfrac{38.5360}{19.1958} \approx 2
\; \mbox{and} \;
\tfrac{t_{\rm sol}(u_h)}{t_{\rm sol}(\flux_h)} \approx \tfrac{3.4323}{0.4727} \approx 7.$
Moreover, in Figure 
\ref{fig:quater-annulus-example-10-theta-40-domains-comparison}, 
we illustrate the evolution of the meshes $\widehat{\mathcal{K}}_h$ and ${\mathcal{K}}_h$ on the 
parametric $\widehat{\Omega}$ and physical $\Omega$ domains, respectively. 

\begin{table}[!t]
\footnotesize
\centering
\newcolumntype{g}{>{\columncolor{gainsboro}}c} 	
\begin{tabular}{c|c|ccc|gc|c}
\quad \# ref. \quad & 
\quad  $\| \nabla e \|_\Omega$ \qquad   & 	  
\quad \quad \quad \quad $\overline{\rm M}$ \quad \quad \quad \quad &    
\quad \quad   $\mdI$ \qquad \quad & 	       
\quad \qquad $\mfI$ \qquad \quad  &  
\qquad $\Ieff (\overline{\rm M})$ \qquad & 
\qquad $\Ieff (\overline{\rm \eta})$ \qquad & 
\qquad \quad e.o.c. \qquad \quad \\
\midrule
  3 &   1.9377e-03 &   2.2682e-03 &   2.1466e-03 &   2.7011e-04 &       1.1706 &       8.0838 &   1.8290 \\
   5 &   2.7731e-04 &   3.9296e-04 &   3.7669e-04 &   3.6143e-05 &       1.4171 &       8.2886 &   2.1330 \\
   7 &   4.7689e-05 &   6.6456e-05 &   5.5448e-05 &   2.4453e-05 &       1.3935 &       8.3850 &   2.1653  \\
\end{tabular}
\caption{\small {\em Ex. \ref{ex:quater-annulus-example-10}}. 
Error, majorant (with dual and reliability terms), 
efficiency indices, error, e.o.c. 
w.r.t. adaptive ref. steps.}
\label{tab:quater-annulus-example-10-error-majorant-v-2-y-4-adaptive-ref}
\end{table}
\begin{table}[!t]
 \footnotesize
\centering
\newcolumntype{g}{>{\columncolor{gainsboro}}c} 	
\begin{tabular}{c|cc|cg|cg|cgc}
\# ref.  & 
\# d.o.f.($u_h$) &  \# d.o.f.($\flux_h$) &  
\; $t_{\rm as}(u_h)$ \; & 
\; $t_{\rm as}(\flux_h)$ \; & 
\; $t_{\rm sol}(u_h)$ \; & 
\; $t_{\rm sol}(\flux_h)$ \; &
$t_{\rm e/w}(\| \nabla e \|)$ & 
$t_{\rm e/w}(\overline{\rm M})$ & 
$t_{\rm e/w}(\overline{\eta})$ \\
\midrule
   1 &        324 &        400 &     0.1059 &     1.7198 &           0.0010 &           0.0245 &       0.1592 &       0.7331 &       0.3351 \\
   3 &       1389 &        400 &     0.4475 &     1.0655 &           0.0081 &           0.0166 &       0.5914 &       0.9428 &       1.1633 \\
   5 &       9125 &        400 &     5.0846 &     1.1225 &           0.1925 &           0.0363 &       5.1422 &       9.5749 &       9.1553 \\
   7 &      50291 &       1623 &    38.5360 &    19.1958 &           3.4323 &           0.4727 &      21.0657 &      96.3344 &      41.5038 \\\end{tabular}
\caption{\small {\em Ex. \ref{ex:quater-annulus-example-10}}. 
Time for assembling and solving the systems generating d.o.f. of $u_h$ and $\flux_h$ 
as well as the time spent on e/w evaluation of error, majorant, and 
residual error estimator w.r.t. adaptive ref. steps.}
\label{tab:quater-annulus-example-10-times-v-2-y-4-adaptive-ref}
\end{table}

\begin{figure}[!t]
	\centering
	%
	%
	{
	\subfloat[ref. \# 3: \quad $\widehat{\Omega}$ and $\widehat{\mathcal{K}}_h$]{
	\includegraphics[width=5.4cm, trim={2cm 2cm 2cm 2cm}, clip]{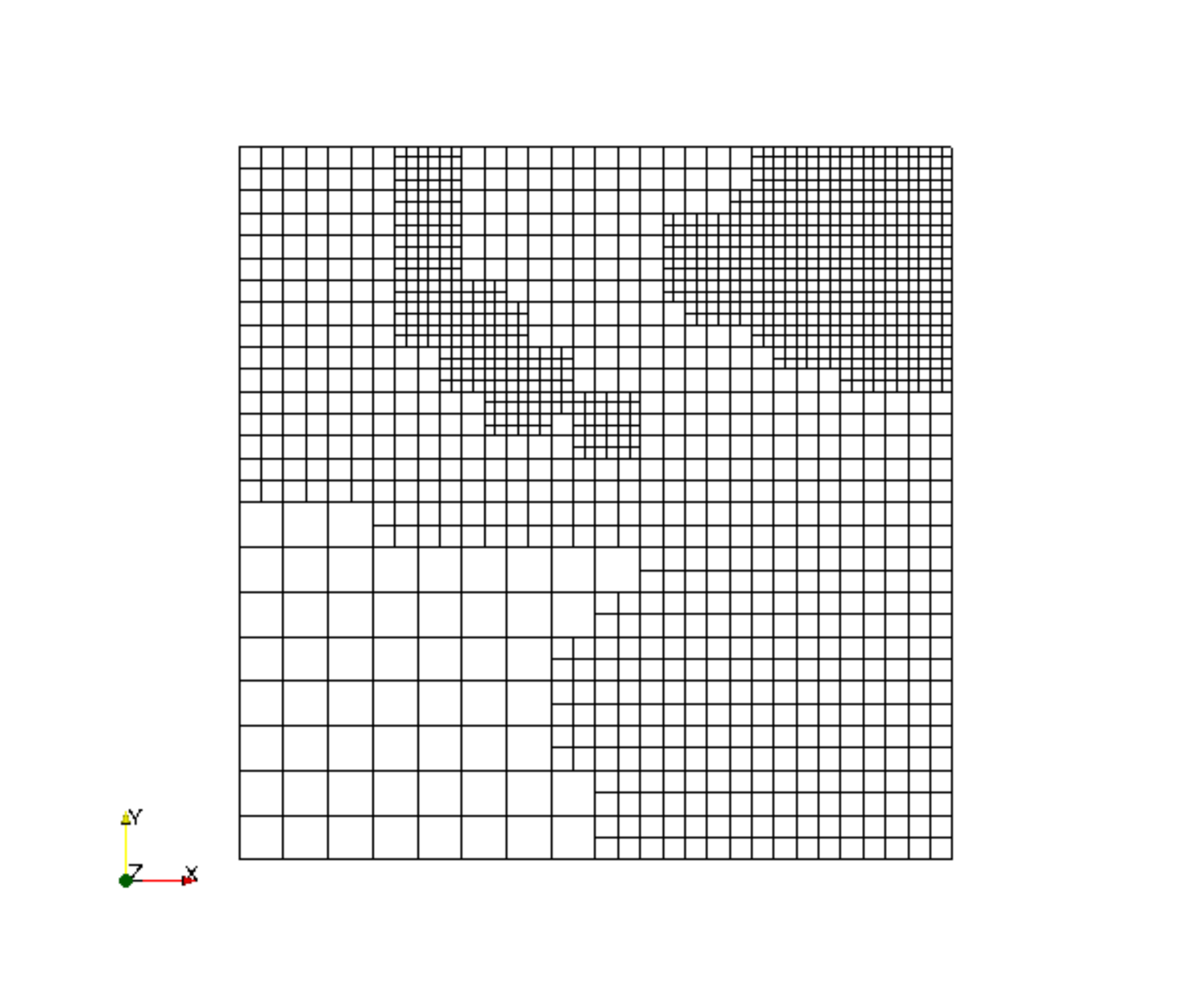}}
	\qquad
	\subfloat[ref. \# 3: \quad $\Omega$ and $\mathcal{K}_h$]{
	\includegraphics[width=5.4cm, trim={2cm 2cm 2cm 2cm}, clip]{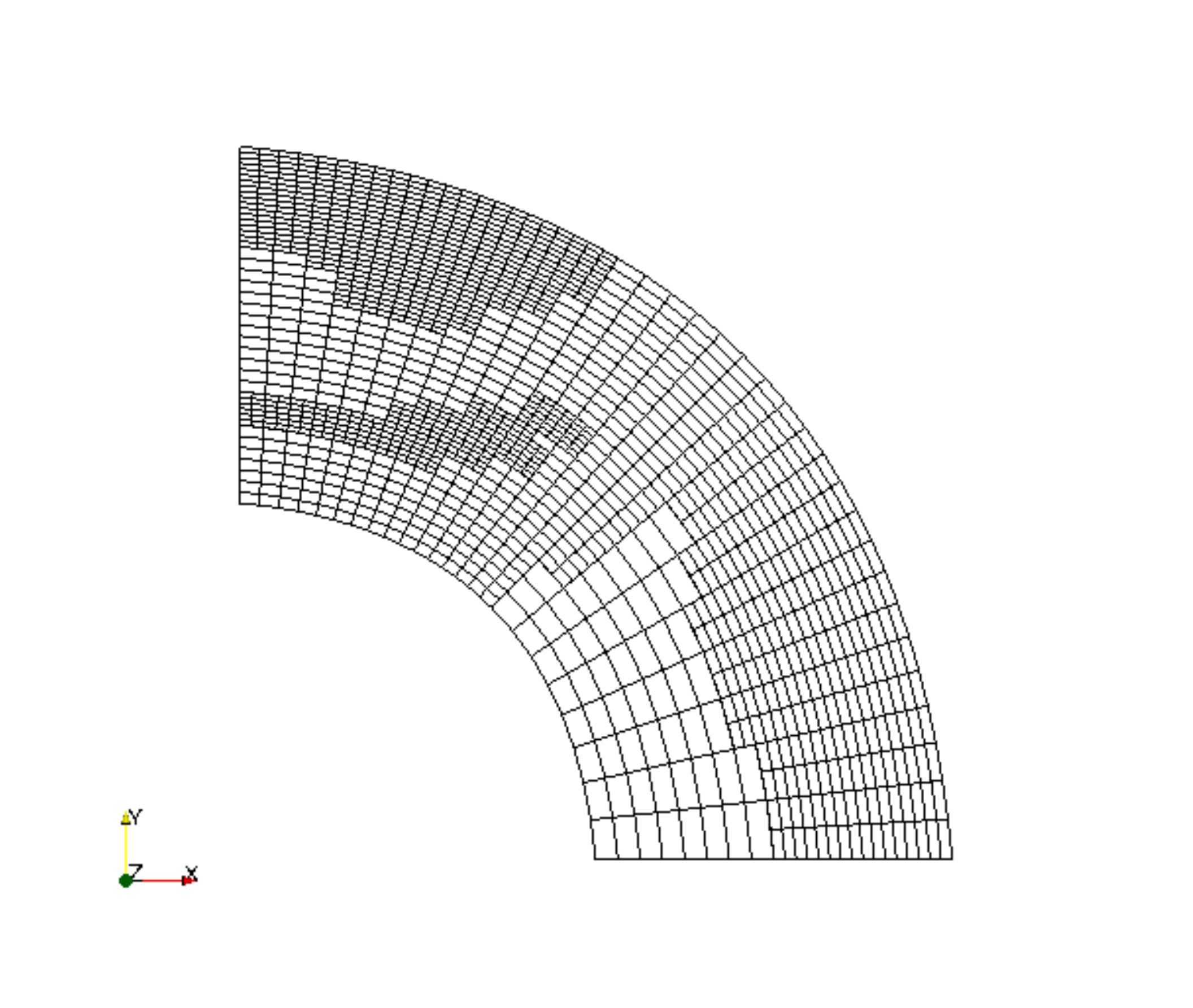}} 
	} \\
	{
	\subfloat[ref. \# 4: \quad $\widehat{\Omega}$ and $\widehat{\mathcal{K}}_h$]{
	\includegraphics[width=5.4cm, trim={2cm 2cm 2cm 2cm}, clip]{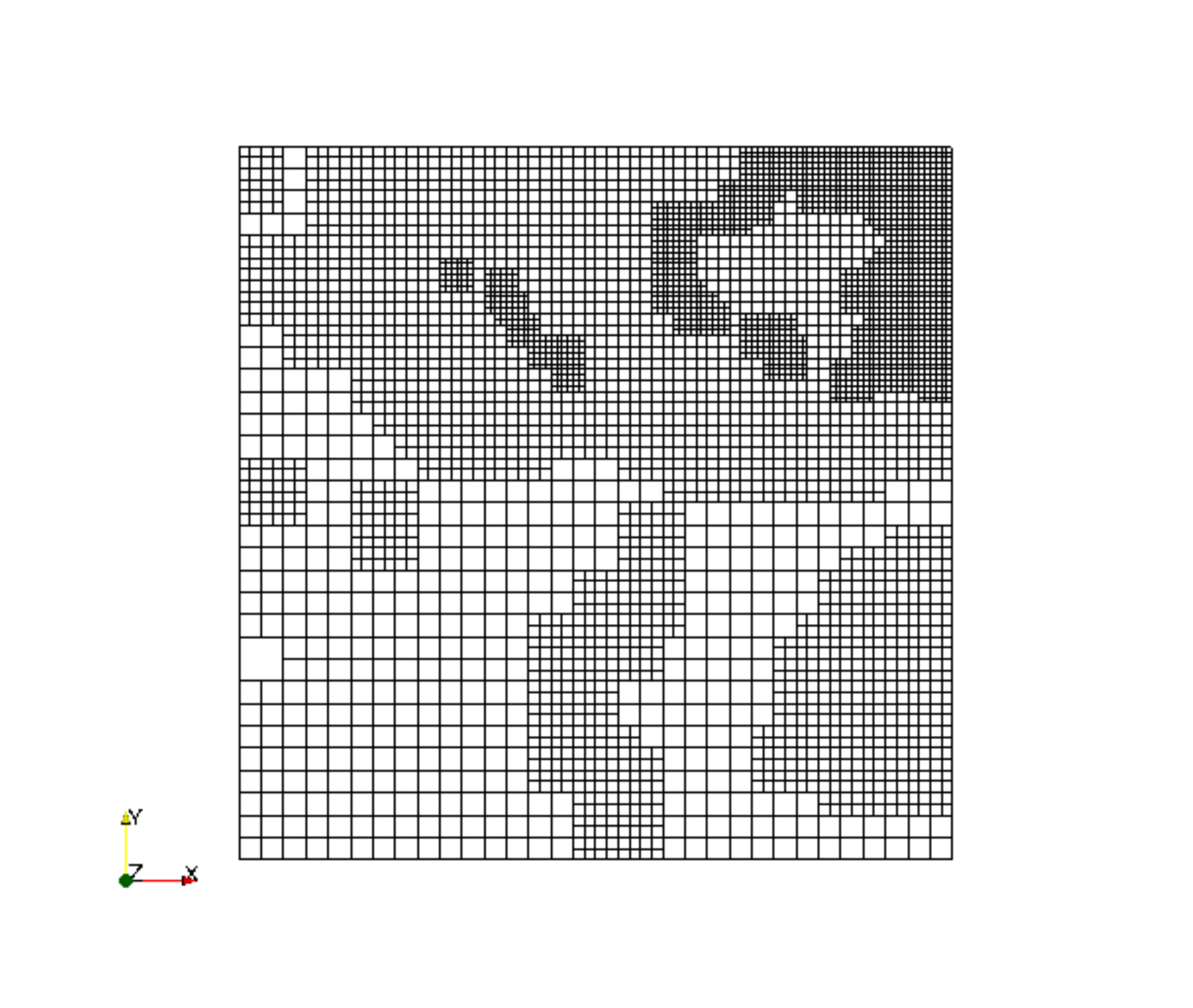}}
	\qquad
	\subfloat[ref. \# 4: \quad ${\Omega}$ and ${\mathcal{K}}_h$]{
	\includegraphics[width=5.4cm, trim={2cm 2cm 2cm 2cm}, clip]{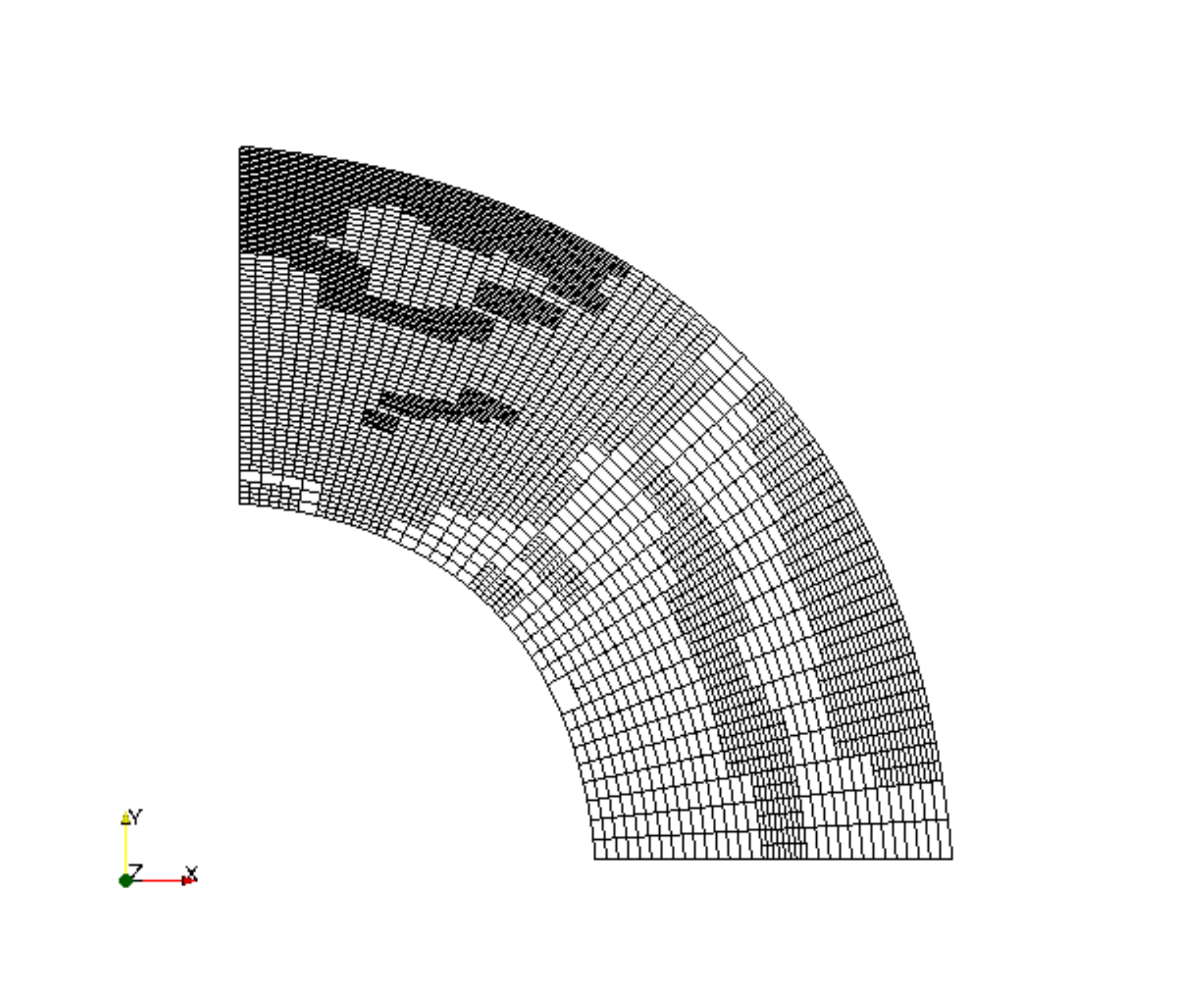}} 
	} \\
	{
	\subfloat[ref. \# 5: \quad $\widehat{\Omega}$ and $\widehat{\mathcal{K}}_h$]{
	\includegraphics[width=5.4cm, trim={2cm 2cm 2cm 2cm}, clip]{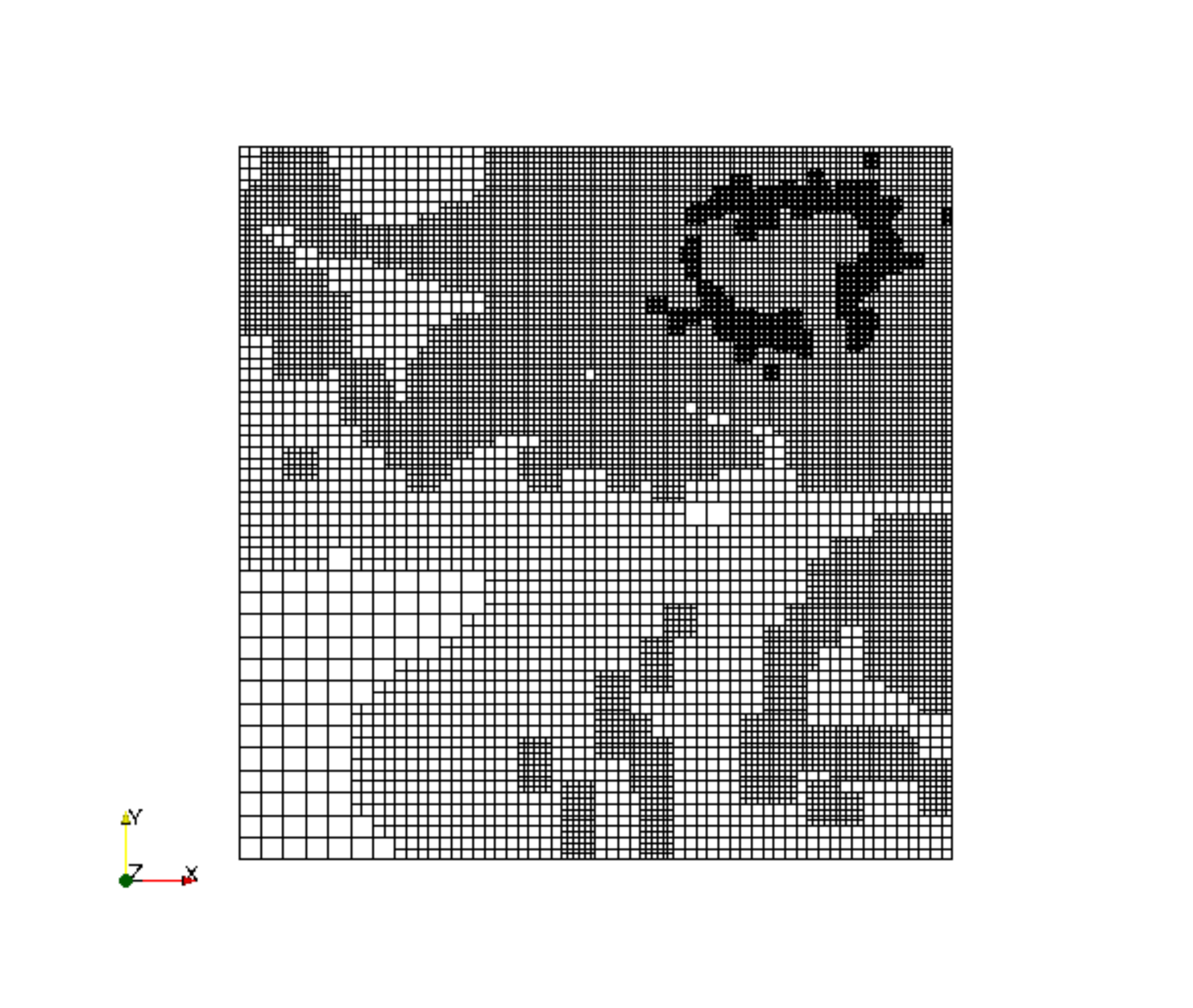}}
	\qquad
	\subfloat[ref. \# 5: \quad ${\Omega}$ and ${\mathcal{K}}_h$]{
	\includegraphics[width=5.4cm, trim={2cm 2cm 2cm 2cm}, clip]{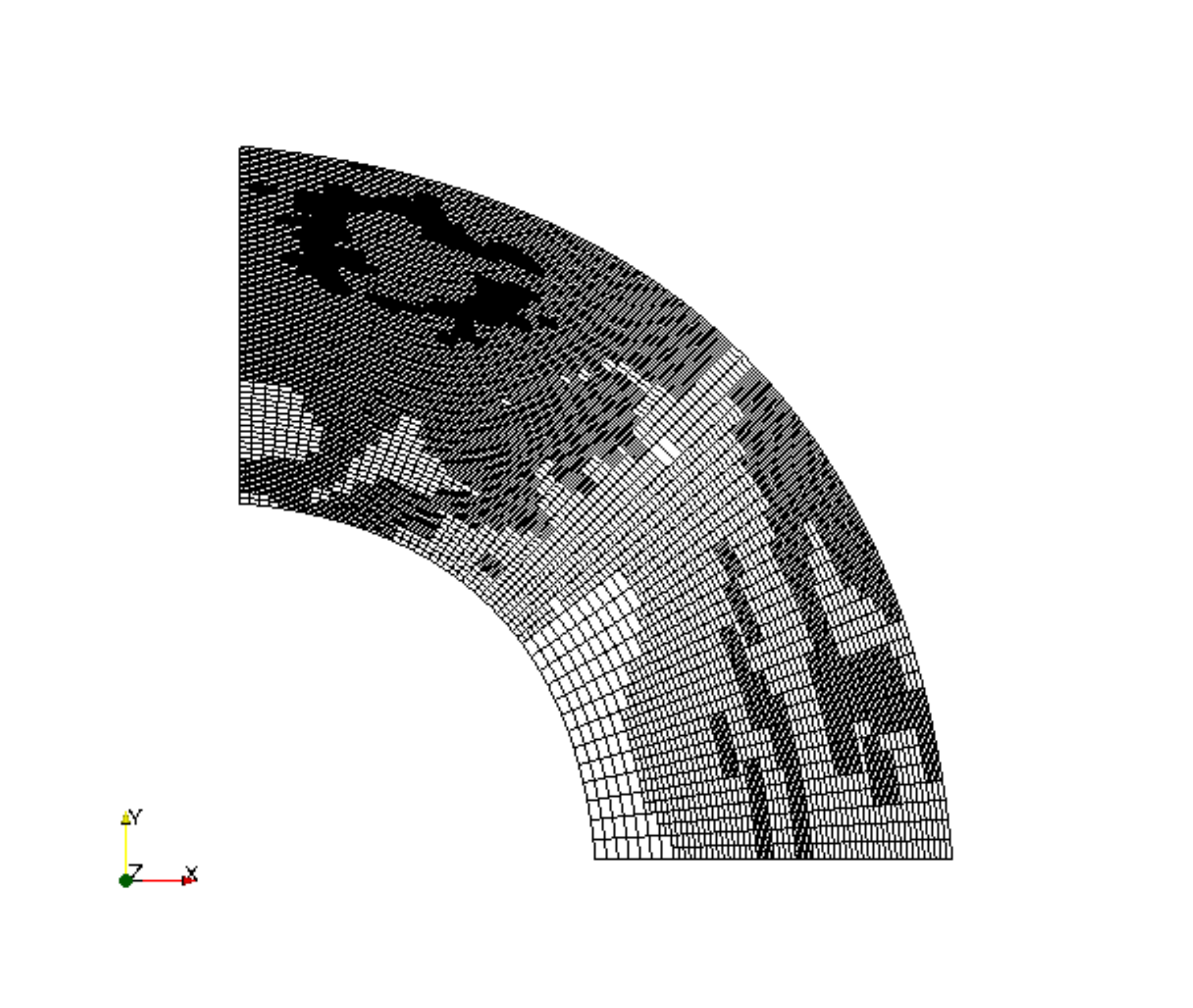}} 
	}
	%
	\caption{\small {\em Ex. \ref{ex:quater-annulus-example-10}}. 
	Comparison of meshes on the physical and parametrical domains w.r.t. adaptive ref. steps.}
	\label{fig:quater-annulus-example-10-theta-40-domains-comparison}
\end{figure}

\end{example}

\begin{example}
\label{ex:unit-cube-example-13}
\rm 
Last two examples are dedicated to three-dimensional problems. Let $\Omega = (0, 1)^3 \subset \mathds{R}^3$,
\begin{alignat*}{3}
u 	& = (1 - x_1) \, x_1^2 \,  (1 - x_2) \, x^2_2 \,  (1 - x_3) \, x^2_3	& \quad  \mbox{in} \quad \Omega, \\
u_D 	& = 0											& \quad \; \mbox{on} \quad \Gamma.
\end{alignat*}
%
The uniform refinement strategy, assuming that $u_h \in S_h^{2, 2}$ and 
$\flux_h \in S_{6h}^{3, 3} \oplus S_{6h}^{3, 3}$, provides numerical results illustrated in Tables 
\ref{tab:unit-cube-example-13-error-majorant-v-2-y-3-uniform-ref} and 
\ref{tab:unit-cube-example-13-times-v-2-y-3-uniform-ref}. 
Comparison of the majorant performance to the accuracy of residual error estimates, i.e., $\Ieff(\maj{})  = 1.2378$ 
and $\Ieff(\overline{\eta}) = 13.4166$, confirms that the latter one always overestimates the error (even for such a 
smooth exact solution). Moreover, the computational costs of majorant generation is a hundred times less than
the computational time for a primal variable, namely, 
%
$
\tfrac{t_{\rm as}(u_h)}{t_{\rm as}(\flux_h)} \approx 578$, and
$\tfrac{t_{\rm sol}(u_h)}{t_{\rm sol}(\flux_h)} \approx 136.$
It is important to note that values of the last three columns of Table 
\ref{tab:unit-cube-example-13-times-v-2-y-3-uniform-ref} illustrate sub-optimal time for the element-wise 
evaluation. As mentioned earlier, this issue is related to the implementation of element-wise iterator currently used 
in G+Smo and will be addressed in the follow-up reports. 
\begin{table}[!t]
\footnotesize
\centering
\newcolumntype{g}{>{\columncolor{gainsboro}}c} 	
\begin{tabular}{c|c|ccc|gc|c}
\quad \# ref. \quad & 
\quad  $\| \nabla e \|_\Omega$ \qquad   & 	  
\quad \quad \quad \quad $\overline{\rm M}$ \quad \quad \quad \quad &    
\quad \quad   $\mdI$ \qquad \quad & 	       
\quad \qquad $\mfI$ \qquad \quad  &  
\qquad $\Ieff (\overline{\rm M})$ \qquad & 
\qquad $\Ieff (\overline{\rm \eta})$ \qquad & 
\qquad \quad e.o.c. \qquad \quad \\
\midrule
  2 &   9.7268e-04 &   1.1328e-03 &   1.0104e-03 &   6.6641e-04 &       1.1646 &      14.2071 &   4.9421 \\
   4 &   5.7843e-05 &   7.9238e-05 &   7.8966e-05 &   1.4807e-06 &       1.3699 &      13.4654 &   2.7348 \\
   6 &   3.6029e-06 &   4.6376e-06 &   4.5126e-06 &   6.8010e-07 &       1.2872 &      13.4195 &   2.1808 \\
   8 &   2.2513e-07 &   2.8772e-07 &   2.7822e-07 &   5.1726e-08 &       1.2780 &      13.4166 &   2.0451 \\
\end{tabular}
\caption{\small {\em Ex. \ref{ex:unit-cube-example-13}}. 
Error, majorant (with dual and reliability terms), efficiency indices, error, e.o.c. 
w.r.t. unif. ref. steps.}
\label{tab:unit-cube-example-13-error-majorant-v-2-y-3-uniform-ref}
\end{table}
\begin{table}[!t]
 \footnotesize
\centering
\newcolumntype{g}{>{\columncolor{gainsboro}}c} 	
\begin{tabular}{c|cc|cg|cg|cgc}
\# ref.  & 
\# d.o.f.($u_h$) &  \# d.o.f.($\flux_h$) &  
\; $t_{\rm as}(u_h)$ \; & 
\; $t_{\rm as}(\flux_h)$ \; & 
\; $t_{\rm sol}(u_h)$ \; & 
\; $t_{\rm sol}(\flux_h)$ \; &
$t_{\rm e/w}(\| \nabla e \|)$ & 
$t_{\rm e/w}(\overline{\rm M})$ & 
$t_{\rm e/w}(\overline{\eta})$ \\
\midrule
   2 &         64 &         64 &     0.0121 &     0.0292 &           0.0000 &           0.0131 &       0.0013 &       0.0079 &       0.0065 \\
   4 &       1000 &         64 &     0.1535 &     0.0477 &           0.0005 &           0.0084 &       0.1913 &       0.2480 &       0.2962 \\
   6 &      39304 &         64 &     6.2179 &     0.0529 &           0.2208 &           0.0106 &      11.7921 &      15.9026 &      17.5110 \\
   8 &    2197000 &        343 &   324.7514 &     0.5641 &          54.6425 &           0.1760 &     560.9245 &     452.5727 &     865.4832 \\
\end{tabular}
\caption{\small {\em Ex. \ref{ex:unit-cube-example-13}}. 
Time for assembling and solving the systems generating d.o.f. of $u_h$ and $\flux_h$ 
as well as the time spent on e/w evaluation of error, majorant, and 
residual error estimator w.r.t. unif. ref. steps.}
\label{tab:unit-cube-example-13-times-v-2-y-3-uniform-ref}
\end{table}
The results obtained while performing an adaptive refinement strategy with the bulk marking criterion 
${\mathds{M}}_{{\rm \bf BULK}}(0.4)$ are summarised in Tables 
\ref{tab:unit-cube-example-13-error-majorant-v-2-y-3-adaptive-ref}--\ref{tab:unit-cube-example-13-times-v-2-y-3-adaptive-ref}.
The obtained ratios between majorants and the true error as well as the computational time it requires to 
be generated, i.e., 
$$
\tfrac{t_{\rm as}(u_h)}{t_{\rm as}(\flux_h)} \approx \tfrac{324.75145}{0.5641} \approx  575  
\quad \mbox{and} \quad
\tfrac{t_{\rm sol}(u_h)}{t_{\rm sol}(\flux_h)} \approx \tfrac{54.6425}{0.1760} \approx 310,$$
confirm the efficiency of the functional approach to the error control for these type of problems. 
The mesh evolution for the considered adaptive refinement strategy is illustrated in Figure 
\ref{fig:unit-cube-example-13-times-v-2-y-3-adaptive-ref} (just for two refinement steps). One can see that 
the main refinement is performed closer to the centre of the computational domain $\Omega$, which is similar 
to the results obtained for two-dimensional case.

\begin{table}[!t]
\footnotesize
\centering
\newcolumntype{g}{>{\columncolor{gainsboro}}c} 	
\begin{tabular}{c|c|ccc|gc|c}
\quad \# ref. \quad & 
\quad  $\| \nabla e \|_\Omega$ \qquad   & 	  
\quad \quad \quad \quad $\overline{\rm M}$ \quad \quad \quad \quad &    
\quad \quad   $\mdI$ \qquad \quad & 	       
\quad \qquad $\mfI$ \qquad \quad  &  
\qquad $\Ieff (\overline{\rm M})$ \qquad & 
\qquad $\Ieff (\overline{\rm \eta})$ \qquad & 
\qquad \quad e.o.c. \qquad \quad \\
\midrule
   3 &   2.3387e-04 &   2.7678e-04 &   2.7288e-04 &   2.1210e-05 &       1.1835 &      13.6135 &   3.5152 \\
   5 &   1.9918e-05 &   2.7263e-05 &   2.3917e-05 &   1.8207e-05 &       1.3687 &      11.7620 &   2.2298 \\
   7 &   2.2272e-06 &   3.1275e-06 &   2.9019e-06 &   1.2276e-06 &       1.4042 &      12.4614 &   2.0710 \\
\end{tabular}
\caption{\small {\em Ex. \ref{ex:unit-cube-example-13}}. 
Error, majorant (with dual and reliability terms), efficiency indices, error, e.o.c. w.r.t. adaptive ref. steps.}
\label{tab:unit-cube-example-13-error-majorant-v-2-y-3-adaptive-ref}
\end{table}
\begin{table}[!t]
 \footnotesize
\centering
\newcolumntype{g}{>{\columncolor{gainsboro}}c} 	
\begin{tabular}{c|cc|cg|cg|cgc}
\# ref.  & 
\# d.o.f.($u_h$) &  \# d.o.f.($\flux_h$) &  
\; $t_{\rm as}(u_h)$ \; & 
\; $t_{\rm as}(\flux_h)$ \; & 
\; $t_{\rm sol}(u_h)$ \; & 
\; $t_{\rm sol}(\flux_h)$ \; &
$t_{\rm e/w}(\| \nabla e \|)$ & 
$t_{\rm e/w}(\overline{\rm M})$ & 
$t_{\rm e/w}(\overline{\eta})$ \\
\midrule
   1 &         27 &         64 &     0.0146 &     0.0718 &           0.0000 &           0.0049 &       0.0066 &       0.0445 &       0.0213 \\
   3 &        216 &         64 &     0.3379 &     0.0773 &           0.0001 &           0.0107 &       0.4165 &       0.8699 &       1.0376 \\
   5 &       4074 &         64 &    21.5658 &     0.0837 &           0.0420 &           0.0074 &      26.7396 &      39.1294 &      50.8124 \\
   7 &     127265 &         64 &  3054.9337 &     0.0838 &           3.2532 &           0.0087 &    1415.2231 &    1720.7416 &    2433.3058 \\
\end{tabular}
\caption{\small {\em Ex. \ref{ex:unit-cube-example-13}}. 
Time for assembling and solving the systems generating d.o.f. of $u_h$ and $\flux_h$ 
as well as the time spent on e/w evaluation of error, majorant, and 
residual error estimator w.r.t. adaptive ref. steps.}
\label{tab:unit-cube-example-13-times-v-2-y-3-adaptive-ref}
\end{table}

\begin{figure}[!t]
	\centering
	%
	%
	{
	\subfloat[ref. \#  2: \quad $\mathcal{K}_h$ with \# d.o.f.($u_h$) = 216]{
	\includegraphics[width=6cm, trim={0cm 0cm 2cm 2cm}, clip]{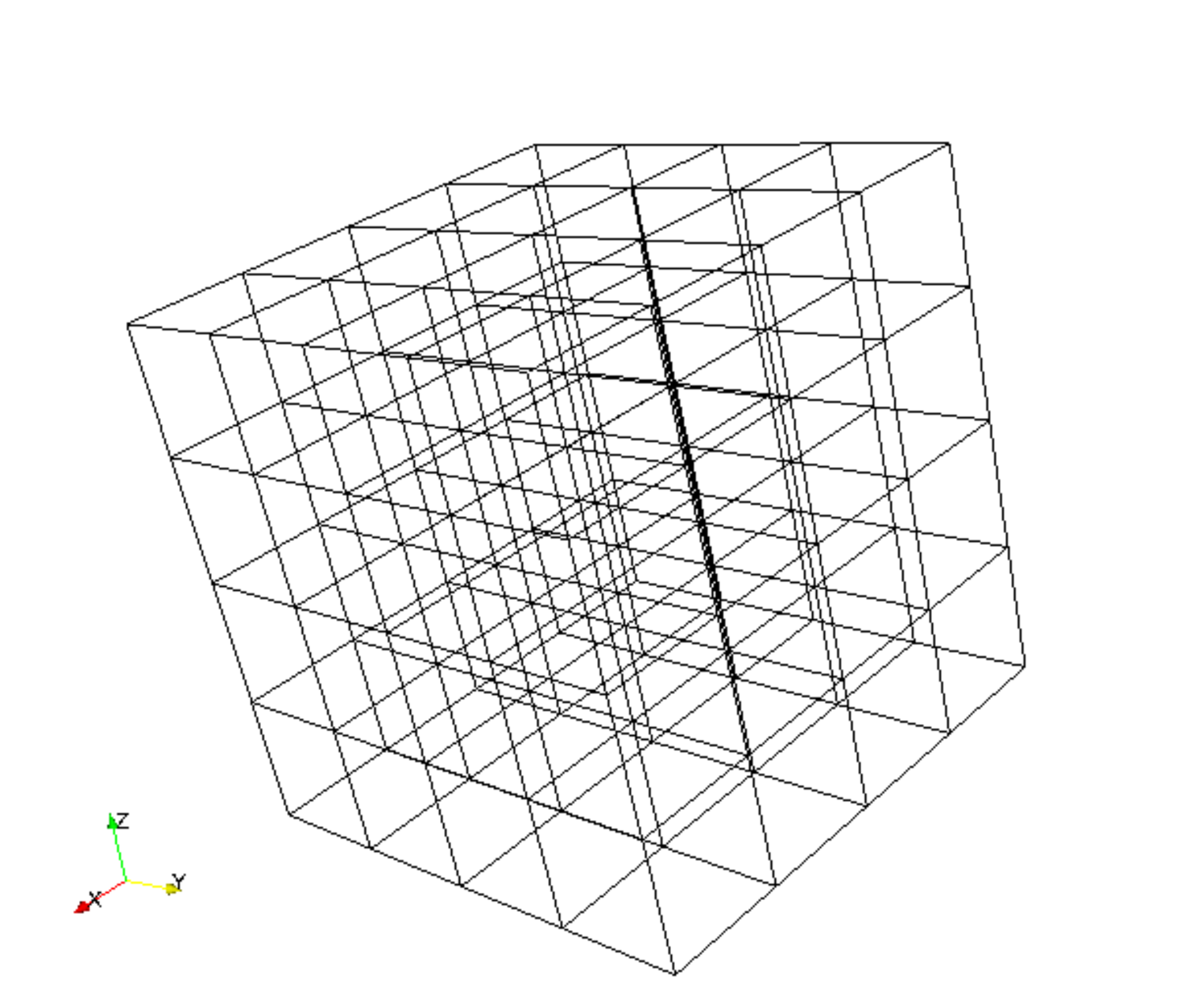}}
	} \quad
	%
	%
	{
	\subfloat[ref. \#  4: \quad $\mathcal{K}_h$ with \# d.o.f.($u_h$) = 4074]{
	\includegraphics[width=6cm, trim={0cm 0cm 2cm 2cm}, clip]{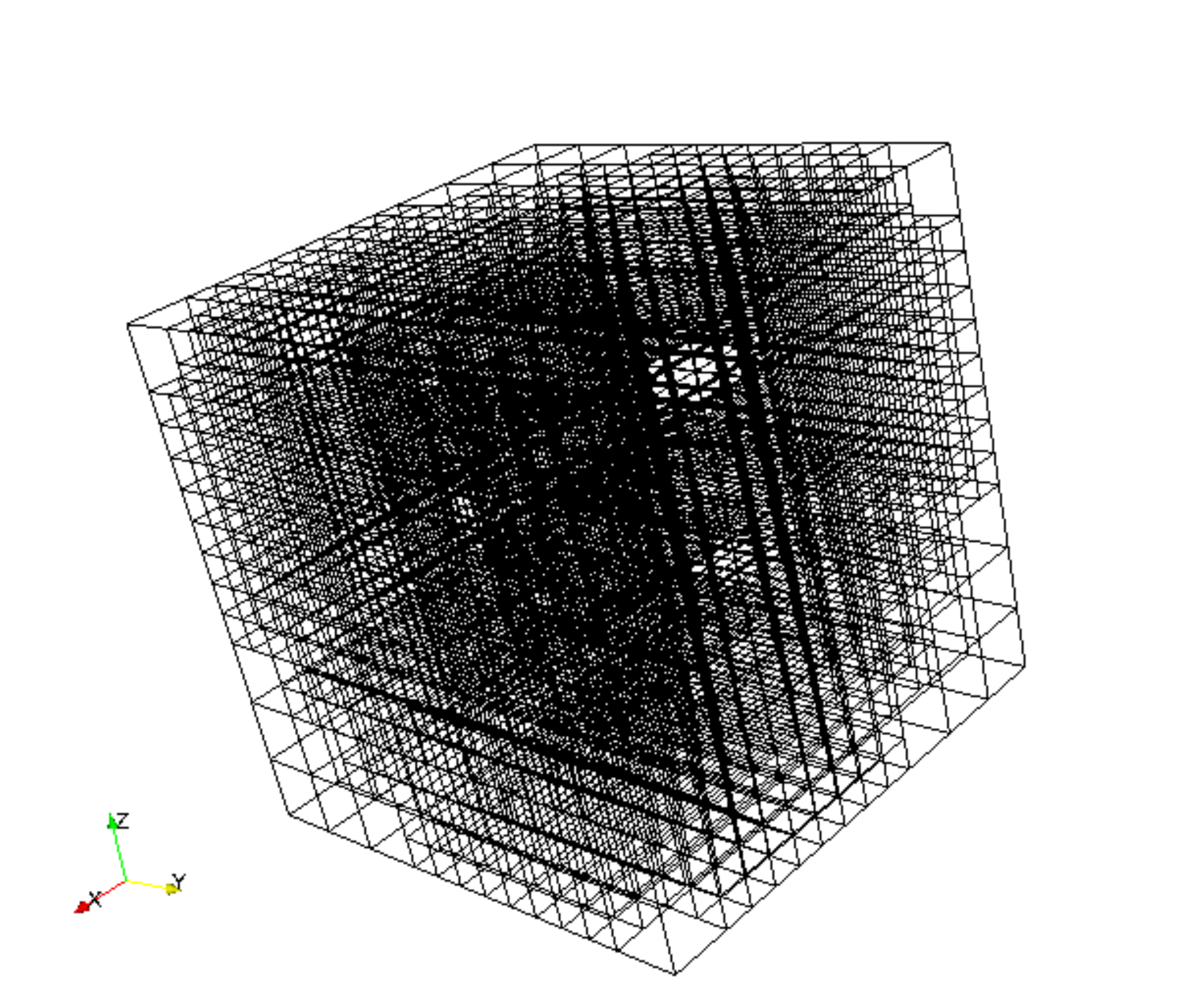}
	}
	}
\caption{\small {\em Ex. \ref{ex:unit-cube-example-13}}. Evolution of meshes on the physical unit cube w.r.t. adaptive ref. steps with bulk
marking criterion ${\mathds{M}}_{\rm BULK}(0.4)$.}
\label{fig:unit-cube-example-13-times-v-2-y-3-adaptive-ref}
\end{figure}
\end{example}

\begin{example}	
\label{ex:g-domain-example-14}
\rm 
Finally, in the last example, we test functional error estimates on a complex three-dimensional geometry of a $G$-shape. 
The exact solution, RHS, and Dirichlet BC are defined as follows:
\begin{alignat*}{3}
u 	& = 10 \, \cos x_1 \, e^{x_2}\, x_3	& \quad  \mbox{in} \quad \Omega, \\
f 	& = 0						& \quad  \mbox{in} \quad \Omega, \\
u_D 	& = 10 \, \cos x_1 \, e^{x_2} \, x_3	& \quad \; \mbox{on} \quad \Gamma.
\end{alignat*}

\begin{figure}[!t]
	\centering
	\includegraphics[width=8cm, trim={0cm 2cm 6cm 2cm}, clip]{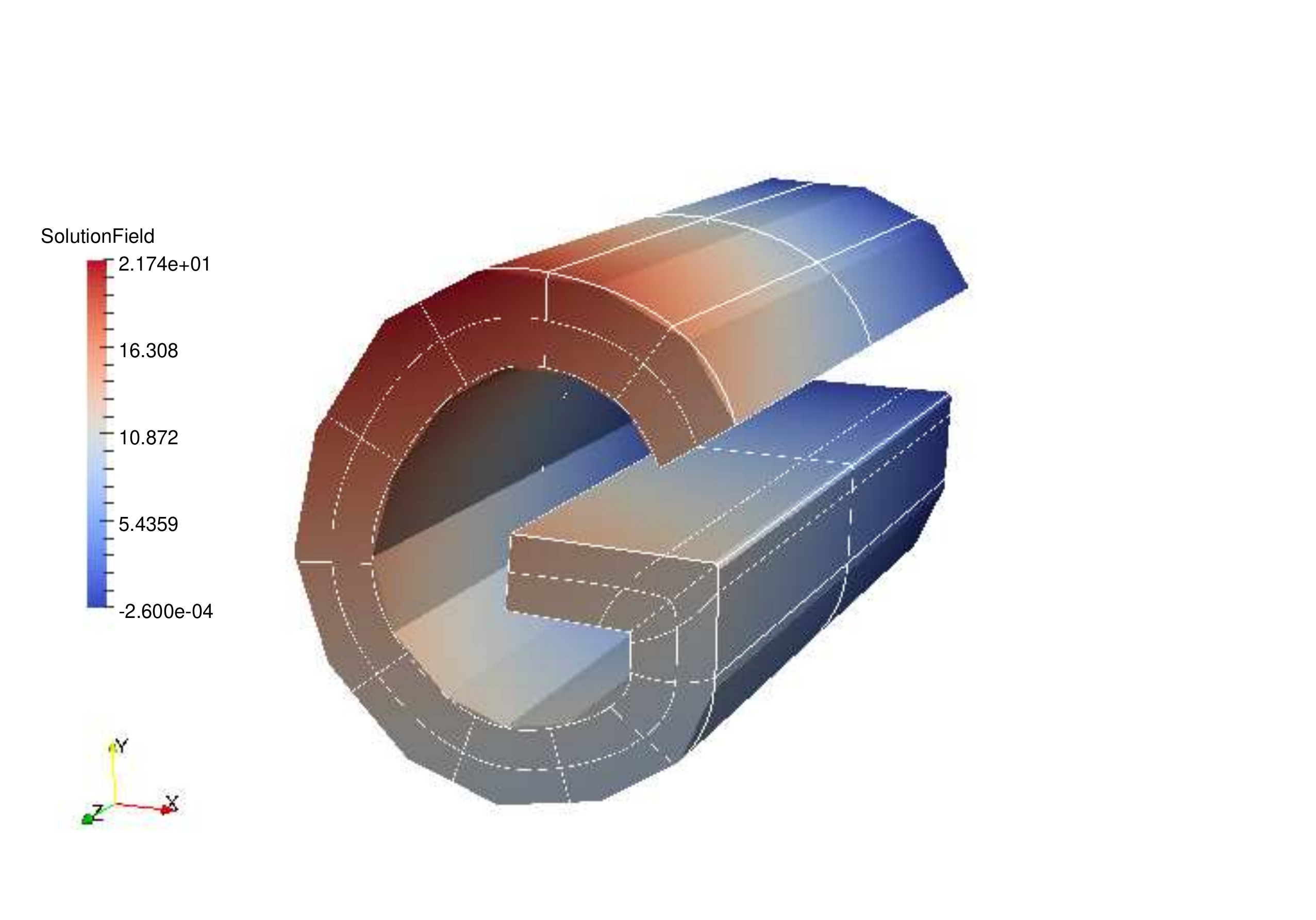}
	\caption{\small {\em Ex. \ref{ex:g-domain-example-14}}. 
	Exact solution $u = 10 \, \cos x_1 \, e^{x_2}\, x_3$.}
	\label{fig:ex-g-domain-example-14-exact-solution}
\end{figure}
%
Let us consider $N_{\rm ref, 0} = 1$ initial unif. ref. step of the mesh assigned to the geometry (illustrated in 
Figure \ref{fig:ex-g-domain-example-14-exact-solution}). The numerical results of 5 adaptive ref. steps are summarised in 
Tables \ref{tab:g-domain-example-14-error-majorant-v-2-y-3-adaptive-ref}--\ref{tab:g-domain-example-14-times-v-2-y-3-adaptive-ref}. 
From Table \ref{tab:g-domain-example-14-error-majorant-v-2-y-3-adaptive-ref} analysing the efficiency of the error 
estimates, it is obvious that $\overline{\rm M}$ is at least 40 times sharper than the residual error indicator. 
If one analyses the computation costs, it is easy to notice that the assembling time of the system \eqref{eq:system-fluxh} 
is better than the assembling time of \eqref{eq:system-uh}, i.e., $\tfrac{t_{\rm as}(u_h)}{t_{\rm as}(\flux_h)} \approx 14$. 
However, the costs for solving the system \eqref{eq:system-fluxh} is higher, i.e., 
$\tfrac{t_{\rm sol}(u_h)}{t_{\rm sol}(\flux_h)} \approx \tfrac{1}{12}$. The element-wise evaluation of $\| \nabla e\|$, 
$\overline{\rm M}$, and $\overline{\rm \eta}$ must take at least as long as $t_{\rm as}(u_h)$ (since both are 
performed with $u_h$ element-wise), which is confirmed 
in the last three columns of Table \ref{tab:g-domain-example-14-times-v-2-y-3-adaptive-ref}.

\begin{table}[!t]
\footnotesize
\centering
\newcolumntype{g}{>{\columncolor{gainsboro}}c} 	
\begin{tabular}{c|c|ccc|gc|c}
\quad \# ref. \quad & 
\quad  $\| \nabla e \|_\Omega$ \qquad   & 	  
\quad \quad \quad \quad $\overline{\rm M}$ \quad \quad \quad \quad &    
\quad \quad   $\mdI$ \qquad \quad & 	       
\quad \qquad $\mfI$ \qquad \quad  &  
\qquad $\Ieff (\overline{\rm M})$ \qquad & 
\qquad $\Ieff (\overline{\rm \eta})$ \qquad & 
\qquad \quad e.o.c. \qquad \quad \\
\midrule
   3 &   1.5495e-02 &   1.6552e-02 &   1.6016e-02 &   2.9195e-03 &       1.0682 &      42.1400 &   3.1601 \\
   4 &   5.5530e-03 &   6.7468e-03 &   6.0466e-03 &   3.8103e-03 &       1.2150 &      59.4362 &   3.0093 \\
   5 &   2.3080e-03 &   3.6335e-03 &   2.7857e-03 &   4.6130e-03 &       1.5743 &      47.1796 &   2.5110 \\
\end{tabular}
\caption{\small {\em Ex. \ref{ex:g-domain-example-14}}. 
Error, majorant (with dual and reliability terms), efficiency indices, error, e.o.c. w.r.t. adaptive ref. steps.}
\label{tab:g-domain-example-14-error-majorant-v-2-y-3-adaptive-ref}
\end{table}
\begin{table}[!t]
 \footnotesize
\centering
\newcolumntype{g}{>{\columncolor{gainsboro}}c} 	
\begin{tabular}{c|cc|cg|cg|cgc}
\# ref.  & 
\# d.o.f.($u_h$) &  \# d.o.f.($\flux_h$) &  
\; $t_{\rm as}(u_h)$ \; & 
\; $t_{\rm as}(\flux_h)$ \; & 
\; $t_{\rm sol}(u_h)$ \; & 
\; $t_{\rm sol}(\flux_h)$ \; &
$t_{\rm e/w}(\| \nabla e \|)$ & 
$t_{\rm e/w}(\overline{\rm M})$ & 
$t_{\rm e/w}(\overline{\eta})$ \\
\midrule
   3 &       1108 &        575 &     8.5398 &     3.4636 &     0.0056 &     3.2101 &      11.8510 &       7.5239 &      21.5933 \\
   4 &       3082 &        575 &    31.4262 &     2.9348 &     0.0686 &     3.7503 &      38.9502 &      26.7215 &      77.4521 \\
   5 &       8798 &        575 &    97.6563 &     7.1416 &     0.2659 &     3.1107 &     122.7187 &     101.5395 &     229.4474 \\
 \end{tabular}
\caption{\small {\em Ex. \ref{ex:g-domain-example-14}}. 
Time for assembling and solving the systems generating d.o.f. of $u_h$ and $\flux_h$ 
as well as the time spent on e/w evaluation of error, majorant, and 
residual error estimator w.r.t. adaptive ref. steps.}
\label{tab:g-domain-example-14-times-v-2-y-3-adaptive-ref}
\end{table}

The evolution of the meshes on both physical and parametric domains is illustrated in Figure 
\ref{fig:g-domain-example-14-meshes-adaptive-ref}.
The final solution and meshes obtained in the final step are presented in Figure 
\ref{fig:g-domain-example-14-final-meshes-adaptive-ref}. 
The latter plots are presented from the axis $O_z$ view-point so that the reader can clearly see the finer parts of the 
physical and parametric meshes.

\begin{figure}[!t]
	\centering
	\subfloat[ref. \#  1: \quad $u_h$ on $\Omega$ with mesh $\mathcal{K}_h$]{	
	\includegraphics[width=7cm, trim={0cm 2cm 6cm 2cm}, clip]{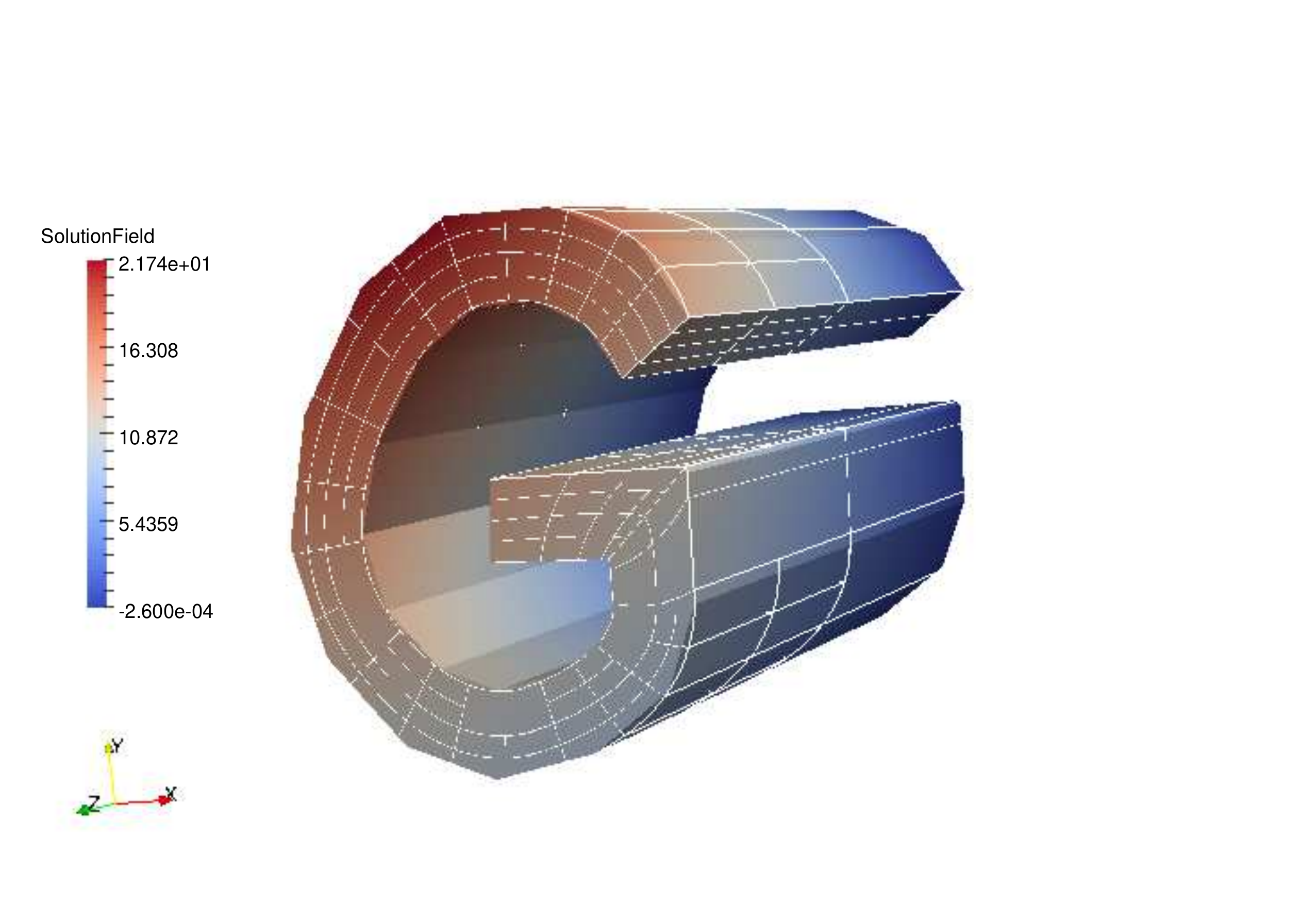}}
	\qquad
	\subfloat[ref. \#  1: \quad $\widehat{\mathcal{K}}_h$]{
	\includegraphics[width=6cm, trim={1cm 0cm 5cm 0cm}, clip]{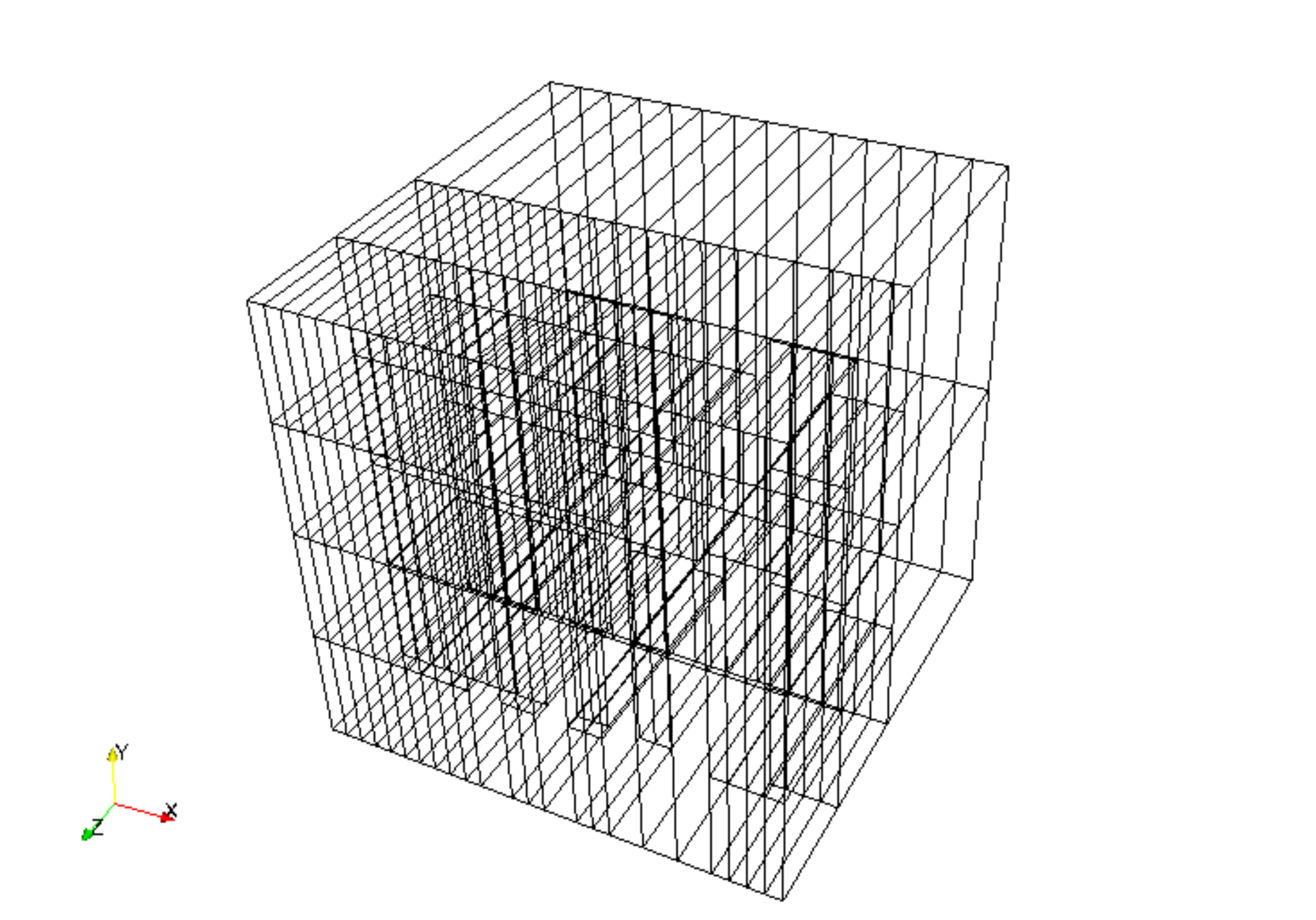}}
	\\
	%
	%
	\subfloat[ref. \#  3: \quad $u_h$ on $\Omega$ with mesh $\mathcal{K}_h$]{	
	\includegraphics[width=7cm, trim={0cm 2cm 6cm 2cm}, clip]{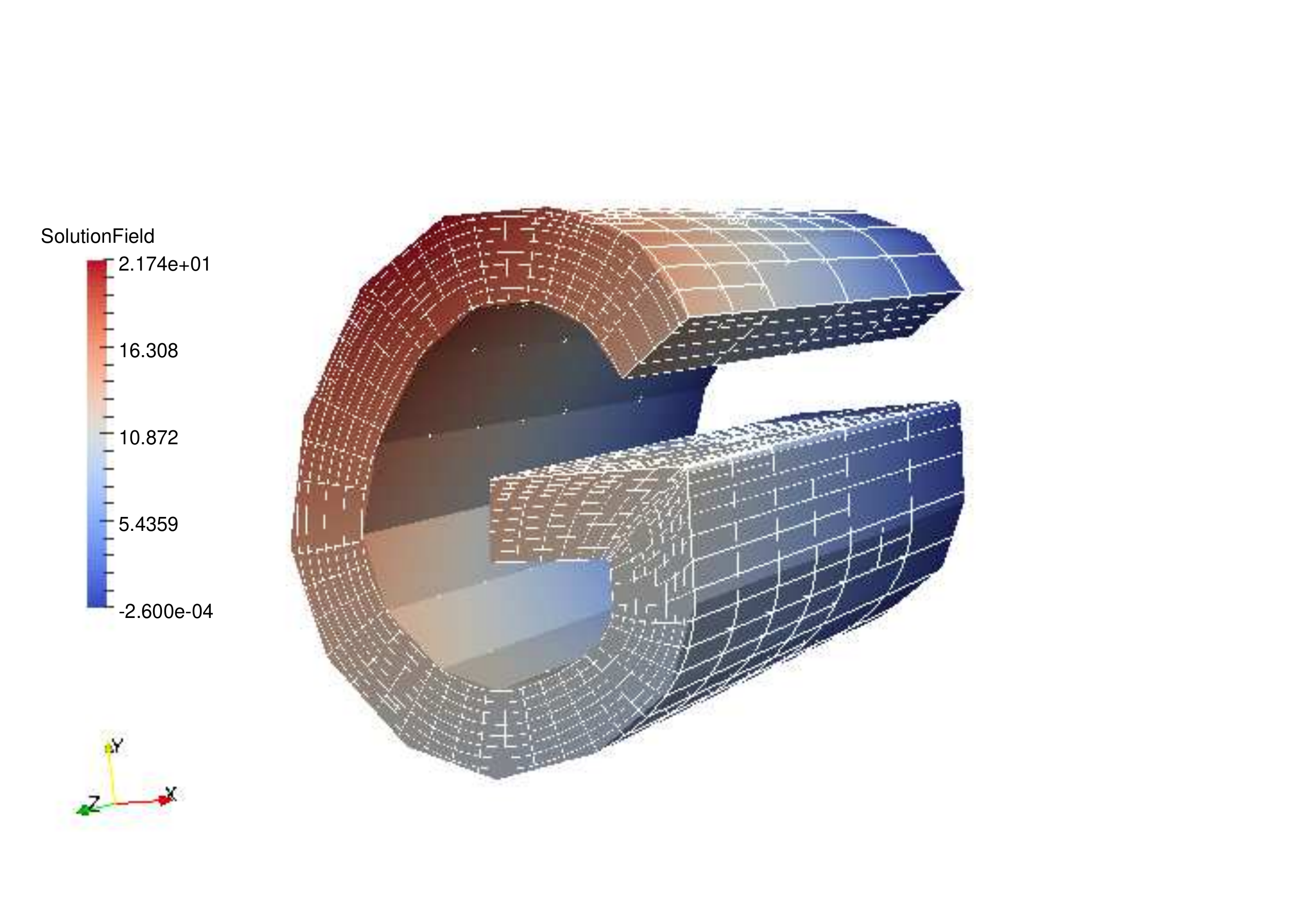}} 
	\qquad
	\subfloat[ref. \#  3: \quad $\widehat{\mathcal{K}}_h$]{
	\includegraphics[width=6cm, trim={1cm 0cm 5cm 0cm}, clip]{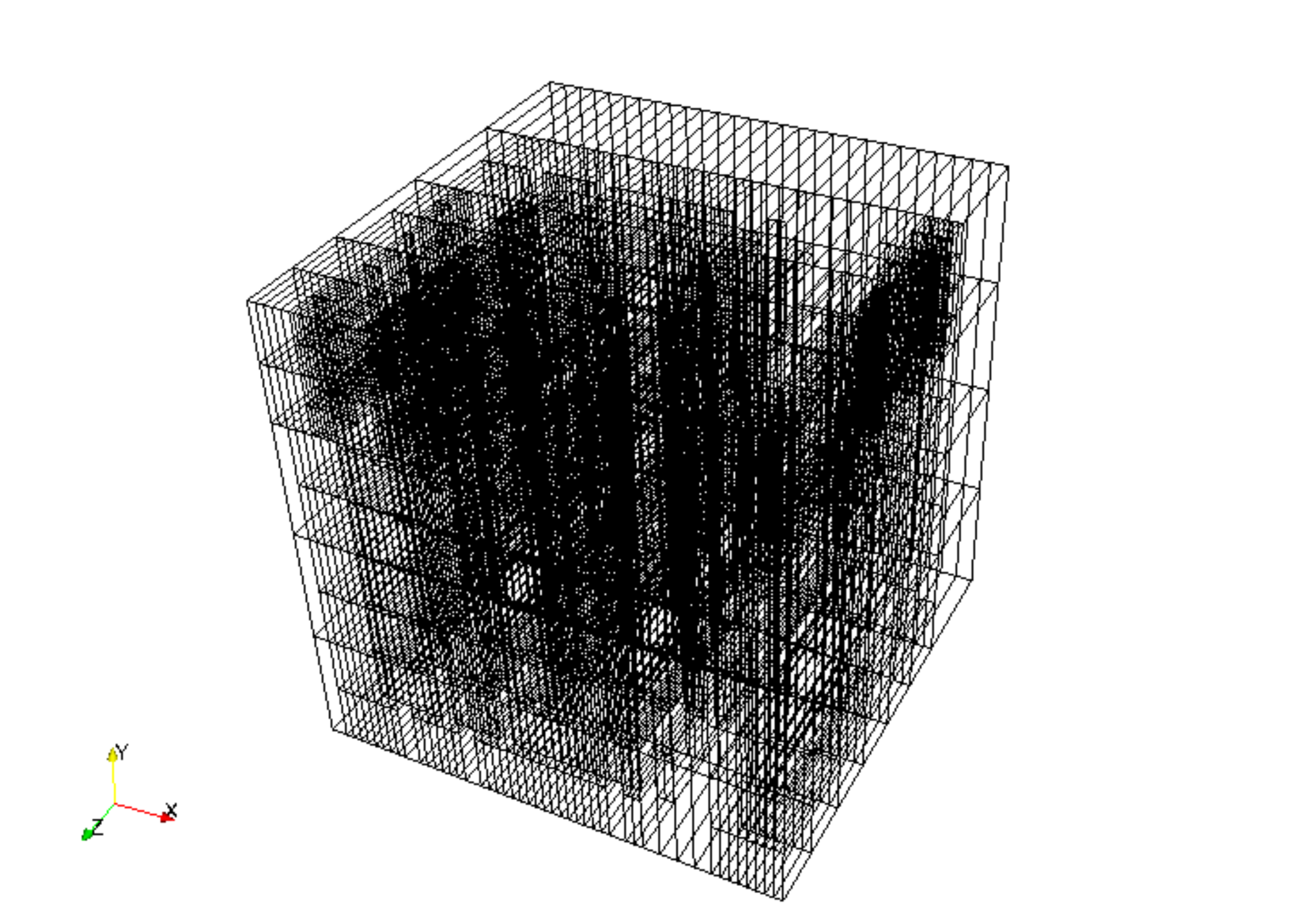}}	
	\caption{\small {\em Ex. \ref{ex:g-domain-example-14}}. 
	Comparison of meshes on the physical and parametrical domains w.r.t. adaptive ref. steps.}
	\label{fig:g-domain-example-14-meshes-adaptive-ref}
\end{figure}

\begin{figure}[!t]
	\centering
	\subfloat[ref. \#  4: \quad $u_h$ on $\Omega$ with mesh $\mathcal{K}_h$]{	
	\includegraphics[width=5cm, trim={1cm 2cm 4cm 2cm}, clip]{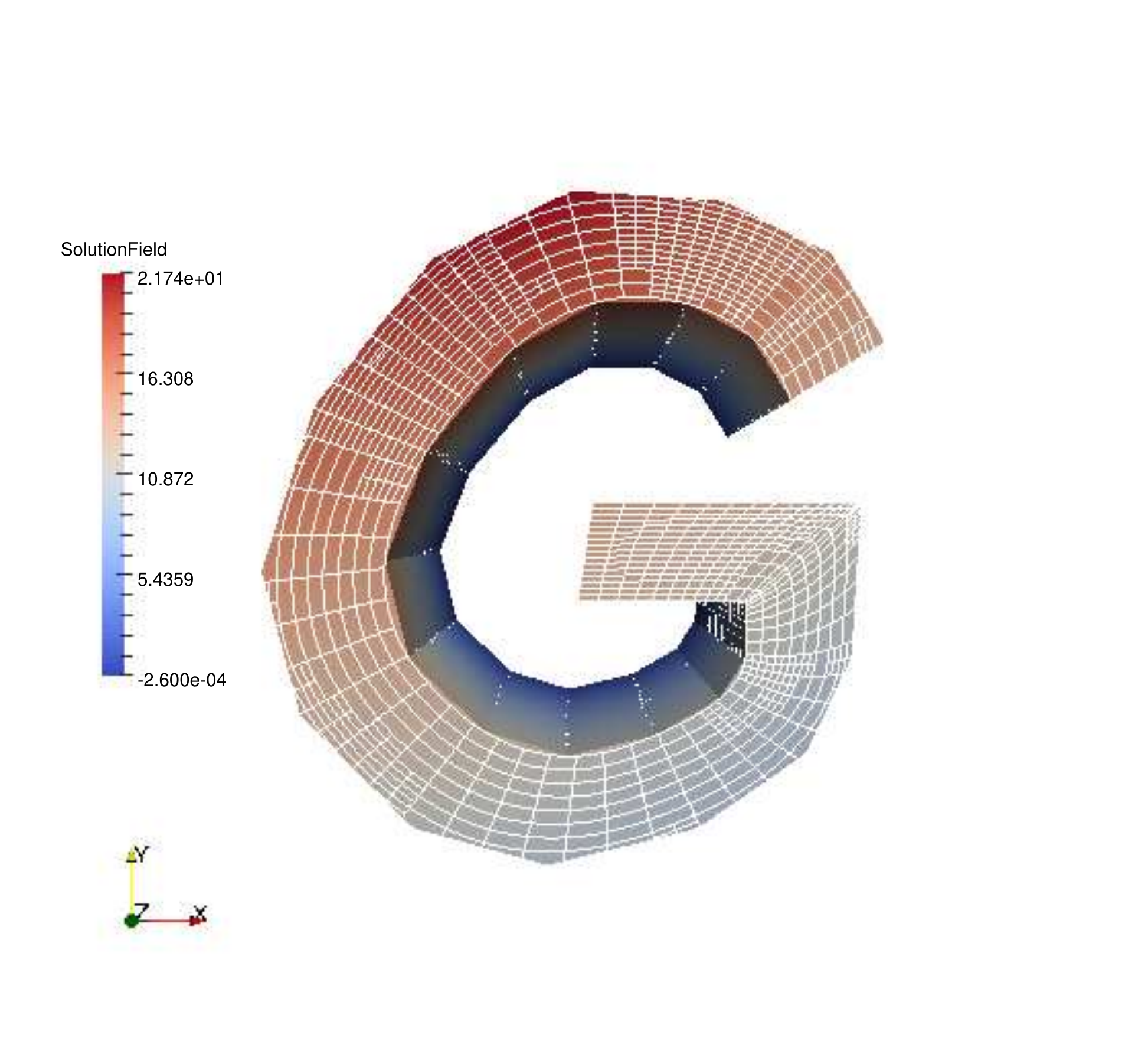}}
	\quad
	\subfloat[ref. \#  4: \quad $t{\mathcal{K}}_h$]{
	\includegraphics[width=5.5cm, trim={1cm 2cm 4cm 2cm}, clip]{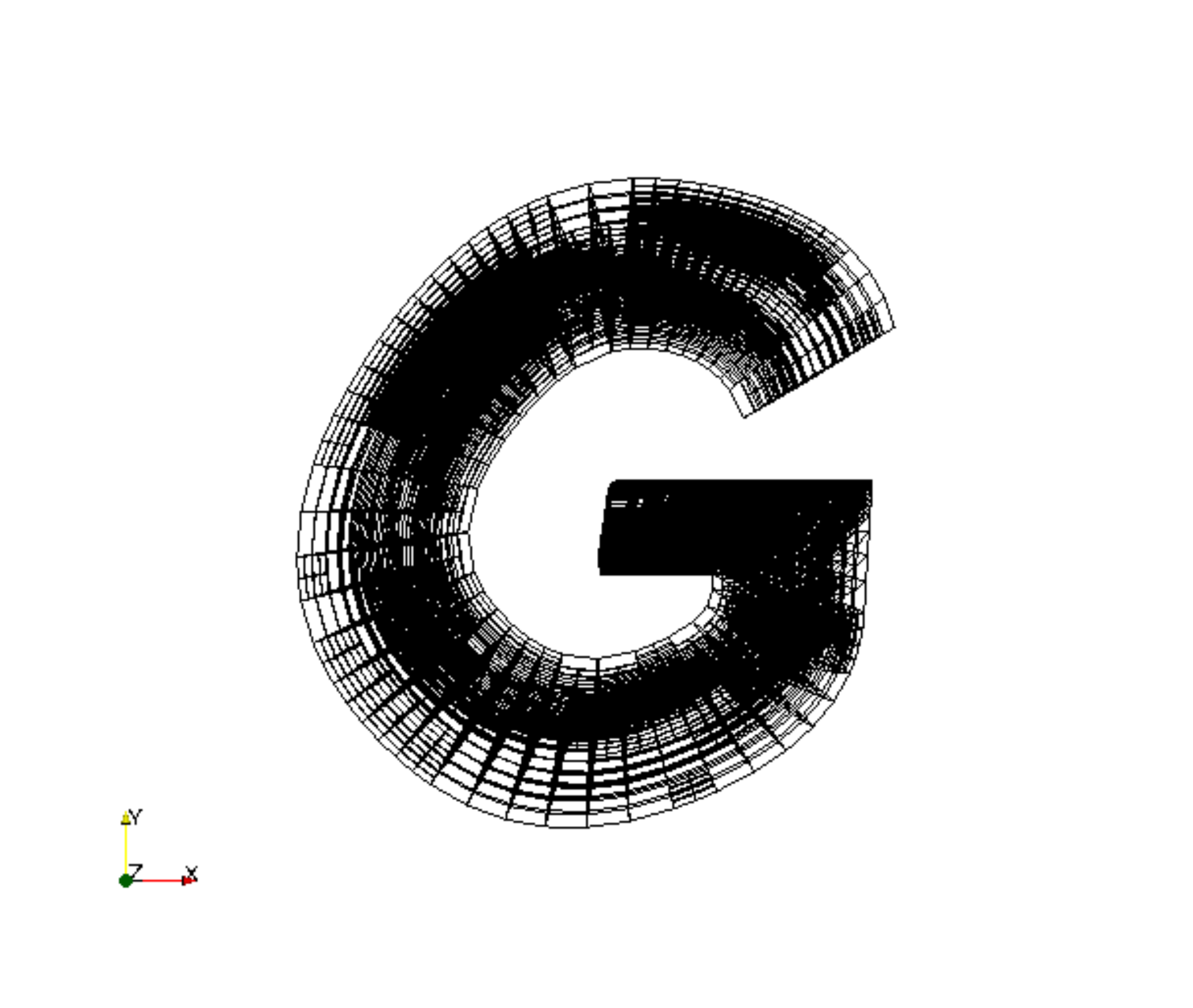}}
	\quad
	\subfloat[ref. \#  4: \quad $\widehat{\mathcal{K}}_h$]{
	\includegraphics[width=5.0cm, trim={1cm 2cm 4cm 2cm}, clip]{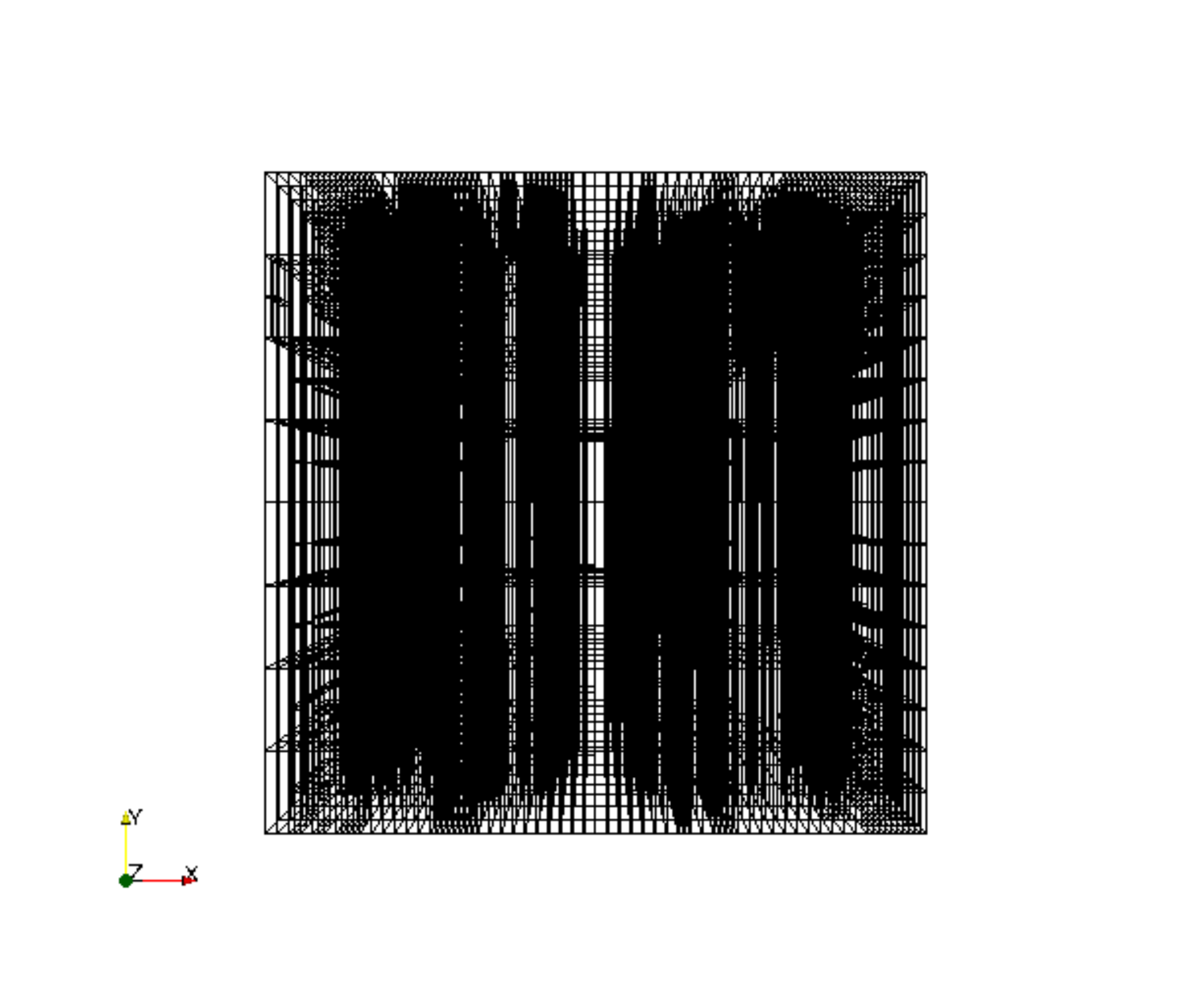}}	
	\caption{\small {\em Ex. \ref{ex:g-domain-example-14}}. 
	Solution on the domain $\Omega$ and corresponding meshes on physical and parametrical domains, 
	view from the axis $0_z$.}
	\label{fig:g-domain-example-14-final-meshes-adaptive-ref}
\end{figure}

\end{example}

{\section*{Conclusion}

In the paper, we have numerically studied two-sided functional error estimates in the framework of IgA schemes. 
The very first derivation of these type of error bounds dates back to 96'-97' (see \cite{LMR:Repin:1997, 
LMR:Repin:1999}). The generalised functional approach applied in these works allows obtaining estimates 
of the error between exact solution and any approximation from the admissible energy space (independently 
on the numerical method used for the reconstruction of this approximation). As the consequence, discussed 
estimates do not depend on any mesh dependent constants and do not pose any requirements (such as Galerkin's 
orthogonality) on the approximate solution. The pioneering study, investigating the application of the majorant
to IgA approximations generated by tensor-product splines in a context of elliptic BVP, can be found in 
\cite{LMR:KleissTomar2015}. Current work extends the study of the fully reliable adaptive IgA schemes using 
THB-splines in combination with the functional two-sided bounds and the local error indicators generated by them.  

In this study, we have highlighted main numerical properties of the functional error estimates 
in the context of both reliable global energy error estimation and efficient local error distribution indication. 
We mainly focused on the algorithmic and implementational parts of two-sided bounds application as well as details 
of their integration into the IgA framework. Presented examples and corresponding numerical results have provided 
an convincing evidence of the global reliability of the majorant and even more importantly local efficiency of the 
error indicator generated by it. Implementation of all the methods was carried out using open  source C++ library 
G+Smo developed primarily by research groups of RICAM and JKU, Linz \cite{gismoweb}.

\newpage
In the majority of the numerical tests presented by Section \ref{sec:numerical-examples}, we were able to
take an advantage of the smoothness of B-(THB-)splines and to reconstruct auxiliary 
functions used by majorant and minorant exploiting higher order splines and coarser meshes (in comparison to the primal 
approximation) in order to achieve the considerable speed-up in the time required for reliable error estimation. 
It naturally allowed to decrease the overall computational time considerably. 
However, for some problems with highly oscillating solutions or solutions possessing singularities, coarsening of the 
mesh used for the auxiliary functions decreases the quality of their reconstruction and, as the consequence  
lowers the effectiveness of the majorant and minorant. Therefore, we had to invest more time to reconstruct 
sufficiently accurate $\flux_h$ and $w_h$ reconstructions. This computational time overhead for assembling and solving 
systems generated from the minimisation of $\overline{\rm M}$ (maximisation of $\underline{\rm M}$) can be however 
minimised by multi-core parallelisation, which stayed out of the focus of the current paper.

Throughout Section \ref{sec:numerical-examples}, majorant has been compared to the residual-type error estimate, 
i.e., its performance confirmed to be always more stable and efficient in comparison to $\overline{\eta}$. Indeed, 
it requires more time for optimising majorant but also reassures that the error is controlled in a reliable way and its 
local distribution is indicated correctly.
Both upper and lower bounds were tested w.r.t. to different marking criterions to make sure that they provide 
an efficient local indication of the error. Meshes produced by the refinement based on the indicator that follows from 
$\overline{\rm M}$ were compared to the meshes obtained by the adaptive strategy using the true error distribution. 
The topology of these meshes appeared to be very close, which confirmed the effectiveness of the majorant application 
in designing a fully adaptive numerical scheme. 
If the exact solution is not provided, the quality of the majorant can be verified by comparison 
of its values to the minorant of the error. That was demonstrated in several numerical examples, where the ration 
between $\overline{\rm M}$ and $\underline{\rm M}$ remained close to one. 

The universality of the functional error estimates allows to easily extend them to parabolic problems. The corresponding
results on the application of these error bounds to space-time IgA techniques in a context of evolutionary problems 
can be found in \cite{LMR:LangerMatculevichRepinArxiv2017}.
}

\noindent
\vskip 10pt
{\bf Acknowledgments}
The research is supported by the Austrian Science Fund (FWF) through the NFN S117-03 project. Implementation 
was carried out using the open-source C++ library G+Smo \cite{gismoweb}.


{\small
\bibliographystyle{plain}
\bibliography{bib/main-thb-splines,bib/lib}
}
\end{document}